\documentclass[11pt]{article}
\usepackage{IEEEtrantools}
\usepackage{mathtools, nccmath}
\usepackage{latexsym}
\usepackage{amsfonts}
\usepackage{amssymb}
\usepackage{amsmath,amsfonts,amsthm}
\usepackage{mathrsfs}
\usepackage{enumerate}
\usepackage{dsfont}
\usepackage{hyperref}
\usepackage{graphicx}
\usepackage{epstopdf}
\usepackage{tikz}
\usetikzlibrary{shapes,arrows.meta,positioning}
\usepackage{cleveref}
\usepackage{tikz-cd}
\usepackage{psfrag}

\usepackage{multirow}
\usepackage{tikz}
\usetikzlibrary{arrows,decorations.markings}
\usepackage{bookmark}
\usepackage{hyperref}
\hypersetup{
     colorlinks   = true,
     citecolor    = blue
}

\newcommand{\newsection}[1]{\setcounter{equation}{0}
\setcounter{dfn}{0}
\section{#1}}

\newtheorem{dfn}{Definition}[section]
\newtheorem{thm}[dfn]{Theorem}
\newtheorem{lmma}[dfn]{Lemma}
\newtheorem{ppsn}[dfn]{Proposition}
\newtheorem{crlre}[dfn]{Corollary}
\newtheorem{xmpl}[dfn]{Example}
\newtheorem{rmrk}[dfn]{Remark}
\newtheorem{notation}[dfn]{Notation}

\DeclareMathOperator*{\dprime}{\prime \prime}

\newcommand{\bbc}{\mathbb{C}}
\newcommand{\bbz}{\mathbb{Z}}
\newcommand{\bbn}{\mathbb{N}}

\def \qed { \mbox{}\hfill
$\Box$\vspace{1ex}}

\title{Sections and Chapters}

\begin{document}

\author{\sc{Keshab Chandra Bakshi\,\footnote{\,Partially supported by DST INSPIRE Faculty grant DST/INSPIRE/04/2019/002754}\,\,, Satyajit Guin}\,\footnote{\,Partially supported by SERB grant MTR/2021/000818}}
\title{Relative position between a pair of spin model subfactors}
\maketitle


\begin{abstract}

Jones pioneered the theory of subfactors, which may be regarded as a quantized version of closed subspaces in Hilbert space. It deals with the relative position of a single factor inside an ambient factor. The level of intricacy increases considerably if there are more than two factors involved. Indeed, Jones proposed the study of two subfactors of a $II_1$ factor as a quantization of two closed subspaces in a Hilbert space. The Pimsner-Popa probabilistic constant, Sano-Watatani angle, interior and exterior angle, and Connes-St{\o}rmer relative entropy (along with a slight variant of it) are a few key invariants for pair of subfactors that analyze their relative position. In practice, however, the explicit computation of these invariants is often difficult.

In this article, we provide an in-depth analysis of  a special class of two subfactors, namely a pair of spin model subfactors of the hyperfinite type $II_1$ factor $R$. We first characterize when two distinct $n\times n$ complex Hadamard matrices give rise to distinct spin model subfactors. Then, a detailed investigation has been carried out for pairs of (Hadamard equivalent) complex Hadamard matrices of order $2\times 2$ as well as Hadamard inequivalent complex Hadamard matrices of order $4\times 4$. To the best of our knowledge, this article is the first instance in the literature where the exact value of the Pimsner-Popa probabilistic constant and the noncommutative relative entropy for pairs of (non-trivial) subfactors have been obtained. Furthermore, we prove the factoriality of the intersection of the corresponding pair of subfactors using the `commuting square technique'. En route, we construct an infinite family of potentially new subfactors of $R$. All these subfactors are irreducible with Jones index $4n,n\geq 2$. As a consequence, the rigidity of the interior angle between the spin model subfactors is established. Last but not least, we explicitly compute the Sano-Watatani angle between the spin model subfactors.

\end{abstract}
\bigskip

{\bf AMS Subject Classification No.:} {\large 46}L{\large 37}, {\large 46}L{\large 55}, {\large 46}L{\large 10}, {\large 37}A{\large 35}.
\smallskip


\hypersetup{linkcolor=blue}
\bigskip

\tableofcontents
\bigskip


\newsection{Introduction}\label{Sec 0}

\subsection{Notation}
\begin{enumerate}
\item Throughout the article underlying field is $\bbc$ and we simply write $M_n$ instead of $M_n(\bbc)$ to denote type $I_n$ factors. $\Delta_n$ will denote the diagonal subalgebra (Masa) in $M_n$.
\item $M_n^{(k)}$ will denote tensor product of $M_n$ with itself $k$-times. For $A\in M_n,\,A^{(k)}$ will denote the matrix $A^{\,\otimes\,k}$.
\item $k\times k$ diagonal matrices will be denoted by $\mbox{diag}\{\mu_1,\ldots,\mu_k\}$, where $\mu_j\in\bbc$.
\item Block diagonal matrix of the form $\left[{\begin{smallmatrix}
A_1 &  &  & \\
 & \ddots &  & \\
 &  &  & A_k\\
\end{smallmatrix}}\right]$ will be denoted by $\mbox{bl-diag}\{A_1,\ldots,A_k\}$.
\item For $A\in M_k,\,I_n\otimes A$ denotes the matrix $\mbox{bl-diag}\{A,\ldots,A\}$ in $M_{nk}$.
\item The following Pauli spin matrices are used throughout the article on several occasions, especially in \Cref{Sec 4},
\[
\sigma_1=\left[{\begin{matrix}
0 & 1\\
1 & 0
\end{matrix}}\right]\quad,\quad\sigma_2=\left[{\begin{matrix}
0 & -i\\
i & 0
\end{matrix}}\right]\quad,\quad\sigma_3=\left[{\begin{matrix}
1 & 0\\
0 & -1
\end{matrix}}\right]\,.
\]
\item The projection $\frac{1}{n}\sum_{i,j=1}^nE_{ij}$ in $M_n(\bbc)$ will be denoted by $J_n\,$.
\item Given an inclusion of finite von Neumann algebras $\mathcal{N}\subset\mathcal{M}$ with a fixed trace $tr$ on $\mathcal{M}$, the $tr$-preserving conditional expectation is denoted by $E^{\mathcal{M}}_{\mathcal{N}}$.
\item For unital inclusion of finite von Neumann algebras $\mathcal{N\subset M}$ with a fixed trace $tr$ on $\mathcal{M}$, $L^2(\mathcal{M})$ (resp. $L^2(\mathcal{N})$) denotes the GNS Hilbert space corresponding to $tr$ (resp., $tr|_{\mathcal{N}}$). The Jones projection, denoted by $e_{\mathcal{N}}$, is the projection onto the closed subspace $L^2(\mathcal{N})$ of $L^2(\mathcal{M})$.
\item We often denote a quadruple of von Neumann algebras
\[
\begin{matrix}
\mathcal{P} & \subset & \mathcal{M}\\
\cup & & \cup\\
\mathcal{N} & \subset & \mathcal{Q}
\end{matrix}
\]
by $(\mathcal{N}\subset\mathcal{P,Q}\subset\mathcal{M})$ for brevity. The Sano-Watatani angle is denoted by $\text{Ang}_M(P,Q)$ and the interior angle (resp. exterior angle) is denoted by $\alpha^N_M(P,Q)$ (resp., $\beta^N_M(P,Q)$).
\item All the logarithms appearing in this paper are with base $e$.
\item Throughout the paper, the notations $\lambda$ and $H$ (resp., $h$) will be reserved for the Pimsner-Popa probabilistic constant and Connes-St{\o}rmer relative entropy (resp., its variant) respectively. 
\end{enumerate}

\subsection{Motivation} 

\begin{figure} 
\centering
\begin{tikzpicture}[font=\large,thick]

\node[draw,
    rounded rectangle,
    align=center,
    minimum width=2.5cm,
    minimum height=1cm] (block1) { Key invariants\\ for two subfactors };

\node[draw,
    align=center,
    left=of block1,
    minimum width=3.5cm,
    minimum height=1cm] (block2) {~Pimsner-Popa constant~};

\node[draw,
    align=center,
    below=of block2,
    minimum width=3.5cm,
    minimum height=1cm] (block3) {~Sano-Watatani angle~};
  
\node[draw,
    align=center,
    right=of block1,
    minimum width=3.5cm,
    minimum height=1cm] (block4) {~Interior and exterior angle~};

\node[draw,
    rectangle,
    below=of block4,
    minimum width=3.5cm,
    minimum height=1cm] (block5) {~Connes-St{\o}rmer entropy~};

\draw[-latex] (block1) edge (block2)
    (block1) edge (block3)
    (block1) edge (block4)
    (block1) edge (block5);
    
    
\end{tikzpicture}
\bigskip

\caption{Few key invariants for two subfactors}\label{key invariants}
\end{figure}
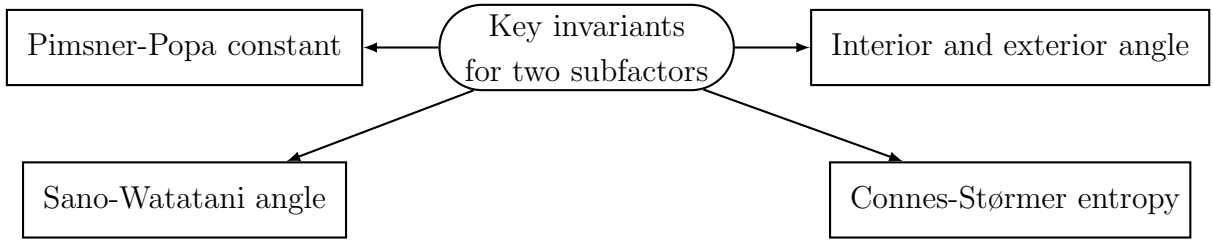

A subfactor $N$ is a unital subalgebra of a type $II_1$-factor $M$, which is itself a  $II_1$ factor, and $1_N=1_M$. An important invariant of a subfactor is the Jones index $[M:N]$ \cite{Jo} which generalizes the subgroup index. The Jones index measures how much bigger $M$ is compared to $N$. In view of the fact that the subfactor theory deals with the relative position of a subfactor $N$ inside an ambient factor $M$, it is a fundamental question to consider the relative positions of multiple subfactors. In the simplest case, one considers an intermediate subfactor $N\subset P\subset M$, and this is relatively well understood. As the next level of difficulty, if one takes a pair of subfactors $P,Q\subset M$, the theory becomes complicated. The study of multiple subfactors was initiated by Ocneanu. He proposed the concept of `maximal atlas' for a compatible family of `finite-index' bimodules arising from the subfactors (see \cite{O}). Ocneanu's work leads Jones to systematically study two subfactors \cite{Jo2}. The first difficulty one encounters is that $P\cap Q$ need not be a factor, and therefore we can not talk about the Jones index of $P\cap Q$ in $M$. However, this can be remedied using the `probabilistic index' $\lambda(M,P\cap Q)^{-1}$ due to Pimsner and Popa \cite{PP}, which is a substitute of the Jones index for a non-factorial inclusion and coincides with the Jones index for subfactor, that is, if $N\subset M$ is a subfactor of type $II_1$ factors, then $\lambda(M,N)={[M:N]}^{-1}$.  We discuss a few key invariants for the theory of `two subfactors' as described in \Cref{key invariants}. Note that we are not claiming here that this list of invariants is exhaustive. Discovering the `complete invariant' for `good class' of two subfactors seems open. Indeed, the formulation of `planar algebra' for a pair of subfactors -- even for commuting square-- remains elusive (see \cite{Jo3}).

We begin with the \textit{first invariant}, which is a generalization of the Pimsner-Popa probabilistic index $\lambda(M,N)^{-1}$ for $N\subset M$. Given a pair of subfactors $P,Q\subset M$ one can associate a number $\lambda(P,Q)\in[0,1]$, called the Pimsner-Popa probabilistic constant, which is an invariant to examine the relative position between the subfactors. Note that the definition makes sense for subalgebras of a finite von Neumann algebra as well. Recently in \cite{B}, given a pair of intermediate subfactors $P$ and $Q$ of a subfactor $N\subset M$ with $[M:N]<\infty$ and $N^{\prime}\cap M=\mathbb{C}$, the first author has given a formula of $\lambda(P,Q)$ in terms of the so-called `biprojections'. In general, given two subfactors $P,Q\subset M$, the formula for $\lambda(P,Q)$ is not known. As a next level of difficulty in the two subfactor theory, one sees that for a pair of subfactors $P,Q\subset M$ with $[M:P], [M:Q]<\infty$ it may very well happen that  $\lambda(M,P\cap Q)^{-1}=\infty$. Indeed, even  if we assume that both $P$ and $Q$ are of index  $2$, it can happen that $P\cap Q$ is of infinite (Pimsner-Popa) index in $M$ (see \cite{Jo2}, for instance). In \cite{JX}, Jones and Xu ask the following question:

\textbf{Question (Jones and Xu).} Consider a pair of finite index subfactors $P,Q\subset M$. Under what condition is the intersection $P\cap Q$ of a finite index?
\smallskip

\noindent They have shown that the finiteness of $\lambda(M,P\cap Q)^{-1}$ is equivalent to the finiteness of the spectrum of the so-called `angle operator', which leads us to our \textit{ second invariant} for two subfactors.

Sano and Watatani introduced a notion of angle \cite{SW}, denoted by $\text{Ang}_\mathcal{M}(\mathcal{P,Q})$, between two subalgebras of a finite von Neumann algebra $\mathcal{M}$, motivated by the angle between two subspaces (projections) in a Hilbert space. The angle between the subfactors $P\subset M$ and $Q\subset M$ determines the degree of non-commutativity of the two subfactors. Finding the possible value of the angle is known to be an important question in the  theory of subfactors. This has been highlighted by Grossman and Jones \cite{GJ}. The computation of the angle, even in simple cases, is not easy. See \cite{SW,GJ,GI} for many interesting applications and explicit computations of the angle operator. In another direction, to understand the relative position between a pair of intermediate subfactors $P$ and $Q$ of a finite index subfactor $N\subset M$ (i.e., $N\subset P,Q\subset M$), the notion of interior angle $\alpha^N_M(P,Q)$ and exterior angle $\beta^N_M(P,Q)$ between $P$ and $Q$ has been introduced in \cite{BDLR}, which is our \textit{third invariant}. This angle has been used to answer an open question by Longo by improving  the existing upper bound for the cardinality of the lattice of intermediate subfactors. Furthermore, a surprising connection between the intermediate subfactor theory and kissing numbers/sphere packing in geometry has been discovered. 

Now we discuss the \textit{fourth and final invariant} for two subfactors, namely the relative entropy of Connes-St{\o}rmer. In recent years, non-commutative entropy has become an active topic of research. Motivated by Shannon's information theory, Kolmogorov introduced (classical) conditional entropy, and subsequently, Sinai improved the definition. Conditional entropy is an important notion in classical ergodic theory. In the noncommutative situation, a notion of conditional (relative) entropy between two finite-dimensional von Neumann subalgebras of a finite von Neumann algebra was introduced by Connes and St{\o}rmer in \cite{CS}. They proved a Kolmogorov-Sinai type theorem using relative entropy as the main technical tool, and as an application, they proved that for $n\neq m$ the $n$-shift of the hyperfinite $II_1$ factor is not conjugate to the $m$-shift using an appropriate definition of  the entropy of an automorphism in the noncommutative framework of operator algebras. Later, Pimsner and Popa observed that relative entropy can be defined for arbitrary von-Neumann subalgebras of type $II_1$ factors (more generally, of tracial von Neumann algebras). More precisely, in analogy with the classical case, one can define $H(P|Q)$ (resp. $H(Q|P)$) which we refer to as the Connes-St{\o}rmer relative entropy. We often use a slight modification of $H(P|Q)$, denoted by $h(P|Q)$, due to Choda (see \cite{choda2}) which is more calculable and provides an alternative invariant for two subfactors. Pimsner and Popa have discovered a surprising connection between relative entropy and the Jones index \cite{PP}. Subsequently, many interesting connections between the Connes-St{\o}rmer entropy and Jones' subfactor theory were found (see \cite{PP2, Po2, Po3, Bi, BH, HI}, just to name a few). The relationship between the minimal index of subfactors (not necessarily of type $II_1$ factors) and entropy was clarified in \cite{Hiai2, Hiai3}. Furthermore, entropy for `canonical shifts' had been discussed in \cite{C, CH, Hiai}. However, note that explicit computation of $H(P|Q)$ in general situations is often difficult, and one can safely say that noncommutative conditional entropy is not yet well-understood. To quote a $1992$ survey article by St{\o}rmer, ``While entropy has for a third of a century been a central concept in ergodic theory, its non-abelian counterpart is still in its adolescent stage with only a few signs of mature strength". The major obstruction here is that in the non-abelian world, it is no longer true that two finite subsystems generate a finite subsystem. Another obstacle in the computation of $H$ is that entropy is not well-behaved with respect to the tensor products (see \cite{HS} for some interesting results along these lines). In Section $6$ \cite{PP}, given a unital inclusion of finite-dimensional $C^*$-algebras, an explicit formula for $H$ has been provided. However, if we consider a {\it pair} of finite-dimensional subalgebras of a finite-dimensional $C^*$-algebra the situation becomes complicated. Indeed, an explicit formula for the Connes-St{\o}rmer relative entropy between a pair of Masas of a matrix algebra, that is $H(\Delta_n|u\Delta_n u^*)$ for $u\in\mathcal{U}(M_n)$, seems to be open even when $n=2$ (see \cite{PSW}). 
 
Given a pair of subfactors $P,Q\subset M$, computation of these invariants is hard. A relatively easier situation is the existence of a subfactor $N\subset M$ such that we have a quadruple of $II_1$ factors
\begin{center}
$\begin{matrix}
P & \subset & M\\
\cup & & \cup\\
N & \subset & Q\\
\end{matrix}$
\end{center}
If the above quadruple is a commuting square, we know that the angle $\mathrm{Ang}_M(P,Q)$ is the singleton set $\{\pi/2\}$ and vice versa. In this situation, as a generalization of the notion of the non-degenerate commuting square, we have introduced the notion of the non-degenerate commuting cube (see \Cref{Sec 2}), a natural and relevant concept to deal with two subfactors. The basic construction of the non-degenerate commuting cube has been discussed. We have investigated a few of its properties. One of the applications of this notion is that if a quadruple of $II_1$ factors is obtained as an iterated basic construction of a non-degenerate commuting cube (see Sections \ref{Sec 4} and \ref{Sec 5} for examples of such), it becomes easy to check whether the quadruple we are dealing with is far from being a commuting square. Using this notion, the angle between the subfactors satisfying certain hypotheses has been computed. In this paper, this notion has been crucially used on several occasions in Sections \ref{Sec 4} and \ref{Sec 5}. 
\smallskip

In the absence of a general theory for two subfactors, and observing the level of difficulty for its invariants, it is natural that one should start investigating some particular important class of subfactors to hope for building a general theory for two subfactors. Motivated by this goal, in this article, our focus is on a pair of `spin model subfactors'. Spin model subfactors were introduced by Jones (\cite{Jo3}, see also \cite{JS}) as an important class of subfactors of the hyperfinite type $II_1$ factor $R$. Given a complex Hadamard matrix $u$ in $M_n$, one obtains a `spin model commuting square' and iterating the basic construction of such a square, one gets a hyperfinite subfactor $R_u\subset R$ with $[R:R_u]=n$ and $R_u^{\,\prime}\cap R=\bbc$. The importance of spin model subfactor has been emphasized by Jones (see \cite{Jo3}).  Not much is known about this subfactor. Indeed, computation of the relative commutants/standard invariant is very hard, and to compute even the first few relative commutants, computer assistance is required. The study of two spin model subfactors is even more intricate. If we start with a pair of distinct complex Hadamard matrices $u\mbox{ and }v$ of order $n\times n$, we (possibly) obtain a `pair' of spin model subfactors $R_u, R_v\subset R$. However, a number of caveats exist. The first one is that it is not guaranteed whether $R_u\neq R_v$ even if $u\neq v$. The second one is the factoriality of $R_u\cap R_v$, and finally, even if $R_u\cap R_v$ is a factor, it can be of infinite index in $R$. The starting point of this paper is the quest for the following problems:

\begin{enumerate}[]
\item \textbf{Problem 1}. Given a pair of (distinct) $n\times n$ complex Hadamard matrices $u$ and $v$, characterize when the spin model subfactors $R_u$ and $R_v$ of the hyperfinite type $II_1$ factor $R$ are distinct.
\item \textbf{Problem 2}. Suppose that we have a pair of spin model subfactors $R_u\subset R\mbox{ and }R_v\subset R$. Is $R_u\cap R_v$ a factor? Under what condition is it a finite index?
\item \textbf{Problem 3.} Find the value of the following invariants for the spin model subfactors $R_u\subset R\mbox{ and }R_v\subset R$.

\begin{enumerate}[$(i)$]
\item Pimsner-Popa probabilistic constant $\lambda(R_u,R_v).$  
\item Sano-Watatani angle $\text{Ang}_R(R_u,R_v)$.
\item The interior (resp., exterior) angle ${\alpha}^{R_u\cap R_v}_R(R_u,R_v)$ \big(resp., ${\beta}^{R_u\cap R_v}_R(R_u,R_v)$\big).
\item Connes-St{\o}rmer relative entropy $H(R_u|R_v)$ (and its slight modification $h(R_u|R_v)$).
\end{enumerate}
\end{enumerate}

We have a complete answer to Problem $1$. We classify pairs $(u,v)$ of complex Hadamard matrices of order $n\times n$ such that they give rise to distinct spin model subfactors $(R_u\neq R_v)$. The classification is obtained in terms of an equivalence relation that is finer than the Hadamard equivalence. The Pimsner-Popa probabilistic constant and the Hamming numbers of the rows of the unitary matrix $u^*v$ play the lead roles. As a consequence of this, it turns out that inequivalent complex Hadamard matrices always give rise to pairs of (distinct) spin model subfactors. We have considered the cases of $n=2\mbox{ and }n=4$ for detailed investigation to incorporate both the instances of Hadamard equivalent and inequivalent matrices. Note that up to Hadamard equivalence, there is a single family of complex Hadamard matrices of order $2\times 2$, whereas, there exists a continuous, one-parameter family of inequivalent complex Hadamard matrices of order $4\times 4$. We take the liberty to remind our readers that the complete classification of complex Hadamard matrices is extremely difficult and unsolved in the literature. For example, even the case of $6\times 6$ is conjectured only. Problem $2$ and Problem $3$ have been solved completely in the cases of $n=2\mbox{ and }n=4$. We observe that while the factoriality and finite index of $R_u\cap R_v$ both hold in the case of $n=2$, they fail for certain situations in the $n=4$ case, and this case behaves rather erratically, which convinces us that Problem $2$ and Problem $3$ (fully) for the general $n\times n$ case seem beyond reach at the moment. 
\smallskip

We now discuss our findings in this paper. First, we briefly discuss the case of $n=2$. We have overcome the first hurdle by proving that $R_u\cap R_v$ is a (non-irreducible) factor. More precisely, we have proved that $R_u\cap R_v\subset R$ is a diagonal subfactor with Jones index $4$ (see the above question due to Jones and Xu). It is indeed a curious fact that the intersection $R_u\cap R_v$ of the spin model subfactors $R_u\mbox{ and }R_v$ in this case turns out to be a vertex model subfactor of index $4$. We have also computed the value of $\lambda(R_u,R_v)$, the Sano-Watatani angle $\mathrm{Ang}_R(R_u,R_v)$, the interior angle $\alpha^{R_u\cap R_v}_R(R_u,R_v)$ and the exterior angle $\beta^{R_u\cap R_v}_R(R_u,R_v)$, and the modified Connes-St{\o}rmer relative entropy $h(R_u|R_v)$, which provides a legitimate bound for the Connes-St{\o}rmer relative entropy $H(R_u|R_v)$. It turns out that $h(R_u|R_v)$ is not constant; rather, it takes value depending on $u\mbox{ and }v$. We have characterized the quadruple of $II_1$ factors $(R_u\cap R_v\subset R_u,R_v\subset R)$ in terms of a single bi-unitary matrix in $M_4$ and completely determined when it is a commuting square. Thus, the relative position between $R_u\mbox{ and }R_v$ has been mostly understood in the case of $n=2$.

Now, we briefly discuss the case of $n=4$. Here, the complex Hadamard matrices $u\mbox{ and }v$ are Hadamard inequivalent. This case is more delicate and very interesting. It turns out that factoriality and the finite index of $R_u\cap R_v$ depend on whether certain rotations naturally arising from the circle parameters associated with $u\mbox{ and }v$ are rational or irrational. In the case of rational rotation, we construct en route an infinite family of potentially new subfactors of $R$. All these subfactors are irreducible with the Jones index $4n,n\in\bbn$ and $n\geq 2$. This has enabled us to conclude factoriality for the intersection $R_u\cap R_v$ for the spin model subfactors $R_u\mbox{ and }R_v$ in the case of rational rotation, and we observe that the intersection is not of the fixed index in $R$, unlike the situation of $2\times 2$. Moreover, $R_u\cap R_v$ is irreducible, and $[R:R_u\cap R_v]\in\{4n:n\in\bbn,\,n\geq 2\}$. Now, if the rotation is irrational, we observe that $R_u\cap R_v$ is of infinite Pimsner-Popa index in $R$. Unfortunately, we could not determine whether this is a factor in this case. Therefore, we see that the case of $4\times 4$ (inequivalent Hadamard matrices) is in sharp contrast to the case of $2\times 2$ (equivalent Hadamard matrices). To understand how the two subfactors `interact' we have computed explicitly the value of $\lambda(R_u,R_v)$, and in the case of rational rotation, the Sano-Watatani angle $\mathrm{Ang}_R(R_u,R_v)$, the interior angle $\alpha^{R_u\cap R_v}_R(R_u,R_v)$ and the exterior angle $\beta^{R_u\cap R_v}_R(R_u,R_v)$, and the exact value of the Connes-St{\o}rmer relative entropy $H(R_u|R_v)$. The rigidity of the interior angle has also been established, namely, that it is always greater than $\pi/3$. In the case of irrational rotation, a legitimate bound for the Connes-St{\o}rmer relative entropy $H(R_u|R_v)$ has been obtained. Finally, as an application of the commuting cube, we have completely characterized when the quadruple of $II_1$ factors $(R_u\cap R_v\subset R_u,R_v\subset R)$ for given $u\mbox{ and }v$ forms a commuting square. The major highlight of the $n=4$ case (rational rotation case) is that starting with two (Hadamard inequivalent complex Hadamard) $4\times 4$ matrices, one produces an infinite family of integer indices (in particular, $4n,\,n\in\bbn$ and $n\geq 2$) and irreducible subfactors of the hyperfinite type $II_1$ factor $R$. However, in the irrational rotation case, many questions remain unanswered.

We draw a table in \Cref{contrast} to highlight the contrast between the $n=2\mbox{ and }n=4$ cases. We only mention the rational rotation case in the $n=4$ situation, as the irrational rotation case remains mysterious except for the fact that $R_u\cap R_v$ is of infinite index.
\begin{table}[ht]\label{contrast}
\centering
\begin{tabular}{|c| c| c| }
\hline\hline
$R_u,R_v\subset R$ & $\,\,u,v\mbox{ are }2\times 2\,\,$ & $\,\,\,\,u,v\mbox{ are }4\times 4\,,$ and rational rotation \\ [0.5ex]
\hline\hline
$R_u\neq R_v$ & only when $u\nsim v$ & $\mbox{always}$  \\
\hline
$R_u\cap R_v$  & always a factor & always a factor \\
\hline
$[R:R_u\cap R_v]$ & always $4$ & $\{4n:n\geq 2\,\,\mbox{integer}\}$ \\
\hline
$R_u\cap R_v\subset R$ is vertex model & yes & no \\
\hline
relative commutant & non-irreducible & irreducible \\
\hline
characterization of $R_u\cap R_v$  & diagonal subfactor & unknown! \\
\hline
$\texttt{\#}\mathrm{Ang}_R(R_u,R_v)$ (set of angles) & singleton set & $\lfloor\frac{n}{2}\rfloor$, if $[R:R_u\cap R_v]=4n$ \\
\hline
interior angle $\alpha_R(R_u,R_v)$ & no rigidity & rigidity present, $\alpha>\frac{\pi}{3}$ always \\
\hline
exterior angle $\beta_R(R_u,R_v)$ & different values & fixed value, $\beta=\arccos(\frac{1}{3})$ \\
\hline
relative entropy $h$ & different values & fixed value $\log 2$ \\
\hline\hline
\end{tabular}
\vspace{1mm}
\caption{Contrast between the cases of $n=2\mbox{ and }n=4$}\label{contrast}
\end{table}

\noindent Now it is not difficult to persuade why the general $n\times n$ case seems beyond our reach at the moment, and extensive investigation is needed. A major difficulty is the lack of complete classification of complex Hadamard matrices of order beyond $5\times 5$.

\subsection{Statement of the main results}

To introduce our main results, we briefly recall some well-known facts about complex Hadamard matrices that appear in many areas like quantum teleportation, coding theory, mutually unbiased bases, unitary error bases, operator algebras, Harmonic analysis, etc. In this paper, we provide yet another application of complex Hadamard matrices in von Neumann algebras.
\smallskip

A matrix with entries $\pm 1$ and mutually orthogonal rows and columns is called a Hadamard matrix. If $H$ is a Hadamard matrix and $H^{\intercal}$ denote its transpose, then $H H^{\intercal}=n I_n$. A complex Hadamard matrix is a generalization of Hadamard matrix.
\begin{dfn}
A complex Hadamard matrix is a $n\times n$ matrix with complex entries of same modulus and $HH^*=nI_n$.
\end{dfn}
Note that $\frac{1}{\sqrt{n}}H$ is a unitary matrix. In this paper, by a complex Hadamard matrix we shall always mean the associated unitary matrix $\frac{1}{\sqrt{n}}H$. Observe that for every $n\geq 1$, the Fourier matrix $(F_n)_{i,j}:=\frac{1}{\sqrt{n}}\exp(2\pi \mathbf{i} (i-1)(j-1)/n)$ is a complex Hadamard matrix. Thus,
\[
 F_1=[1]\,\,,\quad F_2=\frac{1}{\sqrt{2}}\begin{bmatrix}
1 & 1\\
1 & -1
\end{bmatrix}\quad \mbox{and}\quad F_3=\frac{1}{\sqrt{3}}\begin{bmatrix}
1 & 1 &1\\
1 & \omega & {\omega}^2\\
1 & {\omega}^2 &\omega
\end{bmatrix}.
\]
Two complex Hadamard matrices are called Hadamard equivalent, denoted by  ${\displaystyle H_{1}\simeq H_{2}}$, if there exist diagonal unitary matrices ${\displaystyle D_{1},D_{2}}$ and permutation matrices ${\displaystyle P_{1},P_{2}}$ such that
\[
{\displaystyle H_{1}=D_{1}P_{1}H_{2}P_{2}D_{2}\,.}
\]
It is known that for $n=2,3,5$, all complex Hadamard matrices are equivalent to the Fourier matrix $F_{n}\,$. We refer the reader to \cite{Ha} for proof. In particular, any $2\times 2$ complex Hadamard matrix is of the form  $\frac{1}{\sqrt{2}}\begin{bmatrix}
e^{i\alpha_1} & e^{i (\alpha_1+\alpha_3)}\\
e^{i\alpha_2} & -e^{i(\alpha_2+\alpha_3)}
\end{bmatrix}$ for $\alpha_j\in[0,2\pi)$. However, it is known that there exists a continuous, one-parameter family of Hadamard inequivalent $4\times 4$ complex Hadamard matrices, and any element of this family is of the following form
\begin{IEEEeqnarray}{lCl}\label{semicircle}
u(z)= \frac{1}{2}\,\left[{\begin{matrix}
1 & 1 & 1 & 1\\
1 & iz & -1 & -iz\\
1 & -1 & 1 & -1\\
1 & -iz & -1 & iz\\
\end{matrix}}\right]
\end{IEEEeqnarray}
where $\,z=e^{i\alpha}\in\mathbb{S}^1$ with $\alpha\in[0,\pi)$. It is known that the classification of complex Hadamard matrix in higher dimensions is extremely hard.

Our first major theorem in this regard is the following. Recall that given a pair of complex Hadamard matrices $u$ and $v$, one obtains a pair of spin model subfactors $R_u,R_v\subset R$ (see \Cref{pairspin} for details). We have introduced an equivalence relation `$\sim$' finer than the Hadamard equivalence `$\simeq$' which completely characterizes when spin model subfactors $R_u,R_v$ arising from two $n\times n$ distinct complex Hadamard matrices $u,v$ are distinct $(R_u\neq R_v)$.

\begin{thm}[\Cref{general thm}]\label{general thm1}
\begin{enumerate}[$(i)$]
\item For distinct $n\times n$ complex Hadamard matrices $u\mbox{ and }v$, the pair of spin model subfactors $R_u\mbox{ and }R_v$ of the hyperfinite type $II_1$ factor $R$ are distinct $(R_u\neq R_v)$ if and only if $\,u\nsim v$.
\item If two $n\times n$ complex Hadamard matrices $u\mbox{ and }v$ are Hadamard inequivalent, then the corresponding spin model subfactors $R_u\mbox{ and }R_v$ of $R$ are always distinct $(R_u\neq R_v)$.
\end{enumerate}
\end{thm}

\subsubsection{The case of $2\times 2$ complex Hadamard matrices}

If $u$ and $v$ are two distinct $2\times 2$ complex Hadamard matrices, by \Cref{general thm1} we have $R_u\neq R_v$ if and only if $v\neq uD$ (that is, $u\nsim v$), where $D$ is a unitary matrix of the form $\left[{\begin{matrix}
z & 0\\
0 & w\\
\end{matrix}}\right]$ or $\left[{\begin{matrix}
0 & z\\
w & 0\\
\end{matrix}}\right]\mbox{ for }z,w\in\mathbb{S}^1$. It turns out that without loss of generality, we can assume $u,v$ are of the following form
\begin{IEEEeqnarray}{lCl}\label{for intro}
u=\frac{1}{\sqrt{2}}\begin{bmatrix}
1 & 1\\
e^{i\alpha} & -e^{i\alpha}\\
\end{bmatrix}\qquad\mbox{and}\qquad v=\frac{1}{\sqrt{2}}\begin{bmatrix}
1 & 1\\
e^{i\beta} & -e^{i\beta}\\
\end{bmatrix}\,,
\end{IEEEeqnarray}
where $\alpha,\beta\in[0,2\pi)$ and $\alpha\neq\beta$. This situation gives rise to pair of (distinct) spin model subfactors $R_u,R_v\subset R$, and we have {\it a priori} a quadruple of von Neumann algebras
\[
\begin{matrix}
R_v & \subset & R\\
\cup & & \cup\\
R_u\cap R_v &\subset & R_u\\
\end{matrix}
\]

\begin{thm}[\Cref{PP constant final}]
If $u$ and $v$ are as above, then the Pimsner-Popa constant $\lambda(R_u,R_v)$ is $\frac{1}{2}$.
\end{thm}

\begin{thm}[\Cref{factoriality}]
The von Neumann subalgebra $R_u\cap R_v$ of $R$ is a $II_1$ subfactor with $[R:R_u\cap R_v]=4$ and $[R_u:R_u\cap R_v]=[R_v:R_u\cap R_v]=2$.
\end{thm}

\begin{thm}[\Cref{kl}]
The interior and the exterior angle for the quadruple of $II_1$ factors $(R_u\cap R_v\subset R_u,R_v\subset R)$ are given by the following,
\[
\cos\big(\alpha^{R_u\cap R_v}_R(R_u,R_v)\big)=\cos\big(\beta^{R_u\cap R_v}_R(R_u,R_v)\big)=\cos^2(\alpha-\beta)\,.
\]
The quadruple $(R_u\cap R_v\subset R_u,R_v\subset R)$ is a commuting square (and consequently co-commuting square) if and only if $\,\alpha-\beta=\pm\frac{\pi}{2}$. Here, $\alpha\mbox{ and }\beta$ are associated with the complex Hadamard matrices $u\mbox{ and }v$ as described in \Cref{for intro}.
\end{thm}

In this regard, we obtain a concrete example of a quadruple of $II_1$ factors (namely, $(R_u\cap R_v\subset R_u,R_v\subset R)$) which is non-degenerate (that is, $\overline{R_uR_v}=\overline{R_vR_u}=R$) but not co-commuting square.

\begin{thm}[\Cref{relativecommutant}]
The relative commutant ${(R_u\cap R_v)}^\prime\cap R$ is $\bbc\oplus\bbc$.
\end{thm}

The following results completely characterize the quadruple of $II_1$ factors $(R_u\cap R_v\subset R_u,R_v\subset R)$.
\begin{thm}[\Cref{vertex}, \Cref{fullchar}, \Cref{diag}]
The subfactor $R_u\cap R_v\subset R$ is a vertex model subfactor of Jones index $4$. More explicitly, we have the following.
\begin{enumerate}[$(i)$]
\item The pair of subfactors $R_u$ and $R_v$ of the hyperfinite type $II_1$ factor $R$ are conjugate to each other via a unitary in ${(R_u\cap R_v)}^\prime\cap R$.
\item The bi-unitary matrix $\,u_2W_2V_2$ in $M_4$ generates (in the sense of basic construction) the composition of subfactors $R_u\cap R_v\subset R_u\subset R$, and consequently the quadruple $(R_u\cap R_v\subset R_u,R_v\subset R)$.
\item The subfactor $R_u\cap R_v\subset R$ is isomorphic to the diagonal subfactor 
\[
\left\{\begin{bmatrix}
x & 0\\
0 & \theta(x)
\end{bmatrix}:x\in R~\mathrm{and}~\theta\in\mathrm{Out}(R)~\mathrm{such~ that}~{\theta}^2=\mathrm{id}\right\}\subset M_2(\mathbb{C})\otimes R\,.
\]
Under this isomorphism, $R_u$ becomes $\big(\mathrm{diag}\{1,e^{i(\alpha-\beta)}\}\big)R_v\big(\mathrm{diag}\{1,e^{i(\beta-\alpha)}\}\big)$, where $R_v$ is given by the following
\[
\Bigg\{\begin{bmatrix}
x & y\\
\theta(y) & \theta(x)
\end{bmatrix} : x,y \in R\Bigg\}\,.
\]
\end{enumerate}
\end{thm}

The relative position of $R_u\mbox{ and }R_v$ are described by the following results.

\begin{thm}[\Cref{major 2by2}]
We have the following.
\begin{enumerate}[$(i)$]
\item $H(R|R_u\cap R_v)=2\log 2\mbox{ and }H(R_u|R_u\cap R_v)=H(R_v|R_u\cap R_v)=\log 2$.
\item $h(R_u|R_v)=\eta\big(\cos^2\big(\frac{\alpha-\beta}{2}\big)\big)+\eta\big(\sin^2\big(\frac{\alpha-\beta}{2}\big)\big)\,.$
\item $\eta\big(\cos^2\big(\frac{\alpha-\beta}{2}\big)\big)+\eta\big(\sin^2\big(\frac{\alpha-\beta}{2}\big)\big)\leq H(R_u|R_v)\leq\log 2\,,$
and if $\,\alpha-\beta=\pm\frac{\pi}{2}$, then $H(R_u|R_v)=h(R_u|R_v)=\log 2$. Here, $\alpha\mbox{ and }\beta$ are associated with the complex Hadamard matrices $u\mbox{ and }v$ as described in \Cref{for intro}.
\end{enumerate}
\end{thm}

\begin{thm}[\Cref{sanoangle}]
The Sano-Watatani angle between the subfactors $R_u\mbox{ and }R_v$ of $R$ is the singleton set $\{\arccos|\cos(\alpha-\beta)|\}$.
\end{thm}

\subsubsection{The case of $4\times 4$ complex Hadamard matrices}

The case of Hadamard inequivalent $4\times 4$ (distinct) complex Hadamard matrices is more delicate and very interesting. Recall that any such complex Hadamard matrix is of the form described in \Cref{semicircle}. By \Cref{general thm1}, any two distinct elements $u\mbox{ and }v$ of this family give rise to a pair of spin model subfactors. That is, we have a quadruple of von Neumann algebras
\[
\begin{matrix}
R_v & \subset & R\\
\cup & & \cup\\
R_u\cap R_v &\subset & R_u\\
\end{matrix}
\]

\begin{thm}[\Cref{Popa for 4 by 4}]
The Pimsner-Popa constant for the pair of subfactors $R_u\mbox{ and }R_v$ of the hyperfinite $II_1$ factor $R$ is given by $\lambda(R_u,R_v)=\frac{1}{2}$. 
\end{thm}

It turns out that in this case, the (Pimsner-Popa) index of the intersection $R_u\cap R_v$ in $R$ can be both finite and infinite, depending on which $u\mbox{ and }v$ we are starting with. Moreover, in the case of finite index, it can be any number in the set $\{4n:n\in\bbn,\,n\geq 2\}$. Conversely, given any $n\geq 2$, there exist $u\mbox{ and }v$ such that index of $R_u\cap R_v$ in $R$ is equal to $4n$. Below we describe the scenario briefly.
\smallskip

Consider the following subgroup of $\mathbb{S}^1$
\begin{IEEEeqnarray*}{lCl}
\Gamma:=\{\omega\in\mathbb{S}^1:\omega^m=1\mbox{ for some }m\in 2\bbn+2\}\cup\{1\}\,.
\end{IEEEeqnarray*}
That is, $\Gamma$ consists of all {\it even} roots of unity. Let $u:=u(a)\mbox{ and }v:=u(b)$ be two $4\times 4$ Hadamard inequivalent complex Hadamard matrices parametrized by the circle parameters $a\mbox{ and }b$ respectively, and we have a pair of subfactors $R_u\mbox{ and }R_v$ of the hyperfinite $II_1$ factor $R$. We define a relation $u(a)\sim u(b)$ if and only if $\theta(a,b):=\overline{a}b\in\Gamma$, that is, $\theta(a,b)^{\,m}=1\mbox{ for some }m\in 2\bbn+2$. Note that we are considering pair $(a,b)$ such that $a\neq b$, and hence $\theta(a,b)\neq 1$. It is easy to check that `$\sim$' is an equivalence relation.
\smallskip

\noindent\textbf{Caution:} This equivalence relation has no relation with that in \Cref{general thm1}.
\smallskip

Our first major theorem in this direction is the following.
\begin{thm}[\Cref{infinite index}]
If $u\nsim v$, then $\lambda(R,R_u\cap R_v)=\lambda(R_u,R_u\cap R_v)=\lambda(R_v,R_u\cap R_v)=0$. That is, $R_u\cap R_v$ is a von Neumann subalgebra of the hyperfinite $II_1$ factor $R$ with infinite Pimsner-Popa index.
\end{thm}
If $u\sim v$, the outcome is completely opposite, and we obtain a series of potentially new integer index subfactors of $R$.

\begin{thm}[\Cref{producing}, \ref{interesting}, \ref{lop}, \ref{555}]
For a pair of Hadamard inequivalent $4\times 4$ complex Hadamard matrices $u=u(a)$ and $v=u(b)$ such that ${(\bar{a}b)}^m=1$, $R_u\cap R_v\subset R$ is an irreducible subfactor of $R$ with $[R:R_u\cap R_v]=2m$. Furthermore, $H(R_u|R_u\cap R_v)=H(R_v|R_u\cap R_v)=\log m-\log 2$ and also, $H(R|R_u\cap R_v)=\log 2+\log m$.
\end{thm}

The following results describe the relative position of $R_u\mbox{ and }R_v$. The first few of them describe the angle between $R_u\mbox{ and }R_v$, whereas the last one gives relative entropy between these subfactors.

\begin{thm}[\Cref{rigidity of angle}]
For the case of $u\sim v$, the interior angle $\alpha^{R_u\cap R_v}_R\big(R_u,R_v\big)$  is strictly  greater than $\pi/3$ and the exterior angle $\beta^{R_u\cap R_v}_R\big(R_u,R_v\big)=\arccos 1/3$. Furthermore, $\alpha^{R_u\cap R_v}_R\big(R_u,R_v\big)$ converges to $\pi/3$ decreasingly as $m$ tends to infinity.
\end{thm}

\begin{thm}[\Cref{final com result}]
The quadruple of von Neumann algebras
\[
\begin{matrix}
R_{u(b)} & \subset & R\\
\cup & & \cup\\
R_{u(a)}\cap R_{u(b)} &\subset & R_{u(a)}\\
\end{matrix}
\] 
is a commuting square if and only if $(\overline{a}b)^4=1$, that is $b=\pm ia$.
\end{thm}

\begin{thm}[\Cref{dihedralangle}]
Let $u=u(a)$ and $v=u(b)$ be distinct $4\times 4$ inequivalent complex Hadamard matrices parametrized by the circle parameters $a\mbox{ and }b$ such that $u\sim v$ with $(\overline{a}b)^m=1,\,m\in 2\bbn+2$. Then, we have the following.
\begin{enumerate}[$(i)$]
\item If $\,m=4$, then $\mathrm{Ang}_R(R_u,R_v)=\{\pi/2\};$
\item If $\,m\geq 6$, then $\mathrm{Ang}_{R}(R_u,R_v)=\left\{\frac{2k\pi}{m}\,:\,k=1,2,\ldots, \lfloor\frac{m-2}{4}\rfloor\right\}$.
\end{enumerate}
\end{thm}

\begin{thm}[\Cref{entropy for 4 by 4}]
For a pair of Hadamard inequivalent $4\times 4$ complex Hadamard matrices $u\mbox{ and }v$,
\begin{enumerate}[$(i)$]
\item if $u\sim v$, then $H(R_u|R_v)=h(R_u|R_v)=\log 2;$
\item if $u\nsim v$, then $0<\log 2+\frac{1}{8}\left(\eta|1+a\overline{b}|^2+\eta|1-a\overline{b}|^2\right)\leq H(R_u|R_v)\leq\log 2\,.$
\end{enumerate}
\end{thm}

The crucial ingredient that we have used to prove the above theorem is the existence of an intermediate subfactor $R_u\cap R_v\subset\mathcal{I}_{u,v}\subset R$ containing both $R_u\mbox{ and }R_v$. For the case of $u\nsim v$, we are not able to obtain the exact value of the relative entropy $H(R_u|R_v)$. Major difficulties are that $R_u\cap R_v$ is of infinite Pimsner-Popa index in $R$, and its factoriality is not known.

\subsubsection{Results in finite dimensions}

It turns out that both the Pimsner-Popa probabilistic constant and the Connes-St{\o}rmer relative entropy for a pair of subfactors behave well under `controlled' limits. This suggests that to achieve results in infinite-dimension, we should look at finite-dimensional situations. The same theme has been highlighted in Section $6$ \cite{PP} for the case of a single subfactor. In this paper, prior to the above results on infinite factors, we have en route a couple of results in finite dimensions that are instrumental in proving our results. These can also be treated as results of independent interest.

\begin{thm}[\Cref{finitemains}]
If $\Delta_n$ and $u\Delta_n u^*$ are two Masas in $M_n(\mathbb{C})$, $u$ being a unitary matrix in $M_n(\mathbb{C})$, then the Pimsner-Popa constant is given by the following:
\[
\lambda(\Delta_n\,,\,u\Delta_n u^*)=\min_{1\leq i\leq n}\,\big(\mbox{$h\left(u^*\right)_i$}\big)^{-1}
\]
where $\left(u^*\right)_i$ is the $i$-th column of $u^*$ and $h$ denotes the Hamming number of a vector in $\bbc^n$ (see \Cref{hamdef}).
\end{thm}

This result is instrumental in determining when spin model subfactors arising from distinct $n\times n$ complex Hadamard matrices are distinct (\Cref{general thm1}).
 
\begin{thm}[\Cref{2nd lambda}]
For the subalgebras $M_n(\bbc)$ and $uM_n(\bbc)u^*$ of $M_n(\bbc)\oplus M_n(\bbc)$, where $\,u\in\Delta_2\otimes M_n(\bbc)$ is a unitary matrix given by $\,u=\left[{\begin{smallmatrix}
u_1 & 0\\
0 & u_2\\
\end{smallmatrix}}\right]$ with unitary matrices $u_1,u_2\in M_n(\bbc)$, we have the following~:
\begin{IEEEeqnarray*}{lCl}
\lambda\big(M_n(\bbc)\,,\,uM_n(\bbc)u^*\big) = \begin{cases}
                            1 & \mbox{ iff } \,\,u_1^*u_2 \mbox{ is a scalar matrix}, \cr
                            \frac{1}{2} & \mbox{ if } \,\,u_1^*u_2 \mbox{ is not a diagonal matrix}. \cr
                           \end{cases}
\end{IEEEeqnarray*}
\end{thm}

This result is instrumental in computing the value of the Pimsner-Popa constant between the pairs of spin model subfactors (both $2\times 2$ and $4\times 4$ situations in Sections \ref{Sec 4} and \ref{Sec 5}).

\subsection{Outline of the paper}

The organization of the paper is as follows: In \Cref{Sec 1}, the one following this Introduction, we discuss the key invariants for ``two subfactors" introduced earlier. In \Cref{Sec 2}, we introduce the notion of ``commuting cube", which is instrumental for this paper, as well as a general study of  two subfactors. In \Cref{Sec 3}, we concentrate on pairs of spin model subfactors, a particular situation of two subfactors. \Cref{Sec 3.2} contains a key result concerning the Pimsner-Popa constant in a finite-dimensional situation. Then, in \Cref{Sec 4}, a pair of subfactors arising from $2\times 2$ complex Hadamard matrices have been investigated. This is followed by the case of $4\times 4$ inequivalent complex Hadamard matrices in \Cref{Sec 5}, which is one of the major highlights of this paper. Finally, in \Cref{Sec 6}, we discuss a few open problems and possible directions for future work. In the Appendix (\Cref{Sec 7}), we have provided detailed proofs for the construction of the tower of basic constructions for the interested readers.


\section{A few key invariants for two subfactors}\label{Sec 1}

In this section, we discuss a few key invariants for ``two subfactors", namely, the Sano-Watatani angle, the interior and exterior angle, the Pimsner-Popa probabilistic constant, and the Connes-St{\o}rmer relative entropy along with a modified version of it. We also recall the basic construction of a non-degenerate symmetric commuting square, and a few essential results to be used throughout the article.

\subsection{Angle between two subfactors and commuting square}\label{angletams}

In \cite{P}, Popa introduced a notion of orthogonality for a pair of von Neumann subalgebras $\mathcal{B}_1$ and $\mathcal{B}_2$ of a finite von Neumann algebra $\mathcal{M}$. Suppose that $\mathcal{M}$ is commutative, and let $(X,\mathcal{B},\mu)$ be a probability space such that $\mathcal{M}$ is isomorphic to $L^{\infty}(X,\mu)$. If $\mathcal{B}_1$ and $\mathcal{B}_2$ are von Neumann subalgebras of $\mathcal{M}$, then $\mathcal{B}_1$ is orthogonal to $\mathcal{B}_2$ if and only if the corresponding $\sigma$-subalgebras of $\mathcal{B}$ are independent. Thus, orthogonality is the non-commutative version of classical independence. As a marginal generalization of the notion of orthogonality, Popa also introduced the so-called `commuting square', which proves to be an indispensable tool in subfactor theory. Consider an inclusion of finite von Neumann algebras $\mathcal{N}\subset \mathcal{M}$ with a fixed trace $tr$ on $\mathcal{M}$ and intermediate von Neumann subalgebras $\mathcal{P}$ and $\mathcal{Q}$. Thus, we obtain a quadruple of von Neumann algebras
\[
\begin{matrix}
\mathcal{Q} &\subset & \mathcal{M} \cr \cup &\ &\cup\cr \mathcal{N} &\subset & \mathcal{P}\,.
\end{matrix}
\]
If $\mathcal{P}\vee \mathcal{Q}=\mathcal{M}$ and $\mathcal{P}\wedge \mathcal{Q}=\mathcal{N}$, then a quadruple is called a {\it quadrilateral}.

\begin{dfn}[\cite{Po2},\cite{GHJ},\cite{JS}]\label{def of comm}
A quadruple
\[
\begin{matrix}
\mathcal{Q} &\subset & \mathcal{M} \cr \cup &\ &\cup\cr \mathcal{N} &\subset & \mathcal{P}
\end{matrix}
\]
of finite von Neumann algebras is called a commuting square if $E^{\mathcal{M}}_{\mathcal{P}}E^{\mathcal{M}}_{\mathcal{Q}}=E^{\mathcal{M}}_{\mathcal{Q}}E^{\mathcal{M}}_{\mathcal{P}}=E^{\mathcal{M}}_{\mathcal{N}}$. The quadruple is said to be non-degenrate if ${\overline{\mathcal{PQ}}}^{\,\tiny \mbox{SOT}}={\overline{\mathcal{QP}}}^{\,\tiny \mbox{SOT}}=\mathcal{M}$. A quadruple is called a non-degenerate commuting square or symmetric commuting square if it is a commuting square and non-degenerate.
\end{dfn}

For brevity, we sometimes write $(\mathcal{N}\subset\mathcal{P,Q}\subset \mathcal{M})$ to mean the quadruple
\[
\begin{matrix}
\mathcal{Q} &\subset & \mathcal{M} \cr \cup &\ &\cup\cr \mathcal{N} &\subset & \mathcal{P}
\end{matrix}
\]
Note that for a quadruple $(\mathcal{N}\subset \mathcal{P,Q}\subset \mathcal{M})$, if $E^\mathcal{M}_\mathcal{Q}(\mathcal{P})\subseteq \mathcal{N}$ (or $E^\mathcal{M}_\mathcal{P}(\mathcal{Q})\subseteq \mathcal{N}$), then it is a commuting square. Indeed, observe that for any $x\in M$, using the hypothesis $E^M_Q(P)\subseteq N$, we have the following
\begin{IEEEeqnarray*}{lCl}
E_Q^ME^M_P(x)=E^Q_NE_Q^ME^M_P(x)=E^M_NE^M_P(x)=E^P_NE^M_P(x)=E^M_N(x)\,.
\end{IEEEeqnarray*}
Moreover, the converse is also true.

\begin{dfn}[\cite{Po2}, \cite{SW},\cite{Wa}]\label{def of cocomm}
A quadruple
\[
\begin{matrix}
\mathcal{Q} &\subset & \mathcal{M} \cr \cup &\ &\cup\cr \mathcal{N} &\subset & \mathcal{P}
\end{matrix}
\]
of finite von Neumann algebras with a fixed (finite, faithful, normal) trace $\tau^\prime$ on $\mathcal{N}^{\,\prime}$ is called a co-commuting square if the quadruple
\[
\begin{matrix}
\mathcal{P}^\prime &\subset & \mathcal{N}^\prime \cr \cup &\ &\cup\cr \mathcal{M}^\prime &\subset & \mathcal{Q}^\prime
\end{matrix}
\]
is a commuting square. 
\end{dfn}

Since $\mathcal{M}^\prime=\mathcal{P}^\prime\cap\mathcal{Q}^\prime$, it is necessary for $\mathcal{M}=\mathcal{P}\vee\mathcal{Q}$ to hold. We refer to \cite{EK, JS, Po2, BDLR} for a comprehensive treatment on commuting and co-commuting squares.

\begin{rmrk}\rm
In general, co-commuting and non-degenerate need not be same. If a quadruple $(N\subset P,Q\subset M)$ of $II_1$ factors is a commuting square with $[M:N]<\infty$, then it is a co-commuting square if and only if non-degenerate (Theorem $7.1$ and Corollary $7.1$ in \cite{SW}). If the commuting square condition is absent, but irreducibility and finite-index of $N\subset M$ is given, then also this holds (Theorem $3.21$ in \cite{GJ}). In the same spirit, if $(N\subset P,Q\subset M)$ is a commuting square of finite-dimensional $C^*$-algebras with $||\Lambda^M_Q||^2=[M:Q]=[P:N]=||\Lambda^P_N||^2$ (Watatani index), then the quadruple is non-degenerate (Lemma $3.10$ in \cite{BG}).
\end{rmrk}

A quadruple need not always be a commuting square as there might have nontrivial `angle' between the subalgebras. Indeed, as a generalization of commuting square, Sano and Watatani  introduced a notion of ``angle" between a pair of subalgebras of a given finite von Neumann algebra as the spectrum of certain angle operator \cite{SW}. This was motivated by relative position between two different (closed) subspaces $\mathcal{K}$ and $\mathcal{L}$ in a Hilbert space $\mathcal{H}$. Recall the angle operator $\Theta(p,q)$, where $p$ (resp. $q$) is the orthogonal projection onto $\mathcal{K}$ (resp. $\mathcal{L}$). The set $\text{Ang}(p,q)$ of angles between $p$ and $q$ is the subset of $[0,\pi/2]$ defined by the following (see Definition $2.1$ in \cite{SW}),
\begin{equation} 
\mathrm{Ang}(p,q)=\begin{cases}
\text{sp}~\Theta(p,q), & \text{if}~~ pq\neq qp.\\
\{\pi/2\},& \text{otherwise.}
\end{cases}
\end{equation}

Note that $\Theta$ is a positive operator and the spectrum of $\Theta$ is contained in $[0,\frac{\pi}{2}]$, but $0\mbox{ and }\frac{\pi}{2}$ are not eigenvalues.
\begin{dfn}[\cite{SW}]
Let $\mathcal{M}$ be a finite von Neumann algebra with a faithful normal tracial state $\mathrm{tr}$ and $\mathcal{P,Q}$ be von Neumann subalgebras of $\mathcal{M}$. The trace $\mathrm{tr}$ determines the normal faithful conditional expectations $E^{\mathcal{M}}_{\mathcal{P}}:\mathcal{M}\to\mathcal{P}$ and $E^{\mathcal{M}}_{\mathcal{Q}}:\mathcal{M}\to\mathcal{Q}$. They extend to the orthogonal projections $e_{\mathcal{P}}\mbox{ and }e_{\mathcal{Q}}$ on the GNS Hilbert space $L^2(\mathcal{M})$. The angle $\mathrm{Ang}_{\mathcal{M}}(\mathcal{P,Q})$ between $\mathcal{P}$ and $\mathcal{Q}$  is defined as follows~:
\[{\mathrm{Ang}}_{\mathcal{M}}(\mathcal{P,Q}):=\mathrm{Ang}(e_{\mathcal{P}},e_{\mathcal{Q}}).\]
\end{dfn}

\noindent The following facts are well-known (see \cite{SW}).
\begin{enumerate}[$(i)$]
\item For $pq\neq qp$, $\text{Ang}(p,q)\setminus\{\frac{\pi}{2}\}=\sigma\big(\mbox{arccos}\sqrt{pqp-p\wedge q}\big)\setminus\{\frac{\pi}{2}\}$.
\item ${\mathrm{Ang}}_{\mathcal{M}}(\mathcal{P,Q})={\mathrm{Ang}}_{\mathcal{M}}(\mathcal{Q,P})$ and ${\mathrm{Ang}}_{\widetilde{\mathcal{M}}}(\mathcal{P,Q})={\mathrm{Ang}}_{\mathcal{M}}(\mathcal{P,Q})$ for any finite von Neumann algebra $\widetilde{\mathcal{M}}\supset\mathcal{M}$ with trace $\widetilde{\mathrm{tr}}$ such that $\widetilde{\mathrm{tr}}|_{\mathcal{M}}=\mathrm{tr}$.
\item For quadrilateral of $II_1$ factors $(N\subset P,Q\subset M),\,\mathrm{Ang}_M(P,Q)$ is a finite set which contains neither $0$ nor $\frac{\pi}{2}$. Moreover, $\mathrm{Ang}_M(P,Q)=\sigma\big(\mbox{arccos}\sqrt{e_Pe_Qe_P-e_N}\big)\setminus\{\frac{\pi}{2}\}$.
\end{enumerate}
In \cite{JX}, Jones and Xu proved that finiteness of the angle (as a substet of $[0,\pi/2]$) is equivalent to finiteness of the (Pimsner-Popa)index of $P\cap Q$ in $M$. In \cite{GJ}, Grossman and Jones have considered a slight variant of the angle between two subfactors $P,Q\subset M$, namely the spectrum of the positive self-adjoint operator $E_PE_QE_P$ (on $L^2(M)$).
\smallskip

In another direction, to understand the relative position between a pair of intermediate subfactors $P$ and $Q$ of a finite index subfactor $N\subset M$ (i.e., $N\subset P,Q\subset M$), the first author along with Das, Liu and Ren have introduced a new notion of angle $\alpha^N_M(P,Q)$ between $P$ and $Q$ \cite{BDLR}. This angle has been used to answer an open question by Longo. More explicitly, the existing upper bound for the cardinality of the lattice of intermediate subfactors have been improved. Furthermore, a surprising connection between the intermediate subfactor theory and kissing numbers/sphere packing in geometry have been discovered.

\begin{dfn}[\cite{BDLR}](\textit{Angle between intermediate subfactors})
Let $P$ and $Q$ be two intermediate subfactors of a subfactor $N\subset M$. Then, the interior angle $\alpha^N_M(P,Q)$ between $P$ and $Q$ is given by
\[
\alpha^N_M(P, Q)=\cos^{-1}{\langle v_P,v_Q\rangle}_2\,,
\]
where $v_P:=\frac{e_P-e_N}{{\lVert e_P-e_N\rVert}_2}$ (and similarly $v_Q$),$\,{\langle x,y\rangle}_2:=tr(y^*x)$ and ${\lVert x\rVert}_2:=(tr(x^*x))^{1/2}$. The exterior angle between $P$ and $Q$ is given by $\beta^N_M(P, Q)=\alpha^M_{M_1}(P_1, Q_1)$.
\end{dfn}

\begin{figure}[!h]
\begin{center} \resizebox{9 cm}{!}{\includegraphics{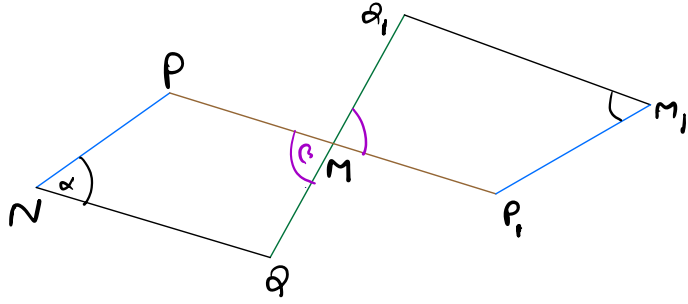}}\end{center}
\caption{Interior and Exterior angles}\label{int and ext ang}
\end{figure}

The interior and exterior angle can be depicted pictorially as in \Cref{int and ext ang}. Here, $P_1,Q_1,M_1$ denote the basic constructions of $P\subset M,\,Q\subset M,\mbox{ and }N\subset M$ respectively.

The following results are well-known.
\begin{ppsn}[\cite{SW,BDLR}]\label{basicresultoncomm}
Consider a quadruple
\begin{equation}\label{quadruple}
\begin{matrix}
\mathcal{Q} &\subset & \mathcal{M} \cr \cup &\ &\cup\cr \mathcal{N} &\subset & \mathcal{P}
\end{matrix}
\end{equation}
of finite von Neumann algebras with a fixed trace $tr$ on $\mathcal{M}$ and suppose that $E^{\mathcal{M}}_{\mathcal{Q}}, E^{\mathcal{M}}_{\mathcal{P}}$ and $E^{\mathcal{M}}_{\mathcal{N}}$ are the corresponding $tr$-preserving conditional expectations. The quadruple $(\ref{quadruple})$ is a commuting square if and only if $\mathrm{Ang}_M(P,Q)=\pi/2.$ Furthermore, if the quadruple $(\ref{quadruple})$ is a quadruple of $II_1$ factors with $[M:N]<\infty$, then it is a commuting square if and only if $\alpha^N_M(P,Q)=\pi/2$ and co-commuting square if and only if $\beta^N_M(P,Q)=\pi/2.$ 
\end{ppsn}

\begin{ppsn}[\cite{BDLR}]\label{to refer}
Consider intermediate subfactors $P\mbox{ and }Q$ of $II_1$ factors $N\subset M$. Let $\{\lambda_i\}$ (resp., $\{\mu_j\}$) be (right) basis for $P/N$ (resp., $Q/N$). Then,
\[
\cos(\alpha^N_M(P,Q))=\frac{\sum_{i,j}\mathrm{tr}\big(E^M_N(\lambda_i^*\mu_j)\mu_j^*\lambda_i\big)-1}{\sqrt{[P:N]-1}\,\sqrt{[Q:N]-1}}\,.
\]
\end{ppsn}

Below we mention a few well-known examples of commuting squares.

\begin{xmpl}
\begin{enumerate}
\item[$(1)$] If $R$ is the hyperfinite $II_1$ factor and $G$ is a finite group with subgroups $H$ and $K$, then the quadruple
$$\begin{matrix}
R\rtimes K &\subset & R\rtimes G \cr \cup &\ &\cup\cr R &\subset & R\rtimes H
\end{matrix}$$ is a commuting square. Thus, $\mathrm{Ang}_R(R\rtimes H, R\rtimes K)=\{\pi/2\}$. 
\item[$(2)$] Suppose that $K\subset L$ is an inclusion of finite von Neumann algebras with a  fixed trace $tr_L$ on $L$. Then, $M_n(L)$ is again a finite von Neumann algebra with the trace $tr_{M_n(\mathbb{C})}\otimes tr_L$, and it follows that the following quadruple
\[\begin{matrix}
M_n(K) &\subset & M_n(L) \cr \cup &\ &\cup\cr K &\subset & L\end{matrix}\]
is a commuting square with respect to the $\,tr_{M_n(\mathbb{C})}\otimes tr_L$-preserving conditional expectations.
\end{enumerate}
\end{xmpl}

Complex Hadamard matrices provide examples of  commuting and co-commuting squares (see \cite{JS}, for instance), which will play a central role in this paper.

\begin{ppsn}
Consider the subalgebra $\Delta_n$ of diagonal matrices in $M_n(\mathbb{C})$, and let $u$ be a unitary matrix in $M_n(\mathbb{C})$. Then, the following conditions are equivalent.
\begin{enumerate}[$(i)$]
\item The quadruple $$\begin{matrix}
 u\Delta_n u^* &\subset & M_n(\mathbb{C}) \cr \cup &\ &\cup\cr \mathbb{C} &\subset & \Delta_n
\end{matrix}$$ is a commuting and co-commuting square.
\item $u$ is a complex Hadamard matrix.
\end{enumerate}
\end{ppsn}

\subsection{Pimsner-Popa probabilistic constant}

Given a subfactor $N\subset M$ of type $II_1$ factors, Jones \cite{Jo} discovered a notion of index $[M:N]$ as the Murray-von Neumann’s coupling constant $\text{dim}_{N}(L^2(M))$. For von Neumann subalgebras $\mathcal{Q \subset P\subset M}$ of a finite von Neumann algebra $\mathcal{M}$, Pimsner and Popa \cite{PP} defined the following probabilistic constant 
\begin{equation}\label{probabilistic index}
\lambda(\mathcal{P,Q})=\text{sup}\{\lambda\geq 0\,:\,E_\mathcal{Q}^\mathcal{P}(x)\geq \lambda x\,\,\forall\,x\in {\mathcal{P}}_{+}\}.
\end{equation}
Moreover, if $M$ is a type $II_1$ factor and $N\subset M$ is a subfactor, then $E_N(x)\geq [M:N]^{-1}x$ for all $x\in M_+$ (with the convention $\frac{1}{\infty}=0$) and
\begin{equation}\label{pimsnerpopa2}
\lambda(M,N)=[M:N]^{-1}.
\end{equation}
Thus, the Pimsner-Popa probabilistic constant serves as a substitute of the Jones index.

\begin{dfn}
For von Neumann algebras $\mathcal{Q \subset P}$, we call $\lambda(\mathcal{P},\mathcal{Q})^{-1}$ the Pimsner-Popa index.
\end{dfn}

We would like to mention the following important properties of $\lambda$.
\smallskip

\noindent\textbf{Property 1:~} $\lambda$ does not respect tensor product. The simplest example is that $\lambda(M_2,\bbc)=\frac{1}{2}$, whereas $\lambda(M_2\otimes M_2,\bbc\otimes M_2)=\lambda(M_4,M_2)=\frac{1}{4}$ (see Sec. $6$ in \cite{PP}).
\smallskip

\noindent\textbf{Property 2:~} For $\mathcal{N\subset Q\subset P\subset M}$ we have $\lambda(\mathcal{P,N})\geq\lambda(\mathcal{P,Q})\lambda(\mathcal{Q,N})$. To see this, observe that $E^\mathcal{Q}_\mathcal{N}(x)\geq\lambda(\mathcal{Q,N})x$ for all $x\in \mathcal{Q}_+$ implies that $E^\mathcal{Q}_\mathcal{N}E^\mathcal{P}_\mathcal{Q}(y)\geq\lambda(\mathcal{Q,N})E^\mathcal{P}_\mathcal{Q}(y)$ for all $y\in \mathcal{P}_+$. This says that $E^\mathcal{P}_\mathcal{N}(y)\geq\lambda(\mathcal{Q,N})E^\mathcal{P}_\mathcal{Q}(y)\geq\lambda(\mathcal{Q,N})\lambda(\mathcal{P,Q})y$ for all $y\in \mathcal{P}_+$, and hence the result follows. We call it the submultiplicativity property of $\lambda$.
\smallskip

\noindent\textbf{Property 3:~} In general, $\lambda$ is not multiplicative. For example, take $\bbc\subset M_2\subset M_4$ and apply Theorem $6.1$ in Sec. $6$, \cite{PP} to observe that the inequality in property $2$ is strict.
\smallskip

\noindent\textbf{Property 4:~} $\lambda$ is not well-behaved with respect to (symmetric) commuting square. For example, consider the quadruple $\big(\bbc\otimes\bbc\subset\bbc\otimes M_2,M_2\otimes\bbc\subset M_2\otimes M_2\big)$ which is a (symmetric) commuting square. Although $\lambda(\bbc\otimes M_2,\bbc\otimes\bbc)=\lambda(M_2,\bbc)=1/2$, we have $\lambda(M_2\otimes M_2,M_2\otimes\bbc)=\lambda(M_4,M_2)=1/4$ by Theorem $6.1$ in Sec. $6$, \cite{PP}.
\medskip

The definition of $\lambda(\mathcal{P,Q})$ as in \Cref{probabilistic index} works for general von Neumann subalgebras $\mathcal{P}$ and $\mathcal{Q}$ (not necessarily $\mathcal{Q \subset P})$ of a finite von Neumann algebra $\mathcal{M}$ as well.

\begin{dfn}[Pimsner-Popa]\label{pimsner-popa}
Consider a pair of von Neumann subalgebras $\mathcal{P}$ and $\mathcal{Q}$ of a finite von Neumann algebra $\mathcal{M}$. The Pimsner-Popa probabilistic constant for the ordered pair $(\mathcal{P,Q})$  is defined as follows~:
\[
\lambda(\mathcal{P,Q})=\sup\{t\geq 0\,:\,E_\mathcal{Q}^\mathcal{M}(x)\geq t x\,\,\forall\,x\in {\mathcal{P}}_{+}\}.
\]
\end{dfn}

Note that in general, $\lambda(\mathcal{P,Q}))\neq\lambda(\mathcal{Q,P})$. For a quadruple $(N\subset P,Q\subset M)$ of type $II_1$ factors with $[M:N]<\infty$, if the quadruple is a commuting square, then $\lambda(P,Q)=[P:N]^{-1}$ and $\lambda(Q,P)=[Q:N]^{-1}$. Moreover, if the quadruple is a co-commuting square, then $\lambda(P,Q)=[M:Q]^{-1}$ and $\lambda(Q,P)=[M:P]^{-1}$ (see Theorem $3.12$ in \cite{B}). We mention that we always work with unital algebras, and we have $0\leq \lambda(\mathcal{P,Q})\leq 1$. Indeed, $E_\mathcal{Q}^\mathcal{M}(x)\geq t x$ implies that $||x||\geq ||E_\mathcal{Q}^\mathcal{M}(x)||\geq t ||x||$, and hence $\lambda(\mathcal{P,Q})\leq 1$. This justifies the term `{\it probabilistic constant}'. Henceforth, we drop the word `probabilistic' for brevity. For a type $II_1$ factor $M$ and subfactors $N_1,N_2\subset M$, we have $E_{N_2}^M(x)\geq [M:N_2]^{-1}x\,\,\forall\,x\in M_+$, and hence $\lambda(M,N_2)\leq\lambda(N_1,N_2)$. Thus, we always have a lower bound for $\lambda(N_1,N_2)$, namely $\lambda(N_1,N_2)\geq[M:N_2]^{-1}$. Computation of $\lambda(\mathcal{P,Q})$ is in general hard. If $(N,Q,P,M)$ is a quadruple of type $II_1$ factors with $[M:N]< \infty$ and $N^{\prime}\cap M=\mathbb{C}$, then the following formula
\begin{center}
$\lambda(P,Q)=\displaystyle \frac{\mbox{tr}(e_Pe_Q)}{\mbox{tr}(e_P)}$
\end{center}
is known due to Theorem $3.3$ in \cite{B}. However, if $N\subset M$ is not irreducible, then the situation is unclear. It seems to be a proper place to mention that we have some results in this direction in Sections \ref{Sec 4} and \ref{Sec 5}.

\begin{lmma}\label{sk}
Consider a pair of von Neumann subalgebras $\mathcal{P}$ and $\mathcal{Q}$ of a finite von Neumann algebra $\mathcal{M}$. Then, $\lambda(\mathcal{P,Q})=1$ if and only if $\mathcal{P}\subseteq\mathcal{Q}$.
\end{lmma}
\begin{prf}
Assume that $\lambda(\mathcal{P,Q})=1$. Then, $E_\mathcal{Q}^\mathcal{M}(x)\geq x$ for all $x\in\mathcal{P}_+$. Since the trace `$\mbox{tr}$' on $\mathcal{M}$ is faithful, we have $\mbox{tr}(E_\mathcal{Q}^\mathcal{M}(x)-x)\geq 0$. However, $\mbox{tr}\circ E_\mathcal{Q}^\mathcal{M}=\mbox{tr}$ implies that $E_\mathcal{Q}^\mathcal{M}(x)=x$, and hence $x\in\mathcal{Q}$. Thus $\mathcal{P}_+\subseteq\mathcal{Q}$, and hence $\mathcal{P}\subseteq\mathcal{Q}$. The converse direction is obvious from the definition.\qed
\end{prf}

Adapting the proof of Proposition $2.6$ in \cite{PP}, we get the following result which is very helpful in finding the value of $\lambda$.

\begin{ppsn}\label{popaadaptation}
\begin{enumerate}[$(i)$]
\item Let $\{M_n\}, \{A_n\}$ and $\{B_n\}$ be increasing sequences of von Neumann subalgebras of a finite von Neumann algebra $M$ such that $\{A_n\}, \{B_n\}\subset M$ and $M=\big(\bigcup_{n=1}^{\infty}M_n\big)^{\dprime}.$ If $A=\big(\bigcup_{n=1}^{\infty} A_n\big)^{\dprime}$ and $B=\big(\bigcup_{n=1}^{\infty} B_n\big)^{\dprime}$, then $\lambda(B,A)\geq \displaystyle \limsup\,\lambda(B_n,A_n).$
\item If in adition, $E_{A_{n+1}}E_{M_n}=E_{A_n}$ and $E_{B_{n+1}}E_{M_n}=E_{B_n}$ for $n\in \mathbb{N}$, then $\lambda(B,A)=\displaystyle \lim\,\lambda(B_n,A_n),$ decreasingly.
\end{enumerate}
\end{ppsn}
\begin{prf}
First observe that for $x\in M$, $E_A(x)=\displaystyle \lim_{n \to \infty} E_{A_n}(x)$ and $E_BE_A(x)=\displaystyle \lim_{n\to \infty} E_{B_n} E_{A_n}(x)$ in strong operator topology. Put $\lambda_n=\lambda(B_n,A_n)$ and let $\lambda=\limsup \lambda_n$. Thus, for arbitrarily fixed $\epsilon>0$, we must obtain a subsequence $\{\lambda_{k_n}\}$ such that $\lambda_{k_n}\geq \lambda-\epsilon$, for $n\geq 1.$ Now, for any $x\in M_{+}$, we have
\[
E_{A_{k_n}}\big(E_{B_{k_n}}(x)\big)\geq \lambda_{k_n} E_{B_{k_n}}(x)\geq (\lambda-\epsilon) E_{B_{k_n}}(x)
\]
and letting $n\to \infty$, we obtain $E_A\big(E_B(x)\big)\geq (\lambda-\epsilon) E_B(x)$ for all $x \in M_{+}$. Therefore, $\lambda(B,A)\geq \lambda-\epsilon$. As $\epsilon>0$ is arbitrary, we conclude the first part.

Now suppose that $E_{A_{n+1}}E_{M_n}=E_{A_n}$ and $E_{B_{n+1}}E_{M_n}=E_{B_n}$ for $n\in \mathbb{N}$. To see that the sequence $\{\lambda_n\}$ is decreasing, consider $x\in (B_n)_+$ and observe that $E_{A_n}(x)=E_{A_{n+1}}(x)\geq \lambda_{n+1} x$, and hence $\lambda_n\geq \lambda_{n+1}$. Furthermore, we see that for any $n\geq 1, E_AE_{M_n}=E_{A_n}$ and $E_BE_{M_n}=E_{B_n}$. So, for any $x\in(B_n)_+$ we have $E_{A_n}(x)= E_A(x)\geq \lambda(B,A) x$. Thus $\lambda_n\geq \lambda(B,A)$, and hence we get $\lambda(B,A)\leq \displaystyle \liminf\,\lambda_n\,$ which concludes the second part.\qed
\end{prf}

Note that if $A_n\subset B_n$ for all $n\geq 1$ in the above, then we recover Proposition $2.6$ in \cite{PP}.

\subsection{Connes-St{\o}rmer relative entropy}

To obtain an appropriate definition of the entropy of an automorphism in the non-commutative framework of operator algebras, Connes and St{\o}rmer \cite{CS} defined relative entropy $H(\mathcal{P|Q})$ between a pair of finite-dimensional von Neumann-subalgebras $\mathcal{P}$ and $\mathcal{Q}$ of a finite von Neumann algebra $\mathcal{M}$ equipped with a fixed faithful normal trace. This generalizes the classical notion of conditional entropy from ergodic theory. Using the relative entropy as the main technical tool, they have proved a non-commutative version of the Kolmogorov-Sinai type theorem (see \cite{NS} also). Pimsner and Popa \cite{PP} observed that the definition of the Connes-St{\o}rmer relative entropy does not depend on $\mathcal{P,Q}$ being finite-dimensional, so that one may also consider the relative entropy $H(\mathcal{P|Q})$ for arbitrary von Neumann subalgebras $\mathcal{P,Q \subset M}.$
\begin{dfn}[\cite{CS}]
Let $(\mathcal{M},\tau)$ be a finite von Neumann algebra and $\mathcal{P,Q\subseteq M}$ are von Neumann subalgebras. Let
\begin{align*}
\gamma &= \{x_j\in \mathcal{M}_+:\sum x_j=1,\,j=1,\ldots,n\}\,\,\,\mbox{ be a finite partition of unity},\cr
\eta &: [0,\infty)\longrightarrow\mathbb{R}\,\mbox{ be the continuous function }\,\,t\longmapsto -t\log t,\cr
H_\gamma(\mathcal{P|Q}) &:= \sum_{j=1}^n\left(\tau\circ\eta\,E_{\mathcal{Q}}^{\mathcal{M}}(x_j)-\tau\circ\eta\,E_{\mathcal{P}}^{\mathcal{M}}(x_j)\right).
\end{align*}
Then, $H(\mathcal{P|Q}):=\sup_{\gamma}\,H_\gamma(\mathcal{P|Q})$ is the Connes-St{\o}rmer relative entropy between $\mathcal{P}\mbox{ and }\mathcal{Q}$.
\end{dfn}

The relative entropy has the following important properties (see \cite{CS},\cite{PP},\cite{NS}).
\smallskip

\noindent\textbf{Property 1:~} $H(\mathcal{P|Q})\geq 0$ and $H(\mathcal{P|Q})=0$ if and only if $\mathcal{P}\subseteq\mathcal{Q}$.
\smallskip

\noindent\textbf{Property 2:~} $H(\mathcal{P|Q})\leq H(\mathcal{P|N})+H(\mathcal{N|Q})$ for $\mathcal{P,Q,N\subseteq M}$.
\smallskip

\noindent\textbf{Property 3:~} $H(\mathcal{P|Q})$ is increasing in $\mathcal{P}$ and decreasing in $\mathcal{Q}$.
\smallskip

\noindent\textbf{Property 4:~} In general, $H$ is not additive. For example, take $\bbc\subset M_2\subset M_4$ and apply ($6.5$ in Sec. $6$, \cite{PP}) to observe that $H(M_4|\bbc)=\log 4<H(M_4|M_2)+H(M_2|\bbc)=\log 4+\log 2$.
\smallskip

By Property $2$, for a pair of von Neumann subalgebras $\mathcal{P,Q\subseteq M}$ we always have an upper bound for the relative entropy $H(\mathcal{P|Q})$, namely $H(\mathcal{P|Q})\leq H(\mathcal{P|M})+H(\mathcal{M|Q})=H(\mathcal{M|Q})$ (using Property $1$ we have $H(\mathcal{P|M})=0$). In fact, $\mathcal{M}$ can be replaced by $\mathcal{P}\vee\mathcal{Q}$, that is, $H(\mathcal{P|Q})\leq H(\mathcal{P}\vee\mathcal{Q}|\mathcal{Q})$. We record the following useful result which shows that the relative entropy behaves well with respect to `controlled' limit, the proof of which follows from a minor modification of Proposition $3.4$ in \cite{PP}, and hence we omit it.

\begin{ppsn}\label{popaadaptation2}
\begin{enumerate}[$(i)$]
\item Let $\{M_n\}, \{A_n\}$ and $\{B_n\}$ be increasing sequences of von Neumann subalgebras of a finite von Neumann algebra $M$ such that $\{A_n\}, \{B_n\}\subset M$ and $M=\big(\bigcup_{n=1}^{\infty}M_n\big)^{\dprime}.$ If $A=\big(\bigcup_{n=1}^{\infty} A_n\big)^{\dprime}$ and $B=\big(\bigcup_{n=1}^{\infty} B_n\big)^{\dprime}$, then $H(B|A)\leq \displaystyle \liminf\,H(B_n|A_n).$
\item If in addition, $E_{A_{n+1}}E_{M_n}=E_{A_n}$ and $E_{B_{n+1}}E_{M_n}=E_{B_n}$ for $n\in \mathbb{N}$, then $H(B|A)=\displaystyle \lim\,H(B_n|A_n)$ increasingly.
\end{enumerate}
\end{ppsn}

In most of the situations, finding exact value of $H(\mathcal{P|Q})$ is difficult. For example, if we consider the following simple situation
\[
\begin{matrix}
\Delta_n &\subset & M_n \cr
 & & \cup\cr
 &  & U\Delta_n U^*
\end{matrix}
\]
where $U\in M_n$ is a unitary and $\Delta_n\subset M_n$ is the diagonal subalgebra (Masa), then the formula for $H(\Delta_n|U\Delta_n U^*)$, even when $n=2$, seems to be hard and is open in the literature as mentioned in \cite{PSW}.

In literature, there is a modified version of the Connes-St{\o}rmer relative entropy which is `more' computable. We will take help of it at a place towards the end of the article to compute the value of the Connes-St{\o}rmer relative entropy, and hence discuss it briefly.

\begin{dfn}[\cite{choda},\cite{choda2}]
Suppose that $(\mathcal{M},\tau)$ be a finite von Neumann algebra and $\mathcal{P,Q\subseteq M}$ are von Neumann subalgebras. Let
\begin{align*}
\gamma &= \{x_j\in \mathcal{P}_+:\sum x_j=1,\,j=1,\ldots,n\}\,\,\,\mbox{ be a finite partition of unity},\cr
\eta &: [0,\infty)\longrightarrow\mathbb{R}\,\mbox{ be the continuous function }\,\,t\longmapsto -t\log t,\cr
h_\gamma(\mathcal{P|Q}) &:= \sum_{j=1}^n\left(\tau\circ\eta\,E_{\mathcal{Q}}^{\mathcal{M}}(x_j)-\tau\circ\eta\,(x_j)\right).
\end{align*}
Then, $h(\mathcal{P|Q}):=\sup_{\gamma}\,h_\gamma(\mathcal{P|Q})$ is the modified Connes-St{\o}rmer relative entropy between $\mathcal{P}\mbox{ and }\mathcal{Q}$.
\end{dfn}

If $\mathcal{M}$ is abelian, then $H(\mathcal{P|Q})=h(\mathcal{P|Q})$. Thus, $h$ also generalizes the classical relative entropy. Moreover, if $\mathcal{P}\subset\mathcal{M}$ then $H(\mathcal{M|P})=h(\mathcal{M|P})$. Also, it is known that $0\leq h(\mathcal{P|Q})\leq H(\mathcal{P|Q})$. Moreover, the following result holds (see \cite{choda2}) to which we give a proof for the sake of completeness.

\begin{ppsn}[\cite{choda2}]\label{df}
Let $(\mathcal{M},\tau)$ be a finite von Neumann algebra and the quadruple
\[
\begin{matrix}
P &\subset & M\\
\cup & &\cup\\
N &\subset & Q\\
\end{matrix}
\]
be a commuting square. Then, $H(\mathcal{P|Q})=h(\mathcal{P|Q})=H(\mathcal{P|N})=h(\mathcal{P|N})$.
\end{ppsn}
\begin{prf}
Since $H(\mathcal{P|N})=h(\mathcal{P|N})$, it is enough to show that $h(\mathcal{P|Q})=h(\mathcal{P|N})$. This is because by Theorem $6$ in \cite{Wa}, it is known that $H(\mathcal{P|Q})=H(\mathcal{P|N})$ in the case of commuting squares, and we also have $H(\mathcal{P|N})=h(\mathcal{P|N})$ because $\mathcal{N}\subset\mathcal{P}$. Let $\mathscr{P}$ be the set of all finite partitions $\{(x_i)_{1\leq i\leq n}:\sum_{i=1}^nx_i=1\mbox{ and each }x_i\in\mathcal{P}_+\}$. Similar to the proof of Theorem $6$ in \cite{Wa}, we have the following,
\begin{IEEEeqnarray*}{lCl}
h(\mathcal{P|N}) &=& \sup_{(x_i)\in\mathscr{P}}\sum_i\tau\eta E^\mathcal{M}_\mathcal{N}(x_i)-\tau\eta(x_i)\\
&=& \sup_{(x_i)\in\mathscr{P}}\sum_i\tau\eta E^\mathcal{M}_\mathcal{Q}E^\mathcal{M}_\mathcal{P}(x_i)-\tau\eta(x_i)\\
&=& \sup_{(x_i)\in\mathscr{P}}\sum_i\tau\eta E^\mathcal{M}_\mathcal{Q}(x_i)-\tau\eta(x_i)\\
&=& h(\mathcal{P|Q})
\end{IEEEeqnarray*}
which completes the proof.\qed
\end{prf}

Therefore, we may safely conclude that the modifed Connes-St{\o}rmer relative entropy $h$ is closely related to the Connes-St{\o}rmer relative entropy $H$. For the situation $\Delta_n,U\Delta_n U^*\subset M_n$ discussed above, the value of $h(\Delta_n|U\Delta_n U^*)$ is known to be $\frac{1}{n}\sum_{i,j}\eta(|u_{ij}|^2)$ for $U=(u_{ij})$ \cite{choda}. This shows that the importance of $h$ lies in the fact that it is a more computable quantity than $H$, and it always provides a lower bound for $H$. Moreover, for certain `nice' situation as indicated in \Cref{df}, it determines the exact value of $H$. Note that we are not claiming here computation of $h$ is easy, rather, we only want to emphasize the fact that $h$ is indeed a `good' invariant for pair of von Neumann algebras, which is closely related to $H$ and theoretically more computable. The Pimsner-Popa constant $\lambda(\mathcal{P,Q})$ defined earlier always provides an upper-bound for $h(\mathcal{P|Q})$.

At this stage, we would like to add a comment here. In Theorem $4.3$, \cite{B} it is claimed that if $N\subset M$ is an irreducible subfactor with finite Jones index and $P$ and $Q$ are two intermediate subfactors, then $H(P|Q)=-\log\lambda(P,Q)$. However, there is a mistake in obtaining the upper bound in the first part of the proof (namely, $H(P|Q)\leq-\log\lambda(P,Q)$). The proof works for $h(P|Q)$ and not $H(P|Q)$. The correct replacement of Theorem $4.3$ in \cite{B} are the following two statements.

\begin{ppsn}\label{upbound}
If $\mathcal{P}$ and $\mathcal{Q}$ are von Neumann subalgebras of a finite von Neumann algebra $\mathcal{M}$, then $\,h(\mathcal{P|Q})\leq -\log \lambda(\mathcal{P},\mathcal{Q})$.
\end{ppsn}

This will follow by the first part of the proof of Theorem $4.3$ in \cite{B}.

\begin{ppsn}\label{B}
Suppose $N\subset M$ is an irreducible subfactor with finite Jones index and $P$ and $Q$ are two intermediate subfactors. Then $h(P|Q)=-\log \lambda(P,Q)$, and consequently $H(P|Q)\geq-\log \lambda(P,Q)$.
\end{ppsn}

This will follow by the second part of the proof of Theorem $4.3$ in \cite{B} together with \Cref{upbound}.

Remarkably, Pimsner and Popa in \cite{PP} had discovered that for a subfactor $N\subset M$ of type $II_1$ factors, $H(M|N)(=h(M|N))$ depends on both the Jones index and the relative commutant. In particular, for a finite index subfactor $N\subset M$ with $N^{\prime}\cap M=\mathbb{C}$ (such a subfactor is called irreducible) they proved that $H(M|N)= \log\,[M:N]$. More generally, the subfactor $N\subset M$ is extremal if and only if $H(M|N)=\log\,[M:N]$. For a finite index irreducible subfactor $N\subset M$ with intermediate subfactors $P$ and $Q$, the Pimsner-Popa constant $\lambda(P,Q)$ is closely related to both $h(P|Q)$ and $H(P|Q)$ (\Cref{upbound} and \ref{B}).
\smallskip

As an application of \Cref{B}, we now show that possible values of $h(P|Q)$ has certain gap in the irreducible situation. This is reminiscent of the Jones index rigidity.

\begin{thm}\label{small index rigidity}
Consider a quadrilateral of finite index subfactors
\[
\begin{matrix}
{Q} &\subset &{M} \cr \cup &\ &\cup\cr {N} &\subset &{P}
\end{matrix}
\]
with $N\subset M$ irreducible and suppose that $[P:N]$ and $[Q:N]$ are less than $4$. Then, either
\begin{center}
$h(P|Q)\in\big\{\log\big(4~{\cos}^2(\pi/n)\big)\,:\,n\geq 3\big\}$
\end{center}
or
\begin{center}
$h(P|Q)=h(Q|P)\in\{\log\big({4~{\cos}^2 (\pi/2k)-1}\big)\,:\,k\geq 3\}.$
\end{center}
\end{thm}
\begin{prf}
Consider the dual factors $\bar{P}$  and $\bar{Q}$  from the basic constructions $P\subset M$ and $Q\subset M$, respectively (with the corresponding Jones projections $e_{\bar{P}}$ and $e_{\bar{Q}}$). Then, we have
\[
tr\big(e_{\bar{P}}e_{\bar{Q}}\big)=\frac{[M:P]}{[Q:N]} tr\big(e_Pe_Q\big).
\]
Furthermore, if $N_{-1}\subset N\subset M, P_{-1}\subset N\subset P$ and $Q_{-1}\subset N\subset Q$ are instances of downward basic constructions with the corresponding Jones projection $e_{N_{-1}}, e_{P_{-1}}$ and $e_{Q_{-1}}$ respectively, then we also have 
\[
tr\big(e_{P_{-1}}e_{Q_{-1}}\big)=\frac{[M:Q]}{[P:N]} tr(e_Pe_Q).
\]
By Lemma 3.16 in \cite{GJ}, we get that $\,tr(e_Pe_Q)=\big(\mbox{dim}_M L^2(\bar{P}\bar{Q})\big)^{-1}$, and hence $\,tr(e_{P_{-1}}e_{Q_{-1}})=\frac{1}{\text{dim}_N L^2(PQ)}$. Therefore, we conclude the following,
\[
\text{dim}_M \big(L^2(\bar{P}\bar{Q})\big)= \frac{[M:Q]}{[P:N]}\text{dim}_N \big(L^2(PQ)\big).
\]
By \cite{B}, we see that
\[
\lambda(P,Q)=\frac{tr(e_Pe_Q)}{tr(e_P)}=\frac{tr(e_Q)}{tr(e_{PQ})}=\frac{[P:N]}{\text{dim}_N \big(L^2(PQ)\big)}.
\]
As $[P:N],[Q:N]<4$, two cases may arise. It may happen that $(N\subset P,Q\subset M)$ forms a commuting square. In this case, thanks to \cite{W}, we have $H(P|Q)=\log [P:N]$ and the result follows from Jones index rigidity. In the other case, by Theorem $4.1.3$ in \cite{G} (see also \cite{GI}) we have $\text{dim}_N\big(L^2(PQ)\big)=[P:N]\big([P:N]-1\big).$ Therefore, $\lambda(P,Q)=\frac{1}{[P:N]-1}$. The rest follows from Corollary  5.2.5 in \cite{G}.\qed
\end{prf}

\begin{crlre}
Consider a quadrilateral of finite index subfactors
\[
\begin{matrix}
{Q} &\subset &{M} \cr \cup &\ &\cup\cr {N} &\subset &{P}
\end{matrix}
\]
with $N\subset M$ irreducible and suppose that $[M:P]$ and $[M:Q]$ are less than $4$. Then, either
\begin{center}
$h(P|Q)\in\big\{\log\big(4~{\cos}^2 (\pi/n)\big)\,:\,n\geq 3\big\}$
\end{center}
or
\begin{center}
$h(P|Q)=h(Q|P)\in\{\log \big({4~{\cos}^2(\pi/2k) -1}\big)\,:\,k\geq 3\}.$
\end{center}
\end{crlre}
\begin{prf}
This follows from Corollary $4.4$ in \cite{B} and \Cref{small index rigidity}.\qed
\end{prf}

\subsection{Basic construction of non-degenerate commuting square}

Let
\[\begin{matrix}
\mathcal{Q} &\subset & (\mathcal{M}, tr) \cr \cup &\ &\cup\cr \mathcal{N} &\subset & \mathcal{P}
\end{matrix}\]
be a quadruple of finite von Neumann algebras which is a commuting and co-commuting (non-degenerate) square with respect to the $tr$-preserving conditional expectations $E^{\mathcal{M}}_{\mathcal{Q}}, E^{\mathcal{M}}_{\mathcal{P}}$ and $E^{\mathcal{M}}_{\mathcal{N}}$. Following (\cite{Po2},\cite{JS}), we discuss the basic construction for this quadruple. To avoid technicalities, we confine ourselves into the following two cases.
\smallskip

\textit{Case I:} Each of the elementary inclusions $\mathcal{N\subset P}$, $\mathcal{Q\subset M}$, $\mathcal{N\subset Q}$, and $\mathcal{P\subset M}$ is an inclusion of type $II_1$ factors with finite Jones index.

\textit{Case II:} Each of the elementary inclusions $\mathcal{N\subset P}$, $\mathcal{Q\subset M}$, $\mathcal{N\subset Q}$, and $\mathcal{P\subset M}$ is a connected inclusion of finite-dimensional $C^*$-algebras (see \cite{GHJ, JS}).

It is known that, in each case, we have a unique Markov trace (for the inclusion $\mathcal{Q\subset M}$) $tr$ on $\mathcal{M}$ (see \cite{JS, GHJ}). Let $e_{\mathcal{Q}}$ be the Jones projection (with respect to $tr$) for the inclusion $\mathcal{Q\subset M}$, and ${\mathcal{Q}}_1:= \langle \mathcal{M}, e_{\mathcal{Q}} \rangle$ be the von Neumann algebra generated in $\mathcal{B}(L^2(\mathcal{M}))$ by $\mathcal{M}$ and $e_{\mathcal{Q}}$, so that $\mathcal{Q\subset M}\subset \mathcal{Q}_1$ is an instance of basic construction. Denote by $\langle \mathcal{P},e^{\mathcal{M}}_{\mathcal{Q}}\rangle$ the von Neumann algebra  generated in $\mathcal{B}(L^2(\mathcal{M}))$ by $\mathcal{P}$ and $e_{\mathcal{Q}}$. Then, $\mathcal{N\subset P}\subset {\mathcal{P}}_1:= \langle \mathcal{P}, e_{\mathcal{Q}}\rangle$ is an instance of basic construction for the inclusion $\mathcal{N\subset P}$ with the Jones projection $e_{\mathcal{Q}}$. It follows that $\mathcal{P}_1\subset \mathcal{Q}_1$, and each of the inclusion $\mathcal{M}\subset {\mathcal{Q}}_1$ and $\mathcal{P}\subset \mathcal{P}_1$ is a subfactor with finite Jones index (resp. connected inclusion of finite-dimensional $C^*$-algebras) in \textit{Case I} (resp. \textit{Case II}). Therefore, we have a unique Markov trace $tr_{\mathcal{Q}_1}$ on $\mathcal{Q}_1$ (which restricts to $tr_{\mathcal{M}}$ on $\mathcal{M}$). The construction of the quadruples~:
\[
\begin{matrix}
\mathcal{Q} &\subset & \mathcal{M} & \subset & \mathcal{Q}_1\cr
\rotatebox{90}{$\subset$} &\ &\rotatebox{90}{$\subset$} &\ & \rotatebox{90}{$\subset$}\cr
\mathcal{N} &\subset & \mathcal{P} & \subset & \mathcal{P}_1\end{matrix}
\]
is called the \textit{basic construction} for the non-degenerate commuting square
\begin{IEEEeqnarray}{lCl}\label{lcomm}
\begin{matrix}
\mathcal{Q} &\subset & \mathcal{M} \cr \cup &\ &\cup\cr \mathcal{N} &\subset & \mathcal{P}
\end{matrix}
\end{IEEEeqnarray}
and the commuting square (with respect to the $tr_{\mathcal{Q}_1}$-preserving conditional expectations $E^{\mathcal{Q}_1}_{\mathcal{M}}, E^{\mathcal{Q}_1}_{\mathcal{P}_1}$ and $E^{\mathcal{Q}_1}_{\mathcal{P}}$) 
\[
\begin{matrix}
\mathcal{M} &\subset & \mathcal{Q}_1\cr \cup &\ &\cup\cr \mathcal{P} &\subset & \mathcal{P}_1
\end{matrix}\]
is called the extension of the commuting square (\ref{lcomm}). 

Non-degenerate commuting squares of finite-dimensional $C^*$-algebras and their basic construction play a central role in the abstract subfactor theory (see \cite{POPA}). They are also instrumental in constructing hyperfinite subfactors with finite Jones index. Below we sketch the construction. We fix a non-degenerate commuting square of finite-dimensional $C^*$-algebras
\[
\begin{matrix}
A_{10} &\subset &A_{11} \cr \cup &\ &\cup\cr A_{00} &\subset & A_{01}
\end{matrix}
\]
with respect to the unique trace $tr$ on $A_{11}$, which is Markov for the inclusion $A_{10}\subset A_{11}$ as in \textit{Case II}. The basic construction of the above square is again a non-degenerate commuting square, and is denoted by 
\[
\begin{matrix}
A_{11} &\subset &A_{12} \cr \cup &\ &\cup\cr A_{01} &\subset & A_{02}\,\,.
\end{matrix}
\]
We can iterate the basic construction and obtain the following ladder of non-degenerate commuting squares
\[
\begin{matrix}
A_{10} &\subset & A_{11} & \subset & A_{12} & \subset &  \cdots \cr
\rotatebox{90}{$\subset $} &\ &\rotatebox{90}{$\subset$} &\ & \rotatebox{90}{$\subset$} \cr
A_{00} &\subset & A_{01} & \subset & A_{02} & \subset & \cdots \,.
\end{matrix}
\]
Then, $tr$ extends to a faithful trace $\tau$ on $\bigcup_k A_{1k}$, and setting $A_{1,\infty}$ (resp.  $A_{0,\infty}$) as the GNS-completion of $\bigcup_k A_{1k}$ (resp. $\bigcup_k A_{0k}$) with respect to $\tau$, we obtain the hyperfinite subfactor $A_{0,\infty}\subset A_{1,\infty}$.  We record here two important well-known facts about this subfactor which will be useful later.

\begin{ppsn}[\cite{PP,GHJ,JS}]\label{index}
Let $A_{01},A_{10}, A_{0,\infty}$, and $A_{1,\infty}$ be as above. If $\Lambda$ denotes the inclusion matrix for the vertical inclusion $A_{00}\subset A_{10}$, then $[A_{1,\infty}:  A_{0,\infty}]={\lVert \Lambda \rVert}^2.$
\end{ppsn}

\begin{ppsn}[\cite{JS}]\label{ocneanucompactness}(Ocneanu compactness)
Let $A_{01},A_{10}, A_{0,\infty}$, and $A_{1,\infty}$ be as above. Then, 
$\big(A_{0,\infty}\big)^{\prime}\cap A_{1,\infty}=\big(A_{01})^{\prime}\cap A_{10}$
\end{ppsn}

\subsection{Some important results}

We first recall few results from Section $6$ in \cite{PP} concerning the Pimsner-Popa constant and relative entropy for inclusion of finite-dimensional $C^*$-algebras. These will be repeatedly used throughout the article in several places, and hence for the benefit of the reader, we collect them at a place together. We keep the same symbols used in \cite{PP} for reader's convenience.

Let $N:=\oplus_{r\in K}N_r$ and $M:=\oplus_{\ell\in L}M_\ell$ be finite-dimensional von Neumann algebras, i,e., the sets of indices $L\mbox{ and }K$ are finite, where $N_r$ is the algebra of $n_r\times n_r$ matrices and $M_\ell$ is the algebra of $m_\ell\times m_\ell$ matrices. Suppose that we have a unital inclusion $N\subset M$. Denote by $\Lambda=(a_{r\ell})_{r\in K\,,\,\ell\in L}$ the inclusion matrix of $N\subset M$ and by $t_\ell$ (respectively $s_r$) the traces of the minimal projections in $M_\ell$ (respectively $N_r$). Thus, if $m=(m_\ell)_{\ell\in L},\,n=(n_r)_{r\in K},\,t=(t_\ell)_{\ell\in L},\,s=(s_r)_{r\in K}$ are column vectors then $\Lambda t=s$ and $\Lambda^Tn=m$. We have the following.

\begin{thm}[Theorems $6.1, 6.2$ in \cite{PP}]\label{ent formula}
For the unital inclusion $N\subset M$, we have the following formulae,
\begin{enumerate}[$(i)$]
\item $\lambda(M,N)^{-1}=\max_\ell\left(\sum_r\frac{b_{r\ell}s_r}{t_\ell}\right)$, where $b_{r\ell}=\min\{a_{r\ell}\,,n_r\}$.
\item $H(M|N)=\sum m_\ell t_\ell\log\left(\frac{m_\ell}{t_\ell}\right)+\sum n_rs_r\log\left(\frac{s_r}{n_r}\right)+\sum_{r,\ell}n_ra_{r\ell}t_\ell\log(c_{r\ell})$,\\
 where $\,c_{r\ell}=\min\{n_r/a_{r\ell}\,,1\}$.
\end{enumerate}
\end{thm}

Particular cases of the above theorem are the following (Example $6.5$ in \cite{PP}).

\begin{crlre}[\cite{PP}]\label{ent formula1}
Let $N\subset M$ be such that $N$ and $M$ are factors of type $I_n$ and $I_m$ respectively.
\begin{enumerate}[$(i)$]
\item If $m/n>n$, then $\lambda(M,N)=1/m$ and $H(M|N)=\log m$.
\item If $m/n\leq n$, then $\lambda(M,N)=(n/m)^2$ and $H(M|N)=2(\log m-\log n)$.
\end{enumerate}
\end{crlre}

\begin{crlre}[\cite{PP}]\label{ent formula2}
Let $N\subset M$ be such that $M$ is the factor of type $I_m$ and $N=\oplus_{r\in K}N_r$, where each $N_r$ is the type $I_{n_r}$ factor, such that $\sum_{r\in K}n_r=m$. Then, $\lambda(M,N)=\frac{1}{\texttt{\#} K}$ and $H(M|N)=\sum_r\eta(n_r/m)$.
\end{crlre}

Note that these formulae are very useful in the sense that if we know the finite-dimensional tower for an inclusion $N\subset M$ of $II_1$ factors, then these formulae can (possibly) give us the exact value of $\lambda(M,N)$ and $H(M|N)$, since both behave well under controlled limits.

Now, we prove few results which will be used on several occasions, especially in Sections \ref{Sec 4} and \ref{Sec 5}.
\begin{lmma}\label{bakproc}
Suppose that for $j=1,2$ we have two symmetric commuting squares $(B_1^j,A_1,B_2^j,A_2)$ of finite-dimensional algebras (see \Cref{com1} in this regard), and we obtain the following ladder of symmetric commuting squares
\[
\begin{matrix}
A_1 &\subset & A_2 &\subset^{\,e_3} & A_3 &\subset^{\,e_4} &\cdots\\
\cup &  & \cup &   & \cup &  & \\
B_1^j &\subset & B_2^j &\subset & B_3^j &\subset &\cdots
\end{matrix}
\]
by iterated basic construction. That is, $A_k=\mathrm{Alg}\{A_{k-1},e_k\}$ for $k\geq 3$, where $e_k$ is the Jones projection, and $B_k^j=\mathrm{Alg}\{B_{k-1}^j,e_k\}$. Let $C_k=\cap_{j=1}^2B_k^j$ for $k\geq 1$ and suppose that $(C_1,A_1,C_2,A_2)$ is a symmetric commuting square. We let $\,D_3=\mathrm{Alg}\{C_2,e_3\}$ and for $k\geq 4,\,D_k=\mathrm{Alg}\{D_{k-1},e_k\}$. Then, $e_k\in C_k$ and $D_k\subset C_k$ for all $k\geq 3$.
\end{lmma}
\begin{prf}
It is immediate that $e_k\in C_k$ for all $k\geq 3$ and $D_3\subset C_3$. Suppose that $D_m\subset C_m$ for some $m\geq 4$. Then, $D_{m+1}=\mathrm{Alg}\{D_m,e_{m+1}\}\subset\mathrm{Alg}\{C_m,e_{m+1}\}\subset C_{m+1}$. Therefore, by induction on $k$ we get the result.\qed
\end{prf}

\begin{lmma}\label{01}
Let $\mathcal{N}\subset\mathcal{M}$ be an inclusion of von Neumann algebras and $u\in\mathcal{M}$ be a unitary. Then, in the following diagram
\[
\begin{matrix}
\mathcal{N} &\subset & \mathcal{M}\cr
 &\ &\rotatebox{90}{$\subset$}\cr
 & & \mathrm{Ad}_u(\mathcal{N})
\end{matrix}
\]
we have $(\mathrm{Ad}_u(\mathcal{N}))^\prime\cap\mathcal{N}=\mathrm{Ad}_u\big(\mathcal{N}^{\,\prime}\cap\mathrm{Ad}_{u^*}(\mathcal{N})\big)$.
\end{lmma}
\begin{prf}
Let $\xi\in\mathrm{Ad}_u\big(\mathcal{N}^{\,\prime}\cap\mathrm{Ad}_{u^*}(\mathcal{N})\big)$ and $\xi=uxu^*$ with $x\in\mathcal{N}^{\,\prime}\cap\mathrm{Ad}_{u^*}(\mathcal{N})$. Clearly, $\xi\in\mathcal{N}$ as $x\in\mathrm{Ad}_{u^*}(\mathcal{N})$. For any $unu^*\in\mathrm{Ad}_u(\mathcal{N})$ with $n\in\mathcal{N}$, since $x\mbox{ and }n$ commutes, we get that $\xi unu^*=uxu^*unu^*=uxnu^*=unxu^*=unu^*\xi$. This says that $\xi\in(\mathrm{Ad}_u(\mathcal{N}))^\prime\cap\mathcal{N}$, i,e., the inclusion `$\supseteq$' holds. The reverse inclusion `$\subseteq$' is easy to check.\qed
\end{prf}

This lemma can be generalized a bit more, as indicated by the following.
\begin{lmma}\label{com and ad interchange}
Let $\mathcal{N}\subset\mathcal{M}$ be an inclusion of von Neumann algebras, $\mathcal{S}$ be a subset of $\mathcal{N}$, and $u\in\mathcal{M}$ be a unitary. Then, $\mathrm{Ad}_u(\mathcal{S})^\prime\cap\mathcal{N}=\mathrm{Ad}_u\big(\mathcal{S}^\prime\cap\mbox{Ad}_{u^*}(\mathcal{N})\big)$. Thus, if $\mathcal{N}=\mathcal{M}$, then $\mathrm{Ad}_u(\mathcal{S})^\prime\cap\mathcal{M}=\mathrm{Ad}_u\big(\mathcal{S}^\prime\cap\mathcal{M})$.
\end{lmma}
\begin{prf}
Follows along the similar line of \Cref{01}.\qed
\end{prf}


\section{Commuting cube and two subfactors}\label{Sec 2}

The results obtained in this section are generic in nature, and we take help of these on many occasions throughout the article. Here we introduce the notion of `commuting cube'. Similar notion has appeared in different contexts (\cite{BV}, \cite{K}) in the finite-dimensional situation.

\begin{ppsn}\label{proving com}
Suppose that the following quadruple
\[
\begin{matrix}
Q &\subset & M\\
\cup & & \cup\\
N &\subset & P\\
\end{matrix}
\]
is a co-commuting square (of either $II_1$ factors or finite-dimensional $C^*$-algebras) with $[P:N]=[M:Q]$ (either Jones index or Watatani index respectively). Then, the quadruple is a commuting square.
\end{ppsn}
\begin{prf}
Since it is given that
\[
\begin{matrix}
Q &\subset & M\\
\cup & & \cup\\
N &\subset & P\\
\end{matrix}
\]
is a co-commuting square, the dual quadruple
\[
\begin{matrix}
\langle P,e_P\rangle &\subset & \langle M,e_N\rangle\\
\cup & & \cup\\
M &\subset & \langle Q,e_Q\rangle\\
\end{matrix}
\]
is a commuting sqaure by definition. Therefore, $\,E^{\langle M,e_N\rangle}_{\langle P,e_P\rangle}E^{\langle M,e_N\rangle}_{\langle Q,e_Q\rangle}(e_Q)=E^{\langle M,e_N\rangle}_M(e_Q)=E^{\langle Q,e_Q\rangle}_M(e_Q)$. This implies that $E^{\langle M,e_N\rangle}_{\langle P,e_P\rangle}(e_Q)=E^{\langle Q,e_Q\rangle}_M(e_Q)=[M:Q]^{-1}$. Hence, $e_PE^{\langle M,e_N\rangle}_{\langle P,e_P\rangle}(e_Q)=[M:Q]^{-1}e_P$. Since $e_P\in\langle P,e_P\rangle$, we get that $E^{\langle M,e_N\rangle}_{\langle P,e_P\rangle}(e_Pe_Q)=[M:Q]^{-1}e_P$. Therefore, we have the following,
\begin{IEEEeqnarray}{lCl}\label{zp}
&  & E^{\langle M,e_N\rangle}_{\langle P,e_P\rangle}(e_Pe_Q)=[M:Q]^{-1}e_P\nonumber\\
&\Rightarrow& \mbox{tr}(e_Pe_Q)=[M:Q]^{-1}[M:P]^{-1}\nonumber\\
&\Rightarrow& \mbox{tr}(e_Pe_Q)=[M:N]^{-1}\qquad (\mbox{as } [M:Q]=[P:N])\nonumber\\
&\Rightarrow& \mbox{tr}(e_Pe_Qe_P)=\mbox{tr}(e_N)\,.
\end{IEEEeqnarray}
Now, $N\subset Q$ implies that $e_N\leq e_Q$, and consequently $e_Pe_Ne_P\leq e_Pe_Qe_P$. Since $e_Pe_Ne_P=e_N$, we have $e_N\leq e_Pe_Qe_P$, and hence $\mbox{tr}(e_Pe_Qe_P-e_N)=\mbox{tr}(e_Pe_Q-e_N)\geq 0$. By the faithfulness of $\mbox{tr}$, along with \Cref{zp}, we get that $e_Pe_Qe_P=e_N$. Now, using $e_Ne_P=e_N,\,e_Ne_Q=e_N$ together with $e_Pe_Qe_P=e_N$, it is easy to verify that $(e_P-e_N)(e_Q-e_N)(e_P-e_N)=0$. Since $e_P-e_N$ and $e_Q-e_N$ both becomes projection, we get the following
\begin{IEEEeqnarray*}{lCl}
||(e_P-e_N)(e_Q-e_N)||^2 &=& ||(e_P-e_N)(e_Q-e_N)(e_Q-e_N)(e_P-e_N)||\\
&=& ||(e_P-e_N)(e_Q-e_N)(e_P-e_N)||\\
&=& 0\,.
\end{IEEEeqnarray*}
This gives us that $(e_P-e_N)(e_Q-e_N)=0$, and consequently $e_Pe_Q=e_N$. Therefore, the quadruple $(N\subset P,Q\subset M)$ is a commuting square.\qed
\end{prf}

\begin{figure}[!h]
\begin{center} \resizebox{8 cm}{!}{\includegraphics{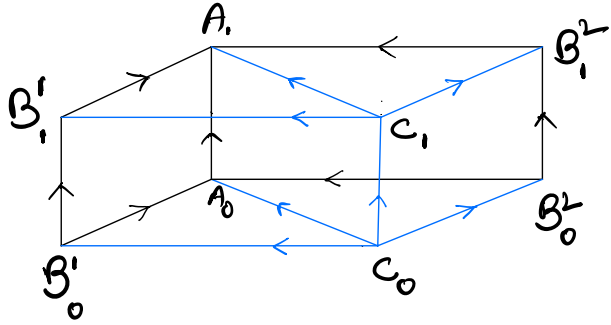}}\end{center}
\caption{Cube of finite von Neumann algebras}\label{com1}
\end{figure}

Consider the cube of finite von Neumann algebras described in \Cref{com1}, where $A_1$ is equipped with a faithful normal tracial state, and $C_0=B_0^1\cap B_0^2,\,C_1=B_1^1\cap B_1^2$. In \Cref{com1}, `$B\rightarrow A$' means $B\subset A$.

\begin{ppsn}\label{comm cube}
Suppose that the adjacent faces in \Cref{com1}
\[
\begin{matrix}
B_1^1 &\subset & A_1\\
\cup & &\cup & \\
B_0^1 &\subset & A_0\\
\end{matrix}\qquad\mbox{and}\qquad\begin{matrix}
A_1 &\supset & B_1^2\\
\cup & &\cup & \\
A_0 &\supset & B_0^2\\
\end{matrix}
\]
are commuting squares. Then, the remaining faces
\[
\begin{matrix}
B_0^j &\subset & B_1^j\\
\cup & &\cup & \\
C_0 &\subset & C_1\\
\end{matrix}
\]
for $j=1,2$ and the slice
\[
\begin{matrix}
C_1 &\subset & A_1\\
\cup & &\cup & \\
C_0 &\subset & A_0\\
\end{matrix}
\]
are also commuting squares. Therefore, $C_1\cap B_0^j=C_0$ for $j=1,2,$ and $C_1\cap A_0=C_0$.
\end{ppsn}
\begin{prf}
Recall that for a quadruple $(N\subset P,Q\subset M)$, if $E^M_Q(P)\subseteq N$ (or $E^M_P(Q)\subseteq N$), then it is a commuting square (\Cref{def of comm}). To show that the slice $(C_0\subset A_0,C_1\subset A_1)$ is a commuting square, in view of the above, it is enough to show that $E^{A_1}_{A_0}(c_1)\in C_0$ for all $c_1\in C_1$. Since $C_1\subset B^1_1$, by hypothesis we have $E^{A_1}_{A_0}(c_1)\in B^1_0$. Similarly, since $C_1\subset B^2_1$ we get $E^{A_1}_{A_0}(c_1)\in B^2_0$. Hence, $E^{A_1}_{A_0}(c_1)\in C_0$ which finishes the proof.

To show that the remaining faces of the cube are also commuting squares, it is enough to show one of these, and the other follows similarly. So we consider one of the remaining face
\[
\begin{matrix}
B_0^1 &\subset & B_1^1\\
\cup & &\cup & \\
C_0 &\subset & C_1\\
\end{matrix}
\]
and claim that $E^{B_1^1}_{B_0^1}(C_1)\subseteq C_0$. Take any $x\in C_1$ and using the hypothesis that $(B_0^1\subset A_0,B_1^1\subset A_1)$ is a commuting square observe the following,
\[
E^{B_1^1}_{B_0^1}(x)=E^{A_1}_{B_0^1}(x)=E^{A_1}_{A_0}E^{A_1}_{B_1^1}(x)=E^{A_1}_{A_0}(x)\,.
\]
Since the slice $(C_0\subset A_0,C_1\subset A_1)$ is already a commuting square, we get that $E^{A_1}_{A_0}(x)\in C_0$, which concludes the proof.\qed
\end{prf}

\begin{dfn}
A cube of finite von Neumann algebras as in \Cref{com1} is called a commuting cube if both the adjacent faces $(B_0^j\subset A_0,B_1^j\subset A_1)$ for $j=1,2$ are commuting squares.
\end{dfn}

Pair of spin model subfactors, to be discussed in the next section, provide natural examples of commuting cube. The following remark is very important.

\begin{rmrk}\rm\label{main difficulty}
\begin{enumerate}[$(i)$]
\item The floor and the roof in a commuting cube need not be commuting square. We will see counterexamples in Sections \ref{Sec 4} and \ref{Sec 5}.
\item If the adjacent faces in a commuting cube are non-degenerate commuting squares, then although the slice $(C_0\subset A_0,C_1\subset A_1)$ is a commuting square, it may fail to become non-degenerate. Same applies to the remaining faces. For example, consider two distinct $2\times 2$ complex Hadamard matrices $u\mbox{ and }v$ of the following form
\[
u=\frac{1}{\sqrt{2}}\,\left[{\begin{matrix}
1 & 1 \\
e^{i\alpha} & -e^{i\alpha} \\
\end{matrix}}\right] \quad\mbox{and}\quad v=\frac{1}{\sqrt{2}}\,\left[{\begin{matrix}
1 & 1 \\
e^{i\beta} & -e^{i\beta} \\
\end{matrix}}\right]
\]
such that $\alpha-\beta\notin\{0,\frac{\pi}{2},-\frac{\pi}{2}\}$. Obtain the cube as illustrated in \Cref{combasic0}.
\begin{figure}[!h]
\begin{center} \resizebox{8 cm}{!}{\includegraphics{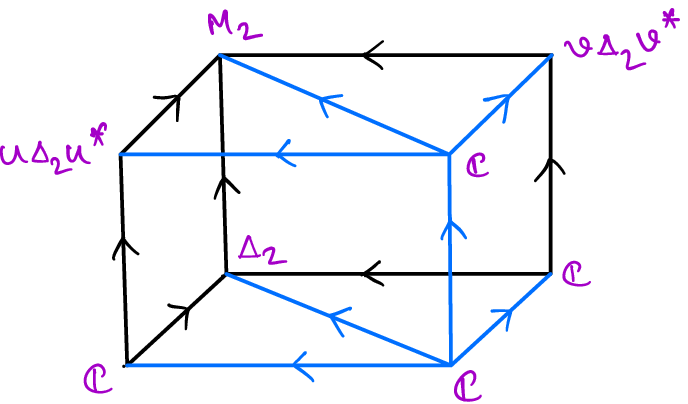}}\end{center}
\caption{A specific commuting cube}\label{combasic0}
\end{figure}
In this cube, both the adjacent faces are non-degenerate commuting squares, since $u\mbox{ and }v$ are complex Hadamard matrices. So, it is a commuting cube. Now, it is a straightforward verification that $\mbox{Ad}_u(\Delta_2)\cap\mbox{Ad}_v(\Delta_2)=\bbc$. Hence, the slice in the cube is the quadruple $(\bbc\subset\bbc,\Delta_2\subset M_2)$, which is clearly not non-degenerate, and so are the remaining faces.
\end{enumerate}
\end{rmrk}

\begin{ppsn}\label{comm cube 2}
Suppose that we have a commuting cube as in \Cref{com1}. If the roof
\[
\begin{matrix}
B_1^2 &\subset & A_1\\
\cup & &\cup & \\
C_1 &\subset & B_1^1\\
\end{matrix}
\]
is a commuting square, then the floor
\[
\begin{matrix}
B_0^2 &\subset & A_0\\
\cup & &\cup & \\
C_0 &\subset & B_0^1\\
\end{matrix}
\]
is also a commuting square.
\end{ppsn}
\begin{prf}
Let $x\in B_0^1$ and observe the following,
\begin{IEEEeqnarray*}{lCl}
E^{A_0}_{B_0^2}(x) &=& E^{A_0}_{B_0^2}E^{A_1}_{A_0}(x)\qquad(\mbox{as }B_0^1\subset A_0)\\
&=& E^{A_1}_{B_0^2}(x)\\
&=& E^{A_1}_{B_1^2}E^{A_1}_{A_0}(x)\qquad(\mbox{as }(B_0^2\subset A_0,B_1^2\subset A_1)\mbox{ is a commuting square})\\
&=& E^{A_1}_{B_1^2}(x)\,\,\,\,\,\quad\qquad(\mbox{as }B_0^1\subset A_0)\\
&=& E^{A_1}_{B_1^2}E^{A_1}_{B_1^1}(x)\qquad(\mbox{as }x\in B_0^1\subset B^1_1)\\
&=& E^{A_1}_{C_1}(x)\,\,\,\,\,\quad\qquad(\mbox{as }(C_1\subset B_1^1,B_1^2\subset A_1)\mbox{ is a commuting square})\,.
\end{IEEEeqnarray*}
Since $x\in B_0^1\subset A_0$ and the slice $(C_0\subset A_0,C_1\subset A_1)$ is a commuting square by \Cref{comm cube}, we get that $E^{A_1}_{C_1}(x)\in C_0$. This shows that $E^{A_0}_{B_0^2}(B_0^1)\subset C_0$, which proves that the floor is also a commuting square.\qed
\end{prf}

The converse of the above result does not hold in general. For example, the floor in \Cref{combasic0} is clearly a commuting square, however the roof is not so, as a non-trivial angle in $\mathrm{Ang}_{M_2}\big(\mbox{Ad}_u(\Delta_2),\mbox{Ad}_v(\Delta_2)\big)$ is present. We postpone the proof until \Cref{Sec 4} (see \Cref{p2}), as it contextually better suits there.

\begin{dfn}\label{defccube}
A commuting cube in \Cref{com1} is called a non-degenerate commuting cube if the adjacent faces and the slice are non-degenerate commuting squares.
\end{dfn}

Therefore, the commuting cube in \Cref{combasic0} is not a non-degenerate commuting cube. We request the reader to wait till Sections \ref{Sec 4} and \ref{Sec 5} for (non-trivial) examples of non-degenerate commuting cube.
\smallskip

\begin{figure}[!h]
\begin{center} \resizebox{7 cm}{!}{\includegraphics{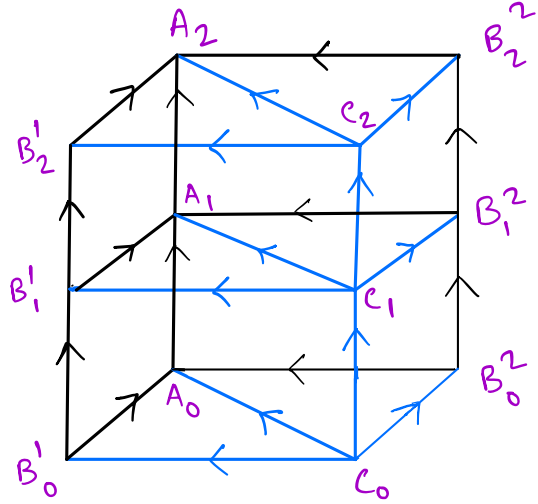}}\end{center}
\caption{Commuting cube and basic construction}\label{combasic1}
\end{figure}

\noindent\textbf{Basic construction of non-degenerate commuting cube~:} Consider the cube in \Cref{com1} and assume that the adjacent faces $(B_0^j\subset A_0,B_1^j\subset A_1)$, for $j=1,2$, are non-degenerate commuting squares. Suppose that the slice $(C_0\subset A_0,C_1\subset A_1)$ is also non-degenerate (see \Cref{main difficulty}). By \Cref{comm cube}, the slice is a commuting square, and such a cube is a non-degenerate commuting cube. Let $A_0\subset A_1\subset^{\,e_2} A_2$ be the basic construction, where $e_2\in A_2$ is the Jones' projection. Then, $C_2=\langle C_1,e_2\rangle$ and $B^j_2=\langle B_j^1,e_2\rangle$ are instances of basic constructions for the inclusions $C_0\subset C_1$ and $B^j_0\subset B^j_1$ respectively for $j=1,2$. Moreover, we obtain a new cube as illustrated in \Cref{combasic1}.
The adjacent faces $(B_1^j\subset A_1,B_2^j\subset A_2)$, for $j=1,2$, and the slice $(C_1\subset A_1,C_2\subset A_2)$ of the {\it upper cube} in \Cref{combasic1} remains non-degenerate commuting squares. Stacking the {\it lower cube} and the {\it upper cube}, we get a cube whose adjacent faces and the resulting slice $(C_0\subset A_0,C_2\subset A_2)$ are non-degenerate commuting squares. Iterating Jones' basic construction, we (possibly) obtain a quadruple of $II_1$ factors $(N\subset P,Q\subset M)$, where $M=\{\cup_{m\geq 0}\,A_m\}^{\prime\prime},\,P=\{\cup_{m\geq 0}\,B^1_m\}^{\prime\prime},\,Q=\{\cup_{m\geq 0}\,B^2_m\}^{\prime\prime}$, and $N=\{\cup_{m\geq 0}\,C_m\}^{\prime\prime}$. However, there is no guarantee that $P\neq Q$ (in \Cref{Sec 3} we shall see an instance of this). If $P\mbox{ and }Q$ are indeed distinct, with a little more effort as in \Cref{bakproc} we can conclude that $N\subset P\cap Q$ (note that $P\cap Q$ is not a factor {\it a priori}). Moreover, $[M:P]=\big|\big|\Lambda^{A_0}_{B^1_0}\big|\big|^2$ (similarly for $Q$) and $[M:N]=\big|\big|\Lambda^{A_0}_{C_0}\big|\big|^2$. By the Ocneanu compactness, we have $N^{\,\prime}\cap M=C_1^{\,\prime}\cap A_0$. Therefore, basic construction of a non-degenerate commuting cube gives rise to (possibly) quadruple of factors. 
\medskip

Instances of this construction will be seen in \Cref{Sec 4} (\Cref{comm cube example}) and \Cref{Sec 5} (\Cref{exmple2}). \Cref{comm cube 2} has far reaching consequences for `two subfactors' $P,Q\subset M$ obtained as iterated basic construction of a non-degenerate commuting cube of finite-dimensional algebras. It says that if there is a non-trivial `angle' present in the cube we start with, which is definitely computable being in a finite-dimensional situation, then the resulting quadruple $(N\subset P,Q\subset M)$ of $II_1$ factors is far from being a commuting square.

\begin{lmma}\label{baksata}
Consider a non-degenerate commuting cube as in \Cref{com1} and define the following operator
\[
S_j:=E^{A_j}_{B^2_j}E^{A_j}_{B^1_j}E^{A_j}_{B^2_j}-E^{A_j}_{C_j}
\]
for $j=0,1$. Suppose that $S_0^{\,2}=\alpha S_0$ for some $\alpha\in\mathbb{R}_+$. Then, $S_1^{\,2}=\alpha S_1$.
\end{lmma}
\begin{prf}
Since the slice is non-degenerate by hypothesis, for any $x\in A_1$, we have $x=\sum a_0c_1$ where $a_0\in A_0$ and $c_1\in C_1$. Then, the following holds :
\begin{IEEEeqnarray*}{lCl}
S_1(x)=\big(E_{B^2_1}E_{B^1_1}E_{B^2_1}-E_{C_1}\big)(x) &=& \big(E_{B^2_1}E_{B^1_1}E_{B^2_1}-E_{C_1}\big)\big(\textstyle\sum a_0c_1\big)\\
&=& E_{B^2_1}E_{B^1_1}\big(\textstyle\sum E_{B^2_1}(a_0)c_1\big)-\textstyle\sum E_{C_1}(a_0)c_1\\
&=& \textstyle\sum\,E_{B^2_1}\big(\big(E_{B^1_1}E_{B^2_1}(a_0)\big)c_1\big)-\textstyle\sum E_{C_1}(a_0)c_1\\
&=& \textstyle\sum\big(E_{B^2_1}E_{B^1_1}E_{B^2_1}-E_{C_1}\big)(a_0)\,c_1\\
&=& \textstyle\sum\big(E_{B^2_0}E_{B^1_0}E_{B^2_0}-E_{C_0}\big)(a_0)\,c_1\,,
\end{IEEEeqnarray*}
where the last line follows since $a_0\in A_0$, and the adjacent faces and the slice are commuting squares. This says that
\[
S_1^{\,2}(x)=\big(E_{B^2_1}E_{B^1_1}E_{B^2_1}-E_{C_1}\big)^2\big(\textstyle\sum a_0c_1\big)=\sum\,\big(E_{B^2_0}E_{B^1_0}E_{B^2_0}-E_{C_1}\big)^2(a_0)\,c_1=\sum\,S_0^{\,2}(a_0)\,c_1\,,
\]
which completes the proof.\qed
\end{prf}

\begin{thm}\label{baksata1}
Suppose that $(N\subset P,Q\subset M)$ is a quadruple of $II_1$ factors obtained as an iterated basic construction of a non-degenerate commuting cube of finite-dimensional algebras in \Cref{com1}. Consider the following nonnegative matrix
\[
S_0:=E^{A_0}_{B^2_0}E^{A_0}_{B^1_0}E^{A_0}_{B^2_0}-E^{A_0}_{C_0}
\]
and suppose that $S_0\neq 0$ with $S_0^{\,2}=\alpha S_0$ for some $\alpha\in\mathbb{R}_+$. Then, $(N\subset P,Q\subset M)$ is not a commuting square, and moreover $\mathrm{Ang}_M(P,Q)$ is the singleton set $\{\mathrm{arccos}\sqrt{\alpha}\}$.
\end{thm}
\begin{prf}
Since $S_0\neq 0$, the quadruple $(C_0\subset B^1_0,B^2_0\subset A_0)$ is not a commuting square. Therefore, none of the roofs, including the final one $(N\subset P,Q\subset M)$, in the iterated basic construction (see \Cref{combasic1}) are commuting squares due to \Cref{comm cube 2}. Since basic construction of non-degenerate commuting square is also so, we get that the resulting final slice $(C_0\subset A_0,N\subset M)$ is a non-degenerate commuting square, and so are the resulting final adjacent faces $(B^1_0\subset A_0,P\subset M)$ and $(B^2_0\subset A_0,Q\subset M)$. Therefore, stacking all these non-degenerate commuting cube, we obtain the final commuting cube which is non-degenerate. By \Cref{baksata}, we obtain that $S_\infty:=E^M_PE^M_QE^M_P-E^M_N$ (and equivalently, $e_Pe_Qe_P-e_N$) satisfies $S_\infty^{\,2}=\alpha S_\infty$. Since $S_\infty\neq 0$ and $\alpha>0$, result now directly follows from Corollary $3.1$ in \cite{SW}.\qed
\end{prf}

In this paper, the notion of `commuting cube' has been crucially used on several occasions throughout Sections \ref{Sec 4} and \ref{Sec 5} (e.g. $6.11,\,6.12,\,6.13,\,6.14,\,6.23,\,7.2,\,7.38,\,7.45$). We end this section with the following question.
\smallskip

\noindent\textbf{Question~:} Under the hypothesis of \Cref{baksata1}, what is the relationship between the spectrum of the angle operators $\Theta(P,Q)$ and $\Theta\big(B^1_0,B^2_0\big)$?


\newsection{Pair of spin model subfactors}\label{pairspin}\label{Sec 3}

A particular instance of `two subfactors' is the pair of spin model subfactors. Recall that for any $n\times n$ complex Hadamard matrix $u$, we have a commuting and co-commuting (non-degenerate commuting) square
\[
\begin{matrix}
A_{10}= \Delta_n &\subset & M_n=A_{11} \cr \cup &\ & \cup\cr A_{00}= \mathbb{C} &\subset &  u\Delta_n u^*=A_{01}\,\,.
\end{matrix}
\]
The basic construction of $\Delta_n\subset M_n$ is $\Delta_n\otimes M_n$ with the Jones' projection $\mbox{bl-diag}\{E_{11},\ldots,E_{nn}\}$, and that of $\Delta_n\otimes M_n$ is $M_n\otimes M_n$ with the Jones' projection $J_n\otimes I_n$, and so on. Here, $J_n=\frac{1}{n}\sum_{i,j=1}^nE_{ij}$. Iterating basic construction we obtain a subfactor $R_u\subset R$, where $R$ is the hyperfinite type $II_1$ factor. This is called a {\it spin model} subfactor. Not much is known about the spin model subfactor, however, the following result is well-known.

\begin{thm}[\cite{Jo3,JS}]
The subfactor $R_u\subset R$ is irreducible and $[R:R_u]=n$.
\end{thm}

Now, if we take another $n\times n$ complex Hadamard matrix $v$ with $u\neq v$, then we obtain a subfactor $R_v\subset R$. Thus, we obtain (possibly) a pair of spin model (irreducible) subfactors
\[
\begin{matrix}
R_u &\subset & R \cr
 & & \cup\cr
 &  & R_v
\end{matrix}\]
Two problems arise here. Firstly, it may happen that $R_u=R_v$ even if $u\neq v$, in which case we fail to obtain a `pair' of subfactors. Secondly, for $R_u\neq R_v$ it is not immediate that $R_u\cap R_v$ is a factor, that is, we do not have a quadruple of factors {\em a priori}. To overcome the first difficulty, that is to obtain ``two subfactors" of $R$, one has to first characterize $R_u=R_v$ in terms of $u\mbox{ and }v$, and then proceed. We have a complete answer for this.
\smallskip

Recall that a matrix in $GL_n(\bbc)$ is said to be monomial (Generalized permutation matrix) if each row and column has exactly one non-zero entry. Let $N$ denote the set of all monomial matrices. It is known that $N$ is a subgroup of $GL_n(\bbc)$ and $N=\mathcal{N}_{GL_n}(\Delta_n)$, the normalizer of $\Delta_n$ in $GL_n(\bbc)$. Consider the subgroup $\mathbb{G}:=N\cap U(n)$ of $GL_n(\bbc)$, where $U(n)$ is the set of all unitary matrices.
\smallskip

\textbf{An equivalence relation~:} For two unitary matrices $A$ and $B$ in $M_n$, consider the equivalence relation $A\sim B$ if $B=AC$, where $C\in\mathbb{G}$. Recall that two complex Hadamard matrices $H_1,H_2$ are called Hadamard equivalent, denoted by $H_1\simeq H_2$, if $H_2=D_1P_1H_1P_2D_2$, where $D_1,D_2$ are unitary diagonal matrices and $P_1,P_2$ are permutation matrices. Hence, the equivalence relation `$\sim$' is a sub-equivalence relation of the Hadamard equivalence `$\simeq$'. This is because for any $C\in N$, there exists an invertible $D\in\Delta_n$ and a permutation matrix $P$ such that $C=PD$. Hence, $C\in\mathbb{G}$ forces $D$ to become a unitary diagonal matrix. It turns out that it is this finer equivalence relation between the complex Hadamard matrices $u$ and $v$ of order $n\times n$ that completely characterize $R_u\neq R_v$. More precisely, we have the following result.

\begin{thm}\label{general thm}
\begin{enumerate}[$(i)$]
\item For distinct $n\times n$ complex Hadamard matrices $u\mbox{ and }v$, the pair of spin model subfactors $R_u\mbox{ and }R_v$ of the hyperfinite type $II_1$ factor $R$ are distinct $(R_u\neq R_v)$ if and only if $\,u\nsim v$.
\item If two $n\times n$ complex Hadamard matrices $u\mbox{ and }v$ are Hadamard inequivalent, then the corresponding spin model subfactors $R_u\mbox{ and }R_v$ of $R$ are always distinct $(R_u\neq R_v)$.
\end{enumerate}
\end{thm}

The proof of \Cref{general thm} requires a bit of work. Although in part $(i)$, the `only if' part is easy, it is the `if' part that requires a bit of work. We will be using the Pimsner-Popa probabilistic constant as a major tool to prove it.

Let $\Delta_n$ denote the diagonal subalgebra of $M_n$ and $\,U\in M_n$ be a unitary matrix. Then, $\Delta_n$ is a Masa (maximal abelian self-adjoint subalgebra) in $M_n$, and in fact any Masa in $M_n$ is of the form $U\Delta_n U^*$ for some unitary $U\in M_n$. Consider the following unital inclusion of subalgebras \[\begin{matrix}
\Delta_n &\subset & M_n\\
 & & \cup\\
 & & U\Delta_n U^*\\
\end{matrix}\]
such that $U$ is not the identity matrix. We are going to determine the value of $\,\lambda(\Delta_n\,,\,U\Delta_n U^*)$. First note that $\lambda(M_n,U\Delta_nU^*)=\lambda(M_n,\Delta_n)=\frac{1}{n}$ (example in $6.5$, Page $94$ in \cite{PP}). That is, $E_{U\Delta_nU^*}(x)\geq \frac{1}{n}x$ for all $x\in(M_n)_+$, and hence $\,\lambda(\Delta_n\,,\,U\Delta_n U^*)\geq \frac{1}{n}$. We shall show that the best constant is determined in terms of the Hamming numbers of the rows of $U$.
\smallskip

The unique trace preserving conditional expectation $\,E:M_n\longrightarrow\Delta_n$ is given by $A\longmapsto (a_{ii})_{ii}$ for any $A=(a_{ij})_{1\leq i,j\leq n}\in M_n$. Here, the notation $(\,.\,)_{ii}$ denotes a diagonal matrix in $M_n$. We are going to use the following well-known facts.
\medskip

\textbf{Weyl's inequality for eigenvalues~:} Let $A\mbox{ and }B$ be $n\times n$ self-adjoint matrices and $\widetilde{A}=A+B$. Suppose that the eigenvalues of $A$ are ordered as $\lambda_1(A)\leq\lambda_2(A)\leq\ldots\leq\lambda_n(A)$, and similarly for the matrices $\widetilde{A}\mbox{ and }B$. Then, one has the following,
\[\lambda_k(A)+\lambda_1(B)\leq\lambda_k(\widetilde{A})\leq\lambda_k(A)+\lambda_n(B)\,.\]

\textbf{Matrix determinant lemma~:} Suppose $A$ is an invertible $n\times n$ matrix and $\mathbf{u},\mathbf{v}$ are two vectors in $\bbc^n$ written in column matrices. Then,
\[
\det\big(A+\mathbf{u}\mathbf{v}^{\textsf{T}}\big)=\big(1+\mathbf{v}^{\textsf {T}}A^{-1}\mathbf{u}\big)\,\det(A)\,.
\]

\noindent By the definition of $\,\lambda$, to find $\lambda(\Delta_n,U\Delta_n U^*)$, we need to find the best constant $t\in[0,1]$ in the following inequality,
\[
U(U^*xU)_{kk}U^*\geq tx\,,\quad\forall\,\,x\in(\Delta_n)_+\,.
\]
This is because $E^{M_n}_{U\Delta_nU^*}=\mbox{Ad}_U\circ E^{M_n}_{\Delta_n}\circ \mbox{Ad}_{U^*}$.

\begin{lmma}\label{minimum}
For $1\leq i\leq n$, consider the minimal projections $p_i:=E_{ii}$ in $M_n$. Then,
\begin{IEEEeqnarray*}{lCl}
\lambda(\Delta\,,\,U\Delta U^*)=\min_{1\leq i\leq n}\,\mathrm{sup}\,\big\{t\geq 0:(U^*p_iU)_{kk}\geq t\,U^*p_iU\big\}\,.
\end{IEEEeqnarray*}
\end{lmma}
\begin{prf}
First observe that $\Delta_n$ is generated by the minimal projections $p_i=E_{ii}$ in $M_n$, and $U(U^*p_iU)_{kk}U^*\geq tp_i$ holds if and only if $(U^*p_iU)_{kk}\geq t\,U^*p_iU$ holds. Let us denote
\[
\widetilde{\lambda}_i:=\sup_{t\geq 0}\,\big\{(U^*p_iU)_{kk}\geq t\,U^*p_iU\big\}\,,
\]
for $\,i\in\{1,\ldots,n\}$. It is obvious from the definition of $\lambda$ that $\lambda(\Delta_n\,,\,U\Delta_nU^*)\leq\widetilde{\lambda}_i$ for each $i$. Thus, we get that
\[
\lambda(\Delta\,,\,U\Delta U^*)\leq\min_{1\leq i\leq n}\,\widetilde{\lambda}_i\,.
\]
Conversely, for any $x\in(\Delta_n)_+$ we have $x=\sum_{i=1}^n\,\alpha_ip_i$ with all $\alpha_i\geq 0$. Thus, for any $t\in[0,1],\,(U^*xU)_{kk}\geq t\,U^*xU$ holds if and only if the following inequality 
\begin{IEEEeqnarray}{lCl}\label{required here}
\sum_{i=1}^n\,\alpha_i(U^*p_iU)_{kk} &\geq& \sum_{i=1}^n\,\alpha_it\,U^*p_iU
\end{IEEEeqnarray}
holds, and the best constant $t$ satisfying the above inequality gives us $\lambda(\Delta\,,\,U\Delta U^*)$. Now, if we set $t_0=\mbox{min}_{1\leq i\leq n}\,\widetilde{\lambda}_i$, then by the definition of $\widetilde{\lambda}_i$ we get that
\[
(U^*p_iU)_{kk}\geq t_0\,U^*p_iU
\]
for all $i$, and consequently \Cref{required here} is satisfied for $t=t_0$. By the definition of $\lambda$, we get that
\[
\lambda(\Delta\,,\,U\Delta U^*)\geq t_0=\min_{1\leq i\leq n}\,\widetilde{\lambda}_i\,,
\]
which completes the proof.\qed
\end{prf}

Let $\sigma=-t$ and for any fixed $i\in\{1,\ldots,n\}$ consider the following one-parameter family of $n\times n$ matrices
\[
\{A_\sigma:=(U^*p_iU)_{kk}+\sigma\,U^*p_iU\,:\,\sigma\in[-1,0]\}\,.
\]
Since $U^*p_iU$ is a rank one matrix, any member of this family is a matrix of the form $D_\sigma+\sigma\,\mathbf{u}\mathbf{u}^*$, where $\mathbf{u}\in\mathbb{C}^n$ is a vector written as a column matrix and $\mathbf{u}^*$ denotes the complex conjugate of the vector $\mathbf{u}$ written as a row matrix. Thus, any matrix $A_\sigma$ is a rank-$1$ perturbation of a diagonal matrix $D_\sigma$. Throughout this subsection, we denote $U^*=(u_{ij}),\,1\leq i,j\leq n$. We are interested in the following inequality
\[(U^*p_iU)_{kk}+\sigma U^*p_iU\geq 0\,,\]
that is,
\begin{IEEEeqnarray}{lCl}\label{best constant matrix}
\mbox{diag}\big(|u_{1i}|^2,\ldots,|u_{ni}|^2\big)+\sigma\,\left[{\begin{matrix}
u_{1i}\\
\vdots\\
u_{ni}\\
\end{matrix}}\right]
\begin{matrix}
\left[{\begin{matrix}
\overline{u_{1i}} & \ldots & \overline{u_{ni}}\\
\end{matrix}}\right]\\[1.0ex] \\ \\ \end{matrix}\,\,\geq 0\,.
\end{IEEEeqnarray}

\begin{lmma}\label{except minimum}
Except possibly the minimum eigenvalue of the following self-adjoint matrix
\[
\mbox{diag}\big(|u_{1i}|^2,\ldots,|u_{ni}|^2\big)+\sigma\,\left[{\begin{matrix}
u_{1i}\\
\vdots\\
u_{ni}\\
\end{matrix}}\right]\begin{matrix}\left[{\begin{matrix}
\overline{u_{1i}} & \ldots & \overline{u_{ni}}\\
\end{matrix}}\right]\\[1.0ex] \\ \\ \end{matrix}\]
all the other $(n-1)$ eigenvalues, counted with multiplicities, are non-negative.
\end{lmma}
\begin{prf}
We write the matrix $A_\sigma$ as $P+D$ where $D$ is the diagonal matrix $\mbox{diag}\big(|u_{1i}|^2,\ldots,|u_{ni}|^2\big)$, and $P$ is the rank one perturbation $A_\sigma-D$. Since $\frac{1}{\sigma}P$ is a rank-$1$ projection, we have the spectrum of $P$ is given by $\sigma(P)=\{0,\sigma\}$. We write the eigenvalues in increasing order as the following,
\[
\lambda_1=\sigma\leq \lambda_2=\lambda_3=\ldots=\lambda_n=0\,.
\]
Similarly, we write the eigenvalues of $D$ in increasing order as the following
\[
0\leq\mu_1\leq\mu_2\leq\ldots\leq\mu_n\,,
\]
where we take $\mu_k=|u_{ki}|^2$ after possible renaming. Now, let the eigenvalues of $A_\sigma$ be $\widetilde{\lambda}_1\leq\ldots\widetilde{\lambda}_n$. By the Weyl's inequality, we get the following
\[
\lambda_k+\mu_1\leq\widetilde{\lambda}_k\leq\lambda_k+\mu_n\,.
\]
Putting $k=2$, we get that $\widetilde{\lambda}_2\geq\lambda_2+\mu_1=\mu_1\geq 0$. Thus, $\widetilde{\lambda}_k\geq 0$ for all $k=2,\ldots,n$, and this completes the proof.\qed
\end{prf}

\begin{lmma}\label{reduced block}
If all the entries of the unitary matrix $U$ are non-zero, then $\lambda\big(\Delta\,,\,U\Delta U^*\big)=\frac{1}{n}\,$.
\end{lmma}
\begin{prf}
Fix any $i\in\{1,\ldots,n\}$, and let $\mathbf{v}^T=\big(\overline{u_{1i}},\ldots,\overline{u_{ni}}\big)$ and $\mathbf{u}^T=\big(\sigma\,u_{1i},\ldots,\sigma\,u_{ni}\big)\,.$ Since, $u_{ij}\neq 0$ is given for all $i,j\in\{1,\ldots,n\}$, using the matrix determinant lemma, we get the following,
\begin{IEEEeqnarray*}{lCl}
&  & \mbox{det}\left(\mbox{diag}\big(|u_{1i}|^2,\ldots,|u_{ni}|^2\big)+\sigma\,\left[{\begin{matrix}
u_{1i}\\
\vdots\\
u_{ni}\\
\end{matrix}}\right]\begin{matrix}\left[{\begin{matrix}
\overline{u_{1i}} & \ldots & \overline{u_{ni}}\\
\end{matrix}}\right]\\[1.0ex] \\ \\ \end{matrix}\right)\\
&=& \left(1+v^T\,\mbox{diag}\Big(\frac{1}{|u_{1i}|^2},\ldots,\frac{1}{|u_{ni}|^2}\Big)\,u\right)\,\mbox{det}\,\Big(\mbox{diag}\big(|u_{1i}|^2,\,\ldots\,,|u_{ni}|^2\big)\Big)\\
&=& \prod_{j=1}^n\,|u_{ji}|^2(1+n\sigma)\,.
\end{IEEEeqnarray*}
The quantity $\prod_{j=1}^n\,|u_{ji}|^2(1+n\sigma)$ is non-negative if and only if $(1+n\sigma)$ is non-negative. The matrix in \Cref{best constant matrix} is positive semi definite if and only if all its eigenvalues are non-negative. Since, all the $(n-1)$ eigenvalues, counted with multiplicities, are non-negative by \Cref{except minimum}, we get that  the best constant $\sigma\in[-1,0]$ in \Cref{best constant matrix} is determined by the inequality $1+n\sigma\geq 0$, since determinant is product of the eigenvalues. That is, the best constant $\sigma$ is given by $\sigma=-1/n\,,$ and consequently the best constant $t$ in the Pimsner-Popa constant is given by $t=1/n$. Since this happens for any $1\leq i\leq n$, by Lemma \ref{minimum} we conclude the proof.\qed
\end{prf}

\begin{dfn}\label{hamdef}
Given a nonzero vector $\mathbf{u}\in\bbc^n$, the Hamming number is given by,
\[
h(\mathbf{u}):=\mbox{number of non-zero entries in }\mathbf{u}\,.
\]
\end{dfn}

\begin{thm}\label{finitemains}
If $\Delta_n$ and $U\Delta_n U^*$ are two Masas in $M_n$, where $U\in U(n)$, then
the Pimsner-Popa constant is given by the following,
\[
\lambda(\Delta_n\,,\,U\Delta_n U^*)=\min_{1\leq i\leq n}\,\big(\mbox{$h\left(U^*\right)_i$}\big)^{-1}
\]
where $\left(U^*\right)_i$ is the $i$-th column of $U^*$.
\end{thm}
\begin{prf}
Consider the $i$-th column of $U^*$ for $1\leq i\leq n$. Let $k_i$ be the number of zero elements in the set $\{u_{1i},\ldots,u_{ni}\}$. That is, $k_i$ many eigenvalues, counted with multiplicities, of the diagonal matrix $\,\mbox{diag}\big(|u_{1i}|^2,\ldots,|u_{ni}|^2\big)$ are $0$. Let $u_{j_1i}=\ldots=u_{j_{k_i}i}=0$. Then, all the $j_\ell$-th row for $1\leq\ell\leq k_i$ of the following $n\times n$ matrix
\[
\left[{\begin{matrix}
u_{1i}\\
\vdots\\
u_{ni}\\
\end{matrix}}\right]\begin{matrix}\left[{\begin{matrix}
\overline{u_{1i}} & \ldots & \overline{u_{ni}}\\
\end{matrix}}\right]\\[1.0ex] \\ \\ \end{matrix}
\]
are identically zero. Since $\overline{u_{j_1i}}=\ldots=\overline{u_{j_{k_i}i}}=0$, all the $j_\ell$-th column for $1\leq\ell\leq k$ of the above matrix are also identically zero. To obtain the Pimsner-Popa constant, we are interested in the best constant $\sigma\in[-1,0]$ for which Eqn. $\ref{best constant matrix}$ holds. That is, all the eigenvalues of the following matrix
\[
A_\sigma=\mbox{diag}\big(|u_{1i}|^2,\ldots,|u_{ni}|^2\big)+\sigma\,\left[{\begin{matrix}
u_{1i}\\
\vdots\\
u_{ni}\\
\end{matrix}}\right]\begin{matrix}\left[{\begin{matrix}
\overline{u_{1i}} & \ldots & \overline{u_{ni}}\\
\end{matrix}}\right]\\[1.0ex] \\ \\ \end{matrix}
\]
must be non-negative. Since swapping two rows or columns only changes the sign of the determinant, the characteristic equation of the matrix $A_\sigma$ becomes the following,
\begin{IEEEeqnarray*}{lCl}
0 &=& \mbox{det }(A_\sigma-xI_n)\\
&=& \mbox{det }\big(\mbox{bl-diag }(-x,\ldots,-x,B_\sigma)\big)
\end{IEEEeqnarray*}
where $B_\sigma$ is a $(n-k_i)\times (n-k_i)$ matrix. This says that positive semi-definiteness of $A_\sigma$ is completely determined by positive semi-definiteness of $B_\sigma$. Now, observe that the matrix $B_\sigma$ is as in the earlier situation of \Cref{reduced block}, where all the entries of the following matrix
\[
\mbox{diag}\big(\,|u_{j_{k_i+1}i}|^2,\ldots,|u_{j_ni}|^2\big)
\]
are non-zero. Hence, the best constant $\sigma$ in \Cref{best constant matrix} for the matrix $B_\sigma$ is given by $\sigma=-\frac{1}{n-k_i}$. The number $n-k_i$ is exactly the number of nonzero entries in the $i$-th column of $U^*$. By Lemma \ref{minimum}, the proof is now completed.\qed
\end{prf}

\begin{crlre}
If $U$ is a unitary matrix in $M_n$ with $U=(u_{ij})$, then $\sum_{i,j=1}^n\eta(|u_{ij}|^2)\leq\log\big(\min_{1\leq i\leq n}\,h(U)_i\big)^n$.
\end{crlre}
\begin{prf}
Let $U=(u_{ij})$ be a unitary matrix in $M_n$ and consider the pair of Masas $(\Delta_n,U\Delta_n U^*)$ in $M_n$. By Cor. $3.2$ in \cite{choda}, $h(\Delta_n|U\Delta_n U^*)=\frac{1}{n}\sum_{i,j=1}^n\eta(|u_{ij}|^2)$. By \Cref{upbound}, we have $h(\Delta_n|U\Delta_n U^*)\leq-\log\lambda(\Delta_n,U\Delta_n U^*)$ and the claim now follows from \Cref{finitemains}.\qed
\end{prf}

\begin{crlre}
The pair $\big(\Delta_n,U\Delta_nU^*\big)$ in the type $I_n$ factor $M_n$ is an orthogonal pair in the sense of Popa $\cite{P}$ only if the unitary $U$ has no zero entry.
\end{crlre}
\begin{prf}
Let $\big(\Delta_n,U\Delta_nU^*\big)$ be an orthogonal pair. Then, $(\bbc\subset \Delta_n,U\Delta_n U^*\subset M_n)$ is a commuting square. By \Cref{df}, $h(\Delta_n|U\Delta_n U^*)=H(\Delta_n|U\Delta_n U^*)=H(\Delta_n|\bbc)=\log n$ (apply Theorem $6.2$ in \cite{PP}). Now by Cor. $3.3$ in \cite{choda}, we know that $\big(\Delta_n,U\Delta_nU^*\big)$ is an orthogonal pair in $M_n$ if and only if $h(\Delta_n|U\Delta_n U^*)=\log n$. Since $h(\Delta_n|U\Delta_n U^*)\leq-\log\lambda(\Delta_n,U\Delta_n U^*)$ by \Cref{upbound}, we have the following inequalities
\[
\log n=h(\Delta_n|U\Delta_n U^*)\leq\log\big(\min_{1\leq i\leq n}\,h(U)_i\big)\leq\log n
\]
by \Cref{finitemains}. Hence, all the entries of the unitary $U$ must be non-zero.\qed
\end{prf}
\medskip

\textbf{Proof of \Cref{general thm}:~} For part $(i)$, we prove the contrapositive, that is, $R_u=R_v$ if and only if $u\sim v$. First observe that $R_u=R_v$ if and only if $\lambda(R_u,R_v)=1$. This is because $\lambda(R_u,R_v)=1$ implies that $R_u\subseteq R_v$ by \Cref{sk}. Since $[R:R_u]=[R:R_v]=n$, we get that $R_u=R_v$. The converse direction is obvious.

Now, assume that $R_u=R_v$. Then $\lambda(R_u,R_v)=1$, and consequently $\lambda(u\Delta_2u^*,v\Delta_2v^*)=1$ because $\lambda(R_u,R_v)$ is the limit of the decreasing sequence of $\lambda$ at each step of the tower of basic constructions (\Cref{popaadaptation}) and $0\leq\lambda\leq 1$. By \Cref{finitemains} we get that $u^*v$ must be a permutation matrix, that is $u\sim v$. Thus, we conclude that $R_u=R_v$ implies $u\sim v$. Conversely, suppose that $u\sim v$. Then $v=uPD$, where $P$ is a permutation matrix in $M_n$ and $D\in U(n)$ is a diagonal unitary matrix. Then $\mbox{Ad}_{u^*v}(\Delta_n)\subseteq\Delta_n$, and consequently $\mbox{Ad}_{v}(\Delta_n)\subseteq\mbox{Ad}_u(\Delta_n)$. Immediately, we get that $R_v\subseteq R_u$. Since $[R:R_u]=[R:R_v]=n$, we get that $[R_u:R_v]=1$, and hence $R_u=R_v$. This completes part $(i)$, and part $(ii)$ immediately follows from part $(i)$, since the equivalence relation `$\sim$' is finer than the Hadamard equivalence relation.\qed

Notice that the `only if' part in \Cref{general thm} is very easy, and it is the `if' part that requires work. \Cref{popaadaptation} and \Cref{finitemains} are the crucial ingredients in proving the `if' part.


\newsection{The Pimsner-Popa constant for a pair of type \texorpdfstring{$I_n$}~ factors}\label{Sec 3.2}

The results in this section is a computational tool and key ingredient to compute the Pimsner-Popa constant in Sections \ref{Sec 4} and \ref{Sec 5}. Here, we determine the Pimsner-Popa constant for the following situation namely,
\[\begin{matrix}
M_n &\subset & M_n\oplus M_n\\
 & & \cup\\
 & & UM_nU^*\\
\end{matrix}\]
where $U\in M_n\oplus M_n$ is a unitary matrix.

Henceforth, we write $M_n\oplus M_n$ as $\Delta_2\otimes M_n$, where $\Delta_2$ denotes the diagonal subalgebra (Masa) of $M_2$. The embedding $M_n\xhookrightarrow{} \Delta_2\otimes M_n$ is the diagonal embedding $x\mapsto I_2\otimes x$. Since $U$ is a unitary in $\Delta_2\otimes M_n$, we have $U=\mbox{bl-diag}\{U_1,U_2\}$ with $U_1,\,U_2\in M_n$ unitary matrices. The unique trace preserving conditional expectation $E:\Delta_2\otimes M_n\longrightarrow M_n$ is given by $\tau\otimes\mbox{id}$, where $\tau:\Delta_2\longrightarrow\bbc$ is the normalized trace induced from $M_2$. For any $x\in (M_n)_+$, we have
\begin{IEEEeqnarray*}{lCl}
E_{UM_nU^*}(x) &=& UE_{M_n}\big(U^*(1\otimes x)U\big)U^*\\
&=& UE_{M_n}\left(\mbox{bl-diag}\{U_1^*xU_1,U_2^*xU_2\}\right)U^*\\
&=& \frac{1}{2}U(U_1^*xU_1+U_2^*xU_2)U^*
\end{IEEEeqnarray*}
Hence, for $0\leq t\leq 1$ the inquality $E_{UM_nU^*}(x)\geq tx$ is equivalent to the following inequality,
\begin{IEEEeqnarray}{lCl}\label{lower bound0}
&  & \big(\frac{1}{2}-t\big)\left[{\begin{matrix}
U_1^*xU_1 & 0\\
0 & U_2^*xU_2\\
\end{matrix}}\right]+\frac{1}{2}\left[{\begin{matrix}
U_2^*xU_2 & 0\\
0 & U_1^*xU_1\\
\end{matrix}}\right]\geq 0\,.
\end{IEEEeqnarray}

\begin{lmma}\label{to be refered here}
We have $E_{UM_nU^*}(x)\geq\frac{1}{2}x$ for all $x\in (M_n)_+$, that is, $\lambda\big(M_n\,,\,UM_nU^*\big)\geq\frac{1}{2}$.
\end{lmma}
\begin{prf}
Follows immediately from \Cref{lower bound0}.\qed
\end{prf}

\begin{crlre}
For $n\geq 3$ and any unitary $U\in M_n\oplus M_n$, the pair $\big(M_n\,,UM_nU^*\big)$ is never an orthogonal pair in $M_n\oplus M_n$ in the sense of Popa $\cite{P}$.
\end{crlre}
\begin{prf}
On contrary assume that $\big(M_n\,,\,UM_nU^*\big)$ is an orthogonal pair. Then the quadruple $(\bbc\subset M_n\,,\,UM_nU^*\subset M_n\oplus M_n)$ is a commuting square. By \Cref{df}, we get that $h\big(M_n|UM_nU^*\big)=H\big(M_n|UM_nU^*\big)=H(M_n|\bbc)=\log\,n$ (see Section $6$ in \cite{PP}). By \Cref{to be refered here}, we have $\lambda\big(M_n\,,\,UM_nU^*\big)^{-1}\leq 2$, and hence $\log n=h\big(M_n|UM_nU^*\big)\leq\log\,2$ using \Cref{upbound}. This contradicts the fact that $n\geq 3$.\qed
\end{prf}

\begin{lmma}\label{ref in thm 1}
If $\,U_1^*U_2$ is not a diagonal matrix in $M_n$, then $\,\lambda\big(M_n\,,\,UM_nU^*\big)=\frac{1}{2}\,.$
\end{lmma}
\begin{prf}
Suppose that
\[
\frac{1}{2}\,\left[{\begin{matrix}
U_1^*xU_1 & 0\\
0 & U_2^*xU_2\\
\end{matrix}}\right]+\frac{1}{2}\,\left[{\begin{matrix}
U_2^*xU_2 & 0\\
0 & U_1^*xU_1\\
\end{matrix}}\right]\geq \big(\frac{1}{2}+\varepsilon\big)\,\left[{\begin{matrix}
U_1^*xU_1 & 0\\
0 & U_2^*xU_2\\
\end{matrix}}\right]
\]
for all $x\in (M_n)_+$ and $\varepsilon>0$. That is, we have the following,
\begin{IEEEeqnarray}{lCl}\label{inequality 1}
\left[{\begin{matrix}
U_2^*xU_2 & 0\\
0 & U_1^*xU_1\\
\end{matrix}}\right] &\geq& 2\varepsilon\,\left[{\begin{matrix}
U_1^*xU_1 & 0\\
0 & U_2^*xU_2\\
\end{matrix}}\right]\quad\forall\,\, x\in (M_n)_+
\end{IEEEeqnarray}
Let $V=U_1^*U_2$, a unitary matrix in $M_n$, and choose $x_j=U_1E_{jj}U_1^*$ for $1\leq j\leq n$. Each $x_j\in (M_n)_+$ and we have for each $1\leq j\leq n$ the following,
\[
\left[{\begin{matrix}
V^*E_{jj}V & 0\\
0 & E_{jj}\\
\end{matrix}}\right]\geq 2\varepsilon\left[{\begin{matrix}
E_{jj} & 0\\
0 & V^*E_{jj}V\\
\end{matrix}}\right]
\]
in $\Delta_2\otimes M_n$. Hence, for any $i\neq j$ in $\{1,\ldots,n\}$ we have
\[\left[{\begin{matrix}
0 & 0\\
0 & E_{1i}\\
\end{matrix}}\right]\left[{\begin{matrix}
V^*E_{jj}V & 0\\
0 & E_{jj}\\
\end{matrix}}\right]\left[{\begin{matrix}
0 & 0\\
0 & E_{i1}\\
\end{matrix}}\right]\geq 2\varepsilon\left[{\begin{matrix}
0 & 0\\
0 & E_{1i}\\
\end{matrix}}\right]\left[{\begin{matrix}
E_{jj} & 0\\
0 & V^*E_{jj}V\\
\end{matrix}}\right]\left[{\begin{matrix}
0 & 0\\
0 & E_{i1}\\
\end{matrix}}\right]\]
as in a $C^*$-algebra $a\geq b$ implies that $x^*ax\geq x^*bx$ for any $x$. Since $\,\varepsilon>0$, letting $V=(v_{pq})_{1\leq p,q\leq n}$ we get the following,
\begin{IEEEeqnarray*}{lCl}
0 &=& E_{1i}V^*E_{jj}VE_{i1}\\
&=& E_{1i}\big(\sum_{p,q}\,\overline{v_{qp}}E_{pq}\big)E_{jj}\big(\sum_{r,s}\,v_{rs}E_{rs}\big)E_{i1}\\
&=& E_{1i}\big(\sum_{p,s}\,\overline{v_{jp}}v_{js}E_{ps}\big)E_{i1}\\
&=& \overline{v_{ji}}v_{ji}E_{11}\,.
\end{IEEEeqnarray*}
This says that $\,v_{ji}=0$ for any $j\neq i$ i,e. $V$ is a diagonal matrix in $M_n$. This is a contradiction, and we get that $\varepsilon=0$.\qed
\end{prf}

\begin{lmma}\label{ref in thm 2}
One has $\lambda\big(M_n\,,\,UM_nU^*\big)=1$ if and only if $U_1^*U_2=\alpha I_n$ for some $\alpha\in\mathbb{S}^1$.
\end{lmma}
\begin{prf}
If $\,U_1^*U_2=\alpha I_n$, then $UM_nU^*=I_2\otimes U_1M_nU_1^*=I_2\otimes M_n$, and hence $\lambda\big(M_n\,,\,UM_nU^*\big)=1$. Conversely, $\lambda\big(M_n\,,\,UM_nU^*\big)=1$ implies that $H(M_n|UM_nU^*)=0$ and consequently, $M_n\subseteq UM_nU^*$, that is, $U^*M_nU\subseteq M_n$. Therefore, $U_1^*AU_1=U_2^*AU_2$ for all $A\in M_n$ and hence $U_1U_2^*=\beta I_n$ for some $\beta\in\bbc$, which completes the proof.
\end{prf}

\begin{thm}\label{2nd lambda}
For the inclusion of algebras $M_n\subset\Delta_2\otimes M_n\supset UM_nU^*$, where $\,U\in\Delta_2\otimes M_n$ is a unitary matrix given by $\,U=\mbox{bl-diag}\{U_1,U_2\}$ with unitary matrices $U_1,U_2\in M_n$, one has the following,
\begin{IEEEeqnarray*}{rCl}
 \lambda\big(M_n\,,\,UM_nU^*\big)=\begin{cases}
                            1 & \mbox{ iff } \,\,U_1^*U_2 \mbox{ is a scalar matrix}, \cr
                            \frac{1}{2} & \mbox{ if } \,\,U_1^*U_2 \mbox{ is not a diagonal matrix}. \cr
                           \end{cases}
\end{IEEEeqnarray*}
\end{thm}
\begin{prf}
Follows from Lemma (\ref{ref in thm 1}, \ref{ref in thm 2}).\qed
\end{prf}

Note that there is one remaining situation when $U_1^*U_2$ is a diagonal unitary matrix. However, this situation does not arise in our case in Sections \ref{Sec 4} and \ref{Sec 5}, and hence we leave it to the interested readers.


\newsection{Subfactors arising from \texorpdfstring{$2\times 2$}~ complex Hadamard matrices}\label{Sec 4}

Recall that $2\times 2$ complex Hadamard matrices form a single family up to Hadamard equivalence, and any element of this family is given by the following
\begin{IEEEeqnarray}{lCl}\label{000001}
F^{\alpha_1}_{\alpha_2}(\alpha_3):=\frac{1}{\sqrt{2}}\begin{bmatrix}
e^{i\alpha_1} & e^{i (\alpha_1+\alpha_3)}\\
e^{i\alpha_2} & -e^{i(\alpha_2+\alpha_3)}
\end{bmatrix}
\end{IEEEeqnarray}
where $\alpha_j\in[0,2\pi),\,j=1,2,3$. Also recall the equivalence relation `$\sim$' from \Cref{general thm}, \Cref{Sec 3}. Consider any two $2\times 2$ complex Hadamard matrices $u$ and $v$ such that $u\nsim v$ and obtain the (distinct) irreducible spin model subfactors each with index 2~:
\[\begin{matrix}
R_u &\subset & R\\
 & & \cup\\
 & & R_v
\end{matrix}
\]
Note that $u\nsim v$ if and only if $v\neq uPD$, where $D$ is a diagonal unitary matrix and $P$ is a permutation matrix. In this $2\times 2$ situation, we have $PD$ is a unitary matrix of the form $\left[{\begin{matrix}
\alpha & 0\\
0 & \beta\\
\end{matrix}}\right]$ or $\left[{\begin{matrix}
0 & \alpha\\
\beta & 0\\
\end{matrix}}\right],$ where $\alpha,\beta\in\mathbb{S}^1$. Given any $u=F^{\alpha_1}_{\alpha_2}(\alpha_3)$ and $v=F^{\beta_1}_{\beta_2}(\beta_3)$ (\Cref{000001}), consider the following matrices
\[
\widetilde{u}=\frac{1}{\sqrt{2}}\begin{bmatrix}
1 &  1\\ 
e^{i(\alpha_2-\alpha_1)} & -e^{i(\alpha_2-\alpha_1)}
\end{bmatrix}\quad\mbox{and}\quad\widetilde{v}=\frac{1}{\sqrt{2}}\begin{bmatrix}
1 &  1\\ 
e^{i(\beta_2-\beta_1)} & -e^{i(\beta_2-\beta_1)}
\end{bmatrix}
\]
associated to them. Since $u\sim\widetilde{u}$, by \Cref{general thm} we have $R_u=R_{\widetilde{u}}$, and similarly for $v$. Hence, the quadruples of von Neumann algebras $(R_u\cap R_v\subset R_u,R_v\subset R)$ and $(R_{\widetilde{u}}\cap R_{\widetilde{v}}\subset R_{\widetilde{u}},R_{\widetilde{v}}\subset R)$ are the same. Therefore, without loss of generality, throughout the \Cref{Sec 4} we will assume that
\begin{IEEEeqnarray}{lCl}\label{reduc}
u=\frac{1}{\sqrt{2}} 
\begin{bmatrix}
1 & 1\\
e^{i\alpha} & -e^{i\alpha}
\end{bmatrix}\quad\mbox{and}\quad v=\frac{1}{\sqrt{2}}
\begin{bmatrix}
1 & 1\\
e^{i\beta} & -e^{i\beta}
\end{bmatrix}
\end{IEEEeqnarray}
for $\alpha,\beta\in[0,2\pi)$. Note that $\beta=\alpha\pm\pi$ implies that $v=u\sigma_1$, which gives $u\sim v$. Therefore, in our situation we have $\alpha,\beta\in[0,2\pi)$ with $\beta\neq\alpha,\alpha\pm\pi$.

In this section, we first compute the Pimsner-Popa constant $\,\lambda(R_u,R_v)$. Recall that as mentioned before, it is not obvious whether intersection of factors is again a factor. We prove that $R_u\cap R_v$ is a $II_1$ subfactor of $R$ with $[R:R_u\cap R_v]=4$. Moreover, we show that $R_u\cap R_v\subset R$ is a vertex model subfactor. Finally, we compute the Sano-Watatani angle $\mbox{Ang}_R(R_u,R_v)$ and relative entropy $h(R_u|R_v)$ between the subfactors $R_u\mbox{ and }R_v$, and obtain legitimate bounds for the relative entropy $H(R_u|R_v)$. We also characterize the quadruple of $II_1$ factors $(R_u\cap R_v\subset R_u,R_v\subset R)$ in terms of a bi-unitary matrix in $M_4$.
\smallskip

\begin{figure}[!h]\begin{center}
		\psfrag{$R$}{ \large  $R$}
		\psfrag{d1}{\large$\Delta_2\otimes M_2\otimes M_2$}
		\psfrag{d2}{\large$M_2\otimes M_2$}
		\psfrag{d3}{\large$\Delta_2\otimes M_2$}
		\psfrag{d4}{\large$M_2$}
		\psfrag{d5}{\large$\Delta_2$}
		\psfrag{R_n}{\large$R_u$}
		\psfrag{R_v}{\large$R_v$}
		\psfrag{f_1}{\large$u_3(M_2\otimes M_2)u_3^*$}
		\psfrag{f_2}{\large$u_2(\Delta_2 \otimes M_2)u_2^*$}
		\psfrag{f_3}{\large$u_1M_2u_1^*$}
		\psfrag{f_4}{\large$u_0\Delta_2u_0^*$}
		\psfrag{f_5}{\large$\mathbb{C}$}
		\psfrag{e_1}{\large$v_3(M_2\otimes M_2)v_3^*$}
		\psfrag{e_2}{\large$v_2(\Delta_2 \otimes M_2)v_2^*$}
		\psfrag{e_3}{\large$v_1M_2v_1^*$}
		\psfrag{e_4}{\large$v_0\Delta_2v_0^*$}
		\psfrag{e_5}{\large$\mathbb{C}$}
	\resizebox{10cm}{!}{\includegraphics{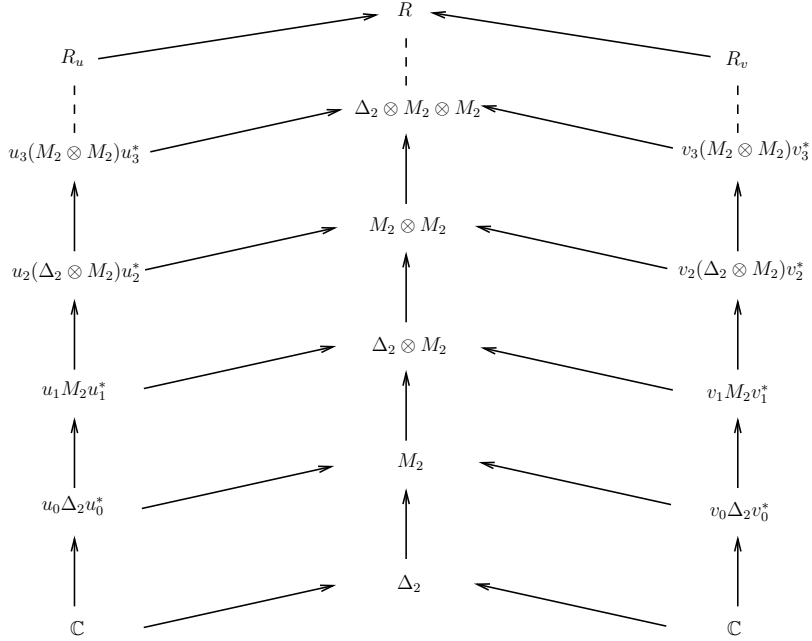}}
	\end{center}
\caption{ A pair of $2\times 2$ Hadamard matrices and basic constructions}\label{fig1}
\end{figure}

Note that
\[
\mathbb{C}\subset \Delta_2\subset M_2\subset \Delta_2 \otimes M_2\subset M_2\otimes M_2\subset \Delta_2\otimes M_2\otimes M_2\subset M_2\otimes M_2\otimes M_2\subset\cdots
\]
is a tower of Jones' basic construction, and thus $R$ is the closure in the SOT topology of the union of these subalgebras. Following \Cref{pairspin}, the ladder of basic constructions of the commuting square
\[
\begin{matrix}
u\Delta_2 u^* &\subset & M_2 \cr \cup &\ & \cup\cr \mathbb{C} &\subset &  \Delta_2
\end{matrix}
\]
is depicted in \Cref{fig1} (note that in our convention $u_0=u$), where the unitary matrices $u_j$ are given in \Cref{1st basic}.
Although the basic construction (in the general $n\times n$ case) is well-known (see \cite{JS, N, Jo3}, for instance), our basic construction is slightly different than the one in the existing literature as it will be more handy in the subsequent computations. We remark that $R_u$ is the closure in the SOT topology of the union of the left vertical algebras as in \Cref{fig1}, and similarly $R_v$ is that of the right vertical algebras. 
\smallskip 

To fix notations, we let $u=(u_{ij})$ be a $2\times 2$ complex Hadamard matrix and consider the following matrices
\[
\eta=\left[{\begin{matrix}
\frac{u_{12}}{u_{11}} & 0\\
0 & \frac{u_{22}}{u_{21}}\\
\end{matrix}}\right]\quad \mbox{and}\quad \xi=\frac{1}{\sqrt{2}}\,\left[{\begin{matrix}
1 & \frac{u_{12}}{u_{11}}\\
1 & \frac{u_{22}}{u_{21}}\\
\end{matrix}}\right]\,.
\]
Since $u$ is unitary, it is easy to check that $\,\eta\mbox{ and }\xi$ are unitary matrices. Let $I_2^{(k)}$ denote the unit element $I_2\otimes\ldots\otimes I_2$ in $(M_2)^{\otimes\,k}$. We have the following tower of basic construction (with the convention that $u_0=u$).

\begin{thm}\label{1st basic}
The tower of the basic construction for $\mathbb{C}\subset u\Delta_2u^*$ is given by
\[
\mathbb{C}\subset u_0\Delta_2u^*_0\subset u_1M_2u^*_1\subset u_2(\Delta_2\otimes M_2)u^*_2\subset u_3(M_2\otimes M_2)u^*_3\subset\cdots
\]
where each $u_i$ is a unitary matrix given by the following prescription~:
\begin{enumerate}[$(i)$]
\item for $k\in\bbn$, we have $\,u_{2k}=u_{2k-1}\left(\xi_k\otimes I_2^{(k)}\right)$, where
\begin{IEEEeqnarray*}{lCl}
\xi_k &=& \begin{cases}
\xi^*\,; & \mbox{ if }\,\,k=1, \cr
F_2\,; & \mbox{ if }\,\,k\geq 2, \cr
\end{cases}
\end{IEEEeqnarray*}
and $\,F_2=\frac{1}{\sqrt{2}}\,
                            \left[{\begin{matrix}
1 & 1\\
1 & -1\\
\end{matrix}}\right]$ is the Fourier matrix in $M_2$.
\item for $k\in\bbn\cup\{0\}$, we have $\,u_{2k+1}=(I_2\otimes u_{2k})\left(E_{11}\otimes I_2^{(k+1)}+E_{22}\otimes\eta_k\right)$, where
\begin{IEEEeqnarray*}{lCl}
\eta_k &=& \begin{cases}
\eta\,; & \mbox{ if }\,\,k=0, \cr
\sigma_3\otimes I_2^{(k)}\,; & \mbox{ if }\,\,k\geq 1. \cr
\end{cases}
\end{IEEEeqnarray*}
and $\sigma_3$ is Pauli spin matrix.
\end{enumerate}
\end{thm}

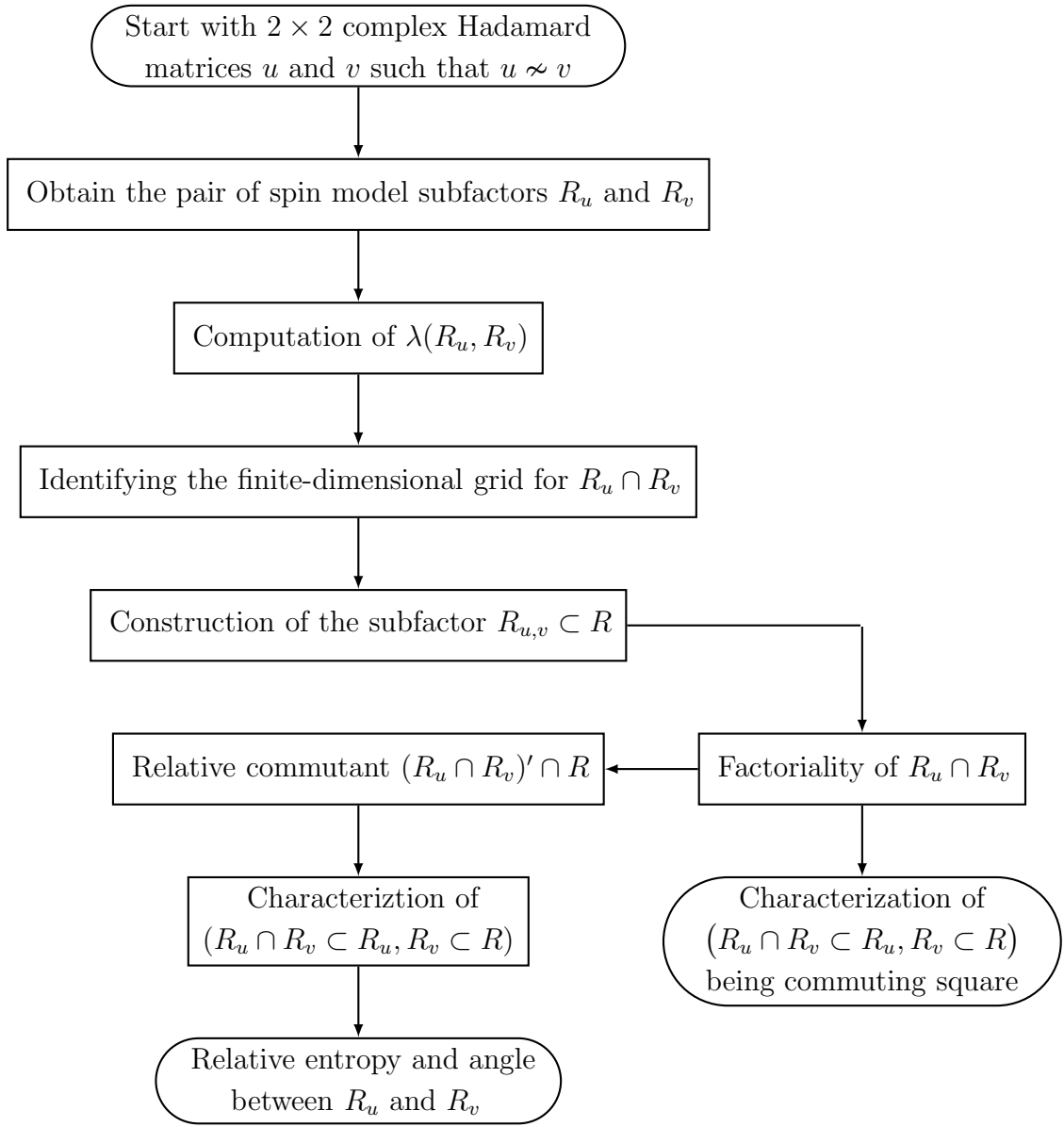
\begin{figure}
\centering
\begin{tikzpicture}[font=\large,thick]

\node[draw,
    rounded rectangle,
    align=center,
    minimum width=2.5cm,
    minimum height=1cm] (block1) {~Start with $2\times 2$ complex Hadamard ~\\
    matrices $u\mbox{ and }v$ such that $u\nsim v$};

\node[draw,
    rectangle, 
    trapezium left angle = 65,
    trapezium right angle = 115,
    trapezium stretches,
    below=of block1,
    minimum width=3.5cm,
    minimum height=1cm] (block2) {~Obtain the pair of spin model subfactors $R_u\mbox{ and }R_v\,$};

\node[draw,
    align=center,
    below=of block2,
    minimum width=3.5cm,
    minimum height=1cm] (block3) {~Computation of $\lambda(R_u,R_v)\,$};

\node[draw,
    align=center,
    below=of block3,
    minimum width=3.5cm,
    minimum height=1cm] (block40) {~Identifying the finite-dimensional grid for $R_u\cap R_v\,$};
  
\node[draw,
    align=center,
    below=of block40,
    minimum width=3.5cm,
    minimum height=1cm] (block4) {~Construction of the subfactor $R_{u,v}\subset R\,$};

\node[draw,
    rectangle,
    below right=of block4,
    minimum width=3.5cm,
    minimum height=1cm] (block6) {~Factoriality of $R_u\cap R_v\,$};

\node[draw,
    rounded rectangle,
    align=center,
    below=of block6,
    minimum width=3.5cm,
    minimum height=1cm] (block50) {Characterization of\\
    $\big(R_u\cap R_v\subset R_u,R_v\subset R\big)$\\
    being commuting square};

\node[draw,
    below=of block4,
    minimum width=3.5cm,
    minimum height=1cm] (block7) {~Relative commutant $(R_u\cap R_v)^\prime\cap R\,$};

\node[draw,
    align=center,
    below=of block7,
    minimum width=2.5cm,
    minimum height=1cm] (block8) {~Characteriztion of\\
     $\,(R_u\cap R_v\subset R_u,R_v\subset R)\,$};
    
\node[draw,
    rounded rectangle,
    align=center,
    below=of block8,
    minimum width=2.5cm,
    minimum height=1cm] (block80) {~Relative entropy and angle\\
     between $R_u\mbox{ and }R_v\,$};

\draw[-latex] (block1) edge (block2)
    (block2) edge (block3)
    (block3) edge (block40)
    (block40) edge (block4)
    (block6) edge (block50)
    (block6) edge (block7)
    (block7) edge (block8)
    (block8) edge (block80);
    
\draw[-latex] (block4) -| (block6)
    node[pos=0.5,fill=white,inner sep=0]{};

\end{tikzpicture}
\bigskip

\caption{Flowchart for our plan of actions in \Cref{Sec 4}}\label{flow chart 1}
\end{figure}

For the proof of this theorem we request the interested reader to visit the Appendix. We draw a flowchart in \Cref{flow chart 1} to display our plan of actions in this section.

\subsection{Computation of the Pimsner-Popa constant}

Let $u\mbox{ and }v$ be two complex Hadamard matrices of the form described in \Cref{reduc}. In this subsection, we compute $\lambda(R_u,R_v)$. Before we begin, we pause for a moment to discuss the following example. Consider the following diagram
\[
\begin{matrix}
u\Delta_2u^* &\subset & M_2 &\supset &v\Delta_2v^*\\
\cup & &\cup & &\cup\\
\bbc &\subset & \Delta_2 &\supset &\bbc
\end{matrix}
\]
where $u$ and $v$ are any two distinct $2\times 2$ complex Hadamard matrices. The quadruple $(\bbc\subset\Delta_2,u\Delta_2u^*\subset M_2)$ (resp., $u$ replaced by $v$) is a symmetric commuting square, since $u$ (resp., $v$) is a complex Hadamard matrix. Now, $\lambda(\bbc,\bbc)=1$ but $\lambda(u\Delta_2u^*,v\Delta_2v^*)=\lambda(\Delta_2,u^*v\Delta_2v^*u)=\min_{1\leq i\leq 2}\,\big(\mbox{$h\left(v^*u\right)_i$}\big)^{-1}$ by \Cref{finitemains}. If we choose $u=\frac{1}{\sqrt{2}}\left(\begin{smallmatrix}
1 & 1\\
i & -i
\end{smallmatrix}\right)$ and $v=\frac{1}{\sqrt{2}}\left(\begin{smallmatrix}
1 & 1\\
1 & -1
\end{smallmatrix}\right)$, we see that $\lambda(u\Delta_2u^*,v\Delta_2v^*)=\frac{1}{2}\neq\lambda(\bbc,\bbc)$. This example illustrates that finding the value of $\lambda(R_u,R_v)$ is not obvious.
\smallskip

\noindent\textbf{Notation~:}~ $A(u_4,v_4):=v_4^*u_4(\sigma_3\otimes I_4)u_4^*v_4(\sigma_3\otimes I_4)$ in $M_8$, and $B(u_2,v_2):=v_2^*u_2(\sigma_3\otimes I_2)u_2^*v_2$ in $M_4$.

\begin{lmma}\label{not diagonal}
Let $k\geq 2$ be a natural number. The following matrix
\[
v_{2k}^*u_{2k}\left(\sigma_3\otimes I_2^{(k)}\right)u_{2k}^*v_{2k}\left(\sigma_3\otimes I_2^{(k)}\right)
\]
in $M_2^{(k+1)}$ is not diagonal if the following matrix
\begin{IEEEeqnarray*}{lCl}
A(u_4,v_4)=v_4^*u_4(\sigma_3\otimes I_4)u_4^*v_4(\sigma_3\otimes I_4)
\end{IEEEeqnarray*}
$($i,e. for $k=2)$ in $M_2^{(3)}$ is not diagonal.
\end{lmma}
\begin{prf}
Recall the tower of basic construction in Theorem \ref{1st basic} for any $2\times 2$ complex Hadamard matrix. For two $2\times 2$ complex Hadamard matrices $u\mbox{ and }v$, we denote by $\xi_k(u)\mbox{ and }\eta_k(u)$ the corresponding matrices $\xi_k\mbox{ and }\eta_k$ in the tower of basic construction for the unitary $u$, and similarly for $v$. Note that we are interested in $k\geq 3$, since the statement is obviously true for $k=2$, and we have $\xi_{n+1}(u)=\xi_{n+1}(v)=F_2\mbox{ and }\eta_n(u)=\eta_n(v)=\sigma_3\otimes I_2^{(n)}$ for any $n\in\bbn$ by \Cref{1st basic}. The matrices $F_2\mbox{ and }\sigma_3$ both are self-adjoint. For $k\geq 3$, we have the following,
\begin{IEEEeqnarray*}{lCl}
&  & v_{2k}^*u_{2k}\left(\sigma_3\otimes I_2^{(k)}\right)u_{2k}^*v_{2k}\left(\sigma_3\otimes I_2^{(k)}\right)\\
&=& \left(F_2\otimes I_2^{(k)}\right)v_{2k-1}^*u_{2k-1}\left(F_2\otimes I_2^{(k)}\right)\left(\sigma_3\otimes I_2^{(k)}\right)\left(F_2\otimes I_2^{(k)}\right)u_{2k-1}^*v_{2k-1}\left(F_2\otimes I_2^{(k)}\right)\left(\sigma_3\otimes I_2^{(k)}\right)\\
&=& \left(F_2\otimes I_2^{(k)}\right)\,\mbox{bl-diag}\left\{I_2^{(k)},\eta_{k-1}(v)\right\}^*\,\left(I_2\otimes v_{2k-2}^*u_{2k-2}\right)\,\mbox{bl-diag}\left\{I_2^{(k)},\eta_{k-1}(u)\right\}(\sigma_1\otimes I_2^{(k)})\\
&  & \mbox{bl-diag}\left\{I_2^{(k)},\eta_{k-1}(u)\right\}^*\left(I_2\otimes u_{2k-2}^*v_{2k-2}\right)\,\mbox{bl-diag}\left\{I_2^{(k)},\eta_{k-1}(v)\right\}\left(\frac{1}{\sqrt{2}}(I_2-i\sigma_2)\otimes I_2^{(k)}\right)\,.
\end{IEEEeqnarray*}
For $k\geq 3$, we have $\eta_{k-1}(u)=\eta_{k-1}(v)=\sigma_3\otimes I_2^{(k-1)}$. Hence,
\begin{IEEEeqnarray*}{lCl}
&  & v_{2k}^*u_{2k}\left(\sigma_3\otimes I_2^{(k)}\right)u_{2k}^*v_{2k}\left(\sigma_3\otimes I_2^{(k)}\right)\\
&=& \left(F_2\otimes I_2^{(k)}\right)\,\mbox{bl-diag}\left\{I_2^{(k)},\sigma_3\otimes I_2^{(k-1)}\right\}\left[{\begin{matrix}
x_{2k-2}\left(\sigma_3\otimes I_2^{(k-1)}\right) & x_{2k-2}\left(\sigma_3\otimes I_2^{(k-1)}\right)\\
x_{2k-2} & -x_{2k-2}\\
\end{matrix}}\right]
\end{IEEEeqnarray*}
with $\,x_{2k-2}=v_{2k-2}^*u_{2k-2}\left(\sigma_3\otimes I_2^{(k-1)}\right)u_{2k-2}^*v_{2k-2}$ in $M_2\otimes M_2^{(k-1)}$. Therefore, we get that
\begin{IEEEeqnarray}{lCl}\label{only basis step}
&  & v_{2k}^*u_{2k}\left(\sigma_3\otimes I_2^{(k)}\right)u_{2k}^*v_{2k}\left(\sigma_3\otimes I_2^{(k)}\right)\nonumber\\
&=& \left[{\begin{matrix}
\left\{x_{2k-2}\,,\,\sigma_3\otimes I_2^{(k-1)}\right\} &   &  \left[x_{2k-2}\,,\,\sigma_3\otimes I_2^{(k-1)}\right]\\
  &   &   \\
\left[x_{2k-2}\,,\,\sigma_3\otimes I_2^{(k-1)}\right] &   &   \left\{x_{2k-2}\,,\,\sigma_3\otimes I_2^{(k-1)}\right\}\\
\end{matrix}}\right]
\end{IEEEeqnarray}
where $[\,,\,]$ and $\{\,,\,\}$ denotes the commutator and the anti-commutator respectively. Now, suppose that for $k=2$, the following matrix
\begin{IEEEeqnarray}{lCl}\label{matrix in relation}
A(u_4,v_4):=v_4^*u_4(\sigma_3\otimes I_4)u_4^*v_4(\sigma_3\otimes I_4)
\end{IEEEeqnarray}
in $M_8$ is not diagonal. Assume that the statement is true up to $(k-1)$-th step. Now for the $k$-th step, first by the induction hypothesis we get that $x_{2k-2}\left(\sigma_3\otimes I_2^{(k-1)}\right)$ is not a diagonal matrix. From \Cref{only basis step}, we see that if the commutator $[x_{2k-2},\sigma_3\otimes I_2^{(k-1)}]\neq 0$, then we are through and if $\,[x_{2k-2},\sigma_3\otimes I_2^{(k-1)}]=0$, then the anticommutator $\{x_{2k-2},\sigma_3\otimes I_2^{(k-1)}\}$ becomes $2x_{2k-2}\left(\sigma_3\otimes I_2^{(k-1)}\right)$, which is not a diagonal matrix by the induction hypothesis. Thus, if the matrix in \Cref{matrix in relation} is not diagonal, then the matrix in \Cref{only basis step} is also not diagonal, which completes the proof.\qed
\end{prf}

\begin{lmma}\label{ef222}
The matrix $A(u_4,v_4)$ is diagonal if and only if the following self-adjoint matrix
\[
B(u_2,v_2)=v_2^*u_2(\sigma_3\otimes I_2)u_2^*v_2
\]
is diagonal in $M_4$.
\end{lmma}
\begin{prf}
Using Theorem \ref{1st basic}, we obtain the following,
\begin{IEEEeqnarray*}{rCl}
u_4^*v_4 &=& \left[{\begin{matrix}
u_2^* & (\sigma_3\otimes I_2)u_2^*\\
u_2^* & -(\sigma_3\otimes I_2)u_2^*\\
\end{matrix}}\right]\left[{\begin{matrix}
v_2 & v_2\\
v_2(\sigma_3\otimes I_2) & -v_2(\sigma_3\otimes I_2)\\
\end{matrix}}\right]\\
&  & \\
&  & \\
&=& \left[{\begin{matrix}
u_2^*v_2+(\sigma_3\otimes I_2)u_2^*v_2(\sigma_3\otimes I_2) &  u_2^*v_2-(\sigma_3\otimes I_2)u_2^*v_2(\sigma_3\otimes I_2)\\
u_2^*v_2-(\sigma_3\otimes I_2)u_2^*v_2(\sigma_3\otimes I_2) &  u_2^*v_2+(\sigma_3\otimes I_2)u_2^*v_2(\sigma_3\otimes I_2)
\end{matrix}}\right]\,.
\end{IEEEeqnarray*}
Hence, we get that
\begin{IEEEeqnarray*}{lCl}
(\sigma_3\otimes I_4)u_4^*v_4(\sigma_3\otimes I_4)
&=& \left[{\begin{matrix}
u_2^*v_2+(\sigma_3\otimes I_2)u_2^*v_2(\sigma_3\otimes I_2) & -u_2^*v_2+(\sigma_3\otimes I_2)u_2^*v_2(\sigma_3\otimes I_2)\\
-u_2^*v_2+(\sigma_3\otimes I_2)u_2^*v_2(\sigma_3\otimes I_2) &  u_2^*v_2+(\sigma_3\otimes I_2)u_2^*v_2(\sigma_3\otimes I_2)
\end{matrix}}\right]\,.
\end{IEEEeqnarray*}
Therefore,
\begin{IEEEeqnarray*}{rCl}
A(u_4,v_4) &=& v_4^*u_4(\sigma_3\otimes I_4)u_4^*v_4(\sigma_3\otimes I_4)\\
&=& \left[{\begin{matrix}
v_2^*u_2+(\sigma_3\otimes I_2)v_2^*u_2(\sigma_3\otimes I_2) &  v_2^*u_2-(\sigma_3\otimes I_2)v_2^*u_2(\sigma_3\otimes I_2)\\
v_2^*u_2-(\sigma_3\otimes I_2)v_2^*u_2(\sigma_3\otimes I_2) &  v_2^*u_2+(\sigma_3\otimes I_2)v_2^*u_2(\sigma_3\otimes I_2)\end{matrix}}\right]\\
&  & \\
& & \left[{\begin{matrix}
u_2^*v_2+(\sigma_3\otimes I_2)u_2^*v_2(\sigma_3\otimes I_2) &  -u_2^*v_2+(\sigma_3\otimes I_2)u_2^*v_2(\sigma_3\otimes I_2)\\
-u_2^*v_2+(\sigma_3\otimes I_2)u_2^*v_2(\sigma_3\otimes I_2) &  u_2^*v_2+(\sigma_3\otimes I_2)u_2^*v_2(\sigma_3\otimes I_2)\end{matrix}}\right]\\
&  &\\
&=& 2\left[{\begin{matrix}
A(u_4,v_4)_{11}  &  A(u_4,v_4)_{12}\\
  &  \\
A(u_4,v_4)_{21}  &  A(u_4,v_4)_{22}
\end{matrix}}\right]
\end{IEEEeqnarray*}
where,
\[
A(u_4,v_4)_{11}=A(u_4,v_4)_{22}=v_2^*u_2(\sigma_3\otimes I_2)u_2^*v_2(\sigma_3\otimes I_2)+(\sigma_3\otimes I_2)v_2^*u_2(\sigma_3\otimes I_2)u_2^*v_2\,,
\]
\[
A(u_4,v_4)_{12}=A(u_4,v_4)_{21}=v_2^*u_2(\sigma_3\otimes I_2)u_2^*v_2(\sigma_3\otimes I_2)-(\sigma_3\otimes I_2)v_2^*u_2(\sigma_3\otimes I_2)u_2^*v_2\,.
\]
Hence, $A(u_4,v_4)=I_2\otimes A(u_4,v_4)_{11}+\sigma_1\otimes A(u_4,v_4)_{12}$. This says that the matrix $A(u_4,v_4)$ is diagonal if and only if $A(u_4,v_4)_{11}$ is diagonal and $A(u_4,v_4)_{12}=0$. That is, $A(u_4,v_4)$ is diagonal if and only if the self-adjoint matrix $v_2^*u_2(\sigma_3\otimes I_2)u_2^*v_2$ commutes with $\sigma_3\otimes I_2$ and $v_2^*u_2(\sigma_3\otimes I_2)u_2^*v_2(\sigma_3\otimes I_2)$ is diagonal, which is further equivalent to $B(u_2,v_2)=v_2^*u_2(\sigma_3\otimes I_2)u_2^*v_2$ is diagonal.\qed
\end{prf}

\begin{lmma}\label{ef22222}
The matrix $B(u_2,v_2)=v_2^*u_2(\sigma_3\otimes I_2)u_2^*v_2$ is not diagonal.
\end{lmma}
\begin{prf}
Using \Cref{1st basic}, a tedious but straightforward computation gives us the following,
\begin{IEEEeqnarray*}{lCl}
v_2^*u_2(\sigma_3\otimes I_2)u_2^*v_2 &=& \frac{1}{2}\left[{\begin{matrix}
\sigma_3v^*u\sigma_3u^*v & & v^*u\sigma_3u^*v\sigma_3\\
-\sigma_3v^*u\sigma_3u^*v &  & v^*u\sigma_3u^*v\sigma_3\\
\end{matrix}}\right]\left(\left[{\begin{matrix}
1 & 1\\
1 & -1\\
\end{matrix}}\right]\otimes I_2\right)
\end{IEEEeqnarray*}
We denote $C=\sigma_3v^*u\sigma_3u^*v\in M_2$. Observe that the matrix $B(u_2,v_2)$ is diagonal if and only if the following holds,
\begin{IEEEeqnarray}{lCl}\label{diagonal conditions1}
C-C^*=0\qquad\mbox{and}\qquad C+C^*\mbox{\,\, is diagonal}\,.
\end{IEEEeqnarray}
Now, another straightforward computation shows that
\[
C=\begin{bmatrix}
\cos(\alpha-\beta) & i\sin(\alpha-\beta)\\
i\sin(\alpha-\beta) & \cos(\alpha-\beta)\\
\end{bmatrix}
\]
using \Cref{reduc}. Since $\alpha-\beta\notin\{0,\pm\pi\}$, we see that $C$ fails to be self-adjoint, and consequently by \Cref{diagonal conditions1} the proof is completed.\qed
\end{prf}
 
\begin{thm}\label{PP constant final}
The Pimsner-Popa constant $\lambda(R_u,R_v)$ for the pair of subfactors $R_u\mbox{ and }R_v$ of the hyperfinite $II_1$ factor $R$ is $1/2$.
\end{thm}
\begin{prf}
Let $k\geq 2$. We claim the following,
\begin{IEEEeqnarray}{lCl}\label{gd}
\lambda\left(u_{2k+1}\big(M_2\otimes M_2^{(k)}\big)u_{2k+1}^{\,*}\,,\,v_{2k+1}\big(M_2\otimes M_2^{(k)}\big)v_{2k+1}^*\right)=\frac{1}{2}\,.
\end{IEEEeqnarray}
Recall from Theorem \ref{1st basic} that
\[
u_{2k+1}=\mbox{bl-diag}\left\{u_{2k}\big(I_2\otimes I_2^{(k)}\big),u_{2k}\big(\sigma_3\otimes I_2^{(k)}\big)\right\}
\]
as an element in $\Delta_2\otimes M_2\otimes M_2^{(k)}$. Therefore,
\[
u_{2k+1}^{\,*}v_{2k+1}=\mbox{bl-diag}\left\{u_{2k}^{\,*}v_{2k}\,,\big(\sigma_3\otimes I_2^{(k)}\big)u_{2k}^{\,*}v_{2k}\big(\sigma_3\otimes I_2^{(k)}\big)\right\}\,.
\]
Now,
\begin{IEEEeqnarray*}{lCl}
&  & \lambda\left(u_{2k+1}\big(M_2\otimes M_2^{(k)}\big)u_{2k+1}^{\,*}\,,\,v_{2k+1}\big(M_2\otimes M_2^{(k)}\big)v_{2k+1}^*\right)\\
&=& \lambda\left(M_2^{(k+1)}\,,\,u_{2k+1}^{\,*}v_{2k+1}M_2^{(k+1)}v_{2k+1}^*u_{2k+1}\right)
\end{IEEEeqnarray*}
with $u_{2k+1}^{\,*}v_{2k+1}\in\Delta_2\otimes M_2\otimes M_2^{(k)}$. We now wish to apply \Cref{2nd lambda}, and for that first consider the following matrix
\begin{IEEEeqnarray}{lCl}\label{8686}
v_{2k}^*u_{2k}\big(\sigma_3\otimes I_2^{(k)}\big)u_{2k}^{\,*}v_{2k}\big(\sigma_3\otimes I_2^{(k)}\big)\,.
\end{IEEEeqnarray}
By \Cref{ef222}, \ref{ef22222} we get that the matrix $A(u_4,v_4)$ is not diagonal. Hence, by \Cref{not diagonal} we get that the matrix in \Cref{8686} is not diagonal. As a consequence of \Cref{2nd lambda}, our claim has been established. By \Cref{popaadaptation}, $\,\lambda(R_u,R_v)$ is the limit of a decreasing sequence involving $\lambda$ at each step of the tower of basic construction starting with the complex Hadamard matrices $u,v$. Since we have produced a constant subsequence in \Cref{gd}, the proof is concluded.\qed
\end{prf}


\subsection{Factoriality of $R_u\cap R_v$ and the quadruple $(R_u\cap R_v\subset R_u,R_v\subset R)$}\label{entropysubsection}

An astute reader must have noticed that  the intersection $R_u\cap R_v$ need not be a factor {\em a priori}. Another difficulty is to determine whether $[R:R_u\cap R_v]$ is finite or not. Our goal in this section is to address both these issues. We prove that $R_u\cap R_v$ is a $II_1$ subfactor of $R$ with $[R:R_u\cap R_v]$=4. We also characterize when the quadruple $(R_u\cap R_v\subset R_u,R_v\subset R)$ of $II_1$ factors is a commuting (and consequently, co-commuting) square. For the sake of brevity, we introduce some notations.

\begin{notation}\rm\label{not}
\noindent \begin{enumerate}
\item For $k\in\bbn\cup\{0\}$, let
\begin{align*}
A_{2k} &= M_2\otimes M_2^{(k)}\,,\cr
B^u_{2k} &= u_{2k}(\Delta_2\otimes M_2^{(k)})u_{2k}^*\,,\cr
B^v_{2k} &= v_{2k}(\Delta_2\otimes M_2^{(k)})v_{2k}^*\,,\cr
C_{2k} &= B^u_{2k}\cap {B}^v_{2k}\,.
\end{align*}
These are the even steps in the tower of basic constructions in \Cref{fig1}. Note that
\[
\overline{\cup_k\,A_{2k}}^{\mbox{ sot}}=R,\,\,\overline{\cup_k\,B^u_{2k}}^{\mbox{ sot}}=R_u,\,\,\overline{\cup_k\,B^v_{2k}}^{\mbox{ sot}}=R_v,\,\mbox{ and }\,\overline{\cup_k\,C_{2k}}^{\mbox{ sot}}=R_u\cap R_v.
\]
\item For $k\in\bbn$, let $W_{2k}:=E_{11}\otimes I_2^{(k)}+E_{22}\otimes\sigma_1^{(k)}$ in $\Delta_2\otimes M_2^{(k)}$. Each $W_{2k}$ is unitary in $M_2^{(k+1)}$. Recall that $A^{(k)}\in M_2^{(k)}$ denotes the matrix $A^{\otimes\,k}=A\otimes A\otimes\cdots\otimes A$ in $M_2^{(k)}$.
\end{enumerate}
\end{notation}

We begin by computing the Sano-Watatani angle between $B^u_0\mbox{ and }B^v_0$ (see Notation \ref{not}), and show that the following quadruple
\begin{IEEEeqnarray}{lCl}\label{spq}
\begin{matrix}
B^u_0 &\subset & A_0\\
 \cup & & \cup\\
C_0 &\subset & B^v_0\\
\end{matrix}
\end{IEEEeqnarray}
may not be a commuting square. This is an example which shows that the floor (and consequently the roof due to \Cref{comm cube 2}) in \Cref{com1} need not be a commuting square. The difficulty starts from here.

\begin{ppsn}\label{p2}
If $u= \frac{1}{\sqrt{2}}\begin{bmatrix}
e^{i\alpha_1} & e^{i (\alpha_1+\alpha_3)}\\
e^{i\alpha_2} & -e^{i(\alpha_2+\alpha_3)}
\end{bmatrix}$ and $v= \frac{1}{\sqrt{2}}\begin{bmatrix}
e^{i\beta_1} & e^{i (\beta_1+\beta_3)}\\
e^{i\beta_2} & -e^{i(\beta_2+\beta_3)}
\end{bmatrix}$, then the cosine of the angle between $u\Delta_2u^*\mbox{ and }v\Delta_2v^*$ is the set $\big\{0,\big|\cos((\alpha_1-\alpha_2)-(\beta_1-\beta_2))\big|\big\}$, and the quadruple in \Cref{spq} is a commuting square if and only if $\alpha_1-\alpha_2=\beta_1-\beta_2\pm\frac{\pi}{2}$.
\end{ppsn}
\begin{prf}
It is clear that $u\Delta_2u^*\cap\,v\Delta_2v^*=\bbc$, as $\mbox{dim}(u\Delta_2u^*)=\mbox{dim}(v\Delta_2v^*)=2$ and $u\Delta_2u^*=v\Delta_2v^*$ implies that $u\sim v$ by \Cref{finitemains}. Consider the non-negative operator $E_{u\Delta_2u^*}E_{v\Delta_2v^*}E_{u\Delta_2u^*}$, and observe that as an element of $\mathcal{B}(M_2,tr)=M_4$, this is same as the operator $\mbox{Ad}_uE_{\Delta_2}\mbox{Ad}_{u^*}\mbox{Ad}_vE_{\Delta_2}\mbox{Ad}_{v^*}\mbox{Ad}_uE_{\Delta_2}\mbox{Ad}_{u^*}$. Put $w=(w_{ij})=u^*v$, a unitary in $M_2$. Then, we have the operator $\mbox{Ad}_uE_{\Delta_2} E_{w\Delta_2w^*}E_{\Delta_2}\mbox{Ad}_{u^*}$ in $\mathcal{B}(M_2)$. First consider $E_{\Delta_2}E_{w\Delta_2w^*}E_{\Delta_2}\in\mathcal{B}(M_2)$. Since $w$ is a $2\times 2$ unitary matrix, we have $|w_{12}|^2=|w_{21}|^2$. We put $c= \sum_{i=1}^2{|w_{1i}|}^4$ and $s= {2}\prod_{i=1}^2 {|w_{i1}|}^2$. It is easy to observe that $\sum_{i=1}^2 {|w_{2i}|}^4=c$ and $2\prod_{i=1}^2 {|w_{i2}|}^2=s$. Therefore, using the fact that $|w_{11}|^2+|w_{12}|^2=1$ (as $w^*$ is a unitary) we get the following,
\[
c+s=|w_{11}|^4+|w_{12}|^4+ 2|w_{11}|^2|w_{21}|^2=|w_{11}|^4+|w_{12}|^4+ 2|w_{11}|^2|w_{12}|^2=1.
\]
Now, it is easy to verify that
\[
E_{\Delta} E_{w\Delta w^*} E_{\Delta}(A)=\mathrm{diag}\{c a_{11}+ (1-c) a_{22}\,,\,(1-c) a_{11}+ c a_{22}\}
\] 
for any matrix $A=(a_{ij})$ in $M_2$. Therefore, with respect to the ordered orthonormal basis $\mathfrak{B}:=\big\{\sqrt{2}E_{11},\sqrt{2}E_{22},\sqrt{2}E_{12},\sqrt{2}E_{21}\big\}$ of $M_2$, the operator $E_{\Delta_2}E_{w\Delta_2 w^*}E_{\Delta_2}:M_2\longrightarrow M_2$, when viewed as a linear map from $\bbc^{4}$ to $\bbc^{4}$, is of the form $\mbox{bl-diag}\{C,0\}$, where
${C}= \begin{bmatrix}
c & 1-c \\
1-c &~ c
\end{bmatrix}$. The conditional expectation from $M_2$ onto $u\Delta_2u^*\cap v\Delta_2v^*$ is simply the trace $tr$. Therefore, with respect to the ordered orthonormal basis $\mathfrak{B}$, the operator $tr:M_2\to\bbc\subset M_2$ given by $A\mapsto tr(A)I_2$ is of the form $\mbox{bl-diag}\{J_2,0\}$, where $J_2$ is the projection $\frac{1}{2}\sum_{i,j=1}^2E_{ij}$ in $M_2$. Hence, $E_{u\Delta_2u^*}E_{v\Delta_2v^*}E_{u\Delta_2u^*}-tr$ has the form $\mbox{bl-diag}\{C-J_2,0\}$. The quadruple in \Cref{spq} is a commuting square if and only if $E_{u\Delta_2u^*}E_{v\Delta_2v^*}E_{u\Delta_2u^*}-tr$, that is, $C=J_2$. Since $J_2=\frac{1}{2}\sum_{i,j=1}^2 E_{ij}$, we see that $C=J_2$ holds if and only if $c=s=\frac{1}{2}$. A straightforward verification leads us that $c=\frac{1}{2}(1+\cos^2(\beta_1-\alpha_1-\beta_2+\alpha_2))$. Therefore, the quadruple in \Cref{spq} is a commuting square if and only if $\alpha_1-\alpha_2=\beta_1-\beta_2\pm\frac{\pi}{2}$. 

Now, suppose that this condition is not satisfied. The spectrum of $\mbox{Ad}_u(E_{\Delta_2}E_{w\Delta_2 w^*}E_{\Delta_2}-tr)\mbox{Ad}_{u^*}$ is simply the spectrum of the matrix $C-J_2\in M_2$. This is because $\mbox{Ad}_u\in\mathcal{B}(M_2)$ is a unitary. Therefore, the spectrum of the operator $E_{u\Delta_2u^*}E_{v\Delta_2v^*}E_{u\Delta_2u^*}-E_{u\Delta_2u^*\,\cap\,v\Delta_2v^*}$ is the spectrum of the matrix $C-J_2$, which is $\{0,2c-1\}$. Observe that $2c-1=\big({|w_{11}|}^2-{|w_{12}|}^2\big)^2$. A straightforward verification shows that $|w_{11}|^2-|w_{12}|^2=\cos\big((\alpha_1-\alpha_2)-(\beta_1-\beta_2)\big)$. Therefore, we obtain that cosine of the angle between $u\Delta_2u^*\mbox{ and }v\Delta_2v^*$ is the set $\big\{0,\big|\cos(\alpha_1-\alpha_2)-(\beta_1-\beta_2))\big|\big\}$.\qed
\end{prf}

\begin{crlre}\label{p23}
Let $\,u\mbox{ and }v$ be of the form as described in \Cref{reduc}. Then, the cosine of the angle between $\mathrm{Ad}_u(\Delta_2)\mbox{ and }\mathrm{Ad}_v(\Delta_2)$ is the set $\{0,|\cos(\alpha-\beta)|\}$.
\end{crlre}

By \Cref{p23}, the quadruple $(C_0\subset B^v_0\,,\,B^u_0\subset A_0)$ in \Cref{spq}, and consequently $(R_u\cap R_v\subset R_u,R_v\subset R)$ due to \Cref{comm cube 2}, is not a commuting square (see \Cref{basicresultoncomm} in this regard) unless $\alpha-\beta=\pm\frac{\pi}{2}$. With the notations as in \Cref{1st basic}, we now identify the finite-dimensional grid $C_{2k}$ for $R_u\cap R_v$.

\begin{lmma}\label{e9}
We have the following
\[
{u}^*_{2k}{v}_{2k} \big(\Delta_2\otimes M^{(k)}_2\big){v}^*_{2k}{u}_{2k} =W_{2k}\big(u^*v\Delta_2v^*u\otimes M^{(k)}_2\big) W^*_{2k}
\]
for any $k\in\bbn$, where $W_{2k}$ is as defined in \Cref{not}.
\end{lmma}
\begin{prf}
We first prove it for $k=1$. We have
\begin{IEEEeqnarray}{lCl}\label{s1}
u^*v &=& \frac{1}{2}(1+e^{i(\beta-\alpha)})I_2+\frac{1}{2}(1-e^{i(\beta-\alpha)})\sigma_1\,.
\end{IEEEeqnarray}
Applying \Cref{1st basic} we see that 
\begin{IEEEeqnarray*}{lCl}
u^*_2v_2=\frac{1}{2}\begin{bmatrix}
u^*v+\sigma_3u^*v\sigma_3 & u^*v-\sigma_3u^*v\sigma_3\\
u^*v+\sigma_3u^*v\sigma_3  & u^*v-\sigma_3u^*v\sigma_3
\end{bmatrix}.
\end{IEEEeqnarray*}
Note that
\begin{center}
$u^*v+\sigma_3u^*v\sigma_3=\big(1+e^{i(\beta-\alpha)}\big)I_2\quad\mbox{ and }\quad u^*v-\sigma_3u^*v\sigma_3=\big(1-e^{i(\beta-\alpha)}\big)\sigma_1\,.$
\end{center}
Therefore, we get that
\begin{IEEEeqnarray}{lCl}\label{s2}
u^*_2v_2=\frac{1}{2}(1+e^{i(\beta-\alpha)})I_2\otimes I_2+\frac{1}{2}(1-e^{i(\beta-\alpha)})\sigma_1\otimes\sigma_1\,.
\end{IEEEeqnarray}
Combining \Cref{s1} and \Cref{s2}, we get the following
\[
u^*_2v_2=\mbox{bl-diag}\{I_2,\sigma_1\}(u^*v\otimes I_2)\mbox{bl-diag}\{I_2,\sigma_1\},
\]
and hence we have
\[
u^*_2v_2(\Delta_2\otimes M_2)v^*_2u_2=\mbox{bl-diag}\{I_2,\sigma_1\}\big(u^*v\Delta_2v^*u\otimes M_2\big) \mbox{bl-diag}\{I_2,\sigma_1\}.
\]
To complete the induction note that for $k\geq 2$, thanks to \Cref{1st basic}, we see that
\[
u^*_{2k}v_{2k}=\frac{1}{2}(1+e^{i(\beta-\alpha)})I_2^{(k)}+\frac{1}{2}(1-e^{i(\beta-\alpha)})\sigma_1^{(k)}\,.
\]
In other words, by \Cref{s1} we obtain the following equality
\[
u^*_{2k}v_{2k}=\mbox{bl-diag}\{I^{(k)}_2,\sigma_1^{(k)}\}\big(u^*v\otimes I^{(k)}_2\big)\mbox{bl-diag}\{I^{(k)}_2,\sigma_1^{(k)}\}\,.
\]
Since $W_{2k}=\mbox{bl-diag}\{I^{(k)}_2,\sigma_1^{(k)}\}\in\Delta_2\otimes M^{(k)}_2$, the proof is completed.\qed
\end{prf}

\begin{ppsn}\label{w}
For any $k\in\bbn$, we have the following
\begin{IEEEeqnarray*}{lCl}
C_{2k}:=u_{2k}\left(\Delta_2\otimes M_2^{(k)}\right)u_{2k}^*\bigcap v_{2k}\left(\Delta_2\otimes M_2^{(k)}\right)v_{2k}^*
&=& \mathrm{Ad}_{u_{2k}W_{2k}}\left(\big(\Delta_2\cap (u^*v\Delta_2v^*u)\big)\otimes M_2^{(k)}\right)\\
&=& \mathrm{Ad}_{u_{2k}W_{2k}}\left(\bbc\otimes M_2^{(k)}\right)\,.
\end{IEEEeqnarray*}
\end{ppsn}
\begin{prf}
We have the following,
\begin{IEEEeqnarray*}{lCl}
& & u_{2k}\left(\Delta_2\otimes M_2^{(k)}\right)u_{2k}^*\bigcap v_{2k}\left(\Delta_2\otimes M_2^{(k)}\right)v_{2k}^*\\
&=& \mbox{Ad}_{u_{2k}}\left(\Delta_2\otimes M_2^{(k)}\bigcap \mbox{Ad}_{u_{2k}^*v_{2k}}\big(\Delta_2\otimes M_2^{(k)}\big)\right)\\
&=& \mbox{Ad}_{u_{2k}}\left(\Delta_2\otimes M_2^{(k)}\bigcap \mbox{Ad}_{W_{2k}}\big(u^*v\Delta_2v^*u\otimes M_2^{(k)}\big)\right)\qquad(\mbox{by \Cref{e9}})\\
&=& \mbox{Ad}_{u_{2k}W_{2k}}\left(\Delta_2\otimes M_2^{(k)}\bigcap\, (u^*v\Delta_2v^*u)\otimes M_2^{(k)}\right)\\
&=& \mbox{Ad}_{u_{2k}W_{2k}}\left(\big(\Delta_2\cap\,(u^*v\Delta_2v^*u)\big)\otimes M_2^{(k)}\right)\,.
\end{IEEEeqnarray*}
This proves the first equality, and for the second, note that since $\beta\neq\alpha,\pi\pm\alpha$ in our case, we have $\Delta_2\cap \left(u^*v\Delta_2v^*u\right)=\bbc\,$.\qed
\end{prf}

We have the following tower of finite von Neumann algebras
\[
\begin{matrix}
A_0 &\subset & A_2 &\subset & A_4 &\subset  &\ldots &\subset & R \cr
\cup &  &\cup &  &\cup & & \ldots &  &\cup \cr
B_0^u &\subset & B_2^u &\subset & B_4^u &\subset  &\ldots &\subset & R_u \cr
\cup &  &\cup &  &\cup &  &\ldots &  &\cup \cr
C_0 &\subset & C_2 &\subset & C_4 &\subset  &\ldots &\subset & R_u\cap R_v
\end{matrix}
\]

\begin{lmma}\label{comm square}
In the following diagram
\[
\begin{matrix}
A_0 &\subset & A_2 &\subset  & A_4 &\subset &\ldots &\subset & R \cr
\cup &  &\cup &  &\cup &  & \ldots &  & \cup \cr
C_0 &\subset & C_2 &\subset & C_4 &\subset &\ldots &\subset & R_u\cap R_v
\end{matrix}
\]
each individual quadruple is a commuting square.
\end{lmma}
\begin{prf}
Choose any $k\in\bbn$ and consider the following quadruple
\[
\begin{matrix}
A_{2k-2} &\subset & A_{2k} \cr
\cup &  &\cup \cr
C_{2k-2} &\subset & C_{2k}\,.
\end{matrix}
\]
Now, draw the cube $\mathscr{C}_k$ as described in \Cref{com10}.
\begin{figure}[!h]
\begin{center} \resizebox{7 cm}{!}{\includegraphics{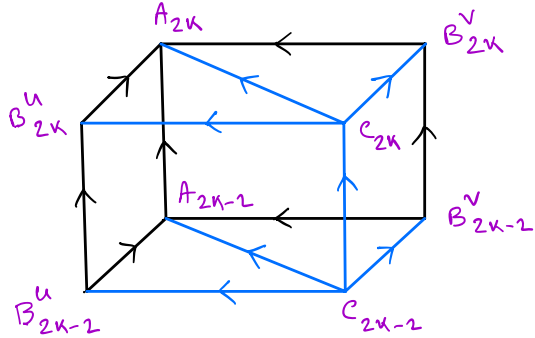}}\end{center}
\caption{Commuting cube $\mathscr{C}_k$}\label{com10}
\end{figure}
Here, the adjacent faces $(B^u_{2k-2}\subset A_{2k-2}\,,\,B^u_{2k}\subset A_{2k})$ and $(B^v_{2k-2}\subset A_{2k-2}\,,\,B^v_{2k}\subset A_{2k})$ are commuting squares by the construction of spin model subfactors $R_u\subset R$ and $R_v\subset R$ respectively. Since the required quadruple is the slice in \Cref{com10}, by \Cref{comm cube} it is a commuting square.\qed
\end{prf}

\begin{thm}\label{factoriality}
The von Neumann subalgebra $R_u\cap R_v$ of $R$ is a $II_1$ subfactor with $[R:R_u\cap R_v]=4$ and $[R_u:R_u\cap R_v]=[R_v:R_u\cap R_v]=2$.
\end{thm}
\begin{prf}
By \Cref{w}, we have $C_0=\bbc$ and $C_2=\mbox{Ad}_{u_2W_2}(\bbc\otimes M_2)$. In \Cref{comm square}, putting $k=1$ we have the following commuting square denoted by the symbol $\mathscr{S}$:
\begin{IEEEeqnarray}{lCl}\label{00aa}
\begin{matrix}
\bbc\otimes M_2 &\subset & M_2\otimes M_2 \cr
\cup &  &\cup \cr
\bbc &\subset & \mbox{Ad}_{u_2W_2}(\bbc\otimes M_2)
\end{matrix}
\end{IEEEeqnarray}
which is the slice in the commuting cube $\mathscr{C}_1$ as described in \Cref{com10}. Observe that the norm of the inclusion matrices of both the vertical embeddings in $\mathscr{S}$ are same and equal to $2$. Therefore, by Lemma 3.10 in \cite{BG}, we see that $\mathscr{S}$ is a non-degenerate commuting square. In other words, the slice in the commuting cube $\mathscr{C}_1$ is non-degenerate, and hence the commuting cube $\mathscr{C}_1$ is a non-degenerate commuting cube (see \Cref{Sec 2}). The tower of the basic construction for the the horizontal inclusion $\mathbb{C}\otimes M_2\subset M_2\otimes M_2$ in $\mathscr{S}$ is given by 
\begin{IEEEeqnarray}{lCl}\label{bakshi}
\mathbb{C}\otimes M_2 \subset M_2\otimes M_2 \subset^{f_1} {M}^{(2)}_2 \otimes M_2\subset^{f_2}  {M}^{(3)}_2\otimes M_2\subset\cdots \subset^{f_k} {M}^{(k)}_2\otimes M_2\subset\cdots\quad
\end{IEEEeqnarray}
where $f_k=\frac{1}{2}\sum_{i,j=1}^{2}E_{ij}\otimes E_{ij}\otimes I^{(k)}_2\in M^{(k+1)}_2\otimes M_2=M^{(k+2)}_2$ is the Jones projection. Indeed, by our construction in \Cref{fig1}, we have the following tower
\[
\bbc\subset\Delta_2\subset^{\,e_1} M_2\subset^{\,e_2}\Delta_2\otimes M_2\subset^{\,e_1\otimes I_2\,} M_2\otimes M_2\subset^{\,e_2\otimes I_2\,}\Delta_2\otimes M_2\otimes M_2\subset\cdots\cdots\subset R
\]
for the hyperfinite factor $R$, where $e_1=J_2,\,e_2=\mbox{diag}\{E_{11},E_{22}\}$. Applying Theorem $3.5$ in \cite{B1}, for instance, we see the above form of $f_k$. We put $L_{-1}=\bbc$ and $L_0= \mbox{Ad}_{u_2W_2}(\bbc\otimes M_2)$ and for $k \geq 1$, suppose that $L_k=\{L_{k-1},f_k\}^{\dprime}$. As $\mathscr{S}$ is a non-degenerate commuting square, it follows that the following tower of algebras is the Jones' basic construction tower:
\[
\mathbb{C}\subset L_0\subset L_1\subset^{f_2} L_2 \subset \cdots \subset^{f_k} L_k\subset\cdots
\]
Define $R_{u,v}=\overline{\bigcup_k L_k}^{\tiny\, \mbox{SOT}}$. Note that $R_{u,v}$ is a $II_1$ subfactor of $R$ with $[R:R_{u,v}]=4$ (see Corollary $5.7.4$ in \cite{JS}, for example). That the tower in \Cref{bakshi} indeed gives us $R$ follows as an application of Theorem $3.5$ in \cite{B1}. Thus, the quadruple $(R_{u,v}\subset R_u,R_v\subset R)$ is obtained as iterated basic construction of the non-degenerate commuting cube (see \Cref{Sec 2}) $\mathscr{C}_1$. Notice that $R_{u,v}\subset R_u\cap R_v$ because $f_k\in R_u\cap R_v$ (by \Cref{bakproc}) and $\mbox{Ad}_{u_2W_2}(\bbc\otimes M_2)\subset R_u\cap R_v\,$. Since $[R:R_{u,v}]=4$, thanks to the multiplicativity of Jones index, we must have $[R_u:R_{u,v}]=2$. Hence, $R_{u,v}\subset R_u$ is irreducible. Since $R_{u,v}\subset R_u\cap R_v\subset R_u$, we conclude that $R_u\cap R_v$ is also a type $II_1$ factor, and $R_{u,v}=R_u\cap R_v$ as $[R_u\cap R_v:R_{u,v}]=1$.\qed
\end{prf}

\begin{crlre}\label{comm cube example}
The quadruple $(R_u\cap R_v\subset R_u,R_v\subset R)$ of $II_1$ factors is obtained as an iterated basic construction of the non-degenerate commuting cube $\mathscr{C}_1$ described in \Cref{com10}. The iterated basic construction of the non-degenerate commuting cube $\mathscr{C}_1$ at the $k$-th step coincides with the commuting cube $\mathscr{C}_k$ for all $k\in\bbn$. 
\end{crlre}
\begin{prf}
In \Cref{factoriality}, the quadruple $(R_{u,v}\subset R_u,R_v\subset R)$ of $II_1$ factors is obtained as an iterated basic construction of the non-degenerate commuting cube $\mathscr{C}_1$, and since $R_{u,v}=R_u\cap R_v$, the first part follows. For the second part, we need to show that the following grid of (finite-dimensional) non-degenerate commuting squares
\[
\begin{matrix}
M_2 &\subset & M_2\otimes M_2 &\subset^{\,f_1} & \cdots &\subset & M_2\otimes M_2^{(k)} & \subset^{\,f_k} & \cdots & \subset R\\
\cup & &\cup & & & &\cup &  & & \cup\\
\bbc &\subset & \mathrm{Ad}_{u_2W_2}(\bbc\otimes M_2) &\subset & \cdots &\subset & \mathrm{Ad}_{u_{2k}W_{2k}}(\bbc\otimes M_2^{(k)}) & \subset & \cdots & \subset R_u\cap R_v
\end{matrix}
\]
is the tower of basic constructions for the subfactor $R_u\cap R_v\subset R$, where $f_k=\frac{1}{2}\sum_{i,j=1}^{2}E_{ij}\otimes E_{ij}\otimes I^{(k)}_2\in M^{(k+1)}_2\otimes M_2=M^{(k+2)}_2$ is the Jones projection. Recall from the proof of \Cref{factoriality} that the tower $L_{-1}\subset L_0\subset L_1\subset\cdots\subset R_{u,v}$ of finite-dimensional algebras is the tower of basic constructions for the subfactor $R_{u,v}$ and $R_{u,v}=R_u\cap R_v$. We show that $L_k=C_{2k+2}$ for all $k\geq 0$. In the following grid of quadruples
\[
\begin{matrix}
A_0 &\subset & A_2 &\subset^{\,f_1} & A_4 &\subset^{\,f_2} & A_6 &\subset^{\,f_3}  \cdots &\subset & R\\
\cup & & \cup & & \cup & & \cup &  & & \cup\\
C_0 &\subset & C_2 &\subset & C_4 &\subset & C_6 & \cdots &\subset & R_u\cap R_v\\
|| & & || & & \rotatebox{90}{$\subseteq$} & & \rotatebox{90}{$\subseteq$} & & & ||\\
L_{-1} &\subset & L_0 &\subset & L_1 &\subset & L_2 & \subset  \cdots &\subset & R_{u,v}\\
\end{matrix}
\]
where $f_k=\frac{1}{2}\sum_{i,j=1}^{2}E_{ij}\otimes E_{ij}\otimes I^{(k)}_2\in A_{2k+2}$, we have the entire quadruple $(L_k\subset A_{2k+2},L_{k+1}\subset A_{2k+4})$ and the upper quadruple $(C_{2k}\subset A_{2k},C_{2k+2}\subset A_{2k+2})$ are commuting squares for each $k\geq 0$. Therefore, the lower quadruple $(L_k\subset C_{2k+2},L_{k+1}\subset C_{2k+4})$ is also a commuting square. Indeed, we see that $E^{C_{2k+4}}_{C_{2k+2}}\big(L_{k+1}\big)=E^{A_{2k+4}}_{A_{2k+2}}\big(L_{k+1}\big)\subseteq L_k$. By \cite{PP}, we have $\lambda(R_u\cap R_v,R_{u,v})$ is the limit of the decreasing sequence $\{\lambda(C_{2k+2},L_k)\}_{k\geq 0}$. Since $R_u\cap R_v=R_{u,v}$, we have $\lambda(R_u\cap R_v,R_{u,v})=1$, and hence $\lambda(C_{2k+2},L_k)=1$ for all $k\geq 0$. Thus, $C_{2k+2}=L_k$ for all $k\geq 0$ by \Cref{sk}, which concludes the proof.\qed
\end{prf}

\begin{thm}[Characterization of commuting square]\label{kl}
The interior and the exterior angle in the quadruple of $II_1$ factors $(R_u\cap R_v\subset R_u,R_v\subset R)$ are equal and given by the following,
\[
\cos\big(\alpha^{R_u\cap R_v}_R(R_u,R_v)\big)=\cos\big(\beta^{R_u\cap R_v}_R(R_u,R_v)\big)=\cos^2(\alpha-\beta)\,.
\]
The quadruple $(R_u\cap R_v\subset R_u,R_v\subset R)$ is a commuting (and consequently, co-commuting) square if and only if $\,\alpha-\beta=\pm\frac{\pi}{2}$. Here, $\alpha\mbox{ and }\beta$ are associated with the complex Hadamard matrices $u\mbox{ and }v$ as described in \Cref{reduc}.
\end{thm}
\begin{prf}
Recall from \Cref{comm cube example} that the quadruple $(R_u\cap R_v\subset R_u,R_v\subset R)$ is obtained as an iterated basic construction of the non-degenerate commuting cube $\mathscr{C}_1$ described in \Cref{com10}, and further the basic construction of $\mathscr{C}_1$ at the $k$-th step coincides with $\mathscr{C}_k$ for all $k\in\bbn$. Stacking all the commuting cubes $\mathscr{C}_k$, we have a commuting cube $\mathscr{C}_\infty$ described in \Cref{com20}.
\begin{figure}[!h]
\begin{center} \resizebox{7 cm}{!}{\includegraphics{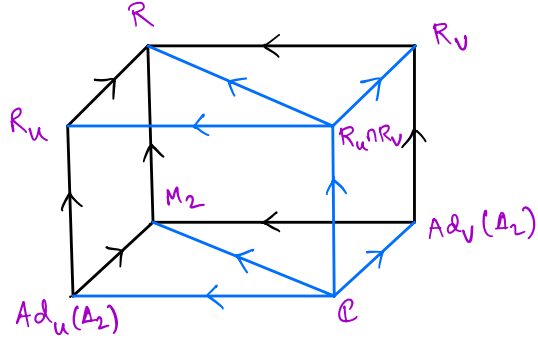}}\end{center}
\caption{Commuting cube $\mathscr{C}_\infty$}\label{com20}
\end{figure}
The adjacent faces in $\mathscr{C}_\infty$ are the commuting squares $(\mbox{Ad}_u(\Delta_2)\subset M_2\,,\,R_u\subset R)$ and $(\mbox{Ad}_v(\Delta_2)\subset M_2\,,\,R_v\subset R)$. Applying \Cref{comm cube}, we see that the quadruples $(\bbc\subset\mbox{Ad}_u(\Delta_2),R_u\cap R_v\subset R_u)$ (which is one of the remaining faces) and $(\bbc\subset M_2,R_u\cap R_v\subset R)$ (which is the slice) are also commuting squares. It can now be easily checked that $\{\lambda_1:=\mbox{Ad}_u(\sqrt{2}E_{11}),\lambda_2:=\mbox{Ad}_u(\sqrt{2}E_{22})\}$ is a basis for $R_u/R_u\cap R_v$, and similarly $\{\mu_1:=\mbox{Ad}_v(\sqrt{2}E_{11}),\mu_2:=\mbox{Ad}_v(\sqrt{2}E_{22})\}$ is a basis for $R_v/R_u\cap R_v$. Since $\lambda_i,\mu_i\in M_2$ for $i=1,2$, we have $E^R_{R_u\cap R_v}(\lambda_i^*\mu_j)=E^{M_2}_\bbc(\lambda_i^*\mu_j)$ for $i,j\in\{1,2\}$. This is because the quadruple $(\bbc\subset M_2,R_u\cap R_v\subset R)$ is a commuting square. Now, by the formula of the interior angle in \Cref{to refer}, we get the following,
\begin{IEEEeqnarray*}{lCl}
\cos\big(\alpha^{R_u\cap R_v}_R(R_u,R_v)\big) &=& \sum_{i,j=1}^2\,tr\big(E^R_{R_u\cap R_v}(\lambda_i^*\mu_j)\mu_j^*\lambda_i\big)-1\\
&=& \sum_{i,j=1}^2\,tr\big(E^{M_2}_\bbc(\lambda_i^*\mu_j)\mu_j^*\lambda_i\big)-1\\
&=& \sum_{i,j=1}^2\,tr(\lambda_i^*\mu_j)\,tr(\mu_j^*\lambda_i)-1\\
&=& \cos^2(\alpha-\beta)\,,
\end{IEEEeqnarray*}
where the last equality is a simple verification using the form of $u\mbox{ and }v$ described in \Cref{reduc}. Finally, recall from \cite{BDLR} that if $(N\subset P,Q\subset M)$ is a quadruple of $II_1$ factors such that $N\subset M$ is extremal with $[M:N]<\infty$ and $[P:N]=[M:Q]$, then $\alpha^N_M(P,Q)=\beta^N_M(P,Q)$, that is, the interior and the exterior angle coincide. In our case, $R_u\cap R_v=R_{u,v}$ (\Cref{factoriality}) is extremal and all the other hypotheses are also satisfied along with it. Hence, we get the first part of the theorem.

Finally, the quadruple $(R_u\cap R_v\subset R_u,R_v\subset R)$ is a commuting (and consequently, co-commuting) square if and only if $\cos\big(\alpha^{R_u\cap R_v}_R(R_u,R_v)\big)=0$, which happens if and only if $\alpha-\beta=\pm\frac{\pi}{2}$.\qed
\end{prf}

\begin{ppsn}\label{comm and cocomm}
The quadruple of $II_1$ factors
\[
\begin{matrix}
R_u &\subset & R\cr
\rotatebox{90}{$\subset$} &\ &\rotatebox{90}{$\subset$}\cr
R_{u}\cap R_v &\subset & R_v
\end{matrix}
\]
is always non-degenerate, that is $\overline{R_uR_v}=\overline{R_vR_u}=R$.
\end{ppsn}
\begin{prf}
For any $k\in\bbn\cup\{0\}$, we claim that the following holds
\begin{IEEEeqnarray}{lCl}\label{bbbb}
u_{2k}\big(\Delta_2\otimes M_2^{(k)}\big)u_{2k}^*v_{2k}\big(\Delta_2\otimes M_2^{(k)}\big)v_{2k}^*=M_2\otimes M_2^{(k)}\,.
\end{IEEEeqnarray}
First consider $k=0$. We have
\[
u^*v=\frac{1}{2}(1+e^{i(\beta-\alpha)})I_2+\frac{1}{2}(1-e^{i(\beta-\alpha)})\sigma_1
\]
where $\beta\neq\alpha,\pi\pm\alpha$ (\Cref{reduc} and discussion after it). Hence, none of the entries of the unitary matrix $u^*v$ are zero. Observe that if any unitary matrix $U$ in $M_2$ is such that none of the entries of it are zero, then $\Delta_2U\Delta_2=M_2$. Therefore, we have $\Delta_2u^*v\Delta_2=M_2$, and consequently $u\Delta_2u^*v\Delta_2v^*=M_2$, which settles Eqn. \ref{bbbb} for $k=0$. For $k\in\bbn$, by \Cref{e9} and \Cref{not} we have the following,
\[
{u}^*_{2k}{v}_{2k}\big(\Delta_2\otimes M^{(k)}_2\big){v}^*_{2k}{u}_{2k} =\mbox{bl-diag}\big\{I^{(k)}_2,\sigma_1^{(k)}\big\}\big(u^*v\Delta_2v^*u\otimes M^{(k)}_2\big)\mbox{bl-diag}\big\{I^{(k)}_2,\sigma_1^{(k)}\big\}.
\]
Therefore,
\begin{IEEEeqnarray*}{lCl}
&  & u_{2k}\big(\Delta_2\otimes M_2^{(k)}\big)u_{2k}^*v_{2k}\big(\Delta_2\otimes M_2^{(k)}\big)v_{2k}^*\\
&=& {u}_{2k}\left(\big(\Delta_2 \otimes M^{(k)}_{2}\big)\mbox{bl-diag}\big\{I^{(k)}_2,\sigma_1^{(k)}\big\}\right)\big(u^*v\Delta_2v^*u\otimes M^{(k)}_2\big)\mbox{bl-diag}\big\{I^{(k)}_2,\sigma_1^{(k)}\big\}{u}_{2k}^*\\
&=& \mbox{Ad}_{{u}_{2k}}\left(\big(\Delta_2 \otimes M^{(k)}_{2}\big)\big(u^*v\Delta_2 v^*u\otimes M^{(k)}_2\big)\mbox{bl-diag}\big\{I^{(k)}_2,\sigma_1^{(k)}\big\}\right)\\
&=& \mbox{Ad}_{{u}_{2k}}\left(\big(\Delta_2u^*v\Delta_2 v^*u\otimes M^{(k)}_2\big)\mbox{bl-diag}\big\{I^{(k)}_2,\sigma_1^{(k)}\big\}\right)\\
&=& \mbox{Ad}_{{u}_{2k}}\left(\big(M_2\otimes M^{(k)}_2\big)\mbox{bl-diag}\big\{I^{(k)}_2,\sigma_1^{(k)}\big\}\right)\qquad(\mbox{the $k=0$ case})\\
&=& \mbox{Ad}_{{u}_{2k}}\big(M_2\otimes M^{(k)}_2\big).
\end{IEEEeqnarray*}
Here in the above, the second last equality follows from the $k=0$-step, that is, $u\Delta_2u^*v\Delta_2v^*=M_2$. Since ${u}_{2k}\in M_2\otimes M^{(k)}_2$, we finally have the following,
\[
u_{2k}\big(\Delta_2\otimes M_2^{(k)}\big)u_{2k}^*v_{2k}\big(\Delta_2\otimes M_2^{(k)}\big)v_{2k}^*=M_2\otimes M^{(k)}_2
\]
for all $k\in\bbn\cup\{0\}$. Hence, we conclude  that $\overline{R_uR_v}^{\tiny\,\mbox{SOT}}=R$. Indeed, since $R_u=\overline{\bigcup_i\,B_{2i}^u}^{\tiny\,\mbox{SOT}}$ and $R_v=\overline{\bigcup_j\,B_{2j}^v}^{\tiny\,\mbox{SOT}},$ we have 
\[
R=\overline{\cup_i\,B_{2i}^uB_{2i}^v}\subseteq\overline{\cup_{i,j}\,B_{2i}^uB_{2j}^v}\subseteq\overline{(\cup_i\,B_{2i}^u)(\cup_j\,B_{2j}^v)}\subseteq\overline{\overline{(\cup_i\,B_{2i}^u)}\,\overline{(\cup_j\,B_{2j}^v)}}=\overline{R_uR_v}\subseteq R\,,
\]
where the closure being taken in the SOT topology. Therefore, the given quadruple is always non-degenerate.\qed
\end{prf}

\begin{rmrk}\rm\label{to next}
If the condition of commuting square is given for a quadruple of $II_1$ factors $(N\subset P,Q\subset M)$ with $[M:N]<\infty$, then the quadruple is non-degenerate if and only if it is a co-commuting square (see \cite{SW}). In the absence of commuting square, if irreducibility of $N\subset M$ is given, then also this holds (see \cite{GJ}). In \Cref{comm and cocomm}, we have a family of quadruples of subfactors indexed by $2\times 2$ complex Hadamard matrices $u\mbox{ and }v$ such that $u\nsim v$, where `non-degeneracy' and `co-commuting' are not the same. To the best of our knowledge, this is the first concrete example in the literature (see Remark $2.3$ in this regard).
\end{rmrk}


\subsection{Characterization of the quadruple $(R_u\cap R_v\subset R_u,R_v\subset R)$}

Due to \Cref{kl}, \Cref{comm and cocomm}, and \Cref{to next}, it is clear that the subfactor $R_u\cap R_v\subset R$ is not irreducible, at least when $\alpha-\beta\neq\pm\frac{\pi}{2}$. We obtain its commutant in $R$, and characterize the subfactor $R_u\cap R_v\subset R$. As a curious fact, we obtain that although we have started with two spin model subfactors $R_u,\,R_v\subset R$, their intersection is a vertex model subfactor. Furthermore, we characterize the quadruple $(R_u\cap R_v\subset R_u,R_v\subset R)$ and show that it is generated (in the sense of iterated basic construction) by a single bi-unitary matrix in $M_2\otimes M_2$.

\begin{thm}\label{vertex}
For any $2\times 2$ complex Hadamard matrices $u\mbox{ and }v$ such that $u\nsim v$, the subfactor $R_u\cap R_v\subset R$ is a vertex model subfactor of Jones index $4$.
\end{thm}
\begin{prf}
We have the symmtric commuting square (\Cref{factoriality} and \Cref{comm cube example}) $\mathscr{S}$:
\[
\begin{matrix}
\mbox{Ad}_{u_2W_2}(\bbc\otimes M_2) &\subset & M_2\otimes M_2 \cr
\cup &  &\cup \cr
\bbc &\subset & \bbc\otimes M_2
\end{matrix}
\]
which is same as the following commuting square
\[
\begin{matrix}
\mbox{Ad}_{u_2W_2V_2}(M_2\otimes\bbc) &\subset & M_2\otimes M_2 \cr
\cup &  &\cup \cr
\bbc &\subset & \bbc\otimes M_2
\end{matrix}
\]
where $\mbox{Ad}_{V_2}:M_2\otimes\bbc\longmapsto\bbc\otimes M_2$ is the flip given by the following permutation matrix
\[
V_2:=\begin{bmatrix}
1 & 0 & 0 & 0\\
0 & 0 & 1 & 0\\
0 & 1 & 0 & 0\\
0 & 0 & 0 & 1
\end{bmatrix}=\begin{bmatrix}
E_{11} & E_{21}\\
E_{12} & E_{22}
\end{bmatrix}
\]
in $M_4$. It is now clear that $u_2W_2V_2$ is a bi-unitary in $M_4=M_2\otimes M_2$ (can be checked directly also) and $R_u\cap R_v\subset R$ is a vertex model subfactor of index $4$. A simple verification will show the following,
\begin{IEEEeqnarray}{lCl}\label{qqz}
u_2W_2V_2=\frac{1}{\sqrt{2}}\begin{bmatrix}
u & u\sigma_1\\
u & -u\sigma_1
\end{bmatrix}
\end{IEEEeqnarray}
where $\sigma_1$ is Pauli spin matrix.\qed
\end{prf}

\begin{thm}\label{relativecommutant}
The relative commutant ${(R_u\cap R_v)}^\prime\cap R$ is $\bbc\oplus\bbc$.
\end{thm}
\begin{prf}
Recall from \Cref{comm cube example} that the following grid of (finite-dimensional) non-degenerate commuting squares
\[
\begin{matrix}
M_2 &\subset & M_2\otimes M_2 &\subset^{\,f_1} & \cdots &\subset & M_2\otimes M_2^{(k)} & \subset^{\,f_k} & \cdots & \subset R\\
\cup & &\cup & & & &\cup &  & & \cup\\
\bbc &\subset & \mathrm{Ad}_{u_2W_2}(\bbc\otimes M_2) &\subset & \cdots &\subset & \mathrm{Ad}_{u_{2k}W_{2k}}(\bbc\otimes M_2^{(k)}) & \subset & \cdots & \subset R_u\cap R_v
\end{matrix}
\]
is the tower of basic constructions for the subfactor $R_u\cap R_v\subset R$. By the Ocneanu compactness theorem (see \Cref{ocneanucompactness}), we have the following
\begin{IEEEeqnarray}{lCl}\label{s7}
{(R_u\cap R_v)}^\prime\cap R=\left(\mbox{Ad}_{u_2W_2}(\bbc\otimes M_2)\right)^\prime\cap\left(\bbc\otimes M_2\right)\,.
\end{IEEEeqnarray}
Then, \Cref{01} implies the following
\begin{IEEEeqnarray}{lCl}\label{88}
{(R_u\cap R_v)}^\prime\cap R=\mbox{Ad}_{u_2W_2}\left((\bbc\otimes M_2)^\prime\cap\mbox{Ad}_{W_2^*u_2^*}(\bbc\otimes M_2)\right)\,,
\end{IEEEeqnarray}
where $W_2$ is as defined in \Cref{not}. Now, for any $A\in M_2$ we have the following using \Cref{1st basic},
\begin{IEEEeqnarray*}{lCl}
\mbox{Ad}_{W_2^*u_2^*}(I_2\otimes A)
&=& \frac{1}{2}\left[{\begin{matrix}
(\mbox{Ad}_{u^*}+\mbox{Ad}_{\sigma_3u^*})(A) & \big((\mbox{Ad}_{u^*}-\mbox{Ad}_{\sigma_3u^*})(A)\big)\sigma_1\\
\sigma_1\big((\mbox{Ad}_{u^*}-\mbox{Ad}_{\sigma_3u^*})(A)\big) & \mbox{Ad}_{\sigma_1}\big((\mbox{Ad}_{u^*}+\mbox{Ad}_{\sigma_3u^*})(A)\big)\\
\end{matrix}}\right]\,.
\end{IEEEeqnarray*}
Since $\left(\bbc\otimes M_2\right)^\prime\cap (M_2\otimes M_2)=M_2\otimes\bbc$, we immediately see that if any $B\otimes I_2$ in $M_2\otimes\bbc\subseteq M_2\otimes M_2$ has to lie in $\mbox{Ad}_{W_2^*u_2^*}(\bbc\otimes M_2)$, then it is necessary that $B=\left[{\begin{smallmatrix}
z & w\\
w & z\\
\end{smallmatrix}}\right]$ for some $z,w\in\bbc$. This proves the following,
\begin{IEEEeqnarray}{lCl}\label{a0a}
(\bbc\otimes M_2)^\prime\cap\mbox{Ad}_{W_2^*u_2^*}(\bbc\otimes M_2)\subseteq\left\{\left[{\begin{smallmatrix}
z & w\\
w & z\\
\end{smallmatrix}}\right]\otimes I_2\,:\,z,w\in\bbc\right\}\,.
\end{IEEEeqnarray}
Conversely, for any $B=\left[{\begin{smallmatrix}
z & w\\
w & z\\
\end{smallmatrix}}\right]\otimes I_2$ consider $C=\mbox{Ad}_u\left(\left[{\begin{smallmatrix}
z & w\\
w & z\\
\end{smallmatrix}}\right]\right)\in M_2$ and observe that $\mbox{Ad}_{W_2^*u_2^*}(I_2\otimes C)=B$. Hence, the inclusion in \Cref{a0a} is in fact equality. By \Cref{88}, using \Cref{1st basic} we get that
\begin{IEEEeqnarray}{lCl}\label{ez}
{(R_u\cap R_v)}^\prime\cap R &=& \left\{\mbox{Ad}_{u_2W_2}\left(\left[{\begin{smallmatrix}
z & w\\
w & z\\
\end{smallmatrix}}\right]\otimes I_2\right)\,:\,z,w\in\bbc\right\}\nonumber\\
&=& \left\{\mbox{Ad}_{u_2}\left(\mbox{diag}\{I_2,\sigma_1\}\,\left(\left[{\begin{smallmatrix}
z & w\\
w & z\\
\end{smallmatrix}}\right]\otimes I_2\right)\mbox{diag}\{I_2,\sigma_1\}\right)\,:\,z,w\in\bbc\right\}\nonumber\\
&=& \left\{I_2\otimes\mbox{Ad}_u\left(\left[{\begin{smallmatrix}
z & w\\
w & z\\
\end{smallmatrix}}\right]\right)\,:\,z,w\in\bbc\right\}\\
&=& \left\{I_2\otimes\mbox{Ad}_{uF_2}\left(\left[{\begin{smallmatrix}
z+w & \\
 & z-w\\
\end{smallmatrix}}\right]\right)\,:\,z,w\in\bbc\right\}\,,\nonumber
\end{IEEEeqnarray}
where $F_2$ is the $2\times 2$ Fourier matrix. Therefore, ${(R_u\cap R_v)}^\prime\cap R=\bbc\oplus\bbc$, where the embedding $\bbc\oplus\bbc\xhookrightarrow{} M_2$ is given by $(\gamma,\delta)\mapsto\mbox{Ad}_{uF_2}(\mbox{diag}\{\gamma,\delta\})$.\qed
\end{prf}

\begin{thm}[Characterization of the quadruple $(R_u\cap R_v\subset R_u,R_v\subset R)$]\label{fullchar}
\begin{enumerate}[$(i)$]
\item[] \mbox{}
\item The pair of subfactors $R_u$ and $R_v$ of the hyperfinite type $II_1$ factor $R$ are conjugate to each other via a unitary in ${(R_u\cap R_v)}^\prime\cap R$.
\item The bi-unitary matrix $\,u_2W_2V_2$ in $M_4$ generates (in the sense of iterated basic construction of finite-dimensional grids) the composition of subfactors $R_u\cap R_v\subset R_u\subset R$, and consequently the quadruple $(R_u\cap R_v\subset R_u,R_v\subset R)$.
\end{enumerate}
\end{thm}
\begin{prf}
Since $u=uv^*v$, we have $B^u_0=\mbox{Ad}_u(\Delta_2)=\mbox{Ad}_{uv^*}\mbox{Ad}_v(\Delta_2)=\mbox{Ad}_{uv^*}(B^v_0)$. We claim that $u_{2k}=(I_2^{(k)}\otimes uv^*)v_{2k}$ for all $k\geq 0$. This follows by induction on $k$, since by \Cref{1st basic}, we have the following,
\begin{IEEEeqnarray*}{lCl}
u_{2k+2} &=& (I_2\otimes u_{2k})\,\mbox{bl-diag}\{I_2^{(k)},\eta_k\}\,\big(\xi_k\otimes I_2^{(k)}\big)\\
&=& (I_2\otimes I_2^{(k)}\otimes uv^*)(I_2\otimes v_{2k})\,\mbox{bl-diag}\{I_2^{(k)},\eta_k\}\,\big(\xi_k\otimes I_2^{(k)}\big)\\
&=& (I_2^{(k+1)}\otimes uv^*)v_{2k+2}\,.
\end{IEEEeqnarray*}
Therefore, we have $B^u_{2k}=\mbox{Ad}_{uv^*}(B^v_{2k})$, where the unitary $uv^*\in M_2$ is identified with $I_2^{(k)}\otimes uv^*\in M_2\otimes M_2^{(k)}=A_{2k}$, and consequently $\cup_{k\geq 0}B^u_{2k}=\mbox{Ad}_{uv^*}\big(\cup_{k\geq 0}B^v_{2k}\big)$. Since $R_u=\overline{\cup_k\,B^u_{2k}}^{\mbox{ sot}}$ and $R_v=\overline{\cup_k\,B^v_{2k}}^{\mbox{ sot}}$,  we have $R_u=\mbox{Ad}_{uv^*}(R_v)$, where the unitary $uv^*$ is the diagonal matrix $\mbox{diag}\{1,e^{i(\alpha-\beta)}\}$ (see \Cref{reduc}) in $M_2$. Now, $uv^*=\mbox{Ad}_u(v^*u)$ and the matrix $v^*u$ is the following
\[
v^*u=\frac{1}{2}\begin{bmatrix}
1+e^{i(\alpha-\beta)} & 1-e^{i(\alpha-\beta)}\\
1-e^{i(\alpha-\beta)} & 1+e^{i(\alpha-\beta)}
\end{bmatrix}\,.
\]
Hence by \Cref{relativecommutant} (in particularly, \Cref{ez}), we see that $uv^*\in{(R_u\cap R_v)}^\prime\cap R$. This completes part $(i)$.

For part $(ii)$, using the fact that $W_2=E_{11}\otimes I_2+E_{22}\otimes\sigma_1$ (see \Cref{not}) lies in $\Delta_2\otimes M_2$ we have $B^u_2=\mbox{Ad}_{u_2}(\Delta_2\otimes M_2)=\mbox{Ad}_{u_2W_2}(\Delta_2\otimes M_2)$. Observe that $\mbox{Ad}_{V_2}:M_2\otimes\Delta_2\mapsto\Delta_2\otimes M_2$, where $V_2$ is the flip as in \Cref{vertex}, and hence $B^u_2=\mbox{Ad}_{u_2W_2V_2}(M_2\otimes \Delta_2)$. By induction, it is easy to verify that for all $k\in\bbn$, if we consider the following unitary matrix
\[
V_{2k}:=\prod_{j=0}^{k-1}\big(I_2^{(j)}\otimes V_2\otimes I_2^{(k-j-1)}\big)
\]
in $M_2^{(k+1)}=M_2\otimes M_2^{(k)}$, then $\mbox{Ad}_{V_{2k}}:M_2^{(k)}\otimes\Delta_2\longmapsto\Delta_2\otimes M_2^{(k)}$. Since $W_{2k}=E_{11}\otimes I_2^{(k)}+E_{22}\otimes\sigma_1^{(k)}\in\Delta_2\otimes M_2^{(k)}$, we have
\[
B^u_{2k}=\mbox{Ad}_{u_{2k}}(\Delta_2\otimes M_2^{(k)})=\mbox{Ad}_{u_{2k}W_{2k}}(\Delta_2\otimes M_2^{(k)})=\mbox{Ad}_{u_{2k}W_{2k}V_{2k}}(M_2^{(k)}\otimes\Delta_2).
\]
for all $k\in\bbn$. Therefore, we have the following grid of non-degenerate commuting squares
\[
\begin{matrix}
M_2 &\subset & M_2\otimes M_2 &\subset^{\,f_1} & \cdots &\subset & M_2\otimes M_2^{(k)} & \subset^{\,f_k} & \cdots & \subset R\\
\cup & &\cup & & & &\cup &  & & \cup\\
\mbox{Ad}_u(\Delta_2) &\subset & \mathrm{Ad}_{u_2W_2V_2}(M_2\otimes\Delta_2) &\subset & \cdots &\subset & \mathrm{Ad}_{u_{2k}W_{2k}V_{2k}}(M_2^{(k)}\otimes\Delta_2) & \subset & \cdots & \subset R_u\\
\cup & &\cup & & & &\cup &  & & \cup\\
\bbc &\subset & \mathrm{Ad}_{u_2W_2V_2}(M_2\otimes\bbc) &\subset & \cdots &\subset & \mathrm{Ad}_{u_{2k}W_{2k}V_{2k}}(M_2^{(k)}\otimes\bbc) & \subset & \cdots & \subset R_u\cap R_v\\
\end{matrix}
\]
where $u_2W_2V_2\in M_4$ (see \Cref{qqz}) is the bi-unitary obtained in \Cref{vertex}, and $f_k=\frac{1}{2}\sum_{i,j=1}^{2}E_{ij}\otimes E_{ij}\otimes I^{(k)}_2\in M^{(k+2)}_2,\,k\in\bbn,$ are the Jones' projections. Finally, since $R_u\mbox{ and }R_v$ are unitary conjugate, and the conjugation is implemented by a unitary belonging at the very first stage (i,e. $M_2$), we see that the single bi-unitary matrix $u_2W_2V_2$ in $M_4$ generates (in the sense of iterated basic construction) the quadruple $(R_u\cap R_v\subset R_u\,,\,R_v\subset R)$ of $II_1$ factors, and we have an explicit grid of finite-dimensional algebras for the entire quadruple as described above.\qed
\end{prf}

Before moving further, we pause for an important remark.

\begin{rmrk}\rm\label{sunder}
It is known that (see \cite{KSV}) any bi-unitary matrix in $M_4$ is equivalent to $W(\omega)=\mbox{diag}\{1,1,1,\omega\}$ for some $\omega\in\mathbb{S}^1$. Moreover, $W(\omega)$ is equivalent to $W(\omega^\prime)$ if and only if $\Re(\omega)=\Re(\omega^\prime)$. In our situation, using \Cref{1st basic} one can see that $u_2W_2V_2=(F_2\otimes u)\mbox{diag}\{I_2,\sigma_1\}$. Since $F_2\sigma_1F_2=\mbox{diag}\{1,-1\}=\sigma_3$, we have $u_2W_2V_2=(F_2\otimes u)(I_2\otimes F_2)\mbox{diag}\{I_2,\sigma_3\}(I_2\otimes F_2)=(F_2\otimes uF_2)W(-1)(I_2\otimes F_2)$, that is, the bi-unitary $u_2W_2V_2$ is equivalent to $W(-1)$.
\end{rmrk}

By \Cref{fullchar} (see also \Cref{sunder}) and Section $4$ in \cite{KSV}, we conclude that $R_u\cap R_v\subset R$ is of depth 2, and by \cite{Po2} we have the following characterization.
\begin{thm}\label{diag}
The subfactor $R_u\cap R_v\subset R$ is isomorphic to the diagonal subfactor 
\[
\left\{\begin{bmatrix}
x & 0\\
0 & \theta(x)
\end{bmatrix}:x\in R\,,\,\theta\in\mathrm{Out}(R)~\mathrm{such~ that}~{\theta}^2=id \right\}\subset M_2(\mathbb{C})\otimes R\,.
\]
Under this isomorphism  $R_u$ becomes $\big(\mathrm{diag}\{1,e^{i(\alpha-\beta)}\}\big)R_v\big(\mathrm{diag}\{1,e^{i(\beta-\alpha)}\}\big)$, where $R_v$ is given by
\[
\Bigg\{\begin{bmatrix}
x & y\\
\theta(y) & \theta(x)
\end{bmatrix} : x,y \in R\Bigg\}\,.
\]
Here, $\alpha\mbox{ and }\beta$ are associated with the complex Hadamard matrices $u\mbox{ and }v$ as described in \Cref{reduc}.
\end{thm}


\subsection{Relative entropy and Sano-Watatani angle between $R_u\mbox{ and }R_v$}

We investigate relative entropy and the Sano-Watatani angle between the pair of spin model subfactors $R_u\mbox{ and }R_v$.

\begin{thm}\label{major 2by2}
For the pair of spin model subfactors $R_u\mbox{ and }R_v$, we have the following.
\begin{enumerate}[$(i)$]
\item $H(R|R_u\cap R_v)=2\log 2\mbox{ and }H(R_u|R_u\cap R_v)=H(R_v|R_u\cap R_v)=\log 2$.
\item $h(R_u|R_v)=\eta\big(\cos^2\big(\frac{\alpha-\beta}{2}\big)\big)+\eta\big(\sin^2\big(\frac{\alpha-\beta}{2}\big)\big)\,.$
\item $\eta\big(\cos^2\big(\frac{\alpha-\beta}{2}\big)\big)+\eta\big(\sin^2\big(\frac{\alpha-\beta}{2}\big)\big)\leq H(R_u|R_v)\leq\log 2\,,$
and if $\,\alpha-\beta=\pm\frac{\pi}{2}$, then $H(R_u|R_v)=h(R_u|R_v)=\log 2$. Here, $\alpha\mbox{ and }\beta$ are associated with the complex Hadamard matrices $u\mbox{ and }v$ as in \Cref{reduc}
\end{enumerate}
\end{thm}
\begin{prf}
Corollary $4.6$ in \cite{PP} gives the value of $H(R|R_u\cap R_v),\,H(R_u|R_u\cap R_v)$ and $H(R_v|R_u\cap R_v)$. Since $H(R_u|R_v)\leq H(R_u|R)+H(R|R_v)=H(R|R_v)$, the upper bound for $H(R_u|R_v)$ is clear. Now, recall from \cite{choda} that $h(u\Delta_2u^*|v\Delta_2v^*)=h(\Delta_2|u^*v\Delta_2v^*u)=\frac{1}{2}\sum_{i,j=1}^2\eta |(u^*v)_{ij}|^2$. Each $|(u^*v)_{ij}|^2$ can easily be obtained from \Cref{reduc}, and since
\[
h(u\Delta_2u^*|v\Delta_2v^*)\leq H(u\Delta_2u^*|v\Delta_2v^*)\leq H(R_u|R_v)
\]
due to \Cref{popaadaptation2}, the lower bound for $H(R_u|R_v)$ is obtained. We also have
\[ 
h(R_u|R_v)\geq\eta\left(\cos^2\left(\frac{\alpha-\beta}{2}\right)\right)+\eta\left(\sin^2\left(\frac{\alpha-\beta}{2}\right)\right),
\]
as $h(R_u|R_v)\geq h(u\Delta_2u^*|v\Delta_2v^*)$. To obtain the reverse inequality, recall that by \Cref{fullchar} we have $R_u=\mbox{Ad}_{uv^*}(R_v)$, where $uv^*=\mathrm{diag}\{1, e^{i(\alpha-\beta)}\}$. Therefore, thanks to Theorem $3.4$ in \cite{choda2}, it is now a straightforward verification that
\[
 h(R_u|R_v)\leq\eta\left(\cos^2\left(\frac{\alpha-\beta}{2}\right)\right)+\eta\left(\sin^2\left(\frac{\alpha-\beta}{2}\right)\right).
 \]
Finally, for $\alpha-\beta=\pm\frac{\pi}{2}$ the lower bound for $H(R_u|R_v)$ also becomes $\log 2$, which completes the proof.\qed
\end{prf}

\begin{thm}\label{sanoangle}
The Sano-Watatani angle between the subfactors $R_u\mbox{ and }R_v$ is the singleton set $\{\arccos|\cos(\alpha-\beta)|\}$, where $\alpha\mbox{ and }\beta$ are associated with the complex Hadamard matrices $u\mbox{ and }v$ as described in \Cref{reduc}.
\end{thm}
\begin{prf}
For $\alpha-\beta=\pm\frac{\pi}{2}$, the result holds as $(R_u\cap R_v\subset R_u,R_v\subset R)$ is a commuting square by \Cref{kl}. Consider now $\alpha-\beta\neq\pm\frac{\pi}{2}$. Recall the commuting cube $\mathscr{C}_k$ introduced in \Cref{comm square}, in particularly, see \Cref{com10}. Consider the commuting cube $\mathscr{C}_1$ (for $k=1$). This is a non-degenerate commuting cube as shown in \Cref{factoriality}. The quadruple of $II_1$ factors $(R_u\cap R_v\subset R_u,R_v\subset R)$ is obtained as an iterated basic construction of this non-degenerate commuting cube (\Cref{comm cube example}), as illustrated in \Cref{combasic2} (recall \Cref{not} in this regard).
\begin{figure}[!h]
\begin{center} \resizebox{7 cm}{!}{\includegraphics{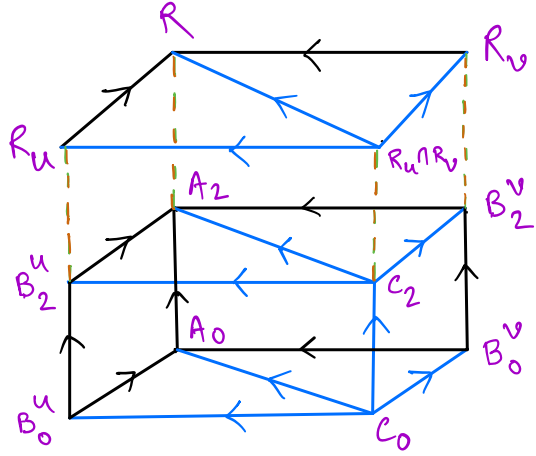}}\end{center}
\caption{Iterated basic construction of non-degenerate commuting cube}\label{combasic2}
\end{figure}
With the notation as in \Cref{reduc}, using \Cref{p23}, we get that the non-negative matrix $E_{B^u_0}^{A_0}E_{B^v_0}^{A_0}E_{B^u_0}^{A_0}-E_{C_0}^{A_0}$ satisfies the following,
\[
\big(E_{B^u_0}E_{B^v_0}E_{B^u_0}-E_{C_0}\big)^2=\cos^2(\alpha-\beta)\big(E_{B^u_0}E_{B^v_0}E_{B^u_0}-E_{C_0}\big)\,.
\]
Since we have the non-degenerate commuting cube $\mathscr{C}_\infty$ in \Cref{com20} (which is the final resulting commuting cube in \Cref{combasic2}), as an application of \Cref{baksata1} we get that the spectrum of the cosine of the angle operator $\sqrt{e^R_{R_u}e^R_{R_v}e^R_{R_u}-e^R_{R_u\cap R_v}}$ for the pair of subfactors $R_u,R_v\subset R$ is the singleton set $\{|\cos(\alpha-\beta)|\}$.\qed
\end{prf}


\newsection{Subfactors arising from \texorpdfstring{$4\times 4$}~ complex Hadamard matrices}\label{Sec 5}

\begin{figure}[!h]
	\begin{center}
		\psfrag{$R$}{ \large  $R$}
		\psfrag{d1}{\large$\Delta_4\otimes M_4\otimes M_4$}
		\psfrag{d2}{\large$M_4\otimes M_4$}
		\psfrag{d3}{\large$\Delta_4\otimes M_4$}
		\psfrag{d4}{\large$M_4$}
		\psfrag{d5}{\large$\Delta_4$}
		\psfrag{R_n}{\large$R_u$}
		\psfrag{R_v}{\large$R_v$}
		\psfrag{f_1}{\large$u_3(M_4\otimes M_4)u_3^*$}
		\psfrag{f_2}{\large$u_2(\Delta_4 \otimes M_4)u_2^*$}
		\psfrag{f_3}{\large$u_1M_4u_1^*$}
		\psfrag{f_4}{\large$u_0\Delta_4u_0^*$}
		\psfrag{f_5}{\large$\mathbb{C}$}
		\psfrag{e_1}{\large$v_3(M_4\otimes M_4)v_3^*$}
		\psfrag{e_2}{\large$v_2(\Delta_4 \otimes M_4)v_2^*$}
		\psfrag{e_3}{\large$v_1M_4v_1^*$}
		\psfrag{e_4}{\large$v_0\Delta_4v_0^*$}
		\psfrag{e_5}{\large$\mathbb{C}$}
		\resizebox{10cm}{!}{\includegraphics{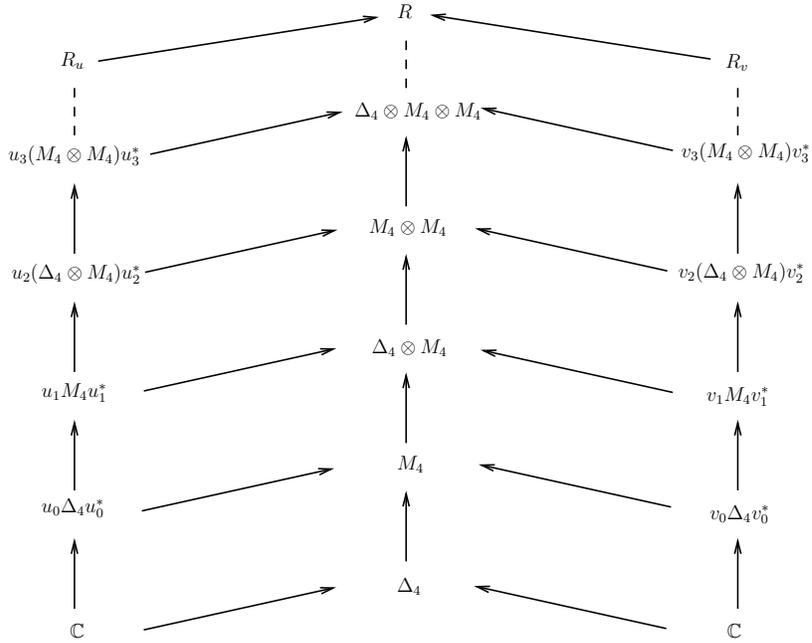}}
	\end{center}
\caption{ A pair of $4\times 4$ Hadamard matrices and basic constructions}\label{fig2}
\end{figure}

Recall that there exists a continuous, one-parameter family of inequivalent $4\times 4$ complex Hadamard matrices, and this family is given by the following unitary
\[
u(z)= \frac{1}{2}\,\left[{\begin{matrix}
1 & 1 & 1 & 1\\
1 & iz & -1 & -iz\\
1 & -1 & 1 & -1\\
1 & -iz & -1 & iz\\
\end{matrix}}\right]
\]
where $\,z=e^{i\phi}\in\mathbb{S}^1$ with $\phi\in[0,\pi)$. Consider two distinct elements of this family
\[
u=\frac{1}{2}\,\left[{\begin{matrix}
1 & 1 & 1 & 1\\
1 & ia & -1 & -ia\\
1 & -1 & 1 & -1\\
1 & -ia & -1 & ia\\
\end{matrix}}\right] \quad\mbox{and}\quad v=\frac{1}{2}\,\left[{\begin{matrix}
1 & 1 & 1 & 1\\
1 & ib & -1 & -ib\\
1 & -1 & 1 & -1\\
1 & -ib & -1 & ib\\
\end{matrix}}\right]
\]
where $a,b\in\mathbb{S}^1$ (semicircle to be precise) and $a\neq b$. Throughout the section, remember that $b\neq\pm a$, since $\phi\in[0,\pi)$. We obtain the pair of spin model subfactors $R_u\subset R$ and $R_v\subset R$ as described in \Cref{pairspin}. Since $u\mbox{ and }v$ are Hadamard inequivalent, we have $R_u\neq R_v$ due to \Cref{general thm}.  That is, we have a pair of irreducible subfactors each with index 4~:
\[\begin{matrix}
R_u &\subset & R\\
 & & \cup\\
 & & R_v\\
\end{matrix}
\]
The following
\[
\mathbb{C}\subset \Delta_4\subset M_4\subset \Delta_4 \otimes M_4\subset M_4\otimes M_4\subset \Delta_4\otimes M_4\otimes M_4\subset M_4\otimes M_4\otimes M_4\subset\cdots
\]
is a tower of Jones' basic construction, and $R$ is the closure in the SOT topology of the union of these subalgebras. Following \Cref{pairspin}, the ladder of basic constructions of the following `spin model commuting square'.
\[
\begin{matrix}
u\Delta_4 u^* &\subset & M_4 \cr \cup &\ & \cup\cr \mathbb{C} &\subset &  \Delta_4
\end{matrix}\]
is depicted in \Cref{fig2} (note that in our convention $u_0=u$), where the unitary matrices $u_j$ are given in \Cref{tower in four by four}.
Although the basic construction (in the general $n\times n$ case) is well-known (see \cite{JS, N, Jo3}, for instance),  our basic construction is slightly different than the one in the existing literature to make computations easier. We remark that $R_u$ is the closure in the SOT topology of the union of the left vertical algebras in \Cref{fig2}, and similarly $R_v$ is that of the right vertical algebras. 

We fix the following one-parameter family of unitary matrices throughout this section. For $\,\alpha\in\mathbb{S}^1$,
\begin{IEEEeqnarray}{lCl}\label{main unitaries}
\xi_\alpha=\frac{1}{2}\left[{\begin{matrix}
1 & 1 & 1 & 1\\
-i & \alpha & i & -\alpha\\
1 & -1 & 1 & -1\\
i & \alpha & -i & -\alpha\\
\end{matrix}}\right]\qquad &,&\qquad W_1^{(\alpha)}=\left[{\begin{matrix}
i & 0 & 0 & 0\\
0 & \overline{\alpha} & 0 & 0\\
0 & 0 & -i & 0\\
0 & 0 & 0 & -\overline{\alpha}\\
\end{matrix}}\right]\,,\nonumber\\
W_2:=W_2^{(\alpha)}=\left[{\begin{matrix}
1 & 0 & 0 & 0\\
0 & -1 & 0 & 0\\
0 & 0 & 1 & 0\\
0 & 0 & 0 & -1\\
\end{matrix}}\right]\qquad &,&\qquad W_3^{(\alpha)}=\left[{\begin{matrix}
-i & 0 & 0 & 0\\
0 & \overline{\alpha} & 0 & 0\\
0 & 0 & i & 0\\
0 & 0 & 0 & -\overline{\alpha}\\
\end{matrix}}\right]\,.
\end{IEEEeqnarray}
These unitary matrices satisfy $\big(W_1^{(\alpha)}\big)^*=W_3^{(\overline{\alpha})},\,W_2^*=W_2,\,\big(W_3^{(\alpha)}\big)^*=W_1^{(\overline{\alpha})}$. Let $I_4^{(k)}$ denote the unit element $I_4\otimes\ldots\otimes I_4$ of $M_4^{(k)}:=M_4^{\,\otimes\,k}$. We have the following tower of basic construction (with the convention $u_0=u$).

\begin{thm}\label{tower in four by four}
The tower of the basic construction for $\mathbb{C}\subset u\Delta_4u^*$ is given by
\[
\mathbb{C}\subset u_0\Delta_4u^*_0\subset u_1M_4u^*_1\subset u_2(\Delta_4\otimes M_4)u^*_2\subset u_3(M_4\otimes M_4)u^*_3\subset\cdots
\]
where each unitary matrix $u_i$ is given by the following prescription~:
\begin{enumerate}[$(i)$]
\item for $k\in\bbn\cup\{0\}$, we have $\,u_{2k+1}=E_{11}\otimes u_{2k}+\sum_{j=1}^3E_{j+1,j+1}\otimes u_{2k}\big(W_j^{(a)}\otimes I_4^{(k)}\big)$;
\item for $k\in\bbn$, we have $\,u_{2k}=u_{2k-1}\big(\xi_a\otimes I_4^{(k)}\big)$.
\end{enumerate}
\end{thm}

For the proof of this theorem we request the interested reader to visit the Appendix. We use the following notations throughout this section which are at par with that in part $(1)$ of Notation \ref{not}. For $k\in\bbn\cup\{0\}$, let
\begin{align*}
A_{2k+1} &= \Delta_4\otimes M_4^{(k+1)}\,,\cr
B^u_{2k+1} &= u_{2k+1}(\bbc\otimes M_4^{(k+1)})u_{2k+1}^*\,,\cr
B^v_{2k+1} &= v_{2k+1}(\bbc\otimes M_4^{(k+1)})v_{2k+1}^*\,,\cr
C_{2k+1} &= B^u_{2k+1}\cap {B}^v_{2k+1}\,.
\end{align*}

\noindent These are the odd steps in the tower of basic constructions in \Cref{fig2}. Note that
\[
\overline{\cup A_{2k+1}}^{\mbox{ sot}}=R,\,\overline{\cup B^u_{2k+1}}^{\mbox{ sot}}=R_u,\,\overline{\cup B^v_{2k+1}}^{\mbox{ sot}}=R_v\mbox{ and } \overline{\cup\, C_{2k+1}}^{\mbox{ sot}}=R_u\cap R_v.
\]
We shall see latter in \Cref{even intersection} and \Cref{odd intersection} the merit for considering odd steps in the tower described in \Cref{fig2}, as against the even steps considered in \Cref{Sec 4}. The following lemma is useful for the subsequent subsections.
 
\begin{lmma}\label{zpq}
In the following diagram
\[
\begin{matrix}
A_1 &\subset & A_3 &\subset  & A_5 &\subset &\ldots &\subset & R \cr
\cup &  &\cup &  &\cup &  &  &  & \cup \cr
C_1 &\subset & C_3 &\subset & C_5 &\subset &\ldots &\subset & R_u\cap R_v
\end{matrix}
\]
each individual quadruple is a commuting square.
\end{lmma}
\begin{prf}
Fix a $k\in\bbn$ and consider the cube in \Cref{com30}.
\begin{figure}[!h]
\begin{center} \resizebox{7 cm}{!}{\includegraphics{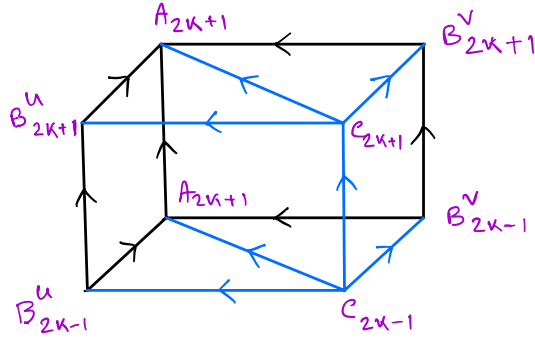}}\end{center}
\caption{A commuting cube in $4\times 4$}\label{com30}
\end{figure}
The adjacent faces are clearly commuting squares by the construction of $R_u,R_v\subset R$. The following quadruple
\[
\begin{matrix}
C_{2k+1} &\subset & A_{2k+1}\\
\cup &  & \cup\\
C_{2k-1} &\subset & A_{2k-1}\\
\end{matrix}
\]
is the slice in the cube. The required proof now follows from \Cref{comm cube}.\qed
\end{prf}

\begin{figure}
\centering
\begin{tikzpicture}[font=\large,thick]

\node[draw,
    align=center,
    rounded rectangle,
    minimum width=2.5cm,
    minimum height=1cm] (block1) {~Start with distinct $4\times 4$ inequivalent ~\\
    complex Hadamard matrices $u\mbox{ and }v$};

\node[draw,
    rectangle, 
    trapezium left angle = 65,
    trapezium right angle = 115,
    trapezium stretches,
    below=of block1,
    minimum width=3.5cm,
    minimum height=1cm] (block2) {~Obtain the pair of spin model subfactors $R_u\mbox{ and }R_v\,$};

\node[draw,
    align=center,
    below=of block2,
    minimum width=3.5cm,
    minimum height=1cm] (block3) {~Computation of $\lambda(R_u,R_v)\,$};

\node[draw,
    align=center,
    below=of block3,
    minimum width=3.5cm,
    minimum height=1cm] (block40) {~Spectral decomposition for the family of unitary~\\
    matrices in the tower of basic constructions~};
  
\node[draw,
    align=center,
    below=of block40,
    minimum width=3.5cm,
    minimum height=1cm] (block4) {~Identifying two situations, namely\\
    $\,u\sim v\mbox{ and }u\nsim v$ depending on certain ~\\
    rational or irrational rotations ~};

\node[draw,
    align=center,
    below=of block4,
    minimum width=3.5cm,
    minimum height=1cm] (block5) {~When $u\sim v$, construction of a family ~\\
     of finite index subfactors $R^{(m)}_{u,v}$ of $R$};

\node[draw,
    align=center,
    below=of block5,
    minimum width=3.5cm,
    minimum height=1cm] (block60) {~The Ocneanu compactness and ~\\
    irreducibility of $R^{(m)}_{u,v}$ in $R$ };
    
\node[draw,
    align=center,
    below=of block60,
    minimum width=3.5cm,
    minimum height=1cm] (block61) {~Factoriality of $R_u\cap R_v\,$\\
     and rigidity of the\\
     ~interior angle $\alpha_R(R_u,R_v)$ ~};
    
\node[draw,
    align=center,
    below=of block61,
    minimum width=3.5cm,
    minimum height=1cm] (block62) {~Intermediate subfactor containing $R_u,R_v\,$};

\node[draw,
    rounded rectangle,
    align=center,
    right=of block62,
    minimum width=3.5cm,
    minimum height=1cm] (block63) {~Angle between $R_u,R_v$ and\\
     relative entropy $H(R_u|R_v)\,$};
  
\node[draw,
    align=center,
    below right=of block4,
    minimum width=3.5cm,
    minimum height=1cm] (block6) {~When $u\nsim v,\,R_u\cap R_v$ has infinite ~\\
    Pimsner-Popa index in $R$};

\node[draw,
    below=of block6,
    minimum width=3.5cm,
    minimum height=1cm] (block7) {~Bounds for $H(R_u|R_v)\,$};

\node[draw,
    rounded rectangle,
    align=center,
    below=of block7,
    minimum width=3.5cm,
    minimum height=1cm] (block8) {~Characteriztion of the commuting ~\\
    square $\big(R_u\cap R_v\subset R_u,R_v\subset R\big)$ \\
    for inequivalent pair of $u\mbox{ and }v$};

\draw[-latex] (block1) edge (block2)
    (block2) edge (block3)
    (block3) edge (block40)
    (block40) edge (block4)
    (block4) edge (block5)
    (block5) edge (block60)
    (block60) edge (block61)
    (block61) edge (block62)
    (block62) edge (block63)
    (block6) edge (block7)
    (block61) edge (block8);
    
\draw[-latex] (block4) -| (block6)
    node[pos=0.5,fill=white,inner sep=0]{};

\end{tikzpicture}
\bigskip

\caption{Flowchart for our plan of actions in \Cref{Sec 5}}\label{flow chart 2}
\end{figure}
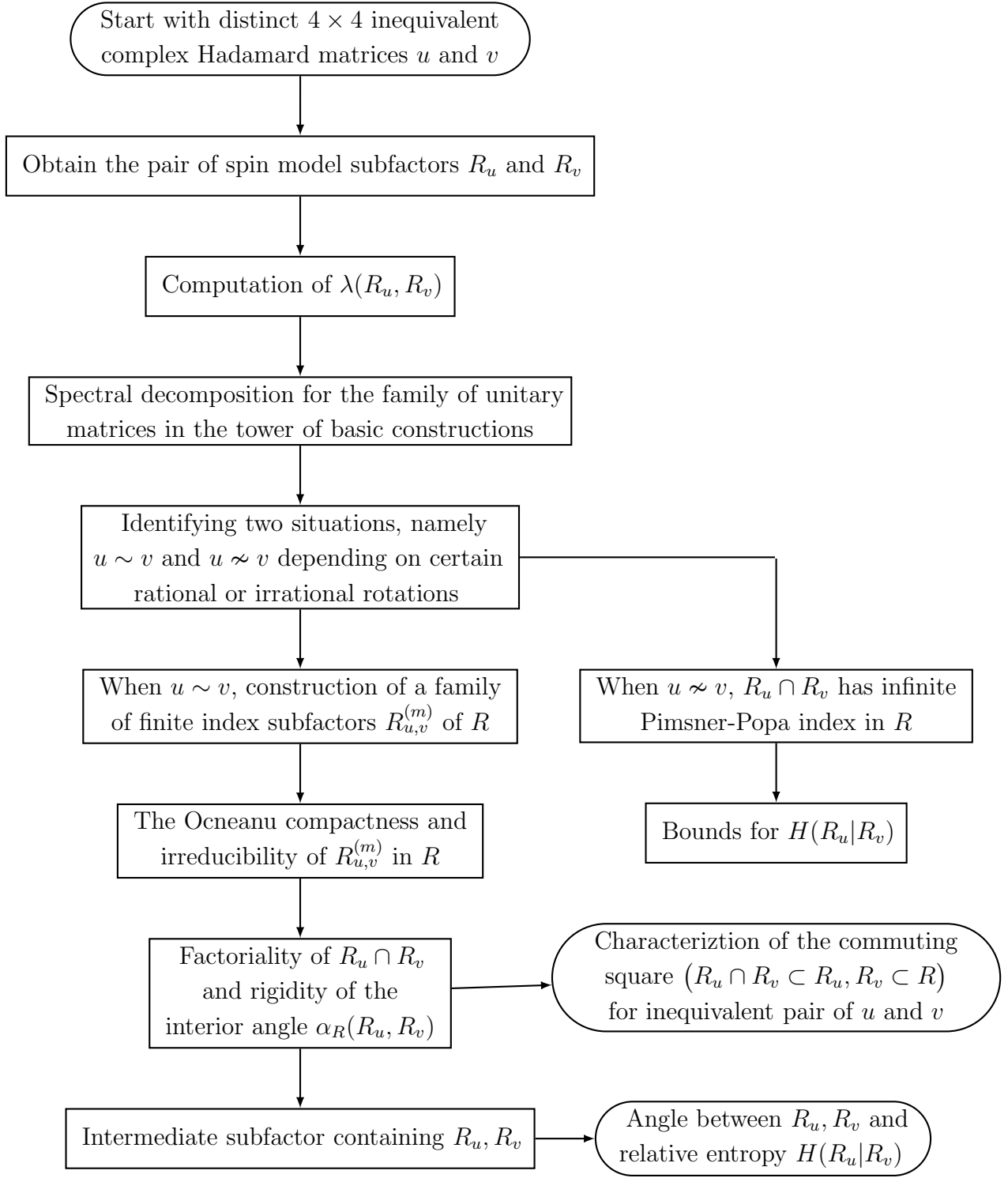

For the sake of reader's benefit, we draw a flowchart in \Cref{flow chart 2} to display our plan of actions in this section.

\subsection{Computation of the Pimsner-Popa constant~}

In this subsection, we compute the Pimsner-Popa probabilistic constant $\lambda(R_u,R_v)$. We introduce the following matrices
\begin{IEEEeqnarray}{lCl}\label{p and q}
p=\frac{1}{2}\left[{\begin{matrix}
1 & 0 & 1 & 0\\
0 & 1 & 0 & 1\\
1 & 0 & 1 & 0\\
0 & 1 & 0 & 1\\
\end{matrix}}\right]\quad,\quad q=\frac{1}{2}\left[{\begin{matrix}
1 & 0 & -1 & 0\\
0 & \overline{a}b & 0 & -\overline{a}b\\
-1 & 0 & 1 & 0\\
0 & -\overline{a}b & 0 & \overline{a}b\\
\end{matrix}}\right]
\end{IEEEeqnarray}
which will play pivotal role throughout the section. Observe that $p$ is a projection (but $q$ is not) and $pq=qp=0$. Also note that $u^*v=p+q$.

\begin{lmma}\label{a cute fact in 4 by 4}
One has the following,
\begin{enumerate}[$(i)$]
\item $\big(W_j^{(a)}\big)^*p\,W_j^{(b)}=q^*\mbox{ and } \big(W_j^{(a)}\big)^*q\,W_j^{(b)}=p~\mbox{ for }j=1,3\,;$
\item $W_2\,p\,W_2=p$ and $W_2\,q\,W_2=q$;
\item $\xi_a^*(E_{11}+E_{33})\xi_b=p\,\mbox{ and } \,\xi_a^*(E_{22}+E_{44})\xi_b=q\,.$
\end{enumerate}
\end{lmma}
\begin{prf}
Straightforward verification using \Cref{main unitaries,p and q}, and hence we leave it to the reader.\qed
\end{prf}

\begin{ppsn}\label{cute lemma for 4 by 4}
For any $k\in\bbn$, we have the following identity
\begin{center}
$u_{2k}^*v_{2k}=p\otimes u_{2k-2}^*v_{2k-2}+q\otimes v_{2k-2}^*u_{2k-2}$
\end{center}
in $M_4\otimes M_4^{(k)}$.
\end{ppsn}
\begin{prf}
We use induction on $k$. For $k=1$, it is a straightforward verification using Theorem \ref{tower in four by four}. Assume that the result holds up to $k-1$. Then, using part (i) of Lemma \ref{a cute fact in 4 by 4} we get the following for $j=1,3$,
\begin{IEEEeqnarray}{lCl}\label{asalp}
&  & \big(W_j^{(a)}\otimes I_4^{(k-1)}\big)^*u_{2k-2}^*v_{2k-2}\big(W_j^{(b)}\otimes I_4^{(k-1)}\big)\nonumber\\
&=& \big(W_j^{(a)}\big)^*pW_j^{(b)}\otimes u_{2k-4}^*v_{2k-4}+\big(W_j^{(a)}\big)^*qW_j^{(b)}\otimes v_{2k-4}^*u_{2k-4}\nonumber\\
&=& q^*\otimes u_{2k-4}^*v_{2k-4}+p\otimes v_{2k-4}^*u_{2k-4}\nonumber\\
&=& v_{2k-2}^*u_{2k-2}
\end{IEEEeqnarray}
Moreover, by part $(ii)$ of \Cref{a cute fact in 4 by 4}, we see that
\begin{IEEEeqnarray}{lCl}\label{asalq}
\big(W_2\otimes I_4^{(k-1)}\big)^*u_{2k-2}^*v_{2k-2}\big(W_2\otimes I_4^{(k-1)}\big)=u_{2k-2}^*v_{2k-2}\,.
\end{IEEEeqnarray}
Now, using part $(iii)$ of Lemma \ref{a cute fact in 4 by 4}, along with Theorem \ref{tower in four by four}, we finally get the following,
\begin{IEEEeqnarray*}{lCl}
&  & u_{2k}^*v_{2k}\\
&=& \big(\xi_a^*\otimes I_4^{(k)}\big)\,\mbox{bl-diag}\Big\{u_{2k-2}^*v_{2k-2},v_{2k-2}^*u_{2k-2},u_{2k-2}^*v_{2k-2},v_{2k-2}^*u_{2k-2}\Big\}\,\big(\xi_b\otimes I_4^{(k)}\big)\\
&=& \big(\xi_a^*\otimes I_4^{(k)}\big)\,\big(\mbox{bl-diag}\{p\otimes u_{2k-4}^*v_{2k-4},q^*\otimes u_{2k-4}^*v_{2k-4},p\otimes u_{2k-4}^*v_{2k-4},q^*\otimes u_{2k-4}^*v_{2k-4}\}\\
&  & +\mbox{bl-diag}\{q\otimes v_{2k-4}^*u_{2k-4},p\otimes v_{2k-4}^*u_{2k-4},q\otimes v_{2k-4}^*u_{2k-4},p\otimes v_{2k-4}^*u_{2k-4}\}\big)\big(\xi_b\otimes I_4^{(k)}\big)\\
&=& \big(\xi_a^*\otimes I_4\otimes I_4^{(k-1)}\big)\Big((E_{11}+E_{33})\otimes p\otimes u_{2k-4}^*v_{2k-4}+(E_{22}+E_{44})\otimes q^*\otimes u_{2k-4}^*v_{2k-4}\\
&  & +(E_{11}+E_{33})\otimes q\otimes v_{2k-4}^*u_{2k-4}+(E_{22}+E_{44})\otimes p\otimes v_{2k-4}^*u_{2k-4}\Big)\big(\xi_b\otimes I_4\otimes I_4^{(k-1)}\big)\\
&=& p\otimes p\otimes u_{2k-4}^*v_{2k-4}+q\otimes q^*\otimes u_{2k-4}^*v_{2k-4}+p\otimes q\otimes v_{2k-4}^*u_{2k-4}+q\otimes p\otimes v_{2k-4}^*u_{2k-4}\\
&=& p\otimes u_{2k-2}^*v_{2k-2}+q\otimes v_{2k-2}^*u_{2k-2}\,,
\end{IEEEeqnarray*}
where the last equality follows from the induction hypothesis.\qed
\end{prf}

\begin{lmma}\label{main reduction in 4 by 4}
For any $k\geq 0$, the unitary $u_{2k+1}^*v_{2k+1}$ is of the following form
\[
\mbox{bl-diag}\left\{u_{2k}^*v_{2k}\,,\,v_{2k}^*u_{2k}\,,\,u_{2k}^*v_{2k}\,,\,v_{2k}^*u_{2k}\right\}\,.
\]
\end{lmma}
\begin{prf}
We have the following identities
\begin{center}
$\big(W_1^{(a)}\otimes I_4^{(k)}\big)^*u_{2k}^*v_{2k}\big(W_1^{(b)}\otimes I_4^{(k)}\big)=\big(W_3^{(a)}\otimes I_4^{(k)}\big)^*u_{2k}^*v_{2k}\big(W_3^{(b)}\otimes I_4^{(k)}\big)=v_{2k}^*u_{2k}\,,$
\end{center}
and
\begin{center}
$\big(W_2\otimes I_4^{(k)}\big)^*u_{2k}^*v_{2k}\big(W_2\otimes I_4^{(k)}\big)=u_{2k}^*v_{2k}\,,$
\end{center}
due to \Cref{asalp,asalq}. By Theorem \ref{tower in four by four}, we have $u_{2k+1}=\mbox{bl-diag}\{u_{2k},u_{2k}\big(W_1^{(a)}\otimes I_4^{(k)}\big),u_{2k}\big(W_2\otimes I_4^{(k)}\big),u_{2k}\big(W_3^{(a)}\otimes I_4^{(k)}\big)\}$, and hence the result follows.\qed
\end{prf}

\begin{lmma}\label{lambda reduction}
Let $\,U_1,U_2\in M_n$ be two unitary matrices and consider the unital inclusions $M_n\subseteq \Delta_4\otimes M_n$ and $UM_nU^*\subseteq\Delta_4\otimes M_n$, where $U\in\Delta_4\otimes M_n$ is the block diagonal unitary matrix $\mbox{\em{bl-diag}}\{U_1,U_2,U_1,U_2\}$. Then, the Pimsner-Popa constant satisfies the following,
\[
\lambda^{\Delta_4\otimes M_n}\left(M_n,UM_nU^*\right)=\lambda^{\Delta_2\otimes M_n}\left(M_n,VM_nV^*\right)\,,
\]
where $V=\mbox{\em{bl-diag}}\{U_1,U_2\}$ is a block diagonal unitary matrix in $\Delta_2\otimes M_n$.
\end{lmma}
\begin{prf}
For all $x\in\left(M_n\right)_+\subseteq\left(\Delta_4\otimes M_n\right)_+$, observe that $E^{\Delta_4\otimes M_n}_{UM_nU^*}(x)\geq tx$ for $t\in[0,1]$ if and only if both the following inequalities
\[
\frac{1}{2}\left(U_1^*xU_1+U_2^*xU_2\right)\geq tU_1^*xU_1
\]
\[
\frac{1}{2}\left(U_1^*xU_1+U_2^*xU_2\right)\geq tU_2^*xU_2
\]
hold simultaneously. Now, observe that for all $y\in\left(M_n\right)_+\subseteq\left(\Delta_2\otimes M_n\right)_+$, the inequality $E^{\Delta_2\otimes M_n}_{VM_nV^*}(y)\geq ty$ for $t\in[0,1]$ holds if and only if the following inequality
\[
\frac{1}{2}\left(U_1^*yU_1+U_2^*yU_2\right)\geq t\,\mbox{bl-diag}\{U_1^*yU_1,U_2^*yU_2\}
\]
holds, which completes the proof.\qed
\end{prf}

We are going to invoke Theorem \ref{2nd lambda} to finally determine the value of the Pimsner-Popa constant for the $k$-th step of the tower of basic construction. For that we need to determine whether the matrix $\big(v_{2k}^*u_{2k}\big)^2$ in $M_4\otimes M_4^{(k)}$ is not diagonal or scalar. Recall that $u\mbox{ and }v$ are Hadamard inequivalent matrices. That is, if we write $a=e^{i\alpha}$ and $b=e^{i\beta}$, then both $\alpha,\beta\in[0,\pi)$. Since $\alpha-\beta\in(-\pi,\pi)$, we have $b\neq -a$.

\begin{lmma}\label{to refer in 4 by 4 PP constant 1}
For any $k\geq 0$, the matrix $\big(v_{2k}^*u_{2k}\big)^2$ is not diagonal in $M_4\otimes M_4^{(k)}$.
\end{lmma}
\begin{prf}
First take $k=0$. We have $v^*u=p+q^*$, and hence $\left(v^*u\right)^2=p+(q^*)^2$ which is not a diagonal matrix as $b\neq -a$. Now for $k\geq 1$, using Proposition \ref{cute lemma for 4 by 4} we have the following,
\begin{IEEEeqnarray*}{lCl}
\big(v_{2k}^*u_{2k}\big)^2 &=& \big(p\otimes v_{2k-2}^*u_{2k-2}+q^*\otimes u_{2k-2}^*v_{2k-2}\big)^2\\
&=& p\otimes\big(v_{2k-2}^*u_{2k-2}\big)^2+(q^*)^2\otimes\big(u_{2k-2}^*v_{2k-2}\big)^2\,.
\end{IEEEeqnarray*}
From here, we see that if $\big(v_{2k-2}^*u_{2k-2}\big)^2$ is not self-adjoint then we are done, and even if $\big(v_{2k-2}^*u_{2k-2}\big)^2$ is self-adjoint, then the fact that $(a\overline{b})^2\neq 1$ (otherwise we get $b=-a$) proves the claim.\qed
\end{prf}

\begin{thm}\label{Popa for 4 by 4}
The Pimsner-Popa constant for the pair of subfactors $R_u\mbox{ and }R_v$ of the hyperfinite $II_1$ factor $R$ is given by $\lambda(R_u,R_v)=1/2$.
\end{thm}
\begin{prf}
Combining \Cref{main reduction in 4 by 4}, \ref{lambda reduction}, \ref{to refer in 4 by 4 PP constant 1}, together with \Cref{2nd lambda}, we get the following
\begin{IEEEeqnarray}{lCl}\label{92}
\lambda\big(u_{2k+1}M_4^{(k+1)}u_{2k+1}^*\,,\,v_{2k+1}M_4^{(k+1)}v_{2k+1}^*\big)=\frac{1}{2}
\end{IEEEeqnarray}
for any $k\geq 0$. This immediately gives us that $\lambda(R_u,R_v)=\frac{1}{2}$, because by \Cref{popaadaptation}, we have $\lambda(R_u,R_v)$ is the limit of a decreasing sequence of $\lambda$ at each step of the tower of basic construction, and we have produced a constant subsequence in \Cref{92}.\qed
\end{prf}


\subsection{Construction of a family of subfactors of hyperfinite $II_1$ factor $R$}

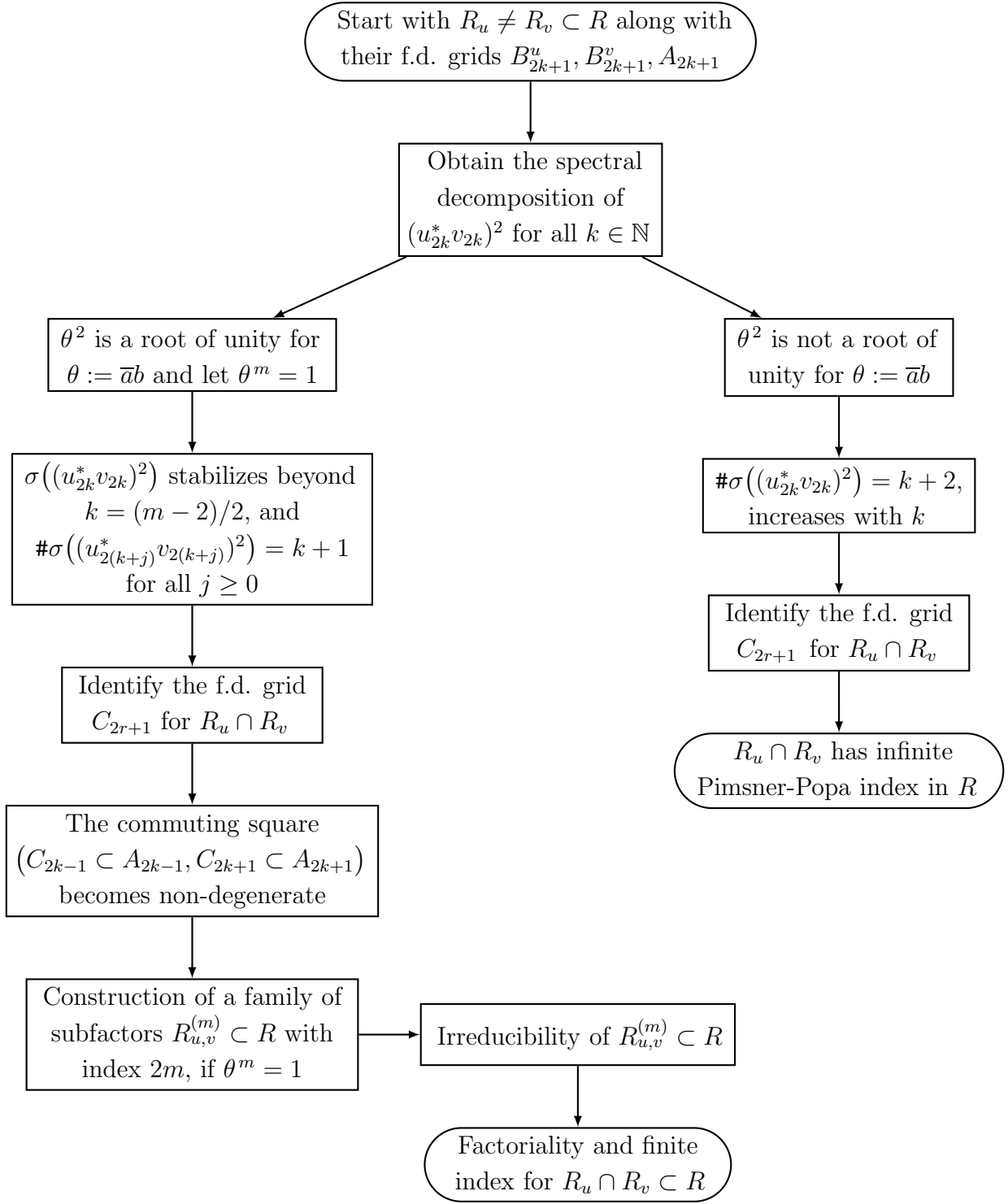
\begin{figure}
\centering
\begin{tikzpicture}[font=\large,thick]

\node[draw,
    align=center,
    rounded rectangle,
    minimum width=2.5cm,
    minimum height=1cm] (block0) {~Start with $R_u\neq R_v\subset R$ along with \\
    their f.d. grids $B^u_{2k+1},B^v_{2k+1},A_{2k+1}$~};
    
\node[draw,
    align=center,
     below=of block0,
    minimum width=2.5cm,
    minimum height=1cm] (block1) {~Obtain the spectral \\
    decomposition of\\
    $(u_{2k}^*v_{2k})^2$ for all $k\in\bbn\,$};

\node[draw,
    align=center,
    rectangle, 
    trapezium left angle = 65,
    trapezium right angle = 115,
    trapezium stretches,
    below left=of block1,
    minimum width=3.5cm,
    minimum height=1cm] (block2) {$\,\theta^{\,2}$ is a root of unity for ~\\
    $\theta:=\overline{a}b$ and let $\theta^{\,m}=1$~};

\node[draw,
    align=center,
    below right=of block1,
    minimum width=3.5cm,
    minimum height=1cm] (block3) {$\,\theta^{\,2}$ is not a root of ~\\
    unity for $\theta:=\overline{a}b$~};

\node[draw,
    align=center,
    below=of block3,
    minimum width=3.5cm,
    minimum height=1cm] (block4) {$\,\texttt{\#}\sigma\big((u_{2k}^*v_{2k})^2\big)=k+2,\,$\\
    increases with $k\,$};
    
  \node[draw,
    align=center,
    below=of block4,
    minimum width=3.5cm,
    minimum height=1cm] (block50) {~Identify the f.d. grid ~\\
    $C_{2r+1}\,$ for $R_u\cap R_v\,$};
  
\node[draw,
    rounded rectangle,
    align=center,
    below=of block50,
    minimum width=3.5cm,
    minimum height=1cm] (block5) {~$R_u\cap R_v$ has infinite \\
    Pimsner-Popa index in $R\,$};

\node[draw,
    align=center,
    below=of block2,
    minimum width=3.5cm,
    minimum height=1cm] (block20) {$\,\sigma\big((u_{2k}^*v_{2k})^2\big)$ stabilizes beyond ~\\
    $k=(m-2)/2$, and\\
    $\texttt{\#}\sigma\big((u_{2(k+j)}^*v_{2(k+j)})^2\big)=k+1$\\
    for all $j\geq 0$~};
     
\node[draw,
    align=center,
    below=of block20,
    minimum width=3.5cm,
    minimum height=1cm] (block60) {~Identify the f.d. grid ~\\
    $C_{2r+1}$ for $R_u\cap R_v\,$};
    
\node[draw,
    align=center,
    below=of block60,
    minimum width=3.5cm,
    minimum height=1cm] (block30) {~The commuting square~ \\
    $\big(C_{2k-1}\subset A_{2k-1},C_{2k+1}\subset A_{2k+1}\big)\,$ \\
    becomes non-degenerate};
     
\node[draw,
    align=center,
    below=of block30,
    minimum width=3.5cm,
    minimum height=1cm] (block7) {~Construction of a family of ~\\
    subfactors $R^{(m)}_{u,v}\subset R$ with\\
    index $2m$, if $\theta^{\,m}=1$~};
  
\node[draw,
    rectangle,
    right=of block7,
    minimum width=3.5cm,
    minimum height=1cm] (block8) {~Irreducibility of $R^{(m)}_{u,v}\subset R\,$};

\node[draw,
    rounded rectangle,
    align=center,
    below=of block8,
    minimum width=3.5cm,
    minimum height=1cm] (block9) {~Factoriality and finite ~\\
    ~index for $R_u\cap R_v\subset R\,$};

\draw[-latex] (block0) edge (block1)
    (block1) edge (block2)
    (block1) edge (block3)
    (block3) edge (block4)
    (block4) edge (block50)
    (block50) edge (block5)
    (block2) edge (block20)
    (block20) edge (block60)
    (block60) edge (block30)
    (block30) edge (block7)
    (block7) edge (block8)
    (block8) edge (block9);

\end{tikzpicture}
\bigskip

\caption{Technical roadmap for subsections $7.2$ and $7.3$}\label{road map 1}
\end{figure}

This subsection and its intimately connected sequel contains one of the most important result of the paper. We will construct a family of potentially new integer index subfactors of $R$ in this subsection, and latter in the sequel we show irreducibility of these family of subfactors. For the sake of reader's convenience, we draw a {\it technical} roadmap for this subsection and its immediate sequel in \Cref{road map 1}.

To fulfil our goal, we need to focus on the von Neumann subalgebra $R_u\cap R_v$ of the hyperfinite type $II_1$ factor $R$. To understand the structure of $R_u\cap R_v$, since we already have the finite-dimensional grids for the spin model subfactors $R_u\mbox{ and }R_v$, our target is to obtain the same for $R_u\cap R_v$. It turns out that in doing so a surprising phenomenon occurs. A family of finite index potentially new subfactors of $R$ arises. Moreover, with the help of it, we can also conclude factoriality of $R_u\cap R_v$ along with its index in $R$ in certain situations. This also helps to conclude infinite Pimsner-Popa index of $R_u\cap R_v$ in certain other situations. We begin by fixing some notations to be used throughout this subsection and latter.
\medskip

\noindent\textbf{Notation:}~~ (i) For complex Hadamard matrices $u\mbox{ and }v$ parametrized by $a\mbox{ and }b$, let $\theta=\overline{a}b\in\mathbb{S}^1$. We replace the matrix $q$ introduced in \Cref{p and q} by the notation $q(\theta)$ to stress its dependence on $\theta$. This is required in what follows next.\\
(ii) Consider the following matrices
\begin{IEEEeqnarray}{lCl}\label{W}
W:=\frac{1}{\sqrt{2}}\left[{\begin{matrix}
1 & 0 & 1 & 0\\
0 & 1 & 0 & 1\\
1 & 0 & -1 & 0\\
0 & 1 & 0 & -1\\
\end{matrix}}\right]\qquad\mbox{and}\qquad C(\theta)=\left[{\begin{matrix}
1 & 0\\
0 & \theta\\
\end{matrix}}\right]\,.
\end{IEEEeqnarray}
Note that $W$ is a self-adjoint unitary matrix in $M_4$.\\
(iii) For a matrix $T$, $\sigma(T)$ will denote its spectrum.

\begin{lmma}\label{even intersection}
For any $k\in\bbn$, we have the following,
\begin{enumerate}[$(i)$]
\item $\mathrm{Ad}_u(\Delta_4)\cap\mathrm{Ad}_v(\Delta_4)=\mathrm{Ad}_u(\bbc^3)$, where the emebedding $\bbc^3\xhookrightarrow{}\Delta_4\subseteq M_4$ is given by $(\lambda_1,\lambda_2,\lambda_3)\mapsto(\lambda_1,\lambda_2,\lambda_3,\lambda_2)\,,$
\item $u_{2k}^*v_{2k}\big(\Delta_4\otimes M_4^{(k)}\big)v_{2k}^*u_{2k}\subseteq\left(u^*v\Delta_4v^*u\right)\otimes M_4^{(k)}\,$,
\item $\mathrm{Ad}_{u_{2k}}\big(\Delta_4\otimes M_4^{(k)}\big)\cap\mathrm{Ad}_{v_{2k}}\big(\Delta_4\otimes M_4^{(k)}\big)\subseteq\mathrm{Ad}_{u_{2k}}\big(\bbc^3\otimes M_4^{(k)}\big)\subseteq M_4\otimes M_4^{(k)}$.
\end{enumerate}
\end{lmma}
\begin{prf}
We have $\mathrm{Ad}_u(\Delta_4)\cap\mathrm{Ad}_v(\Delta_4)=\mathrm{Ad}_u(\Delta_4\cap\mathrm{Ad}_{u^*v}(\Delta_4))$. Recall that $u^*v=p+q(\theta)$, where $p\mbox{ and }q(\theta)$ are as in \Cref{p and q}. Hence, for any $\mbox{bl-diag}\{\mu_1,\mu_2,\mu_3,\mu_4\}\in\Delta_4$, we have the following
\begin{IEEEeqnarray*}{lCl}
\mathrm{Ad}_{u^*v}\big(\mbox{bl-diag}\{\mu_1,\mu_2,\mu_3,\mu_4\}\big)=\begin{bmatrix}
\mu_1 & 0 & 0 & 0\\
0 & \mu_2\frac{(2+\theta+\overline{\theta})}{4}+ \mu_4\frac{(2-\theta-\overline{\theta})}{4} & 0 &  \mu_2\frac{(\theta-\overline{\theta})}{4}+ \mu_4\frac{(\overline{\theta}-\theta)}{4}\\
0 & 0 & \mu_3 & 0\\
0 & \mu_2\frac{(\overline{\theta}-\theta)}{4}+ \mu_4\frac{(\theta-\overline{\theta})}{4} & 0 & \mu_2\frac{(2-\theta-\overline{\theta})}{4}+ \mu_4\frac{(2+\theta+\overline{\theta})}{4}\\
\end{bmatrix}\,.
\end{IEEEeqnarray*}
Observe that $\theta\neq\overline{\theta}$. This is because $\theta=\overline{\theta}$ implies that $b=\pm a$, since $\theta=\overline{a}b$. However, in our situation this never happens as already discussed at the beginning of \Cref{Sec 5}. Therefore, the intersection $\Delta_4\cap\mathrm{Ad}_{u^*v}(\Delta_4)$ consists of the diagonal matrices $\mbox{bl-diag}\{\mu_1,\mu_2,\mu_3,\mu_4\}$ such that $\mu_1,\mu_3$ are arbitrary complex numbers, and $\mu_2,\mu_4$ are equal. In other words,
\[
\Delta_4\cap\mathrm{Ad}_{u^*v}(\Delta_4)=\big\{\mbox{bl-diag}\{\mu_1,\mu_2,\mu_3,\mu_2\}\,:\,\mu_1,\mu_2,\mu_3\in\bbc\big\}\,,
\]
which completes part $(i)$. For part $(ii)$, by \Cref{cute lemma for 4 by 4} we get the following,
\begin{IEEEeqnarray*}{lCl}
&  & u_{2k}^*v_{2k}\left(\Delta_4\otimes M_4^{(k)}\right)v_{2k}^*u_{2k}\\
&=& \left(p\otimes u_{2k-2}^*v_{2k-2}+q(\theta)\otimes v_{2k-2}^*u_{2k-2}\right)\left(\Delta_4\otimes M_4^{(k)}\right)\left(p\otimes v_{2k-2}^*u_{2k-2}+q(\theta)^*\otimes u_{2k-2}^*v_{2k-2}\right)\\
&\subseteq& \left(p\Delta_4 p+p\Delta_4 q(\theta)^*+q(\theta)\Delta_4 p+q\Delta_4 q(\theta)^*\right)\otimes M_4^{(k)}\\
&=& \left((p+q(\theta))\Delta_4(p+q(\theta)^*)\right)\otimes M_4^{(k)}\,.
\end{IEEEeqnarray*}
Since $p+q(\theta)=u^*v\,$, and consequently $p+q(\theta)^*=v^*u\,$, part $(ii)$ is done. Part $(iii)$ follows from part $(i)$ and $(ii)$ together.\qed
\end{prf}

\begin{ppsn}\label{odd intersection}
For any $k\geq 0$, the intersection
\[
C_{2k+1}=\mathrm{Ad}_{u_{2k+1}}\big(\bbc\otimes M_4^{(k+1)}\big)\cap\mathrm{Ad}_{v_{2k+1}}\big(\bbc\otimes M_4^{(k+1)}\big)
\]
in $A_{2k+1}=\Delta_4\otimes M_4^{(k+1)}$ is given by the following,
\[
\left\{\mathrm{Ad}_{u_{2k+1}}\big(I_4\otimes A\big):A\in\big((u_{2k}^*v_{2k})^2\big)^\prime\cap M_4^{(k+1)}\right\}\,.
\]
\end{ppsn}
\begin{prf}
First observe the following,
\begin{IEEEeqnarray*}{lCl}
& & \mathrm{Ad}_{u_{2k+1}}\big(\bbc\otimes M_4^{(k+1)}\big)\cap\mathrm{Ad}_{v_{2k+1}}\big(\bbc\otimes M_4^{(k+1)}\big)\\
&=& \mathrm{Ad}_{u_{2k+1}}\big(\bbc\otimes M_4^{(k+1)}\cap\mathrm{Ad}_{u_{2k+1}^*v_{2k+1}}\big(\bbc\otimes M_4^{(k+1)}\big)\,.
\end{IEEEeqnarray*}
Recall from \Cref{main reduction in 4 by 4} that $u_{2k+1}^*v_{2k+1}=\mbox{bl-diag}\left\{u_{2k}^*v_{2k},v_{2k}^*u_{2k},u_{2k}^*v_{2k},v_{2k}^*u_{2k}\right\}$. An element $\xi$ is in $\bbc\otimes M_4^{(k+1)}\cap \mbox{Ad}_{u_{2k+1}^*v_{2k+1}}\big(\bbc\otimes M_4^{(k+1)}\big)$ if and only if
\[
\xi=I_4\otimes A=\mbox{Ad}_{u_{2k+1}^*v_{2k+1}}(I_4\otimes B)
\]
for some $A,B\in M_4^{(k+1)}$. This implies that $A=\mbox{Ad}_{u_{2k}^*v_{2k}}(B)=\mbox{Ad}_{v_{2k}^*u_{2k}}(B)$, which means $A$ remains invariant under the conjugation action by the unitary $\big(u_{2k}^*v_{2k}\big)^2\in M_4^{(k+1)}$, and this completes the proof.\qed
\end{prf}

Both the inclusions in part $(ii)$ and $(iii)$ of \Cref{even intersection} is strict for $k>0$. To obtain $\mathrm{Ad}_{u_{2k}}\big(\Delta_4\otimes M_4^{(k)}\big)\cap\mathrm{Ad}_{v_{2k}}\big(\Delta_4\otimes M_4^{(k)}\big)$, one has to further identify the intersection $C_{2k+1}\cap\big(\bbc\otimes M_4^{(k+1)}\big)$ concretely, as $(C_{2k}\subset A_{2k}\,,\,C_{2k+1}\subset A_{2k+1})$ forms a commuting square, where $C_{2k}$ and $A_{2k}$ are in the even steps of the tower of basic constructions in \Cref{fig2}. Thus, the intersections in the even steps of the tower are little harder to identify as compared to that in the odd steps obtained in \Cref{odd intersection}. This is where we differ from \Cref{Sec 4}.
\smallskip

\textbf{A scheme to identify the finite-dimensional grid for $R_u\cap R_v$:~} In order to concretely identify the finite-dimensional grid $\{C_{2k+1}:k\geq 0\}$ of the von Neumann subalgebra $R_u\cap R_v$ of $R$, it is now clear from \Cref{odd intersection} that we need to find the Wedderburn-Artin decomposition of the algebra $\big((u_{2k}^*v_{2k})^2\big)^\prime\cap M_4^{(k+1)}$ for each $k$. The following lemma is crucial in this regard, and although it is an elementary fact from Linear algebra, for the sake of completeness we provide a short proof.

\begin{lmma}\label{com}
Suppose that we have a block-diagonal matrix $A=[A_{i,j}]=A_{1,1}\oplus\cdots\oplus A_{n,n}$ such that each block $A_{i,i},\,1\leq i\leq n,$ has exactly one eigenvalue $\lambda_i$ and $\lambda_i\neq\lambda_j$ for $i\neq j$. Then, the only matrices which commute with $A$ are block-diagonal with the same block sizes as $A$.
\end{lmma}
\begin{prf}
Suppose that $A$ is a $k\times k$ matrix and consider a matrix $B$ in $M_k$. Putting $B$ in block form $(B_{i,j})$, let $A_iB_{i,j}=B_{i,j}A_j$ for all $i,j$. We claim that $i\neq j$ implies $B_{i,j}=0$. Note that $A_iB_{i,j}=B_{i,j}A_j$ implies $\mathcal{P}(A_i)B_{i,j}=B_{i,j}\mathcal{P}(A_j)$ for any polynomial $\mathcal{P}$. There exists $r_j\in\bbn$ such that $(A_j-\lambda_j)^{r_j}=0$, and hence $(A_i-\lambda_j)^{r_j}B_{i,j}=0$. Since $A_i-\lambda_j$ is invertible for $i\neq j$ by hypothesis, we get that $B_{i,j}=0$ for $i\neq j$.\qed
\end{prf}

We describe the scheme that we are now going to follow. Suppose that we have a normal matrix $x\in M_n$. Let $\{\lambda_1,\ldots,\lambda_r\},\,1\leq r\leq n,$ be the distinct eigenvalues of $x$ with geometric multiplicities $m(\lambda_i),\,1\leq i\leq r$. By the spectral theorem, $n=\sum_{i=1}^rm(\lambda_i)$. That is, there exists a unitary matrix $U\in M_n$ such that $x=U\mbox{diag}\{\lambda_11_{m(\lambda_1)},\ldots,\lambda_r1_{m(\lambda_r)}\}U^*$. By \Cref{com}, we have the Wedderburn-Artin decomposition of the algebra $\{x\}^\prime\cap M_n$ as $\oplus_{j=1}^rM_{m(\lambda_j)}$ with $n=\sum_{j=1}^rm(\lambda_j)$.

In view of the above, it is now clear that our first job is to find the spectral decomposition of the unitary matrices $(u_{2k}^*v_{2k})^2\in M_4^{(k+1)}$ for each $k\in\bbn\cup\{0\}$, and the next subsection deals with precisely this.


\subsubsection{Spectral decomposition of $\big(u_{2k}^*v_{2k}\big)^2$ for $k\in\bbn\cup\{0\}$}

We obtain the spectral decomposition for the unitary matrices $\big(u_{2k}^*v_{2k}\big)^2\in M_4^{(k+1)}$ for all $k\in\bbn\cup\{0\}$. Due to the analysis carried out here, we will be able to identify the intersections $C_{2k+1}$ concretely for all $k$ as discussed above. In this regard, recall the tower of basic construction obtained in \Cref{tower in four by four}, the matrices $W$ and $C(\theta)$ from \Cref{W}, and finally the matrices $p,q(\theta)$ from \Cref{p and q}.

\begin{lmma}\label{Sa}
\begin{enumerate}[$(i)$]
\item For any $n\in\bbn,\,C(\theta)^n=C(\theta^n)$ and $C(\theta)^*=C(\overline{\theta})$.
\item We have $WpW=\mathrm{diag}\{I_2,0\}$ and $Wq(\theta)W=\mathrm{diag}\{0,C(\theta)\}$.
\item For any $n\in\bbn,\,q(\theta)^n=q(\theta^{\,n})$.
\end{enumerate}
\end{lmma}
\begin{prf}
These are straightforward verification, and hence left to the reader.\qed
\end{prf}

\begin{lmma}\label{nn}
For any $k\in\bbn$, we have the following,
\[
\sigma\big(\big(u_{2k}^*v_{2k}\big)^2\big)=\sigma\big(\big(u_{2k-2}^*v_{2k-2}\big)^2\big)\cup\sigma\big(\big(v_{2k-2}^*u_{2k-2}\big)^2\big)\cup\sigma\big(\big(\theta\,v_{2k-2}^*u_{2k-2}\big)^2\big).
\]
\end{lmma}
\begin{prf}
By \Cref{cute lemma for 4 by 4}, we have the identity $u_{2k}^*v_{2k}=p\otimes u_{2k-2}^*v_{2k-2}+q(\theta)\otimes v_{2k-2}^*u_{2k-2}$ in $M_4\otimes M_4^{(k)}$. Since $pq(\theta)=q(\theta)p=0$ and $p$ is a projection, we get that
\[
\big(u_{2k}^*v_{2k}\big)^2=p\otimes \big(u_{2k-2}^*v_{2k-2}\big)^2+q(\theta^2)\otimes\big(v_{2k-2}^*u_{2k-2}\big)^2
\]
in $M_4\otimes M_4^{(k)}$ by part $(iii)$ of \Cref{Sa}. Using part $(i),(ii)$ of \Cref{Sa}, we have the following identity in $M_4\otimes M_4^{(k)}$,
\begin{IEEEeqnarray}{lCl}\label{mm}
&  & \mbox{Ad}_{W\otimes W^{(k)}}\big(\big(u_{2k}^*v_{2k}\big)^2\big)\\
&=& E_{11}\otimes I_2\otimes\mbox{Ad}_{W\otimes W^{(k-1)}}\big(\big(u_{2k-2}^*v_{2k-2}\big)^2\big)+E_{22}\otimes C(\theta^{\,2})\otimes\mbox{Ad}_{W\otimes W^{(k-1)}}\big(\big(v_{2k-2}^*u_{2k-2}\big)^2\big).\nonumber
\end{IEEEeqnarray}
This immediately says the following,
\begin{IEEEeqnarray*}{lCl}
&  & \sigma\left(\mbox{Ad}_{W\otimes W^{(k)}}\big(\big(u_{2k}^*v_{2k}\big)^2\big)\right)\\
&=& \sigma\left(\mbox{Ad}_{W^{(k)}}\big(\big(u_{2k-2}^*v_{2k-2}\big)^2\big)\right)\cup\sigma\left(C(\theta^{\,2})\otimes\mbox{Ad}_{W^{(k)}}\big(\big(v_{2k-2}^*u_{2k-2}\big)^2\big)\right)\,,
\end{IEEEeqnarray*}
which concludes the proof since $W$ is a unitary matrix in $M_4$ and $C(\theta^{\,2})$ is a diagonal matrix in $M_2$.\qed
\end{prf}

\begin{lmma}\label{mq}
For any $k\in\bbn\cup\{0\}$, the set of all eigenvalues (including repeated ones possibly) of the unitary matrix $(u_{2k}^*v_{2k})^2$ is $\big\{1,\theta^{\,2},\overline{\theta}^{\,2},\theta^{\,4},\overline{\theta}^{\,4},\ldots,\theta^{\,k},\overline{\theta}^{\,k},\theta^{\,k+2}\big\}$ if $k$ is even, and $\big\{1,\theta^{\,2},\overline{\theta}^{\,2},\theta^{\,4},\overline{\theta}^{\,4},\ldots,\theta^{\,k-1},\overline{\theta}^{\,k-1},\theta^{\,k+1},\overline{\theta}^{\,k+1}\big\}$ if $k$ is odd.
\end{lmma}
\begin{prf}
We use induction on $k$. First suppose that $k$ is even. The basis step of $k=0$ follows from the following fact (using \Cref{Sa})
\[
\mbox{Ad}_W\big((u^*v)^2\big)=\mbox{Ad}_W\big((p+q(\theta))^2\big)=\mbox{Ad}_W\big(p+q(\theta^{\,2})\big)=\mbox{bl-diag}\{I_2,C(\theta^{\,2})\}\,,
\]
and hence $\sigma\big((u^*v)^2\big)=\{1,\theta^{\,2}\}$. Assume that the statement holds for $k-2$. That is, we have $\sigma\big(\big(u_{2k-4}^*v_{2k-4}\big)^2\big)=\{1,\theta^{\,2},\overline{\theta}^{\,2},\theta^{\,4},\overline{\theta}^{\,4},\ldots,\theta^{\,k-2},\overline{\theta}^{\,k-2},\theta^{\,k}\}$. Hence, $\sigma\big(\big(v_{2k-4}^*u_{2k-4}\big)^2\big)=\{1,\theta^{\,2},\overline{\theta}^{\,2},\theta^{\,4},\overline{\theta}^{\,4},\ldots,\theta^{\,k-2},\overline{\theta}^{\,k-2},\overline{\theta}^{\,k}\}$. Now, observe that when each element of this set is multiplied by $\theta^{\,2}$, we get back the set $\sigma\big(\big(u_{2k-4}^*v_{2k-4}\big)^2\big)$. That is, we have $\,\sigma\big(\big(\theta\,v_{2k-4}^*u_{2k-4}\big)^2\big)=\sigma\big(\big(u_{2k-4}^*v_{2k-4}\big)^2\big)$. By \Cref{nn} we now have the following,
\begin{IEEEeqnarray*}{lCl}
\sigma\big(\big(u_{2k-2}^*v_{2k-2}\big)^2\big) &=& \sigma\big(\big(u_{2k-4}^*v_{2k-4}\big)^2\big)\cup\sigma\big(\big(v_{2k-4}^*u_{2k-4}\big)^2\big)\cup\sigma\big(\big(\theta\,v_{2k-4}^*u_{2k-4}\big)^2\big)\\
&=& \{1,\theta^{\,2},\overline{\theta}^{\,2},\theta^{\,4},\overline{\theta}^{\,4},\ldots,\theta^{\,k-2},\overline{\theta}^{\,k-2},\theta^{\,k}\}\cup\{1,\theta^{\,2},\overline{\theta}^{\,2},\theta^{\,4},\overline{\theta}^{\,4},\ldots,\theta^{\,k-2},\overline{\theta}^{\,k-2},\overline{\theta}^{\,k}\}\\
&=& \{1,\theta^{\,2},\overline{\theta}^{\,2},\theta^{\,4},\overline{\theta}^{\,4},\ldots,\theta^{\,k},\overline{\theta}^{\,k}\},
\end{IEEEeqnarray*}
which proves the statement for $k-1$, which is odd. Similar application of \Cref{nn} shows that
\[
\sigma\big(\big(u_{2k}^*v_{2k}\big)^2\big)=\{1,\theta^{\,2},\overline{\theta}^{\,2},\theta^{\,4},\overline{\theta}^{\,4},\ldots,\theta^{\,k},\overline{\theta}^{\,k},\theta^{\,k+2}\},
\]
which proves the statement for $k$, and finishes the induction. Similarly, starting with $k$ odd and using the following fact that (use \Cref{cute lemma for 4 by 4} and \Cref{Sa})
\[
\mbox{Ad}_{W\otimes W}\big((u_2^*v_2)^2\big)=WpW\otimes\mbox{Ad}_W(u^*v)^2+Wq(\theta^{\,2})W\otimes\mbox{Ad}_W(v^*u)^2\,,
\]
we can use induction on $k$ to conclude the proof.\qed
\end{prf}

Note that we are not claiming in the above lemma that $\sigma\big((u_{2k}^*v_{2k})^2\big)$ consists of {\it distinct} eigenvalues. In fact, we are going to investigate this now.
\smallskip

\textbf{Terminology~:} For a complex number $z$ of modulus $1$, we will say that $z$ is an {\em even root of unity} if $z^n=1$ for some even integer $n$. Note that this is same as saying $z^2$ is a root of unity. Throughout the article, by even root of unity we always mean primitive even root of unity. That is, if $z^n=1$ for some $n\in 2\bbn$, then $z^m\neq 1$ for any $m\in 2\bbn$ with $m<n$.
\medskip

We now investigate whether $\sigma\big((u_{2k}^*v_{2k})^2\big)$ obtained in \Cref{mq} consists of distinct eigenvalues only, and it turns out that it depends on $\,\theta\in\mathbb{S}^1$. If $\,\theta$ is not an even root of unity, then it is easy to observe that all the elements of the following sets
\[
\big\{1,\theta^{\,2},\overline{\theta}^{\,2},\theta^{\,4},\overline{\theta}^{\,4},\ldots,\theta^{\,k},\overline{\theta}^{\,k},\theta^{\,k+2}\big\}\,,\,\mbox{when $k$ is even}
\]
\[
\big\{1,\theta^{\,2},\overline{\theta}^{\,2},\theta^{\,4},\overline{\theta}^{\,4},\ldots,\theta^{\,k-1},\overline{\theta}^{\,k-1},\theta^{\,k+1},\overline{\theta}^{\,k+1}\big\}\,,\,\mbox{when $k$ is odd}
\]
are distinct and we have cardinality of $\sigma\big((u_{2k}^*v_{2k}\big)^2\big)$ is $k+2$. Therefore, when $k$ increases, number of elements in the spectrum of the unitary matrices $\big((u_{2k}^*v_{2k})^2\big)$ also increase, and $\sigma\big((u_{2k}^*v_{2k}\big)^2\big)\subsetneq\sigma\big((u_{2k+2}^*v_{2k+2}\big)^2\big)$. However, as soon as $\,\theta$ is an even root of unity, that is $\,\theta^{\,m}=1$ for some $m\in2\bbn+2$ (and $\theta^{\,r}\neq 1$ for any $r<m$), the situation changes drastically. We notice that for $k=(m-2)/2$, the last eigenvalue in the enumeration described in \Cref{mq} coincides with the second last eigenvalue, and we get that $\sigma\big((u_{2k}^*v_{2k}\big)^2\big)=\sigma\big((u_{2k-2}^*v_{2k-2}\big)^2\big)$ irrespective of whether $k$ is even or odd. Note that we are not considering the case $\theta^{\,2}=1$ (i,e., $m=2$) because otherwise we would get $b=\pm a$, which is not possible in our situation as already discussed at the beginning of \Cref{Sec 5}. These discussions culminate into the following results.

\begin{ppsn}\label{oo}
For any $k\in\bbn\cup\{0\}$, if $\,\theta$ is not an even root of unity, then we have
\begin{IEEEeqnarray*}{lCl}
\sigma\big(\big(u_{2k}^*v_{2k}\big)^2\big) &=& \begin{cases}
\{1,\theta^{\,2},\overline{\theta}^{\,2},\theta^{\,4},\overline{\theta}^{\,4},\ldots,\theta^{\,k},\overline{\theta}^{\,k},\theta^{\,k+2}\} & \mbox{ if }\,\,k\mbox{ is even}\,, \cr
\{1,\theta^{\,2},\overline{\theta}^{\,2},\theta^{\,4},\overline{\theta}^{\,4},\ldots,\theta^{\,k-1},\overline{\theta}^{\,k-1},\theta^{\,k+1},\overline{\theta}^{\,k+1}\} & \mbox{ if }\,\,k\mbox{ is odd}\,, \cr
\end{cases}
\end{IEEEeqnarray*}
and $\texttt{\#}\,\sigma\big(\big(u_{2k}^*v_{2k}\big)^2\big)= k+2$ for all $k\in\bbn\cup\{0\}$.
\end{ppsn}

\begin{ppsn}\label{101}
Let $\theta^{\,m}=1$ for some $m\in 2\bbn+2$ and choose $k=(m-2)/2$. Then, we have
\[
\sigma\big((u_{2k-2}^*v_{2k-2})^2\big)=\sigma\big((u_{2k+2j}^*v_{2k+2j})^2\big)=\begin{cases}
\{1,\theta^{\,2},\overline{\theta}^{\,2},\ldots,\theta^{\,k-1},\overline{\theta}^{\,k-1},\theta^{\,k+1}\}& \mbox{ if }\,\,k\mbox{ is odd}\,,\cr
\{1,\theta^{\,2},\overline{\theta}^{\,2},\ldots,\theta^{\,k-2},\overline{\theta}^{\,k-2},\theta^{\,k},\overline{\theta}^{\,k}\}& \mbox{ if }\,\,k\mbox{ is even},
\end{cases}
\]
for all $j\in\bbn\cup\{0\}$, and $\texttt{\#}\,\sigma\big((u_{2k+2j}^*v_{2k+2j})^2\big)=k+1$.
\end{ppsn}
\begin{prf}
If $\theta^{\,m}=1$ for some $m\in 2\bbn+2$, then choose $k=(m-2)/2$. If $k$ is odd, then by \Cref{mq} we have $\theta^{\,k+1}=\overline{\theta}^{\,k+1}$ and $\texttt{\#}\sigma\big((u_{2k}^*v_{2k})^2\big)=k+1$. If $k$ is even, then again by \Cref{mq} we have $\theta^{\,k+2}=\overline{\theta}^{\,k}$ and $\texttt{\#}\sigma\big((u_{2k}^*v_{2k})^2\big)=k+1$. It now follows by induction on $j$ that $\sigma\big((u_{2k-2}^*v_{2k-2})^2\big)=\sigma\big((u_{2k}^*v_{2k})^2\big)=\sigma\big((u_{2(k+j)}^*v_{2(k+j)})^2\big)$ for all $j\in\bbn$ with cardinality of the spectrum being equal to $k+1$.\qed
\end{prf}

Our first job of finding the distinct eigenvalues for each element in the family of unitary matrices $\{(u_{2k}^*v_{2k})^2\,:\,k\in\bbn\}$ is completed by \Cref{oo} and \ref{101}. Our next job is to find the multiplicities of the eigenvalues. Throughout the rest of the article, we fix the above enumeration of the eigenvalues of $\big(u_{2k}^*v_{2k}\big)^2$, that is, the eigenvalue $1$ is followed by $\theta^{\,2}$, then $\overline{\theta}^{\,2}$, and so on. In other words, we always arrange the eigenvalues according to increasing power of $\theta$ immediately followed by $\overline{\theta}$. The reason for this will be clear in the sequel.

We denote the multiplicity of the eigenvalue $1$ by $m_k(1)$; and that of $\theta^{\,2r}\mbox{ and }\overline{\theta}^{\,2r}$ by $m_k(\theta^{\,2r})\mbox{ and }m_k(\overline{\theta}^{\,2r})$ respectively for $r=1,2,\ldots, \frac{\lfloor k+1\rfloor}{2}$. First we shall deal with the case of $\,\theta$ not being an even root of unity.

\begin{ppsn}\label{recurrence}
Let $\,\theta$ be not an even root of unity. The multiplicities of the eigenvalues $1\mbox{ and }\theta^{\,2}$ of $(u^*v)^2$ are respectively given by $m_0(1)=3$ and $m_0(\theta^{\,2})=1$, and that of $(u_{2k}^*v_{2k})^2$ for $k\in\bbn$ are given by the following chain of recurrence relations~:
\begin{enumerate}[$(i)$]
\item For $k\in\bbn$ odd,
\begin{IEEEeqnarray*}{lCl}
&  & m_1(1)=10,\,m_1(\theta^{\,2})=5,\,m_1(\overline{\theta}^{\,2})=1,\\
&  & m_k(1)=3m_{k-1}(1)+m_{k-1}(\theta^{\,2}),\\
&  & m_k(\theta^{\,2r})=2m_{k-1}(\theta^{\,2r})+\big(1-\delta_{r,\frac{k+1}{2}}\big)m_{k-1}(\overline{\theta}^{\,2r})+m_{k-1}(\overline{\theta}^{\,2r-2})\quad\mbox{ for } 1\leq r\leq\frac{k+1}{2},\\
&  & m_k(\overline{\theta}^{\,2r})=2\big(1-\delta_{r,\frac{k+1}{2}}\big)m_{k-1}(\overline{\theta}^{\,2r})+m_{k-1}(\theta^{\,2r})+\big(1-\delta_{r,\frac{k+1}{2}}\big)m_{k-1}(\theta^{\,2r+2})\,\mbox{ for } 1\leq r\leq\frac{k+1}{2}.
\end{IEEEeqnarray*}
\item For $k\in\bbn$ even,
\begin{IEEEeqnarray*}{lCl}
&  & m_2(1)=35,\,m_2(\theta^{\,2})=21,\,m_2(\overline{\theta}^{\,2})=7,\,m_2(\theta^{\,4})=1,\\
&  & m_k(1)=3m_{k-1}(1)+m_{k-1}(\theta^{\,2}),\\
&  & m_k(\theta^{\,2r})=2m_{k-1}(\theta^{\,2r})+m_{k-1}(\overline{\theta}^{\,2r})+m_{k-1}(\overline{\theta}^{\,2r-2})\quad\mbox{ for } 1\leq r\leq\frac{k}{2},\\
&  & m_k(\overline{\theta}^{\,2r})=2m_{k-1}(\overline{\theta}^{\,2r})+m_{k-1}(\theta^{\,2r})+\big(1-\delta_{r,\frac{k}{2}}\big)m_{k-1}(\theta^{\,2r+2})\quad\mbox{ for } 1\leq r\leq\frac{k}{2},\\
&  & m_k(\theta^{\,k+2})=m_{k-1}(\overline{\theta}^{\,k})\,.
\end{IEEEeqnarray*}
\end{enumerate}
Here, $\delta_{i,j}$ denotes the Kronecker delta.
\end{ppsn}
\begin{prf}
The case of $k=0$ follows immediately from \Cref{Sa}, and we see that $m_0(1)=3$ and $m_0(\theta^{\,2})=1$ because $(u^*v)^2=p+q(\theta^{\,2})$ (recall $p\mbox{ and }q(\theta)$ from \Cref{p and q}). Now, recall from \Cref{mm} the following (also, recall the fact that $C(\theta^{\,2})=\mbox{diag}\{1,\theta^{\,2}\}\in M_2$),
\begin{IEEEeqnarray*}{lCl}
&  & \mbox{Ad}_{W\otimes W^{(k)}}\big(\big(u_{2k}^*v_{2k}\big)^2\big)\\
&=& \begin{bsmallmatrix}
1 &  &  & \\
 & 1 &  & \\
 &  & 0 & \\
 &  &  & 0\\
\end{bsmallmatrix}\otimes\mbox{Ad}_{W\otimes W^{(k-1)}}\big(\big(u_{2k-2}^*v_{2k-2}\big)^2\big)+\begin{bsmallmatrix}
0 &  &  & \\
 & 0 &  & \\
 &  & 1 & \\
 &  &  & \theta^{\,2}\\
\end{bsmallmatrix}\otimes\mbox{Ad}_{W\otimes W^{(k-1)}}\big(\big(v_{2k-2}^*u_{2k-2}\big)^2\big).
\end{IEEEeqnarray*}
The case of $k\in\bbn$ now easily follows from the above equation by induction on $k$ in the same spirit of \Cref{oo}.\qed
\end{prf}

Thus, once we know the multiplicities of the eigenvalues $1,\theta^{\,2},\overline{\theta}^{\,2}$ when $k=1$, we know that for all $k\in\bbn$. The reason for enumerating the eigenvalues of $\big(u_{2k}^*v_{2k}\big)^2$ according to increasing power of $\theta$ immediately followed by $\overline{\theta}$ shall be clear now from the following result.

\begin{lmma}\label{pg}
Let $\,\theta$ be not an even root of unity. The multiplicities $m_k(.)$ of the eigenvalues of $\big(u_{2k}^*v_{2k}\big)^2$ satisfy the ordering $m_k(1)>m_k(\theta^{\,2})>m_k(\overline{\theta}^{\,2})>\ldots$ for any $k\in\bbn\cup\{0\}$.
\end{lmma}
\begin{prf}
We prove by induction on $k$. The basis step of the induction is clear from \Cref{recurrence}. Assume that the statement holds for the $k$-th step and $k$ is an odd integer. By \Cref{oo}, we have
\[
\sigma\big((u_{2k}^*v_{2k})^2\big)=\big\{1,\theta^{\,2},\overline{\theta}^{\,2},\ldots,\theta^{\,k-1},\overline{\theta}^{\,k-1},\theta^{\,k+1},\overline{\theta}^{\,k+1}\big\}\,,
\]
and by the induction hypothesis we have
\[
m_k(1)>m_k(\theta^{\,2})>m_k(\overline{\theta}^{\,2})>\ldots>m_k(\theta^{\,k+1})>m_k(\overline{\theta}^{\,k+1}).
\]
Now to show for the $(k+1)$-th step, first observe that
\[
 m_{k+1}(1)=3m_k(1)+m_k(\theta^{\,2})>2m_k(\theta^{\,2})+m_k(\overline{\theta}^{\,2})+m_k(1)=m_{k+1}(\theta^{\,2})
 \]
 using \Cref{recurrence}. Similarly, it follows that for each $1\leq r\leq\frac{k-1}{2}$, we have $m_{k+1}(\theta^{\,2r})>m_{k+1}(\overline{\theta}^{\,2r})>m_{k+1}(\theta^{\,2r+2})$. It remains only to show that $m_{k+1}(\theta^{\,k+1})>m_{k+1}(\overline{\theta}^{\,k+1})>m_{k+1}(\theta^{\,k+3})$. Again by \Cref{recurrence}, observe that
\begin{IEEEeqnarray*}{lCl}
m_{k+1}(\theta^{\,k+1}) &=& 2m_k(\theta^{\,k+1})+m_k(\overline{\theta}^{\,k+1})+m_k(\overline{\theta}^{\,k-1})\\
&>& 2m_k(\overline{\theta}^{\,k+1})+m_k(\theta^{\,k+1})\\
&=& m_{k+1}(\overline{\theta}^{\,k+1})
\end{IEEEeqnarray*}
and $m_{k+1}(\overline{\theta}^{\,k+1})=2m_k(\overline{\theta}^{\,k+1})+m_k(\theta^{\,k+1})>m_k(\overline{\theta}^{\,k+1})=m_{k+1}(\theta^{\,k+3})$. This completes the induction and we have the following,
\[
m_{k+1}(1)>m_{k+1}(\theta^{\,2})>m_{k+1}(\overline{\theta}^{\,2})>\ldots>m_{k+1}(\theta^{\,k+1})>m_{k+1}(\overline{\theta}^{\,k+1})>m_{k+1}(\theta^{\,k+3}).
\]
Since $k+1$ is an even integer, by \Cref{oo} we have
\[
\sigma\big((u_{2k+2}^*v_{2k+2})^2\big)=\big\{1,\theta^{\,2},\overline{\theta}^{\,2},\ldots,\theta^{\,k+1},\overline{\theta}^{\,k+1},\theta^{\,k+3}\big\},
\]
and hence we are done.

Exact similar argument works if we start with $k$ even in the induction hypothesis, instead of $k$ being odd.\qed
\end{prf}

\begin{lmma}\label{ug}
Let $\,\theta$ be not an even root of unity. For $k\in\bbn$ odd, the multiplicity $m_k(\overline{\theta}^{\,k+1})$ of the last eigenvalue $\overline{\theta}^{\,k+1}$ obtained in $\Cref{oo}$ is $1$, and for $k\in\bbn$ even, the multiplicity $m_k(\theta^{\,k+2})$ of the last eigenvalue $\theta^{\,k+2}$ is also $1$.
\end{lmma}
\begin{prf}
When $k\in\bbn$ is odd, by induction on $k$ it follows from \Cref{recurrence} that
 \begin{IEEEeqnarray*}{lCl}
m_k(\overline{\theta}^{\,k+1}) &=& 2m_{k-1}(\overline{\theta}^{\,k+1})+m_{k-1}(\theta^{\,k+1})\\
&=& m_{k-1}(\theta^{\,k+1})\\
&=& m_{k-2}(\overline{\theta}^{\,k-1})\\
&=& 1\,.\qquad(\mbox{by induction hypothesis})
 \end{IEEEeqnarray*}
When $k\in\bbn$ is even, by \Cref{recurrence} we have $m_k(\theta^{\,k+2})=m_{k-1}(\overline{\theta}^{\,k})$, and the result follows from the first part.\qed
\end{prf}

Combining \Cref{oo} and \ref{recurrence}, we obtain the following spectral decomposition of each element in the family of unitary matrices $\{\big(u_{2k}^*v_{2k}\big)^2\in M_4^{(k+1)}\,:\,k\in\bbn\}$, when $\,\theta$ is not an even root of unity.

\begin{thm}\label{spec decom}
Let $\,\theta$ be not an even root of unity. The spectral decomposition for the unitary $\big(u_{2k}^*v_{2k}\big)^2$ in $M_4^{(k+1)}$ is $\sum_{j=1}^{k+2}\lambda_jE_j$, where $\lambda_1=1,\lambda_2=\theta^{\,2},\lambda_3=\overline{\theta}^{\,2}$ and so on, as described in $\Cref{oo}$. The dimension of the eigenspaces $E_j$, that is, the geometric multiplicities of the eigenvalues $\lambda_j$, are determined by the chain of recurrence relations obtained in $\Cref{recurrence}$, and we have $\sum_{j=1}^{k+2}\mbox{dim}(E_j)=4^{k+1}$.
\end{thm}
\medskip

Now, we shall deal with the case of $\,\theta$ being an even root of unity. Recall the eigenvalues from \Cref{101}. To distinguish the multiplicities of the eigenvalues in this case with that obtained in \Cref{recurrence}, we denote them by $\widetilde{m}_k$ in this case.

\begin{ppsn}\label{recurrence1}
Let $\theta^{\,m}=1$ for some $m\in 2\bbn+2$ and choose $k=(m-2)/2$. Then, up to $\big(u_{2k-2}^*v_{2k-2}\big)^2$, the multiplicities $\widetilde{m}_{k-1}(\lambda_r)$ of the eigenvalues $\lambda_r$ coincide with $m_{k-1}(\lambda_r)$ as defined in \Cref{recurrence}, and from $\big(u_{2k}^*v_{2k}\big)^2$ onward, we have the following new chain of recurrence relations for the multiplicities $\widetilde{m}_{k+j}(.)$ for $0\leq j\leq k-1\,$:
\begin{enumerate}[$(i)$]
\item If $k\in\bbn$ is odd, then for any $0\leq j\leq k-1$ we have
\begin{enumerate}[$(a)$]
\item if $j\geq 0$ even, then
\begin{IEEEeqnarray*}{lCl}
&  & \widetilde{m}_{k+j}(1)=m_{k+j}(1),\,\widetilde{m}_{k+j}(\theta^{\,2})=m_{k+j}(\theta^{\,2}),\ldots,\,\widetilde{m}_{k+j}(\overline{\theta}^{\,k-(j+1)})=m_{k+j}(\overline{\theta}^{\,k-(j+1)})\\
&  & \widetilde{m}_{k+j}(\theta^{\,k+1})=m_{k+j}(\theta^{\,k+1})+m_{k+j}(\overline{\theta}^{\,k+1}),\\
&  & \widetilde{m}_{k+j}(\overline{\theta}^{\,k-r})=m_{k+j}(\overline{\theta}^{\,k-r})+m_{k+j}(\theta^{\,k+r+2})\,\,\,\,\mbox{ for }\,1\leq r\leq j-1\,\mbox{ odd integer},\\
&  & \widetilde{m}_{k+j}(\theta^{\,k-r})=m_{k+j}(\theta^{\,k-r})+m_{k+j}(\overline{\theta}^{\,k+r+2})\,\,\,\mbox{ for }\,1\leq r\leq j-1\,\mbox{ odd integer}.
\end{IEEEeqnarray*}
\item if $j\geq 1$ odd, then
\begin{IEEEeqnarray*}{lCl}
&  & \widetilde{m}_{k+j}(1)=m_{k+j}(1),\,\widetilde{m}_{k+j}(\theta^{\,2})=m_{k+j}(\theta^{\,2}),\ldots,\,\widetilde{m}_{k+j}(\theta^{\,k-j})=m_{k+j}(\theta^{\,k-j})\\
&  & \widetilde{m}_{k+j}(\theta^{\,k+1})=m_{k+j}(\theta^{\,k+1})+m_{k+j}(\overline{\theta}^{\,k+1}),\\
&  & \widetilde{m}_{k+j}(\overline{\theta}^{\,k-r})=m_{k+j}(\overline{\theta}^{\,k-r})+m_{k+j}(\theta^{\,k+r+2})\,\,\,\,\mbox{ for }\,1\leq r\leq j\,\mbox{ odd integer},\\
&  & \widetilde{m}_{k+j}(\theta^{\,k-r})=m_{k+j}(\theta^{\,k-r})+m_{k+j}(\overline{\theta}^{\,k+r+2})\,\,\,\mbox{ for }\,1\leq r\leq j-2\,\mbox{ odd integer}.
\end{IEEEeqnarray*}
\end{enumerate}
\item If $k\in\bbn$ is even, then for any $0\leq j\leq k-1$ we have
\begin{enumerate}[$(a)$]
\item if $j\geq 0$ even, then
\begin{IEEEeqnarray*}{lCl}
&  & \widetilde{m}_{k+j}(1)=m_{k+j}(1),\,\widetilde{m}_{k+j}(\theta^{\,2})=m_{k+j}(\theta^{\,2}),\ldots,\,\widetilde{m}_{k+j}(\theta^{\,k-j})=m_{k+j}(\theta^{\,k-j})\\
&  & \widetilde{m}_{k+j}(\overline{\theta}^{\,k-r})=m_{k+j}(\overline{\theta}^{\,k-r})+m_{k+j}(\theta^{\,k+r+2})\,\,\,\,\mbox{ for }\,0\leq r\leq j\,\mbox{ even integer},\\
&  & \widetilde{m}_{k+j}(\theta^{\,k-r})=m_{k+j}(\theta^{\,k-r})+m_{k+j}(\overline{\theta}^{\,k+r+2})\,\,\,\mbox{ for }\,0\leq r\leq j-2\,\mbox{ even integer}.
\end{IEEEeqnarray*}
\item if $j\geq 1$ odd, then
\begin{IEEEeqnarray*}{lCl}
&  & \widetilde{m}_{k+j}(1)=m_{k+j}(1),\,\widetilde{m}_{k+j}(\theta^{\,2})=m_{k+j}(\theta^{\,2}),\ldots,\,\widetilde{m}_{k+j}(\overline{\theta}^{\,k-j-1})=m_{k+j}(\overline{\theta}^{\,k-j-1})\\
&  & \widetilde{m}_{k+j}(\overline{\theta}^{\,k-r})=m_{k+j}(\overline{\theta}^{\,k-r})+m_{k+j}(\theta^{\,k+r+2})\,\,\,\,\mbox{ for }\,0\leq r\leq j-1\,\mbox{ even integer},\\
&  & \widetilde{m}_{k+j}(\theta^{\,k-r})=m_{k+j}(\theta^{\,k-r})+m_{k+j}(\overline{\theta}^{\,k+r+2})\,\,\,\mbox{ for }\,0\leq r\leq j-1\,\mbox{ even integer}.
\end{IEEEeqnarray*}
\end{enumerate}
\end{enumerate}
Here, $\,m_{k+j}(.)$ is as defined in $\Cref{recurrence}$.
\end{ppsn}
\begin{prf}
Recall \Cref{oo} in this regard. First assume that $k$ is an even integer. If $j$ increases as an even integer starting from $0$, then the eigenvalue $\theta^{\,k+j+2}$ coincides with $\overline{\theta}^{\,k-j}$ for all $j$ and $\overline{\theta}^{\,k+j}$ coincides with $\theta^{\,k-j+2}$ for all $j$. Moreover, multiplicities of the eigenvalues from $1$ up to $\theta^{\,k-j}$ remain unchanged. This is because $\theta^{\,m}=\theta^{\,2k+2}=1$. Similarly, if $j$ increases as an odd integer starting from $1$, then the eigenvalues $\theta^{\,k+j+1}$ and $\overline{\theta}^{\,k+j+1}$ coincide with $\overline{\theta}^{\,k-j+1}$ and $\theta^{\,k-j+1}$ respectively for all $j$ because of similar reason, and multiplicities of the eigenvalues from $1$ up to $\overline{\theta}^{\,k-j-1}$ remain unchanged. When $k$ is an odd integer, exact similar analysis holds. Now, the result follows in the same spirit of \Cref{recurrence} using the following identity (\Cref{mm})
\begin{IEEEeqnarray*}{lCl}
&  & \mbox{Ad}_{W\otimes W^{(k)}}\big(\big(u_{2k}^*v_{2k}\big)^2\big)\\
&=& \begin{bsmallmatrix}
1 &  &  & \\
 & 1 &  & \\
 &  & 0 & \\
 &  &  & 0\\
\end{bsmallmatrix}\otimes\mbox{Ad}_{W\otimes W^{(k-1)}}\big(\big(u_{2k-2}^*v_{2k-2}\big)^2\big)+\begin{bsmallmatrix}
0 &  &  & \\
 & 0 &  & \\
 &  & 1 & \\
 &  &  & \theta^{\,2}\\
\end{bsmallmatrix}\otimes\mbox{Ad}_{W\otimes W^{(k-1)}}\big(\big(v_{2k-2}^*u_{2k-2}\big)^2\big).
\end{IEEEeqnarray*}
as the basic tool.\qed
\end{prf}

\begin{rmrk}\rm\label{guin}
If we consider the unitary $\big(u_{2k}^*v_{2k}\big)^2$ and put the condition $\theta^{\,2k+2}=1$ for some $k\in\bbn$, then the multiplicity of only the last eigenvalue changes, irrespective of $k$ being even or odd (compare the situations in \Cref{recurrence} and putting $j=0$ in \Cref{recurrence1}). The multiplicities of all the remaining eigenvalues do not change. This fact will be crucially used in subsequent sections. This is another reason (apart from \Cref{pg}) of choosing the enumeration of the eigenvalues according to increasing power of $\theta$ immediately followed by $\overline{\theta}$, to keep track of what are the changes when $k$ increases.
\end{rmrk}

Note that in the above we have obtained the multiplicities $\widetilde{m}_{k+j}(.)$ of the eigenvalues only for $0\leq j\leq k-1$. In this article, we do not need $j$ to go beyond $k-1$. This will be clear when we apply this in the subsection $7.2.3$. We remark that for $j\geq k$, with a little more effort, it is possible to obtain all the multiplicities $\widetilde{m}_{k+j}(.)$ in terms of recurrence relations, and write down the full spectral decomposition for this case of $\theta^{\,m}=1$ for $m\in 2\bbn+2$. However, this will be completely unnecessary for the purpose of this paper, and to stay focused to our goal, we leave it to the interested readers.

Combining \Cref{101} and \ref{recurrence1}, we obtain the following spectral decomposition when $\theta^{\,m}=1,\,m\in2\bbn+2$ (i,e., when $\theta$ is an even root of unity).

\begin{thm}\label{spec decom1}
Let $\theta^{\,m}=1$ for some $m\in 2\bbn+2$ and choose $k=(m-2)/2$. The spectrum of $(u_{2k+2j}^*v_{2k+2j})^2$ remains the same with cardinality $k+1$ for all $j\in\bbn\cup\{-1,0\}$ as described in \Cref{101}. The spectral decomposition for the unitary $\big(u_{2k+2j}^*v_{2k+2j}\big)^2,\,0\leq j\leq k-1,$ is $\sum_{i=1}^{k+1}\lambda_iE_i^{(j)}$, where the dimension of the eigenspaces $E_i^{(j)}$, that is, the geometric multiplicities of the eigenvalues $\lambda_i$, are determined by the chain of recurrence relations described in $\Cref{recurrence1}$, and we have $\sum_{i=1}^{k+1}\mbox{dim}(E_i^{(j)})=4^{k+1+j}$ for each $0\leq j\leq k-1$.
\end{thm}


\subsubsection{A family of von Neumann subalgebras of $R$ with the infinite\\
 Pimsner-Popa index}

The spectral analysis carried out in the previous subsection says that if we consider the following subgroup of $\mathbb{S}^1$
\begin{IEEEeqnarray}{lCl}\label{main set}
\Gamma:=\{\omega\in\mathbb{S}^1:\omega^m=1,\,m\in 2\bbn+2\}\cup\{1\}\,,
\end{IEEEeqnarray}
that is, $\Gamma$ consists of all even roots of unity, then the spectrum of $\big(u_{2k}^*v_{2k}\big)^2$ depends on whether $\theta\in\Gamma$ or $\theta\in\mathbb{S}^1\setminus\Gamma$, where $\theta=\overline{a}b$ (note that $\theta\neq 1$ as $b\neq a$). It turns out that depending on whether $\theta$ is a primitive even root of unity or not, that is, whether $\theta^{\,2}$ is a rational rotation or an irrational rotation, the intersection $R_u\cap R_v$, which is {\it a priori} only a von Neumann subalgebra of the hyperfinite $II_1$ factor $R$, has finite or infinite Pimsner-Popa index in $R$. The set $\Gamma$ can also be interpreted in terms of an equivalence relation among the Hadamard matrices $u$ and $v$, which we discuss now.
\smallskip

\textbf{An equivalence relation~:} Let $u(a)\mbox{ and }u(b)$ be two distinct $4\times 4$ Hadamard inequivalent complex Hadamard matrices parametrized by circle parameters $a\mbox{ and }b$ respectively (recall that $b\neq\pm a$ in our situation). We define a relation $u(a)\sim u(b)$ if and only if $\theta(a,b):=\overline{a}b\in\Gamma$, that is, $\theta(a,b)^{\,m}=1\mbox{ for some }m\in 2\bbn+2$. Note that $\theta(a,b)\neq 1$. It is easy to check that `$\sim$' is an equivalence relation. Therefore, for $u=u(a)$ and $v=u(b)$, the analysis carried out in previous subsection indicates that the spectrum of $\big(u_{2k}^*v_{2k}\big)^2$ for $k\in\bbn$ depends on whether $u\sim v$ or $u\nsim v$.
\smallskip

\noindent\textbf{Caution:} This equivalence relation has no relation with the one defined in \Cref{Sec 3}.
\smallskip

We shall see that this equivalence relation plays a major role in our analysis here. In this subsection, we will deal with the case of $u\nsim v$ (i,e. $\theta\notin\Gamma$) and request the reader to wait until the next subsection, where $u\sim v$ (i,e. $\theta\in\Gamma$) will be dealt.

Recall once again the unitary matrix $W\in M_4$ from \Cref{W}. Also, note that by $W^{(k)},\,k\in\bbn,$ we mean the matrix $W\otimes W\otimes\ldots\otimes W\in M_4^{(k)}:=M_4^{\,\otimes\,k}$.

\begin{lmma}\label{mb}
We have $\big((u^*v)^2\big)^\prime\cap M_4=\mathrm{Ad}_W(M_3\oplus\bbc)$.
\end{lmma}
\begin{prf}
Since $\mbox{Ad}_W(u^*v)^2=\mbox{bl-diag}\{I_2,C(\theta^{\,2})\}=\mbox{diag}\{I_3,\theta^{\,2}\}$ in $M_4$ due to \Cref{Sa}, an application of \Cref{com} and \ref{com and ad interchange} proves the claim.\qed
\end{prf}

\begin{ppsn}\label{aef}
If $u\nsim v$, then there exist certain unitary matrices $\,U_k\in M_4^{(k+1)}$ for each $k\in\bbn$ such that the following holds,
\begin{IEEEeqnarray*}{lCl}
&  & \big((u_{2k}^*v_{2k})^2\big)^\prime\cap M_4^{(k+1)}\\
&=& \begin{cases}
\mathrm{Ad}_{W^{(k+1)}U_k}\Big(M_{m_k(1)}\bigoplus\bigoplus_{r=1}^{\frac{k+1}{2}}\Big(M_{m_k(\theta^{\,2r})}\oplus M_{m_k\big(\overline{\theta}^{\,2r}\big)}\Big)\Big)\,\quad\qquad\mbox{if k is odd};\\
\mathrm{Ad}_{W^{(k+1)}U_k}\Big(M_{m_k(1)}\bigoplus\bigoplus_{r=1}^{\frac{k}{2}}\Big(M_{m_k(\theta^{\,2r})}\oplus M_{m_k\big(\overline{\theta}^{\,2r}\big)}\Big)\bigoplus\bbc\Big)\quad\mbox{if k is even}.
\end{cases}
\end{IEEEeqnarray*}
\end{ppsn}
\begin{prf}
We use \Cref{com and ad interchange} throughout. For $k=1$, by \Cref{Sa} and \Cref{mm} we have
\begin{IEEEeqnarray*}{lCl}
\mathrm{Ad}_{W\otimes W}\big((u_2^*v_2)^2\big) &=& E_{11}\otimes I_2\otimes \mathrm{Ad}_W\big((u^*v)^2\big)+E_{22}\otimes C(\theta^{\,2})\otimes \mathrm{Ad}_W\big((v^*u)^2\big)\\
&=& E_{11}\otimes I_2\otimes\big(\mbox{diag}\{1,1,1,\theta^{\,2}\}\big)+E_{22}\otimes C(\theta^{\,2})\otimes\big(\mbox{diag}\{1,1,1,\overline{\theta}^{\,2}\}\big)\,.
\end{IEEEeqnarray*}
There exists a unitary matrix (a permutation matrix to be precise) $U_1\in M_4^{(2)}$ which leaves $\Delta_4^{(2)}$ invariant and
\begin{IEEEeqnarray}{lCl}\label{lkrs}
\mathrm{Ad}_{W\otimes W}\big((u_2^*v_2)^2\big) &=& \mathrm{Ad}_{U_1}\big(I_{10}\oplus\theta^{\,2}I_5\oplus\overline{\theta}^{\,2}\big)\,.
\end{IEEEeqnarray}
For example, take $U_1=\big(E_{4,11}+E_{11,4}+\sum_{i\neq 4,11}E_{ii}\big)\big(E_{8,12}+E_{12,8}+\sum_{i\neq 8,12}E_{ii}\big)\big(E_{12,16}+E_{16,12}+\sum_{i\neq 12,16}E_{ii}\big)$. Now, applying \Cref{com and ad interchange}, \Cref{lkrs}, and \Cref{recurrence}, we have the following,
\begin{IEEEeqnarray*}{lCl}
\mathrm{Ad}_{U_1^*(W\otimes W)}\Big(\big((u_2^*v_2)^2\big)^\prime\cap M_4^{(2)}\Big) &=& \big(\mathrm{Ad}_{U_1^*(W\otimes W)}(u_2^*v_2)^2\big)^\prime\cap M_4^{(2)}\\
&=& \big(I_{10}\oplus\theta^{\,2}I_5\oplus\overline{\theta}^{\,2}\big)^\prime\cap M_{16}\\
&=& M_{10}\oplus M_5\oplus\bbc\qquad (\mbox{by } \Cref{com},\mbox{ since }\theta\notin\Gamma)\\
&=& M_{m_1(1)}\oplus M_{m_1(\theta^{\,2})}\oplus M_{m_1(\overline{\theta}^{\,2})},
\end{IEEEeqnarray*}
which finishes the $k=1$ case. Now, suppose that $k$ is odd and we have a unitary matrix $U_k$ in $M_4^{(k+1)}$ such that the following holds
\begin{IEEEeqnarray}{lCl}\label{pa}
\mathrm{Ad}_{W^{(k+1)}}\big((u_{2k}^*v_{2k})^2\big)=\mathrm{Ad}_{U_k}\Big(I_{m_k(1)}\bigoplus\bigoplus_{r=1}^{\frac{k+1}{2}}\Big(\theta^{\,2r}I_{m_k(\theta^{\,2r})}\oplus \overline{\theta}^{\,2r}I_{m_k\big(\overline{\theta}^{\,2r}\big)}\Big)\Big)\,.
\end{IEEEeqnarray}
Using \Cref{mm}, by \Cref{oo} and \Cref{pa} we have the following,
\begin{IEEEeqnarray}{lCl}\label{ppq}
&  & \mathrm{Ad}_{W^{(k+2)}}\big((u_{2k+2}^*v_{2k+2})^2\big)\\
&=& E_{11}\otimes I_2\otimes\mathrm{Ad}_{W^{(k+1)}}(u_{2k}^*v_{2k})^2+E_{22}\otimes C(\theta^{\,2})\otimes\mathrm{Ad}_{W^{(k+1)}}(v_{2k}^*u_{2k})^2\nonumber\\
&=& E_{11}\otimes I_2\otimes\mathrm{Ad}_{U_k}\Big(I_{m_k(1)}\bigoplus\,\bigoplus_{r=1}^{\frac{k+1}{2}}\Big(\theta^{\,2r}\,I_{m_k(\theta^{\,2r})}\oplus\overline{\theta}^{\,2r}\,I_{m_k\big(\overline{\theta}^{\,2r}\big)}\Big)\Big)\nonumber\\
&  & +E_{22}\otimes C(\theta^{\,2})\otimes\mathrm{Ad}_{U_k}\Big(I_{m_k(1)}\bigoplus\,\bigoplus_{r=1}^{\frac{k+1}{2}}\Big(\theta^{\,2r}\,I_{m_k(\theta^{\,2r})}\oplus\overline{\theta}^{\,2r}\,I_{m_k\big(\overline{\theta}^{\,2r}\big)}\Big)\Big)^*\nonumber\\
&=& \mathrm{Ad}_{(I_4\otimes U_k)V_k}\Big(I_{m_{k+1}(1)}\bigoplus \bigoplus_{r=1}^{\frac{k+1}{2}}\Big(\theta^{\,2r}\,I_{m_{k+1}(\theta^{\,2r})}\oplus \overline{\theta}^{\,2r}\,I_{m_{k+1}(\overline{\theta}^{\,2r})}\Big)\bigoplus \theta^{\,k+3}\,I_{m_{k+1}(\theta^{\,k+3})}\Big)\nonumber\,,
\end{IEEEeqnarray}
where the last line follows from \Cref{recurrence}. Here, $V_k$ is a permutation matrix which interchange the diagonal entries suitably to put all the same eigenvalue together at one block which are initially situated in various places at the diagonal. Let $U_{k+1}=(I_4\otimes U_k)V_k$, which is a unitary matrix in $M_4^{(k+2)}$, and we obtain the following by \Cref{ppq},
\begin{IEEEeqnarray*}{lCl}
&  & \mathrm{Ad}_{U_{k+1}^*W^{(k+2)}}\Big(\big((u_{2k+2}^*v_{2k+2})^2\big)^\prime\cap M_4^{(k+2)}\Big)\\
&=& \big(\mathrm{Ad}_{U_{k+1}^*W^{(k+2)}}(u_{2k+2}^*v_{2k+2})^2\big)^\prime\cap M_4^{(k+2)}\\
&=& \Big(I_{m_{k+1}(1)}\bigoplus \bigoplus_{r=1}^{\frac{k+1}{2}}\Big(\theta^{\,2r}\,I_{m_{k+1}(\theta^{\,2r})}\oplus \overline{\theta}^{\,2r}\,I_{m_{k+1}(\overline{\theta}^{\,2r})}\Big)\bigoplus \theta^{\,k+3}\,I_{m_{k+1}(\theta^{\,k+3})}\Big)^\prime\cap M_4^{(k+2)}\\
&=& M_{m_{k+1}(1)}\bigoplus \bigoplus_{r=1}^{\frac{k+1}{2}}\Big(M_{m_{k+1}(\theta^{\,2r})}\oplus M_{m_{k+1}(\overline{\theta}^{\,2r})}\Big)\bigoplus M_{m_{k+1}(\theta^{\,k+3})}\quad (\mbox{by } \Cref{com},\mbox{ since }\theta\notin\Gamma).
\end{IEEEeqnarray*}
Now, $m_{k+1}(\theta^{\,k+3})=1$ by \Cref{ug}, which completes the induction argument. Exactly similar argument works when we start the induction with $k$ even instead of $k$ being odd.\qed
\end{prf}

\begin{ppsn}\label{intersection infinite}
If $\,u\nsim v$, then we have $\,C_1=\mathrm{Ad}_{u_1(I_4\otimes W)}\big(\bbc\otimes(M_3\oplus\bbc)\big)$, and for any $k\in\bbn$
\begin{IEEEeqnarray*}{lCl}
&  & C_{2k+1}\\
&=& \begin{cases}
\mathrm{Ad}_{u_{2k+1}\big(I_4\otimes W^{(k+1)}U_k\big)}\Big(\bbc\otimes\Big(M_{m_k(1)}\bigoplus\bigoplus_{r=1}^{\frac{k+1}{2}}\Big(M_{m_k(\theta^{\,2r})}\oplus M_{m_k\big(\overline{\theta}^{\,2r}\big)}\Big)\Big)\Big)\,\,\,\qquad\mbox{if k is odd};\\
\mathrm{Ad}_{u_{2k+1}\big(I_4\otimes W^{(k+1)}U_k\big)}\Big(\bbc\otimes\Big(M_{m_k(1)}\bigoplus\bigoplus_{r=1}^{\frac{k}{2}}\Big(M_{m_k(\theta^{\,2r})}\oplus M_{m_k\big(\overline{\theta}^{\,2r}\big)}\Big)\bigoplus\bbc\Big)\Big)\,\,\mbox{if k is even};
\end{cases}
\end{IEEEeqnarray*}
where the unitary matrices $U_k$'s are as in \Cref{aef}.
\end{ppsn}
\begin{prf}
By \Cref{odd intersection}, we know that $C_{2k+1}$ for $k\geq 0$ is completely determined by $\big((u_{2k}^*v_{2k})^2\big)^\prime\cap M_4^{(k+1)}$ up to the conjugation by the unitary $u_{2k+1}$. Now, \Cref{mb} proves the case of $C_1$, and \Cref{aef} proves the claim for any $k\in\bbn$.\qed
\end{prf}
\smallskip

For the sake of better understanding, we mention below first few commutants $\big((u_{2k}^*v_{2k})^2\big)^\prime\cap M_4^{(k+1)},\,k\in\bbn\cup\{0\},$ concretely up to unitary matrices~:
\[
\begin{matrix}
\underline{\mbox{for }k=0}\\
M_4\\
\cup\\
M_3\oplus\bbc\\
\end{matrix}\qquad\begin{matrix}
\underline{\mbox{for }k=1}\\
M_4\otimes M_4\\
\cup\\
M_{10}\oplus M_5\oplus\bbc\\
\end{matrix}\qquad\begin{matrix}
\underline{\mbox{for }k=2}\\
M_4\otimes M_4\otimes M_4\\
\cup\\
M_{35}\oplus M_{21}\oplus M_7\oplus\bbc\\
\end{matrix}\qquad\begin{matrix}
\underline{\mbox{for }k=3}\\
M_4\otimes M_4\otimes M_4\otimes M_4\\
\cup\\
M_{126}\oplus M_{84}\oplus M_{36}\oplus M_{9}\oplus\bbc\\
\end{matrix}
\]
The first vertical inclusion for $k=0$ is determined by the unitary $W$ in $M_4$, and all the others are determined by the unitary $W^{(k+1)}U_k$ in $M_4^{(k+1)}$ for $k\in\bbn$, where $U_k$'s are the unitary matrices in $M_4^{(k+1)}$ obtained in \Cref{aef}. Therefore, in view of \Cref{odd intersection}, first few stages of the following quadruple for $k\in\bbn$
\[
\begin{matrix}
A_{2k-1} &\subset & A_{2k+1}\\
|| &  & ||\\
\bbc\otimes\Delta_4\otimes M_4^{(k)} & & \Delta_4\otimes M_4^{(k+1)}\\
\cup &  & \cup\\
\bbc\otimes\mbox{Ad}_{u_{2k-1}}\big(\bbc\otimes(M_4^{(k)})\big)\cap\mbox{Ad}_{v_{2k-1}}\big(\bbc\otimes(M_4^{(k)})\big) & & \mbox{Ad}_{u_{2k+1}}\big(\bbc\otimes(M_4^{(k+1)})\big)\cap\mbox{Ad}_{v_{2k+1}}\big(\bbc\otimes(M_4^{(k+1)})\big)\\
|| &  & ||\\
C_{2k-1} &\subset & C_{2k+1}\\
\end{matrix}
\]
are the following~:
\smallskip

\noindent \textbf{For $k=1:$}
\[
\begin{matrix}
\bbc\otimes\Delta_4\otimes M_4 &\subset & \Delta_4\otimes M_4\otimes M_4\\
\cup &  & \cup\\
\bbc\otimes\mbox{Ad}_{u_1(I_4\otimes W)}(\bbc\otimes(M_3\oplus\bbc)) &\subset & \mbox{Ad}_{u_3\big(I_4\otimes W^{(2)}U_1\big)}(\bbc\otimes(M_{10}\oplus M_5\oplus\bbc))\\
\end{matrix}
\]
\smallskip

\noindent \textbf{For $k=2:$}
\[
\begin{matrix}
\bbc\otimes\Delta_4\otimes M_4\otimes M_4 &\subset & \Delta_4\otimes M_4\otimes M_4\otimes M_4\\
\cup &  & \cup\\
\bbc\otimes\mbox{Ad}_{u_3\big(I_4\otimes W^{(2)}U_1\big)}\big(\bbc\otimes(M_{10}\oplus M_5\oplus\bbc)\big) &\subset & \mbox{Ad}_{u_5\big(I_4\otimes W^{(3)}U_2\big)}(\bbc\otimes(M_{35}\oplus M_{21}\oplus M_7\oplus\bbc))\\
\end{matrix}
\]
\smallskip

\noindent \textbf{For $k=3:$}
\[
\begin{matrix}
\mbox{Ad}_{u_7\big(I_4\otimes W^{(4)}U_3\big)}(\bbc\otimes(M_{126}\oplus M_{84}\oplus M_{36}\oplus M_9\oplus\bbc)) &\subset & \Delta_4\otimes M_4\otimes M_4\otimes M_4\otimes M_4\\
\cup &  & \cup\\
\bbc\otimes\mbox{Ad}_{u_5\big(I_4\otimes W^{(3)}U_2\big)}(\bbc\otimes(M_{35}\oplus M_{21}\oplus M_7\oplus\bbc)) &\subset & \bbc\otimes\Delta_4\otimes M_4\otimes M_4\otimes M_4\\
\end{matrix}
\]

\begin{lmma}\label{bca}
For $u\nsim v$ and $k\in\bbn\cup\{0\}$, the inclusion matrix for the inclusion $\big((u_{2k}^*v_{2k})^2\big)^\prime\cap M_4^{(k+1)}\subset M_4\otimes M_4^{(k)}$ is the column matrix $(1\,\,1\,\ldots\,1)^{T}$ of order $(k+2)\times 1$.
\end{lmma}
\begin{prf}
By \Cref{spec decom}, we have the spectral decomposition of the unitary matrix $(u_{2k}^*v_{2k})^2$ for all $k\in\bbn\cup\{0\}$. Suppose that $k$ is odd. Then, $\{1,\theta^{\,2},\overline{\theta}^{\,2},\ldots,\theta^{\,k+1},\overline{\theta}^{\,k+1}\}$ is the set of distinct eigenvalues of $(u_{2k}^*v_{2k})^2$, and $m_k(1)+\sum_rm_k(\theta^{\,2r})+m_k(\overline{\theta}^{\,2r})=4^{k+1}$, where $m_k(.)$ are the multiplicities of the eigenvalues. By \Cref{intersection infinite} it is now clear that the required inclusion matrix is the column matrix $(1\,\,1\,\ldots\,1)^{T}$ of order $(k+2)\times 1$, as there are $(k+2)$-many eigenvalues of $(u_{2k}^*v_{2k})^2$ by \Cref{oo}. Exactly similar analysis holds for $k$ even.\qed
\end{prf}

\begin{lmma}\label{inf}
For $u\nsim v$ and $k\in\bbn\cup\{0\}$, the inclusion matrix $\Lambda_k$ for the inclusion
\[
C_{2k+1}:=\mathrm{Ad}_{u_{2k+1}}\big(\bbc\otimes M_4^{(k+1)}\big)\cap\mathrm{Ad}_{v_{2k+1}}\big(\bbc\otimes M_4^{(k+1)}\big)\subset \Delta_4\otimes M_4\otimes M_4^{(k)}=:A_{2k+1}
 \]
is the $(k+2)\times 4$ matrix each of whose entry is $1$. Thus, $||\Lambda_k||^2=4(k+2)$.
\end{lmma}
\begin{prf}
Since $\Delta_4\otimes M_4\otimes M_4^{(k)}$ is the direct sum of four copies of $M_4\otimes M_4^{(k)}$, and the unitary matrix $u_{2k+1}(I_4\otimes W^{(k+1)}U_k)$ in \Cref{intersection infinite} lies in $\Delta_4\otimes M_4\otimes M_4^{(k)}$, each entry of the inclusion matrix $\Lambda_k$ equals to $1$ is clear by \Cref{bca}. To determine $||\Lambda_k||^2$, observe that $\Lambda_k\Lambda_k^T$ is the matrix $4(k+2)J_{k+2}$, where $J_{k+2}$ is the projection $\sum_{j=1}^{k+2}(k+2)^{-1}E_{jj}$ in $M_{k+2}$.\qed
\end{prf}

\begin{rmrk}\rm\label{justification of main difficulty}
In the following tower of commuting squares
\[
\begin{matrix}
A_1 &\subset & A_3 &\subset  & A_5 &\subset &\ldots\cr
\cup &  &\cup &  &\cup &  &  & \cr
C_1 &\subset & C_3 &\subset & C_5 &\subset &\ldots 
\end{matrix}
\]
the first quadruple $(C_1\subset A_1,C_3\subset A_3)$ is not non-degenerate. This is because the inclusion matrix $\Lambda_{C_1}^{A_1}=\Lambda_0$ is a $2\times 4$ matrix with $||\Lambda_{C_1}^{A_1}||^2=8$ by \Cref{inf}, whereas the inclusion matrix $\Lambda_{C_3}^{A_3}=\Lambda_1$ is a $3\times 4$ matrix with $||\Lambda_{C_3}^{A_3}||^2=12$. That is, $||\Lambda_{C_1}^{A_1}||^2\neq ||\Lambda_{C_3}^{A_3}||^2$. In fact, none of the quadruples are non-degenerate, which follows from \Cref{inf}. This also justifies part $(ii)$ in \Cref{main difficulty}. Moreover, we can also conclude that $C_1\subset C_3\subset C_5\subset\ldots$ is not a basic construction, as otherwise the quadruple $(C_1\subset A_1,C_3\subset A_3)$ would have been non-degenerate by \cite{BG}.
\end{rmrk}

\begin{thm}\label{infinite index}
If $u\nsim v$, then $\lambda(R,R_u\cap R_v)=\lambda(R_u,R_u\cap R_v)=\lambda(R_v,R_u\cap R_v)=0$. That is, $R_u\cap R_v$ is a von Neumann subalgebra of the hyperfinite $II_1$ factor $R$ with infinite Pimsner-Popa index.
\end{thm}
\begin{prf}
First observe that it is enough to only conclude that $\lambda(R,R_u\cap R_v)=0$. This is because by the submultiplicativity property of $\lambda$ for $R_u\cap R_v\subset R_u\subset R$, we have $\lambda(R,R_u\cap R_v)\geq\lambda(R_u,R_u\cap R_v)\lambda(R,R_u)=\frac{1}{4}\lambda(R_u,R_u\cap R_v)$ (similarly for $v$ in place of $u$) because $\lambda(R,R_u)=[R:R_u]^{-1}=\frac{1}{4}$ by \cite{PP}. Hence if $\lambda(R,R_u\cap R_v)=0$, then automatically we get that $\lambda(R_u,R_u\cap R_v)=0$.

Now by \Cref{ent formula}, we have the following
\[
\lambda\left(\Delta_4\otimes M_4^{(k+1)}\,,\,\mathrm{Ad}_{u_{2k+1}}\big(\bbc\otimes M_4^{(k+1)}\big)\cap\mathrm{Ad}_{v_{2k+1}}\big(\bbc\otimes M_4^{(k+1)}\big)\right)^{-1}=\max_\ell\Big(\sum_{r=1}^{k+2}\frac{b_{r\ell}s_r}{t_\ell}\Big)\,,
\]
where $b_{r\ell}=\min\{a_{r\ell},n_r\}$ and $s_r$ denotes the trace of minimal projections in each summand $N_r$ of $\big((u_{2k}^*v_{2k})^2\big)^\prime\cap M_4^{(k+1)}$ described in \Cref{aef}, and $t_\ell$ denotes trace of minimal projections in $\bbc\otimes M_4^{(k+1)}\subset\Delta_4\otimes M_4^{(k+1)}$ for $\ell=1,\ldots,4$. Hence, $t_\ell=\frac{1}{4^{(k+2)}}$ for each $\ell$ (recall the unique trace on $R$ obtained through the basic construction as in \Cref{fig2}), and since the inclusion matrix $\Lambda_k$ is the $(k+2)\times 4$ matrix each of whose entry is $1$, we get that $b_{r\ell}=1$ for each $r,\ell$. This says that $s_r=\frac{1}{4^{(k+1)}}$ for each $1\leq r\leq k+2$ because $\Lambda_k(t_1\,\,t_2\,\,t_3\,\,t_4)^T=(s_1\,\,\cdots\,\,s_{k+2})^T$. Thus,
\[
\lambda^{-1}=4^{(k+2)}\sum_{r=1}^{k+2}s_r=4^{(k+2)}\sum_{r=1}^{k+2}\frac{1}{4^{(k+1)}}=4(k+2)\,.
\]
Recall that $R_u\cap R_v$ is the SOT-limit of $C_{2k+1}=\mathrm{Ad}_{u_{2k+1}}\big(\bbc\otimes M_4^{(k+1)}\big)\cap\mathrm{Ad}_{v_{2k+1}}\big(\bbc\otimes M_4^{(k+1)}\big)$. Now, using \Cref{popaadaptation} we get that $\lambda(R,R_u\cap R_v)=\lim_{k\to\infty}\frac{1}{4(k+2)}=0$.\qed
\end{prf}

Therefore, we obtain the following family of von Neumann subalgebras of $R$
\[
\left\{R_{u(a)}\cap R_{u(b)}\subset R:\,u(a)\nsim u(b);\,a,b\in\mathbb{S}^1\right\}\,,
\]
such that any member of this family has infinite Pimsner-Popa index in $R$.
\medskip

\noindent\textbf{Question:} Is $R_u\cap R_v$ for $u\nsim v$ a factor?


\subsubsection{Construction of a family of finite index subfactors of $R$}

We now deal with the case of $u\sim v$. Before we jump into the general situation, first let us understand few initial cases to realize where it differs from the corresponding situation in $u\nsim v$. Recall \Cref{intersection infinite} in this regard. Consider again the following tower of commuting squares
\[
\begin{matrix}
A_1 &\subset & A_3 &\subset  & A_5 &\subset &\ldots\cr
\cup &  &\cup &  &\cup &  &  & \cr
C_1 &\subset & C_3 &\subset & C_5 &\subset &\ldots 
\end{matrix}
\]
In the case of $u\nsim v$, for $k=0$, we have seen that the inclusion matrix $\Lambda_0$ for the inclusion $C_1=\mbox{Ad}_{u_1(I_4\otimes W)}(\bbc\otimes(M_3\oplus\bbc))\subset\Delta_4\otimes M_4=A_1$ is the $2\times 4$ matrix each of whose entry is $1$ (\Cref{inf}), and the square of its norm is $||\Lambda_0||^2=8$. Now let $k=1$. The inclusion matrix $\Lambda_1$ for the inclusion $C_3=\mbox{Ad}_{u_3(I_4\otimes W^{(2)}U_1)}(\bbc\otimes(M_{10}\oplus M_5\oplus\bbc))\subset\Delta_4\otimes M_4\otimes M_4=A_3$ is the $3\times 4$ matrix each of whose entry is $1$, and the square of its norm is $||\Lambda_1||^2=12$. Thus, the following quadruple
\[
\begin{matrix}
A_1:=\bbc\otimes\Delta_4\otimes M_4 &\subset & \Delta_4\otimes M_4\otimes M_4=:A_3\\
\cup &  & \cup\\
C_1:=\bbc\otimes\mbox{Ad}_{u_1(I_4\otimes W)}(\bbc\otimes(M_3\oplus\bbc)) &\subset & \mbox{Ad}_{u_3\big(I_4\otimes W^{(2)}U_1\big)}(\bbc\otimes(M_{10}\oplus M_5\oplus\bbc))=:C_3\\
\end{matrix}
\]
is a commuting square which is not non-degenerate if $u\nsim v$. However, if we assume that $u\sim v$ and $\theta=\overline{a}b$ is such that $\theta^{\,4}=1$, then the situation changes drastically. By \Cref{101}, we see that $\sigma\big((u_2^*v_2)^2\big)=\{1,\theta^{\,2}\}$ as $\theta^{\,4}=1$ (putting $k=2$ there). By \Cref{recurrence1}, the multiplicity of the eigenvalue $1$ remains the same, i,e. $10$, but that of $\theta^{\,2}$ is now $6$ (as $\theta^{\,2}=\overline{\theta}^{\,2}$) instead of $5$. In view of \Cref{intersection infinite}, we now have the following modified quadruple,
\[
\begin{matrix}
A_1:=\bbc\otimes\Delta_4\otimes M_4 &\subset & \Delta_4\otimes M_4\otimes M_4=:A_3\\
\cup &  & \cup\\
C_1:=\bbc\otimes\mbox{Ad}_{u_1(I_4\otimes W)}(\bbc\otimes(M_3\oplus\bbc)) &\subset & \mbox{Ad}_{u_3\big(I_4\otimes W^{(2)}U_1\big)}(\bbc\otimes(M_{10}\oplus M_6))=:C_3\\
\end{matrix}
\]
The inclusion matrix $\tilde{\Lambda}_1$ for the inclusion $\mbox{Ad}_{u_3(I_4\otimes W^{(2)}U_1)}(\bbc\otimes(M_{10}\oplus M_6))\subset\Delta_4\otimes M_4\otimes M_4$ is the $2\times 4$ matrix each of whose entry is $1$. That is, $\tilde{\Lambda}_1=\Lambda_0$ and $||\tilde{\Lambda}_1||^2=||\Lambda_0||^2=8$. Therefore, in this case the quadruple $(C_1\subset A_1,C_3\subset A_3)$, which is already a commuting square, becomes non-degenerate also by Lemma $3.10$ in \cite{BG}. Iterating Jones' basic construction on the quadruple $(C_1\subset A_1,C_3\subset A_3)$, we obtain a subfactor $R^{(4)}$ of $R$ such that $R^{(4)}\subseteq R_u\cap R_v$ (by \Cref{bakproc}) and $[R:R^{(4)}]=||\Lambda_0||^2=8$. Since $[R_u:R^{(4)}]=2$, by the irreducibility of $R^{(4)}$ in $R_u$, it follows that $R^{(4)}=R_u\cap R_v$. That is, $R_u\cap R_v$ is a finite index subfactor of $R$ with $[R:R_u\cap R_v]=8$. This is in sharp contrast with \Cref{infinite index}.

Similarly, if we assume that $u\sim v$ and $\theta=\overline{a}b$ is such that $\theta^{\,6}=1$, we see that $\sigma\big((u_2^*v_2)^2\big)=\sigma\big((u_4^*v_4)^2\big)=\{1,\theta^{\,2},\overline{\theta}^{\,2}\}$ by \Cref{101}. In this case, the quadruple $(C_1\subset A_1,C_3\subset A_3)$ remains not non-degenerate as before, however, a similar analysis as above tells us that $\tilde{\Lambda}_2=\Lambda_1$ with the square of their norm being equal to $12$. In other words, the quadruple $(C_3\subset A_3,C_5\subset A_5)$ becomes a non-degenerate commuting square. Therefore, iterating Jones' basic construction on the quadruple $(C_3\subset A_3,C_5\subset A_5)$, we obtain a subfactor $R^{(6)}$ of $R$ such that $R^{(6)}\subseteq R_u\cap R_v$ (by \Cref{bakproc}) and $[R:R^{(6)}]=||\Lambda_1||^2=12$. Since $[R_u:R^{(6)}]=3$, by the irreducibility of $R^{(6)}$ in $R_u$, it follows that $R^{(6)}=R_u\cap R_v$. That is, $R_u\cap R_v$ is a finite index subfactor of $R$ with $[R:R_u\cap R_v]=12$, which is again in sharp contrast with \Cref{infinite index}.
\smallskip

It turns out that this phenomena happens in general, that is, for $\theta^{\,m}=1$ where $m\in 2\bbn+2$. We need the following technical lemma first. Recall \Cref{101} and \Cref{recurrence1} in this regard, along with the notations.

\begin{lmma}\label{mpr}
For $u=u(a),v=u(b)$ let $\theta=\overline{a}b$, and suppose that $\theta$ is a primitive even root of unity. That is, $\,\theta^{\,m}=1$ for some $m\in 2\bbn+2$. Then, for $k=(m-2)/2$ we have the following,
\begin{IEEEeqnarray*}{lCl}
&  & C_{2k+1}:=\mbox{Ad}_{u_{2k+1}}\big(\bbc\otimes M_4^{(k+1)}\big)\cap\mbox{Ad}_{v_{2k-1}}\big(\bbc\otimes M_4^{(k+1)}\big)\\
&=& \begin{cases}
\mathrm{Ad}_{u_{2k+1}\big(I_4\otimes W^{(k+1)}U_k\big)}\Big(\bbc\otimes\Big(M_{\widetilde{m}_k(1)}\bigoplus\displaystyle{\bigoplus_{r=1}^{\frac{k-1}{2}}}\Big(M_{\widetilde{m}_k(\theta^{\,2r})}\oplus M_{\widetilde{m}_k\big(\overline{\theta}^{\,2r}\big)}\Big)\\
\hspace{8.5cm}\bigoplus M_{\widetilde{m}_k(\theta^{\,k+1})}\Big)\Big)\qquad\,\mbox{if $k$ is odd},\\
\mathrm{Ad}_{u_{2k+1}\big(I_4\otimes W^{(k+1)}U_k\big)}\Big(\bbc\otimes\Big(M_{\widetilde{m}_{k}(1)}\bigoplus\displaystyle{\bigoplus_{r=1}^{\frac{k}{2}}}\Big(M_{\widetilde{m}_k(\theta^{\,2r})}\oplus M_{\widetilde{m}_k\big(\overline{\theta}^{\,2r}\big)}\Big)\Big)\Big)\quad\mbox{if $k$ is even}.
\end{cases}
\end{IEEEeqnarray*}
Here, $U_k$ remains the same unitary obtained in \Cref{aef}.
\end{lmma}
\begin{prf}
Consider the unitary $(u_{2k}^*v_{2k})^2$ in $M_4^{(k+1)}$. By \Cref{101}, the set of distinct eigenvalues are either $\{1,\theta^{\,2},\overline{\theta}^{\,2},\ldots,\theta^{\,k-1},\overline{\theta}^{\,k-1},\theta^{\,k+1}\}$ or $\{1,\theta^{\,2},\overline{\theta}^{\,2},\ldots,\theta^{\,k},\overline{\theta}^{\,k}\}$ depending on $k$ being odd or even respectively, and $\sigma\big((u_{2k}^*v_{2k})^2\big)=\sigma\big((u_{2k-2}^*v_{2k-2})^2\big)$ with cardinality $k+1$. The multiplicities $\widetilde{m}_k(.)$ of the eigenvalues are given by \Cref{recurrence1}. The rest of the proof follows along the similar line of \Cref{aef}.\qed
\end{prf}

The above lemma can be thought of as a counterpart to \Cref{intersection infinite} for the case of $u\sim v$. We have crucially used here the enumeration of eigenvalues according to increasing power of $\theta$ followed immediately by $\overline{\theta}$ (see also \Cref{pg} and \Cref{guin}) to keep the unitary $U_k$ same with that in \Cref{aef}. Now we are ready to prove our main theorem.

\begin{thm}\label{producing}
Let $u=u(a)\mbox{ and }v=u(b)$ be $4\times 4$ Hadamard inequivalent complex Hadamard matrices such that $u\sim v$. If $\theta=\overline{a}b$ satisfies $\theta^{\,m}=1$ for $m\in 2\bbn+2$, then there is a finite index subfactor $R^{(m)}_{u,v}$ of the hyperfinite type $II_1$ factor $R$ such that $[R:R^{(m)}_{u,v}]=2m$ and $R^{(m)}_{u,v}\subseteq R_u\cap R_v$.
\end{thm}
\begin{prf}
The cases of $m=4,6$ are already discussed in detail at the beginning of this subsection. For $\theta^{\,m}=1$ with $m\geq 8$, choose $k_m\in\bbn$ such that $m=2k_m+2$. Such a positive integer $k_m$ always exists and $k_m\geq 3$. Consider the following inclusion
\[
C_{2k_m-1}:=\mbox{Ad}_{u_{2k_m-1}}\big(\bbc\otimes M_4^{(k_m)}\big)\cap\mbox{Ad}_{v_{2k_m-1}}\big(\bbc\otimes M_4^{(k_m)}\big)\subset\Delta_4\otimes M_4^{(k_m)}=:A_{2k_m-1},
\]
which is same as $\mbox{Ad}_{u_{2k_m-1}}\big(\bbc\otimes\big((u_{2k_m-2}^*v_{2k_m-2})^2\big)^\prime\cap M_4^{(k_m)}\big)\subset\Delta_4\otimes M_4^{(k_m)}$ due to \Cref{odd intersection}. The set of eigenvalues of $(u_{2k_m-2}^*v_{2k_m-2})^2$ is given by \Cref{101}, and their multiplicities are given by \Cref{recurrence1}. Observe that up to this stage, the eigenvalues, along with their multiplicities, remain the same as in \Cref{oo} and \Cref{recurrence} for the case of $u\nsim v$. This is because we always take $\theta$ to be a primitive even root of unity. Therefore, \Cref{intersection infinite} remains valid for $C_{2k_m-1}$, and the inclusion matrix $\Lambda_{k_m-1}$ for the inclusion $C_{2k_m-1}\subset A_{2k_m-1}$ is a $(k_m+1)\times 4$ matrix each of whose entry is $1$, and the square of the norm of it is $||\Lambda_{k_m-1}||^2=4(k_m+1)$ (see \Cref{inf}). Now, consider the next inclusion
\[
C_{2k_m+1}:=\mbox{Ad}_{u_{2k_m+1}}\big(\bbc\otimes M_4^{(k_m+1)}\big)\cap\mbox{Ad}_{v_{2k_m+1}}\big(\bbc\otimes M_4^{(k_m+1)}\big)\subset\Delta_4\otimes M_4^{(k_m+1)}=:A_{2k_m+1},
\]
which is same as $\mbox{Ad}_{u_{2k_m+1}}\big(\bbc\otimes\big((u_{2k_m}^*v_{2k_m})^2\big)^\prime\cap M_4^{(k_m+1)}\big)\subset\Delta_4\otimes M_4^{(k_m+1)}$ due to \Cref{odd intersection}. That is, we are considering the following quadruple
\[
\begin{matrix}
A_{2k_m-1} &\subset &A_{2k_m+1}\\
|| &   &  ||\\
\bbc\otimes\Delta_4\otimes M_4^{(k_m)} & & \Delta_4\otimes M_4^{(k_m+1)}\\
\cup &  & \cup\\
\bbc\otimes\mbox{Ad}_{u_{2k_m-1}}\big(\bbc\otimes\big((u_{2k_m-2}^*v_{2k_m-2})^2\big)^\prime\cap M_4^{(k_m)}\big) & & \mbox{Ad}_{u_{2k_m+1}}\big(\bbc\otimes\big((u_{2k_m}^*v_{2k_m})^2\big)^\prime\cap M_4^{(k_m+1)}\big)\\
|| &   &  ||\\
C_{2k_m-1} &\subset &C_{2k_m+1}\\
\end{matrix}
\]
Recall that if $u$ were not related to $v$, then this quadruple is not non-degenerate as explained in \Cref{justification of main difficulty}. This is because in the absence of any such condition on $\theta$, the inclusion matrix for $C_{2k_m+1}\subset A_{2k_m+1}$ is a $(k_m+2)\times 4$ matrix each of whose entry is $1$ by \Cref{inf}. However, we show that putting the condition $\theta^{\,2k_m+2}=1$ forces a row-drop in the inclusion matrix. More precisely, we claim that the inclusion matrix $\Lambda_{k_m}$ for the right vertical inclusion $C_{2k_m+1}\subset A_{2k_m+1}$ is also a $(k_m+1)\times 4$ matrix each of whose entry is $1$, and hence $||\Lambda_{k_m-1}||^2=||\Lambda_{k_m}||^2=4(k_m+1)$. Consider the unitary $(u_{2k_m}^*v_{2k_m})^2$ in $M_4^{(k_m+1)}$. By \Cref{mpr}, we see that the number of direct summands in the Wedderburn–Artin decomposition of the algebra $\big((u_{2k_m}^*v_{2k_m})^2\big)^\prime\cap M_4^{(k_m+1)}$ is exactly $k_m+1$. Since $A_{2k_m+1}=\Delta_4\otimes M_4^{(k_m+1)}=\oplus M_4^{(k_m+1)}$, and the unitary $u_{2k+1}(I_4\otimes W^{(k+1)}U_k)$ in \Cref{mpr} belongs to $A_{2k_m+1}$, the inclusion matrix $\Lambda_{k_m}$ for the inclusion $C_{2k_m+1}\subset A_{2k_m+1}$ becomes a $(k_m+1)\times 4$ matrix each of whose entry is $1$. This finishes the claim.

Therefore, the following quadruple
\begin{IEEEeqnarray*}{lCl}
\begin{matrix}
A_{2k_m-1} & \subset & A_{2k_m+1}\\
\cup & & \cup\\
C_{2k_m-1} & \subset & C_{2k_m+1}\\
\end{matrix}
\end{IEEEeqnarray*}
which is already a commuting square, becomes non-degenerate also by Lemma $3.10$ in \cite{BG}. Iterating Jones' basic construction starting with this non-degenerate commuting square, we obtain a subfactor $R^{(m)}_{u,v}$ of the hyperfinite type $II_1$ factor $R$ such that $R^{(m)}_{u,v}\subseteq R_u\cap R_v$ (by \Cref{bakproc}) and $[R:R^{(m)}_{u,v}]=||\Lambda_{k_m-1}||^2=4(k_m+1)=2m$. That the tower $A_1\subset A_3\subset A_5\subset\cdots$ indeed gives us $R$ follows as an application of Theorem $3.5$ in \cite{B1}, for instance.\qed
\end{prf}

Therefore, we obtain the following family of finite index subfactors
\[
\left\{R^{(m)}_{u(a),u(b)}\subset R:\,u(a)\sim u(b)\mbox{ with }(\overline{a}b)^m=1,\,m\in 2\bbn+2\right\}
\]
such that $\big[R:R^{(m)}_{u(a),u(b)}\big]=2m$, where $u(a),u(b)$ are complex Hadamard matrices of order $4\times 4$ parametrized by the circle parameters $a\mbox{ and }b$ respectively such that $u(a)\sim u(b)$ with $(\overline{a}b)^m=1$ for $m\in 2\bbn+2$. The converse of the above result is also true.

\begin{crlre}
Given any $m\in 2\bbn+2$, there exists a pair $(u,v)$ of Hadamard inequivalent $4\times 4$ complex Hadamard matrices and a subfactor $R_0\subset R$ such that $R_0\subset R_u\cap R_v$ and $[R:R_0]=2m$. 
\end{crlre}
\begin{prf}
Given any integer $m\geq 4$, choose a complex Hadamard matrix $u(a)$ of the required form, where $a$ is any arbitrary complex number on the semicircle. Consider $b=\omega a$, where $\omega\neq\pm 1$ is any primitive $m$-th root of unity. If $b$ lies in the lower semicircle, consider $-b$ instead of $b$, which lies in the upper semicircle. Obtain the complex Hadamard matrix $u(b)$. For the pair $(u,v):=(u(a),u(b))$, we have $(\overline{a}b)^m=1$. By \Cref{producing}, we have the required subfactor $R_0:=R^{(m)}_{u,v}\subset R$.\qed
\end{prf}

We shall see in the next subsection that all these subfactors $R^{(m)}_{u(a),u(b)}\subset R$ are irreducible. We end this subsection with the following question.
\medskip

\noindent \textbf{Question:} Can we characterize the following family of finite index subfactors of $R$ \[\left\{R^{(m)}_{u(a),u(b)}\subset R:\,u(a)\sim u(b)\mbox{ with }(\overline{a}b)^m=1,\,m\in 2\bbn+2\right\}?\]


\subsection{Relative commutant and factoriality of $R_u\cap R_v$ when $u\sim v$}

This subsection is intimately connected with its immediate prequel. Goal of this subsection is to prove that all the subfactors $R^{(m)}_{u(a),u(b)}\subset R$ constructed there are irreducible in $R$. As a consequence, we will get that $R_{u(a)}\cap R_{u(b)}$ is factor when $u(a)\sim u(b)$. Major tools are the Ocneanu compactness and the spectral decomspoition obtained in subsection $7.2.1$. We begin by recalling the following well-known fact from Linear algebra.

\begin{lmma}\label{sylvester}
Let $B$ be a block-diagonal matrix, that is $B\equiv B_1\oplus B_2$, and $A$ be a matrix such that $[A,B]=0$. Suppose that $B_1\mbox{ and }B_2$ do not share any eigenvalue over $\bbc$. Then, $A$ has the same block-diagonal structure as $B$, that is, it is block-diagonal with blocks placed at the same positions as those of $B$.
\end{lmma}
\begin{prf}
Let $A=\begin{pmatrix}
X & Y\\
Z & W
\end{pmatrix}$. Since $[A,B]=0$, we obtain two equations $B_1Y-YB_2=0$ and and $B_2Z-ZB_1=0$. Since $B_1\mbox{ and }B_2$ do not share any eigenvalue, by the Sylvester criterion, we get that $Y=0$ and $Z=0$.\qed
\end{prf}

\begin{lmma}\label{z}
The commutant of any finite-dimensional algebra $N=\oplus_{j=1}^r N_{n_j}$, where $N_{n_j}$ is the type $I_{n_j}$ factor, in the type $I_m$ factor $M_m$ such that $\sum_jn_j=m$ is $\bbc^{\,r}$.
\end{lmma}
\begin{prf}
By the previous \Cref{sylvester}, if $A\in M_m$ commutes with every $B=\oplus_{j=1}^rB_j$ where $B_j\in N_{n_j}$, then $A$ must be of the form $\oplus_{j=1}^rA_j$ where $A_j\in N_{n_j}$. This is because it is always possible to choose such a block-diagonal matrix $B$ satisfying the hypothesis of \Cref{sylvester}. Since each $N_{n_j}$ is the type $I_{n_j}$ factor, each $A_j$ must be a scalar matrix of order $n_j$.\qed
\end{prf}

Now we prove our main theorem. The crucial ingredients are \Cref{aef} and \ref{intersection infinite}, and more precisely \Cref{mpr}. We request the readers to recall the corresponding notations also. Recall that the unitary $U_k$ in \Cref{mpr} is defined inductively by $U_k=(I_4\otimes U_{k-1})V_{k-1}$, where $V_{k-1}$ is certain permutation matrix in $M_4^{(k+1)}$ that only interchange the diagonal entries (see the proof of \Cref{aef}). This permutation matrix has crucial role in what follows next.

\begin{thm}\label{interesting}
Let $u=u(a)\mbox{ and }v=u(b)$ be $4\times 4$ Hadamard inequivalent complex Hadamard matrices such that $u\sim v$ with $(\overline{a}b)^m=1$. Then, the subfactor $R^{(m)}_{u,v}\subset R$ is irreducible.
\end{thm}
\begin{prf}
Since $u\sim v$, we have $\theta=\overline{a}b\in\Gamma$. If $\theta^{\,m}=1$ ($\theta$ is a primitive even root of unity), take $k=(m-2)/2$. The proof is a bit technical, and hence we break it into several steps for clarity. We separate two situations, namely $k$ being odd or even. First suppose that $k$ is odd. By \Cref{recurrence1}, we have the multiplicities $\widetilde{m}_k(.)$ of the eigenvalues of $(u_{2k}^*v_{2k})^2$ (and similarly $\widetilde{m}_{k-1}(.)$ for $(u_{2k-2}^*v_{2k-2})^2$). To avoid the cumbersome notation, we denote $\widetilde{m}_k(.)$ simply by $m_k(.)$ specifically for this proof. We hope that this does not create any confusion because by \Cref{recurrence1}, we see that $\widetilde{m}_{k-1}(.)$ is anyway equal to $m_{k-1}(.)$, and so is $\widetilde{m}_k(.)$ equal to $m_k(.)$ except for the last eigenvalue $\theta^{\,k+1}$ (that is, only $\widetilde{m}_k(\theta^{\,k+1})\neq m_k(\theta^{\,k+1})$).
\smallskip

\noindent\textbf{Step 1:~ The action of the permutation matrix} $V_{k-1}$ (\Cref{mpr}, \Cref{aef}).
\smallskip

\noindent Consider the unitary matrix $(u_{2k}^*v_{2k})^2$. Recall from \Cref{cute lemma for 4 by 4} that
\[
(u_{2k}^*v_{2k})^2=p\otimes (u_{2k-2}^*v_{2k-2})^2+q(\theta^{\,2})\otimes(v_{2k-2}^*u_{2k-2})^2
\]
for $k\in\bbn$. Recall from \Cref{ppq} the following,
\begin{IEEEeqnarray*}{lCl}
&  & \mathrm{Ad}_{W^{(k)}}(u_{2k}^*v_{2k})^2\\
&=& E_{11}\otimes I_2\otimes\mathrm{Ad}_{W^{(k-1)}}(u_{2k-2}^*v_{2k-2})^2+E_{22}\otimes C(\theta^{\,2})\otimes\mathrm{Ad}_{W^{(k-1)}}(v_{2k-2}^*u_{2k-2})^2\\
&=& E_{11}\otimes I_2\otimes\mathrm{Ad}_{U_{k-1}}\Big(I_{m_{k-1}(1)}\bigoplus\,\bigoplus_{r=1}^{\frac{k-1}{2}}\Big(\theta^{\,2r}\,I_{m_{k-1}(\theta^{\,2r})}\oplus\overline{\theta}^{\,2r}\,I_{m_{k-1}\big(\overline{\theta}^{\,2r}\big)}\Big)\bigoplus\theta^{\,k+1}\,I_{m_{k-1}(\theta^{\,k+1})}\Big)\\
&  & +E_{22}\otimes C(\theta^{\,2})\otimes\mathrm{Ad}_{U_{k-1}}\Big(I_{m_{k-1}(1)}\bigoplus\,\bigoplus_{r=1}^{\frac{k-1}{2}}\Big(\overline{\theta}^{\,2r}\,I_{m_{k-1}(\theta^{\,2r})}\oplus\theta^{\,2r}\,I_{m_{k-1}\big(\overline{\theta}^{\,2r}\big)}\Big)\bigoplus\overline{\theta}^{\,k+1}\,I_{m_{k-1}(\theta^{\,k+1})}\Big)\\
&=& \mathrm{Ad}_{(I_4\otimes U_{k-1})V_{k-1}}\Big(I_{m_k(1)}\bigoplus \bigoplus_{r=1}^{\frac{k-1}{2}}\Big(\theta^{\,2r}\,I_{m_k(\theta^{\,2r})}\oplus\overline{\theta}^{\,2r}\,I_{m_k(\overline{\theta}^{\,2r})}\Big)\bigoplus \theta^{\,k+1}\,I_{m_k(\theta^{\,k+1})}\Big)\,.
\end{IEEEeqnarray*}
Note that here $\theta^{\,k+1}=\overline{\theta}^{\,k+1}$ as $\theta^{\,m}=1$. Since $C(\theta^{\,2})$ is the diagonal matrix $\mbox{diag}\{1,\theta^{\,2}\}$, we obtain that under the conjugation action by the permutation matrix $V_{k-1}$, the element
\[
\alpha I_{m_k(1)}\bigoplus\bigoplus_{r=1}^{\frac{k-1}{2}}\Big(\beta_{2r}I_{m_k(\theta^{\,2r})}\oplus\gamma_{2r}I_{m_k(\overline{\theta}^{\,2r})}\Big)\bigoplus\delta I_{m_k(\theta^{\,k+1})}
\]
where $\alpha,\beta_{2r},\gamma_{2r},\delta\in\bbc$ for $r=1,\ldots,\frac{k-1}{2}$, becomes the following
\begin{IEEEeqnarray}{lCl}\label{easlr}
&  & (E_{11}+E_{22})\otimes\Big\{\alpha I_{m_{k-1}(1)}\bigoplus\bigoplus_{r=1}^{\frac{k-1}{2}}\Big(\beta_{2r}I_{m_{k-1}(\theta^{\,2r})}\oplus\gamma_{2r}I_{m_{k-1}(\overline{\theta}^{\,2r})}\Big)\bigoplus\delta I_{m_{k-1}(\theta^{\,k+1})}\Big\}\nonumber\\
&  &  +E_{33}\otimes\Big\{\alpha I_{m_{k-1}(1)}\bigoplus\bigoplus_{r=1}^{\frac{k-1}{2}}\Big(\gamma_{2r}I_{m_{k-1}(\theta^{\,2r})}\oplus\beta_{2r}I_{m_{k-1}(\overline{\theta}^{\,2r})}\Big)\bigoplus\delta I_{m_{k-1}(\theta^{\,k+1})}\Big\}\nonumber\\
&  & +E_{44}\otimes\Big\{\big(\beta_2 I_{m_{k-1}(1)}\oplus\alpha I_{m_{k-1}(\theta^{\,2})}\big)\bigoplus\bigoplus_{r=1}^{\frac{k-3}{2}}\Big(\beta_{2r+2}I_{m_{k-1}(\overline{\theta}^{\,2r})}\oplus\gamma_{2r}I_{m_{k-1}(\theta^{\,2r+2})}\Big)\nonumber\\
&  & \qquad\qquad\bigoplus\big(\delta I_{m_{k-1}(\overline{\theta}^{\,k-1})}\oplus\gamma_{k-1}I_{m_{k-1}(\theta^{\,k+1})}\big)\Big\}\,.
\end{IEEEeqnarray}
This finishes our first step.
\smallskip

\noindent\textbf{Step 2:~ The structure of} $\,C_{2k+1}^{\,\prime}\cap A_{2k-1}$.
\smallskip

\noindent We find the general structure of this commutant. With the convention $W_0=I_4$, we have the following,
\begin{IEEEeqnarray}{lCl}\label{pi}
&  & C_{2k+1}^{\,\prime}\cap A_{2k-1}\\
&=& \Big\{\mbox{Ad}_{u_{2k+1}}\Big(\bbc\otimes \big((u_{2k}^*v_{2k})^2\big)^\prime\cap M_4^{(k+1)}\Big)\Big\}^\prime\cap\Big(\bbc\otimes\Delta_4\otimes M_4^{(k)}\Big)\nonumber\\
&=& \Big\{\mbox{Ad}_{u_{2k+1}}\Big(\bbc\otimes \big((u_{2k}^*v_{2k})^2\big)^\prime\cap M_4^{(k+1)}\Big)\Big\}^\prime\cap\mbox{Ad}_{I_4\otimes u_{2k-1}}\Big(\bbc\otimes\Delta_4\otimes M_4^{(k)}\Big)\nonumber\\
&=& \mbox{Ad}_{u_{2k+1}}\Big\{\Big(\bbc\otimes \big((u_{2k}^*v_{2k})^2\big)^\prime\cap M_4^{(k+1)}\Big)^\prime\cap\,\mbox{Ad}_{u_{2k+1}^*(I_4\otimes u_{2k-1})}\Big(\bbc\otimes\Delta_4\otimes M_4^{(k)}\Big)\Big\}\nonumber\\
&=& \mbox{Ad}_{u_{2k+1}}\Big\{\Big(\bbc\otimes \big((u_{2k}^*v_{2k})^2\big)^\prime\cap M_4^{(k+1)}\Big)^\prime\cap\,\mbox{Ad}_{\big(\sum_{j=0}^3E_{j+1,j+1}\otimes W_j^*\xi^*\otimes I_4^{(k)}\big)}\big(\bbc\otimes\Delta_4\otimes M_4^{(k)}\big)\Big\}\nonumber\\
&=& \mbox{Ad}_{u_{2k+1}}\Big\{\Big(\bbc\otimes\mathrm{Ad}_{W^{(k+1)}}\Big(\mathrm{Ad}_{U_k}\Big(M_{m_k(1)}\oplus M_{m_k(\theta^{\,2})}\oplus M_{m_k(\overline{\theta}^{\,2})}\oplus\cdots\oplus M_{m_k(\theta^{\,k+1})}\Big)\Big)\Big)^\prime\nonumber\\
& & \hspace*{1.7cm}\bigcap\,\mbox{Ad}_{\big(\sum_{j=0}^3E_{j+1,j+1}\otimes W_j^*\xi^*\otimes I_4^{(k)}\big)}\big(\bbc\otimes\Delta_4\otimes M_4^{(k)}\big)\Big\}\quad\qquad(\mbox{by }\Cref{mpr})\nonumber\\
&=& \mbox{Ad}_{u_{2k+1}(I_4\otimes W^{(k+1)})}\Big\{\Big(\bbc\otimes\mathrm{Ad}_{U_k}\Big(M_{m_k(1)}\oplus M_{m_k(\theta^{\,2})}\oplus M_{m_k(\overline{\theta}^{\,2})}\oplus\cdots\oplus M_{m_k(\theta^{\,k+1})}\Big)\Big)^\prime\nonumber\\
& & \hspace*{3.5cm}\bigcap\,\mbox{Ad}_{\big(\sum_{j=0}^3E_{j+1,j+1}\otimes WW_j^*\xi^*\otimes I_4^{(k)}\big)}\big(\bbc\otimes\Delta_4\otimes M_4^{(k)}\big)\Big\}\nonumber\\
&=& \mbox{Ad}_{u_{2k+1}\big(I_4\otimes W\otimes W^{(k)}U_{k-1}\big)}\Big\{\Big(\bbc\otimes\mathrm{Ad}_{V_{k-1}}\Big(M_{m_k(1)}\oplus M_{m_k(\theta^{\,2})}\oplus M_{m_k(\overline{\theta}^{\,2})}\oplus\cdots\oplus M_{m_k(\theta^{\,k+1})}\Big)\Big)^\prime\nonumber\\
& & \hspace*{4.5cm}\bigcap\,\mbox{Ad}_{\big(\sum_{j=0}^3E_{j+1,j+1}\otimes WW_j^*\xi^*\otimes I_4^{(k)}\big)}\big(\bbc\otimes\Delta_4\otimes M_4^{(k)}\big)\Big\}\,,\nonumber
\end{IEEEeqnarray}
where the last line follows from the fact that $U_s=(I_4\otimes U_{s-1})V_{s-1}$ with $U_{s-1}\in M_4^{(s)}$. Now, we have
\begin{IEEEeqnarray}{lCl}\label{mz}
& & \mbox{Ad}_{\big(\sum_{j=0}^3E_{j+1,j+1}\otimes WW_j^*\xi^*\otimes I_4^{(k)}\big)}\big(\bbc\otimes\Delta_4\otimes M_4^{(k)}\big)\\
&=& \mbox{bl-diag}\Big\{\mbox{Ad}_{W\xi^*\otimes I_4^{(k)}}(x)\,,\,\mbox{Ad}_{WW_1^*\xi^*\otimes I_4^{(k)}}(x)\,,\,\mbox{Ad}_{WW_2\xi^*\otimes I_4^{(k)}}(x)\,,\,\mbox{Ad}_{WW_3^*\xi^*\otimes I_4^{(k)}}(x)\Big\}\nonumber
\end{IEEEeqnarray}
where $x\in\Delta_4\otimes M_4^{(k)}$. We write $x=\mbox{bl-diag}\{x_1,x_2,x_3,x_4\}$ with $x_j\in M_4^{(k)}$ for $j=1,2,3,4$. Then, using \Cref{tower in four by four} it is easy to verify the following four identities,
\begin{IEEEeqnarray}{lCl}\label{e}
\mbox{Ad}_{W\xi^*\otimes I_4^{(k)}}(x) &=& \frac{1}{2}\Big((E_{11}+E_{22})\otimes(x_1+x_3)+(E_{33}+E_{44})\otimes(x_2+x_4)\nonumber\\
& & \quad+(E_{12}+E_{21})\otimes(x_1-x_3)+(iaE_{34}-i\overline{a}E_{43})\otimes(x_2-x_4)\Big)\nonumber\\
\mbox{Ad}_{WW_1^*\xi^*\otimes I_4^{(k)}}(x) &=& \frac{1}{2}\big((E_{11}+E_{22})\otimes(x_2+x_4)+(E_{33}+E_{44})\otimes(x_1+x_3)\nonumber\\
& & \quad+(E_{12}+E_{21})\otimes(x_2-x_4)+(-i\overline{a}E_{34}+iaE_{43})\otimes(x_1-x_3)\big)\nonumber\\
\mbox{Ad}_{WW_2\xi^*\otimes I_4^{(k)}}(x) &=& \frac{1}{2}\big((E_{11}+E_{22})\otimes(x_1+x_3)+(E_{33}+E_{44})\otimes(x_2+x_4)\nonumber\\
& & \quad+(E_{12}+E_{21})\otimes(x_3-x_1)+(-iaE_{34}+i\overline{a}E_{43})\otimes(x_2-x_4)\big)\nonumber\\
\mbox{Ad}_{WW_3^*\xi^*\otimes I_4^{(k)}}(x) &=& \frac{1}{2}\big((E_{11}+E_{22})\otimes(x_2+x_4)+(E_{33}+E_{44})\otimes(x_1+x_3)\nonumber\\
& & \quad+(E_{12}+E_{21})\otimes(x_4-x_2)+(i\overline{a}E_{34}-iaE_{43})\otimes(x_1-x_3)\big)\,.
\end{IEEEeqnarray}
\Cref{e} describes arbitrary elements of the set in \Cref{mz}.
\smallskip

\noindent\textbf{Step 3:~ Proof of} $\,C_{2k+1}^{\,\prime}\cap A_{2k-1}=\bbc$.
\smallskip

\noindent This is the major step in the proof. In view of \Cref{pi}, it is enough to prove the following,
\begin{IEEEeqnarray}{lCl}\label{f}
\bbc &=& \Big(\bbc\otimes\mathrm{Ad}_{V_{k-1}}\Big(M_{m_k(1)}\oplus M_{m_k(\theta^{\,2})}\oplus M_{m_k(\overline{\theta}^{\,2})}\oplus\cdots\oplus M_{m_k(\theta^{\,k+1})}\Big)\Big)^\prime\nonumber\\
& & \,\bigcap\,\,\mbox{Ad}_{\big(\sum_{j=0}^3E_{j+1,j+1}\otimes WW_j^*\xi^*\otimes I_4^{(k)}\big)}\big(\bbc\otimes\Delta_4\otimes M_4^{(k)}\big)\,.
\end{IEEEeqnarray}
First observe that since the permutation matrix $V_{k-1}$ only interchanges the diagonal elements, we get the following,
\[
\Big(\bbc\otimes\mathrm{Ad}_{V_{k-1}}\Big(M_{m_k(1)}\oplus M_{m_k(\theta^{\,2})}\oplus M_{m_k(\overline{\theta}^{\,2})}\oplus\cdots\oplus M_{m_k(\theta^{\,k+1})}\Big)\Big)^\prime\cap A_{2k+1}\subseteq\Delta_4\otimes\Delta_4^{(k+1)}\,.
\]
To see this, apply \Cref{com and ad interchange} and observe the following,
\begin{IEEEeqnarray*}{lCl}
& & \Big(\mathrm{Ad}_{V_{k-1}}\Big(M_{m_k(1)}\oplus M_{m_k(\theta^{\,2})}\oplus M_{m_k(\overline{\theta}^{\,2})}\oplus\cdots\oplus M_{m_k(\theta^{\,k+1})}\Big)\Big)^\prime\cap M_4\otimes M_4^{(k)}\\
&=& \mathrm{Ad}_{V_{k-1}}\Big(\Big(M_{m_k(1)}\oplus M_{m_k(\theta^{\,2})}\oplus M_{m_k(\overline{\theta}^{\,2})}\oplus\cdots\oplus M_{m_k(\theta^{\,k+1})}\Big)^\prime\cap M_4^{(k+1)}\Big)\\
&=& \mathrm{Ad}_{V_{k-1}}\Big(\bbc I_{m_k(1)}\oplus\bbc I_{m_k(\theta^{\,2})}\oplus\bbc I_{m_k(\overline{\theta}^{\,2})}\oplus\ldots\oplus\bbc I_{m_k(\theta^{\,k+1})}\Big)\,,
\end{IEEEeqnarray*}
where the last line follows from \Cref{z}. Therefore, the intersection in the right hand side of \Cref{f} must consists of purely diagonal matrices because
\[
\mbox{Ad}_{\big(\sum_{j=0}^3E_{j+1,j+1}\otimes WW_j^*\xi^*\otimes I_4^{(k)}\big)}\big(\bbc\otimes\Delta_4\otimes M_4^{(k)}\big)\subseteq A_{2k+1}:=\Delta_4\otimes M_4\otimes M_4^{(k)}\,.
\]
Since
\begin{IEEEeqnarray*}{lCl}
& & \mbox{Ad}_{\big(\sum_{j=0}^3E_{j+1,j+1}\otimes WW_j^*\xi^*\otimes I_4^{(k)}\big)}\big(\bbc\otimes\Delta_4\otimes M_4^{(k)}\big)\\
&=& \mbox{bl-diag}\Big\{\mbox{Ad}_{W\xi^*\otimes I_4^{(k)}}(x)\,,\,\mbox{Ad}_{WW_1^*\xi^*\otimes I_4^{(k)}}(x)\,,\,\mbox{Ad}_{WW_2\xi^*\otimes I_4^{(k)}}(x)\,,\,\mbox{Ad}_{WW_3^*\xi^*\otimes I_4^{(k)}}(x)\Big\}
\end{IEEEeqnarray*}
for $x=\mbox{bl-diag}\{x_1,x_2,x_3,x_4\}\in\Delta_4\otimes M_4^{(k)}$, we see that any element of the intersection in \Cref{f} must satisfy $x_1=x_3$ and $x_2=x_4$, with $x_1,x_2\in\Delta_4^{(k)}$ because of the identities established in \Cref{e} in Step $2$. Thus, to prove that the \Cref{f} holds, it is enough to prove that the following,
\begin{IEEEeqnarray}{lCl}\label{caslr}
& & \Big\{\Delta_4\otimes\mathrm{Ad}_{V_{k-1}}\Big(\alpha I_{m_k(1)}\bigoplus\bigoplus_{r=1}^{\frac{k-1}{2}}\Big(\beta_{2r} I_{m_k(\theta^{\,2})}\oplus\gamma_{2r} I_{m_k(\overline{\theta}^{\,2})}\Big)\bigoplus\delta I_{m_k(\theta^{\,k+1})}\Big)\Big)\,:\,\alpha,\beta_{2r},\gamma_{2r},\delta\in\bbc\Big\}\nonumber\\
& & \bigcap\,\Big\{I_2\otimes\mbox{bl-diag}\{x_1,x_1,x_2,x_2,x_2,x_2,x_1,x_1\}\,:\,x_1,x_2\in\Delta_4^{(k)}\Big\}=\bbc\,.
\end{IEEEeqnarray}
Because of the fact that $\Delta_4=\Delta_2\otimes\Delta_2$, to establish \Cref{caslr}, it is enough to prove the following,
\begin{IEEEeqnarray}{lCl}\label{daslr}
& & \Big\{\Delta_2\otimes\mathrm{Ad}_{V_{k-1}}\Big(\alpha I_{m_k(1)}\bigoplus\bigoplus_{r=1}^{\frac{k-1}{2}}\Big(\beta_{2r} I_{m_k(\theta^{\,2})}\oplus\gamma_{2r} I_{m_k(\overline{\theta}^{\,2})}\Big)\bigoplus\delta I_{m_k(\theta^{\,k+1})}\Big)\Big)\,:\,\alpha,\beta_{2r},\gamma_{2r},\delta\in\bbc\Big\}\nonumber\\
& & \bigcap\,\Big\{\mbox{bl-diag}\{x_1,x_1,x_2,x_2,x_2,x_2,x_1,x_1\}\,:\,x_1,x_2\in\Delta_4^{(k)}\Big\}=\bbc\,.
\end{IEEEeqnarray}
Now, observe that the intersection in \Cref{daslr} is determined by the following two intersections,
\begin{IEEEeqnarray}{lCl}\label{lk1}
& & \Big\{\mathrm{Ad}_{V_{k-1}}\Big(\alpha I_{m_k(1)}\bigoplus\bigoplus_{r=1}^{\frac{k-1}{2}}\Big(\beta_{2r} I_{m_k(\theta^{\,2})}\oplus\gamma_{2r} I_{m_k(\overline{\theta}^{\,2})}\Big)\bigoplus\delta I_{m_k(\theta^{\,k+1})}\Big)\Big)\,:\,\alpha,\beta_{2r},\gamma_{2r},\delta\in\bbc\Big\}\nonumber\\
& & \bigcap\,\,\mbox{bl-diag}\{x_1,x_1,x_2,x_2\}\,,
\end{IEEEeqnarray}
and
\begin{IEEEeqnarray}{lCl}\label{lk2}
& & \Big\{\mathrm{Ad}_{V_{k-1}}\Big(\alpha^\prime I_{m_k(1)}\bigoplus\bigoplus_{r=1}^{\frac{k-1}{2}}\Big(\beta^\prime_{2r} I_{m_k(\theta^{\,2})}\oplus\gamma^\prime_{2r} I_{m_k(\overline{\theta}^{\,2})}\Big)\bigoplus\delta^\prime I_{m_k(\theta^{\,k+1})}\Big)\Big)\,:\,\alpha^\prime,\beta^\prime_{2r},\gamma^\prime_{2r},\delta^\prime\in\bbc\Big\}\nonumber\\
& & \bigcap\,\mbox{bl-diag}\{x_2,x_2,x_1,x_1\}\,.
\end{IEEEeqnarray}
Consider arbitrary element $\xi(\alpha,\beta,\gamma,\delta)$ in \Cref{lk1}
\[
\xi(\alpha,\beta,\gamma,\delta):=\mathrm{Ad}_{V_{k-1}}\Big(\alpha I_{m_k(1)}\bigoplus\bigoplus_{r=1}^{\frac{k-1}{2}}\Big(\beta_{2r}I_{m_k(\theta^{\,2r})}\oplus\gamma_{2r}I_{m_k(\overline{\theta}^{\,2r})}\Big)\bigoplus\delta I_{m_k(\theta^{\,k+1})}\Big)\,,
\]
and similarly $\xi(\alpha^\prime,\beta^\prime,\gamma^\prime,\delta^\prime)$ in \Cref{lk2}. Then $\xi(\alpha,\beta,\gamma,\delta)$ (similarly $\xi(\alpha^\prime,\beta^\prime,\gamma^\prime,\delta^\prime)$) is of the form described in \Cref{easlr} in Step $1$. Therefore, we get the following six equations,
\begin{IEEEeqnarray*}{lCl}
x_1 &=& \alpha I_{m_{k-1}(1)}\bigoplus\bigoplus_{r=1}^{\frac{k-1}{2}}\Big(\beta_{2r} I_{m_{k-1}(\theta^{\,2})}\oplus\gamma_{2r} I_{m_{k-1}(\overline{\theta}^{\,2})}\Big)\bigoplus\delta I_{m_{k-1}(\theta^{\,k+1})}\\
x_1 &=& \alpha^\prime I_{m_{k-1}(1)}\bigoplus\bigoplus_{r=1}^{\frac{k-1}{2}}\Big(\gamma^\prime_{2r} I_{m_{k-1}(\theta^{\,2})}\oplus\beta^\prime_{2r} I_{m_{k-1}(\overline{\theta}^{\,2})}\Big)\bigoplus\delta^\prime I_{m_{k-1}(\theta^{\,k+1})}\\
x_1 &=& \big(\beta^\prime_2 I_{m_{k-1}(1)}\oplus\alpha^\prime I_{m_{k-1}(\theta^{\,2})}\big)\bigoplus\bigoplus_{r=1}^{\frac{k-3}{2}}\Big(\beta^\prime_{2r+2}I_{m_{k-1}(\overline{\theta}^{\,2r})}\oplus\gamma^\prime_{2r}I_{m_{k-1}(\theta^{\,2r+2})}\Big)\\
&  & \hspace*{4.5cm}\bigoplus\big(\delta^\prime I_{m_{k-1}(\theta^{\,k+1})}\oplus\gamma^\prime_{k-1}I_{m_{k-1}(\theta^{\,k+1})}\big)\,,
\end{IEEEeqnarray*}
and
\begin{IEEEeqnarray*}{lCl}
x_2 &=& \alpha I_{m_{k-1}(1)}\bigoplus\bigoplus_{r=1}^{\frac{k-1}{2}}\Big(\gamma_{2r} I_{m_{k-1}(\theta^{\,2})}\oplus\beta_{2r} I_{m_{k-1}(\overline{\theta}^{\,2})}\Big)\bigoplus\delta I_{m_{k-1}(\theta^{\,k+1})}\\
x_2 &=& \big(\beta_2 I_{m_{k-1}(1)}\oplus\alpha I_{m_{k-1}(\theta^{\,2})}\big)\bigoplus\bigoplus_{r=1}^{\frac{k-3}{2}}\Big(\beta_{2r+2}I_{m_{k-1}(\overline{\theta}^{\,2r})}\oplus\gamma_{2r}I_{m_{k-1}(\theta^{\,2r+2})}\Big)\\
&  & \hspace*{4.4cm}\bigoplus\big(\delta I_{m_{k-1}(\theta^{\,k+1})}\oplus\gamma_{k-1}I_{m_{k-1}(\theta^{\,k+1})}\big)\\
x_2 &=& \alpha^\prime I_{m_{k-1}(1)}\bigoplus\bigoplus_{r=1}^{\frac{k-1}{2}}\Big(\beta^\prime_{2r} I_{m_{k-1}(\theta^{\,2})}\oplus\gamma^\prime_{2r} I_{m_{k-1}(\overline{\theta}^{\,2})}\Big)\bigoplus\delta^\prime I_{m_{k-1}(\theta^{\,k+1})}\,.
\end{IEEEeqnarray*}
If the associated `matrix order' with each complex number is understood and creates no confusion, writing in terms of tuple, these six equations are the following set of equalities,
\begin{IEEEeqnarray*}{lCl}
x_1 &=& (\alpha,\beta_2,\gamma_2,\beta_4,\gamma_4,\ldots,\beta_{k-1},\gamma_{k-1},\delta)\\
&=& (\alpha^\prime,\gamma^\prime_2,\beta^\prime_2,\gamma^\prime_4,\beta^\prime_4,\ldots,\gamma^\prime_{k-1},\beta^\prime_{k-1},\delta^\prime)\\
&=& (\beta_2^\prime,\alpha^\prime,\beta^\prime_4,\gamma^\prime_2,\beta^\prime_6,\gamma^\prime_4,\ldots,\beta^\prime_{k-1},\gamma^\prime_{k-3},\delta^\prime,\gamma^\prime_{k-1})\,,
\end{IEEEeqnarray*}
and
\begin{IEEEeqnarray*}{lCl}
x_2 &=& (\alpha,\gamma_2,\beta_2,\gamma_4,\beta_4,\ldots,\gamma_{k-1},\beta_{k-1},\delta)\\
&=& (\beta_2,\alpha,\beta_4,\gamma_2,\beta_6,\gamma_4,\ldots,\beta_{k-1},\gamma_{k-3},\delta,\gamma_{k-1})\\
&=& (\alpha^\prime,\beta^\prime_2,\gamma^\prime_2,\beta^\prime_4,\gamma^\prime_4,\ldots,\beta^\prime_{k-1},\gamma^\prime_{k-1},\delta^\prime)\,.
\end{IEEEeqnarray*}
Immediately we see that $\alpha=\alpha^\prime$ and $\delta=\delta^\prime$. From the second and third equations associated with $x_1$, we get the following,
\[
\alpha^\prime=\beta^\prime_{2r}=\beta^\prime_{2r+2}=\delta^\prime\,\,\mbox{ and }\,\,\alpha^\prime=\gamma_{2r}^\prime=\gamma^\prime_{2r+2}=\delta^\prime
\]
for $1\leq r\leq\frac{k-3}{2}$. These two equations together with $\alpha=\alpha^\prime$ and $\delta=\delta^\prime$ shows that $x_1=\alpha(1,\ldots,1)$. Moreover, by the third equation associated with $x_2$, we see that $x_2=\alpha(1,\ldots,1)$. Hence, we finally get that $x_1=x_2=\alpha(1,\ldots,1)$, that is, $\xi(\alpha,\beta,\gamma,\delta)=\xi(\alpha^\prime,\beta^\prime,\gamma^\prime,\delta^\prime)=\alpha 1$. This finally establishes the claim in \Cref{f}, and finishes this step.
\smallskip

\noindent\textbf{Step 4:~ The Ocneanu compactness.}
\smallskip

\noindent This is the final step of the proof. First, briefly recall the construction of $R^{(m)}_{u,v}$ from \Cref{producing}. The following quadruple
\[
\begin{matrix}
A_{2k-1} &\subset & A_{2k+1}\\
|| &  & ||\\
\bbc\otimes\Delta_4\otimes M_4^{(k)} & & \Delta_4\otimes M_4\otimes M_4^{(k)}\\
\cup &  & \cup\\
\bbc\otimes\mbox{Ad}_{u_{2k-1}}\big(\bbc\otimes\big((u_{2k-2}^*v_{2k-2})^2\big)^\prime\cap M_4^{(k)}\big) & & \mbox{Ad}_{u_{2k+1}}\big(\bbc\otimes\big((u_{2k}^*v_{2k})^2\big)^\prime\cap M_4^{(k+1)}\big)\\
|| &  & ||\\
C_{2k-1} &\subset & C_{2k+1}
\end{matrix}
\]
becomes a non-degenerate commuting square, and we have started iterating the basic construction from this stage. If $D_{2k+j}$, for $j\geq 3$ odd integer, denote the finite-dimensional subalgebras of $C_{2k+j}$ obtained through the basic construction, then $R^{(m)}_{u,v}$ is constructed as the SOT-limit of $D_{2k+j}$, whereas $R$ is the SOT-limit of $A_{2k+j}$. Note that we have used the convention that $D_{2k+i}=C_{2k+i}$ for $i=\pm 1$. Since for $j\geq -1$ odd integer, we have $(D_{2k+j}\subset A_{2k+j},D_{2k+j+2}\subset A_{2k+j+2})$ is a tower of symmetric commuting squares with limit $R^{(m)}_{u,v}\subset R$, by the Ocneanu compactness we have the following,
\[
\big(R^{(m)}_{u,v}\big)^{\,\prime}\cap R=D_{2k+1}^{\,\prime}\cap A_{2k-1}=C_{2k+1}^{\,\prime}\cap A_{2k-1}\,.
\]
Now, by the Step $3$ we get the irreducibility of $R^{(m)}_{u,v}\subset R$ when $k$ is odd.

Finally, the exact similar analysis holds when $k$ is even in Step $1$, which completes the proof.\qed
\end{prf}

\begin{thm}\label{lop}
For the pair of spin model subfactors $R_{u(a)}, R_{u(b)}\subset R$ such that $u(a)\sim u(b)$, we have $R_{u(a)}\cap R_{u(b)}$ is a factor and $R_{u(a)}\cap R_{u(b)}=R^{(m)}_{u(a),u(b)}$ if $\,(\overline{a}b)^{\,m}=1$.
\end{thm}
\begin{prf}
Let $u=u(a)\mbox{ and }v=u(b)$ be $4\times 4$ Hadamard inequivalent complex Hadamard matrices such that $u\sim v$. If $\theta=\overline{a}b$ satisfies $\,\theta^{\,m}=1$, we take $k=(m-2)/2$. By \Cref{interesting}, we immediately see that $R_u\cap R_v$ is a factor since $R^{(m)}_{u,v}\subset R_u\cap R_v$ and $R^{(m)}_{u,v}\subset R$ is irreducible. For the identification of $R_u\cap R_v$ with $R^{(m)}_{u,v}$, it is enough to prove that $[R_u:R_u\cap R_v]=[R_u:R^{(m)}_{u,v}]$. Note that $[R_u:R^{(m)}_{u,v}]=\frac{m}{2}$.

We claim that $[R_u:R_u\cap R_v]=\lambda^{-1}(R_u,R_u\cap R_v)=\frac{m}{2}$. The first equality is obvious by \cite{PP}, since $R_u\cap R_v$ is a factor. To prove the second equality, recall that $R_u$ is the SOT-limit of $B^u_{2j+1}$ and $R_u\cap R_v$ is that of $C_{2j+1}$. Given $\theta^{\,m}=1$, if we take $k=(m-2)/2$, then $\sigma\big((u_{2k-2}^*v_{2k-2})^2\big)=\sigma\big((u_{2k+2j}^*v_{2k+2j})^2\big)$ for all $j\in\bbn\cup\{0\}$, and cardinality of the spectrum is $k+1$ (see \Cref{101}). Since $(u_{2k+2j}^*v_{2k+2j})^2$ is diagonalizable for all $j\in\bbn\cup\{-1\}$, the sum of the multiplicities of the eigenvalues of $(u_{2k+2j}^*v_{2k+2j})^2$ will add up to the order of $M_4^{(k+j+1)}$ for each $j\in\bbn\cup\{-1\}$, which is equal to $4^{k+j+1}$. Note that the eigenvalues remain the same, but the multiplicities are changing with $j$. Since $C_{2k+2j+1}$ is given by the algebra $\big((u_{2k+2j}^*v_{2k+2j})^2\big)^\prime\cap M_4^{(k+j+1)}$ (\Cref{odd intersection}), we see that at each stage of the inclusion $C_{2k+2j+1}\subset B^u_{2k+2j+1}$, we have $(k+1)$-many Wedderburn-Artin direct summands in $M_{4^{k+j+1}}$, and order of the summands add up to $4^{k+j+1}$.

Now, each quadruple $(C_{2r+1}\subset B^u_{2r+1},C_{2r+3}\subset B^u_{2r+3})$ is a commuting square for $r\in\bbn\cup\{0\}$, which follows as an application of \Cref{comm cube}. By \Cref{ent formula2} and \Cref{popaadaptation}, we see that $\lambda(R_u,R_u\cap R_v)$ is the limit of the eventually constant sequence $\{(k+1)^{-1}\}$. Therefore, $\lambda(R_u,R_u\cap R_v)^{-1}=k+1=\frac{m}{2}$, which proves the claim. Thus, we get that $\lambda(R_u,R_u\cap R_v)^{-1}=[R_u:R_u\cap R_v]=\frac{m}{2}=[R_u:R^{(m)}_{u,v}]$. By the multiplicativity of the Jones' index, we immediately get that $R^{(m)}_{u,v}=R_u\cap R_v$.\qed
\end{prf}

We remark that using the above theorem, it takes only a little more effort to show that in the following tower of finite-dimensional algebras
\[
C_1\subset C_3\subset C_5\subset\ldots\subset C_{2k-1}\subset C_{2k+1}\subset\ldots\subset R_u\cap R_v
\]
if $\theta^{\,m}=1$, then $C_{2k-1}\subset C_{2k+1}\subset\ldots\subset R_u\cap R_v$ is a tower of basic construction, where $k=(m-2)/2$. To see this, fix a $m$ and let $k=(m-2)/2$. Recall that the subfactor $R^{(m)}_{u,v}\subset R$ is the limit of the tower of symmetric commuting squares $(D_{2k+j}\subset A_{2k+j},D_{2k+j+2}\subset A_{2k+j+2})$, for $j\geq -1$ odd integer, where $D_{2k+j}=C_{2k+j}$ for $j=\pm 1$ and $D_{2k+j}\subset C_{2k+j}$ for $j\geq 3$. Since for $j\geq 3,\,(D_{2k+j}\subset A_{2k+j},D_{2k+j+2}\subset A_{2k+j+2})$ and $(C_{2k+j}\subset A_{2k+j},C_{2k+j+2}\subset A_{2k+j+2})$ both are commuting squares in the following diagram
\[
\begin{matrix}
\subset & A_{2k+3} &\subset & A_{2k+5} &\subset & A_{2k+7} &\subset &\cdots &\subset & R\\
 & \cup & & \cup & & \cup & & & & \cup\\
\subset & C_{2k+3} &\subset & C_{2k+5} &\subset & C_{2k+7} &\subset &\cdots &\subset & R_u\cap R_v\\
& || & & \cup & & \cup & & & & \cup\\
& D_{2k+3} &\subset  & D_{2k+5} &\subset & D_{2k+7} &\subset &\cdots &\subset & R^{(m)}_{u,v}\\
\end{matrix}
\]
we have that $(D_{2k+j}\subset C_{2k+j},D_{2k+j+2}\subset C_{2k+j+2})$ is also a commuting square for any $j\geq 3$ odd integer. This follows from the following,
\[
E_{D_{2k+j+2}}^{C_{2k+j+2}}(C_{2k+j})=E_{D_{2k+j+2}}^{B^u_{2k+j+2}}(C_{2k+j})\subset E_{D_{2k+j+2}}^{B^u_{2k+j+2}}(B^u_{2k+j})\subset D_{2k+j}\,.
\]
Now, by Proposition $2.6$ in \cite{PP} and \Cref{lop}, we get that $1=\lambda(R_u\cap R_v,R^{(m)}_{u,v})=\lim_{j\to\infty}\lambda(C_{2k+j},D_{2k+j})$ decreasingly. Since $\lambda$ is always less than or equal to $1$, we get that $\lambda(C_{2k+j},D_{2k+j})=1$ for all $j\geq 3$ odd integer. This says that $C_{2k+j}\subseteq D_{2k+j}$ for all $j$ by \Cref{sk}, and since the reverse inclusion is already there, we finally get that $C_{2k+j}=D_{2k+j}$ for all $j\geq 3$ odd integer. Since the tower $D_{2k+3}\subset D_{2k+7}\subset\cdots$ is the iterated basic construction for the inclusion $C_{2k-1}\subset C_{2k+1}$, our claim follows.

\begin{crlre}
For $a,b,c,d\in\mathbb{S}^1$, if $\big(\overline{a}b\big)^m=\big(\overline{c}d\big)^n=1$ and $m\neq n$, then the corresponding subfactors $R_{u(a)}\cap R_{u(b)}\subset R$ and $R_{u(c)}\cap R_{u(d)}\subset R$ are non-isomorphic. 
\end{crlre}

\begin{crlre}\label{exmple2}
For Hadamard inequivalent complex Hadamard matrices $u\mbox{ and }v$ such that $u\sim v$, the quadruple of $II_1$ factors $(R_u\cap R_v\subset R_u,R_v\subset R)$ is obtained as an iterated basic construction of a non-degenerate commuting cube.
\end{crlre}
\begin{prf}
Follows from \Cref{producing} and \ref{lop}, and the discussion in \Cref{Sec 2}.
\end{prf}


\subsubsection{The interior, exterior angles and its rigidity}

We show that the (interior) angle between $R_u$ and $R_v$ in the sense of \cite{BDLR} (see \Cref{angletams}) has certain rigidity when $u\sim v$. Recall from \Cref{lop} that if $u\sim v$ and $\theta^{\,m}=1$ for some $m\in 2\mathbb{N}+2$, where $u=u(a)$ and $v=u(b)$ with $\theta=\overline{a}b$, then we have $R^{(m)}_{u,v}=R_u\cap R_v$. Therefore, we get  a family of quadruples of $II_1$ factors
\[
\begin{matrix}
R_u & \subset  & R \cr
\cup &  &\cup  \cr
R_u\cap R_v& \subset & R_v&
\end{matrix}
\]
such that $[R: R_u\cap R_v]=2m$ and $[R_u:R_u\cap R_v]=[R_v:R_u\cap R_v]=m/2$. Recall that given any quadruple of finite index subfactors
\[ \begin{matrix}
P & \subset  & M \cr
\cup &  &\cup  \cr
N & \subset & Q&
\end{matrix}
\]
the interior angle $\alpha^N_M(P,Q)$ and the exterior angle $\beta^N_M(P,Q)$ have the following formulae (see \cite{BDLR}):
\begin{align*}
\alpha^N_M(P,Q) &=\displaystyle \frac{[M:N]\,tr(e_Pe_Q)-1}{\sqrt{[P:N]-1}\sqrt{[Q:N]-1}}\,,\\
\beta^N_M(P,Q) &=\displaystyle \frac{tr(e_Pe_Q)-tr(e_{P}) tr(e_{Q})}{\sqrt{tr(e_P)-{tr(e_P)}^2} \sqrt{tr(e_Q)-{tr(e_Q)}^2}}\,.
\end{align*}
We also know that $\lambda(P,Q)=\displaystyle\frac{tr(e_Pe_Q)}{tr(e_P)}$ if $N\subset M$ is irreducible (Theorem $3.3$ in \cite{B}). Since in our case, $R^{(m)}_{u,v}\subset R$ is irreducible by \Cref{interesting}, and $\lambda(R_u,R_v)=1/2$ by \Cref{Popa for 4 by 4}, we have $tr(e_{R_u}e_{R_v})=1/8$. Thus, by the above formulae of angles we see that $\cos \beta (R_u,R_v)= 1/3$ and $\cos \alpha(R_u,R_v)=\frac{m-4}{2(m-2)}<1/2$. Therefore, we conclude that $\alpha(R_u,R_v)>\pi/3.$ We summarize these findings as the following theorem.

\begin{thm}\label{rigidity of angle}
The interior angle $\alpha \big(R_u,R_v\big)$  is strictly  greater than $\pi/3$ and the exterior angle $\beta \big(R_u,R_v\big)=\arccos 1/3$. Furthermore, $\alpha \big(R_u,R_v\big)$ converges to $\pi/3$ decreasingly as $m$ tends to infinity.
\end{thm}

As an immediate corollary, using \Cref{basicresultoncomm} we obtain the following.
\begin{crlre}\label{final cor 0}
For any $4\times 4$ Hadamard inequivalent complex Hadamard matrices $u\mbox{ and }v$ such that $u\sim v$, the following quadruple of subfactors
\[
\begin{matrix}
R_u & \subset  & R \cr
\cup &  &\cup  \cr
R_u\cap R_v& \subset & R_v&
\end{matrix}
\]
is a commuting square if and only if $m=4$, that is, $b=\pm ia$. Moreover, the above quadruple is never a co-commuting square.
\end{crlre}


\subsection{The Sano-Watatani angle and Connes-St{\o}rmer relative entropy}

For inequivalent complex Hadamard matrices $u\mbox{ and }v$, we compute the exact value of relative entropy $H(R_u|R_v)$ (along with its variant $h(R_u|R_v)$) and the Sano-Watatani angle $\mathrm{Ang}_R(R_u,R_v)$ for the case of $u\sim v$. When $u\nsim v$, we give a legitimate upper and lower bounds for $H(R_u|R_v)$. As a crucial tool, we have shown the existence of a proper subfactor of $R$ containing both $R_u$ and $R_v$.

\subsubsection{The subfactor $R_u\vee R_v\subset R$}

This subsection is devoted to  prove the existence of an index $2$ subfactor of $R$ containing both $R_u\mbox{ and }R_v$. Recall that we have pair of irreducible subfactors $R_u$ and $R_v$ of $R$ each with index $4$. Consider the von Neumann subalgebra $R_u\vee R_v\subset R$. By the irreducibility of $R_u\subset R$ we get that $R_u\vee R_v\subset R$ is a factor. We show that $R_u\vee R_v$ is a proper subfactor $R$. We begin with the following notations which are at par with that in the beginning of \Cref{Sec 5}. For $k\in\bbn\cup\{0\}$, let
\begin{align*}
A_{2k} &= M_4\otimes M_4^{(k)}\,,\cr
B^u_{2k} &= u_{2k}(\Delta_4\otimes M_4^{(k)})u_{2k}^*\,,\cr
B^v_{2k} &= v_{2k}(\Delta_4\otimes M_4^{(k)})v_{2k}^*\,,\cr
\mathcal{I}_{2k} &= B^u_{2k}\vee B^v_{2k}\,.
\end{align*}

\noindent Note that these are the even steps in the tower of basic constructions described in \Cref{fig2}.

\begin{lmma}\label{proper}
For all $k\in\bbn\cup\{0\},\,\mathcal{I}_{2k}$ is a proper subalgebra of $A_{2k}$.
\end{lmma}
\begin{prf}
Since each $A_{2k}$ is factor, it is enough to show that $\bbc\oplus\bbc\subseteq\mathcal{I}_{2k}^{\,\prime}\cap A_{2k}$ for any $k\in\bbn\cup\{0\}$. First observe that
\begin{IEEEeqnarray}{lCl}\label{als}
(B^u_{2k})^\prime\cap A_{2k} &=& \mbox{Ad}_{u_{2k}}\big((\Delta_4\otimes M_4^{(k)})^\prime\cap\mbox{Ad}_{u_{2k}^*}(A_{2k})\big)\nonumber\\
&=& \mbox{Ad}_{u_{2k}}\big((\Delta_4\otimes M_4^{(k)})^\prime\cap(M_4\otimes M_4^{(k)})\big)\qquad\big(\mbox{as }u_{2k}\in A_{2k}=M_4^{(k+1)}\big)\nonumber\\
&=& \mbox{Ad}_{u_{2k}}(\Delta_4\otimes\bbc^{(k)})\,,
\end{IEEEeqnarray}
and similarly, $(B^v_{2k})^\prime\cap A_{2k}=\mbox{Ad}_{v_{2k}}(\Delta_4\otimes\bbc^{(k)})$. Now, $\mathcal{I}_{2k}^{\,\prime}\cap A_{2k}=\{B^u_{2k}\cup B^v_{2k}\}^\prime\cap A_{2k}=(B^u_{2k})^\prime\cap(B^v_{2k})^\prime\cap A_{2k}$. Hence, by \Cref{als} we need to show the following,
\begin{IEEEeqnarray}{lCl}\label{bls}
\bbc\oplus\bbc\subseteq\mbox{Ad}_{u_{2k}}(\Delta_4\otimes\bbc^{(k)})\cap\mbox{Ad}_{v_{2k}}(\Delta_4\otimes\bbc^{(k)}),
\end{IEEEeqnarray}
or in other words,
\begin{IEEEeqnarray}{lCl}\label{cls}
\bbc\oplus\bbc\subseteq\mbox{Ad}_{u_{2k}}\big((\Delta_4\otimes\bbc^{(k)})\cap\mbox{Ad}_{u_{2k}^*v_{2k}}(\Delta_4\otimes\bbc^{(k)})\big).
\end{IEEEeqnarray}
Recall from \Cref{cute lemma for 4 by 4} that for any $k\in\bbn$, one has the following identity in $M_4^{(k+1)}$
\[
u_{2k}^*v_{2k}=p\otimes u_{2k-2}^*v_{2k-2}+q(\theta)\otimes v_{2k-2}^*u_{2k-2}\,,
\]
with the convention $u_0=u\mbox{ and }v_0=v$. Applying \Cref{Sa} to this identity, we get the following,
\begin{IEEEeqnarray}{lCl}\label{lk}
& & \mbox{Ad}_{W\otimes I_4^{(k)}}(u_{2k}^*v_{2k})\nonumber\\
&=& \mbox{bl-diag}\{u_{2k-2}^*v_{2k-2}\,,\,u_{2k-2}^*v_{2k-2}\,,\,v_{2k-2}^*u_{2k-2}\,,\,\theta\,v_{2k-2}^*u_{2k-2}\}\,,
\end{IEEEeqnarray}
where $W$ is the self-adjoint unitary defined in \Cref{W}. Now, we have the following identity
\begin{IEEEeqnarray}{lCl}\label{dls}
&  & (\Delta_4\otimes\bbc^{(k)})\cap\mbox{Ad}_{u_{2k}^*v_{2k}}(\Delta_4\otimes\bbc^{(k)})\nonumber\\
&=& \mbox{Ad}_{W\otimes I_4^{(k)}}\Big(\big(\mbox{Ad}_W(\Delta_4)\otimes\bbc^{(k)}\big)\cap\mbox{Ad}_{W\otimes I_4^{(k)}}\mbox{Ad}_{u_{2k}^*v_{2k}}\big(\Delta_4\otimes\bbc^{(k)}\big)\Big)\nonumber\\
&=& \mbox{Ad}_{W\otimes I_4^{(k)}}\Big(\big(\mbox{Ad}_W(\Delta_4)\otimes\bbc^{(k)}\big)\cap\mbox{Ad}_{W\otimes I_4^{(k)}}\mbox{Ad}_{u_{2k}^*v_{2k}}\mbox{Ad}_{W\otimes I_4^{(k)}}\big(\mbox{Ad}_W(\Delta_4)\otimes\bbc^{(k)}\big)\Big)\nonumber\\
&=& \mbox{Ad}_{W\otimes I_4^{(k)}}\Big(\big(\mbox{Ad}_W(\Delta_4)\otimes\bbc^{(k)}\big)\cap\mbox{Ad}_{(W\otimes I_4^{(k)})u_{2k}^*v_{2k}(W\otimes I_4^{(k)})}\big(\mbox{Ad}_W(\Delta_4)\otimes\bbc^{(k)}\big)\Big)\,.
\end{IEEEeqnarray}
Therefore, to prove \Cref{cls}, it is enough to show the following inclusion
\begin{IEEEeqnarray}{lCl}\label{els}
\bbc\oplus\bbc\subseteq\big(\mbox{Ad}_W(\Delta_4)\otimes\bbc^{(k)}\big)\cap\mbox{Ad}_{(W\otimes I_4^{(k)})u_{2k}^*v_{2k}(W\otimes I_4^{(k)})}\big(\mbox{Ad}_W(\Delta_4)\otimes\bbc^{(k)}\big)
\end{IEEEeqnarray}
in view of \Cref{dls}. Consider any arbitrary $(z,w)\in\bbc\oplus\bbc$ and the diagonal matrix $D=\mbox{diag}\{z,w,z,w\}\in\Delta_4$. It is a straightforward verification that $WDW=D$ (recall $W$ from \Cref{W}). Hence, $D\otimes I_4^{(k)}\in\mbox{Ad}_W(\Delta_4)\otimes\bbc^{(k)}$. Now by \Cref{lk}, we have the following,
\begin{IEEEeqnarray*}{lCl}
&  & \mbox{Ad}_{(W\otimes I_4^{(k)})u_{2k}^*v_{2k}(W\otimes I_4^{(k)})}(D\otimes I_4^{(k)})\\
&=& \mbox{Ad}_{\mbox{bl-diag}\{u_{2k-2}^*v_{2k-2}\,,\,u_{2k-2}^*v_{2k-2}\,,\,v_{2k-2}^*u_{2k-2}\,,\,\theta\,v_{2k-2}^*u_{2k-2}\}}(D\otimes I_4^{(k)})\\
&=& D\otimes I_4^{(k)}\,,
\end{IEEEeqnarray*}
since $u_{2k-2}^*v_{2k-2}$ is a unitary and $\theta\in\mathbb{S}^1$. This shows that the matrix $D\otimes I_4^{(k)}$, where $D=\mbox{diag}\{z,w,z,w\}$ for $z,w\in\bbc$, lies in the intersection
\[
\big(\mbox{Ad}_W(\Delta_4)\otimes\bbc^{(k)}\big)\cap\mbox{Ad}_{(W\otimes I_4^{(k)})u_{2k}^*v_{2k}(W\otimes I_4^{(k)})}\big(\mbox{Ad}_W(\Delta_4)\otimes\bbc^{(k)}\big)
\]
which validates \Cref{els}. Therefore, \Cref{cls} holds, and consequently we have
\[
\bbc\oplus\bbc\subseteq\mbox{Ad}_{u_{2k}}(\Delta_4\otimes\bbc^{(k)})\cap\mbox{Ad}_{v_{2k}}(\Delta_4\otimes\bbc^{(k)})=\mathcal{I}_{2k}^{\,\prime}\cap A_{2k}\,,
\]
by \Cref{als,bls}, which concludes the proof.\qed
\end{prf}

\begin{lmma}\label{proper1}
For all $k\in\bbn\cup\{0\}$, both $B^u_{2k}\mbox{ and }B^v_{2k}$ are proper subalgebras of $\,\mathcal{I}_{2k}$.
\end{lmma}
\begin{prf}
On contrary assume that there exists $k_0\in\bbn$ such that $\mathcal{I}_{2k_0}=B^u_{2k_0}$ (or $v$ in place of $u$). Then $B^v_{2k_0}\subseteq B^u_{2k_0}$, and consequently $B^v_{2k}\subseteq B^u_{2k}$ for all $k\geq k_0$ (recall the tower of basic constructions descibed in \Cref{fig2}). Thus, we get that $R_v\subset R_u$. Since $[R:R_u]=[R:R_v]=4$, we get that $R_u=R_v$, and consequently $\lambda(R_u,R_v)=1$. However, this is a contradiction to \Cref{Popa for 4 by 4}, and completes the proof.\qed
\end{prf}

\begin{thm}\label{intermediate}
Let $u$ and $v$ be (distinct) $4\times 4$ inequivalent complex Hadamard matrices. Consider the pair of spin model subfactors $R_u,R_v\subset R$. Then, $R_u\vee R_v$ is a proper subfactor of $R$ and $[R:R_u\vee R_v]=2$.
\end{thm}
\begin{prf}
Given $u\mbox{ and }v$, consider the von Neumann subalgebra $R_u\vee R_v\subset R$. Note that $R_u\vee R_v=\big(\cup_k\mathcal{I}_{2k}\big)^{\prime\prime}\subset R$, where $\mathcal{I}_{2k}=B^u_{2k}\vee B^v_{2k}$. By the irreducibility of the spin model subfactor $R_u\subset R$, we get the factoriality of $R_u\vee R_v\subset R$. Hence, the only nontrivial part is to show that $[R:R_u\vee R_v]=2$. We appeal \Cref{proper}, \ref{proper1} and \Cref{proving com} to prove this.
\smallskip

\noindent\textbf{Step 1~:} For each $k$, the quadruple
\[
\begin{matrix}
B^u_{2k+2}\vee B^v_{2k+2} &\subset & A_{2k+2}\\
\cup & & \cup\\
B^u_{2k}\vee B^v_{2k} &\subset & A_{2k}\\
\end{matrix}
\]
is a co-commuting square.
\smallskip

\noindent We need to show that the following quadruple
\begin{IEEEeqnarray}{lCl}\label{wq}
\begin{matrix}
A_{2k}^\prime &\subset & (B^u_{2k})^\prime\cap(B^v_{2k})^\prime\\
\cup & & \cup\\
A_{2k+2}^\prime &\subset & (B^u_{2k+2})^\prime\cap(B^v_{2k+2})^\prime\\
\end{matrix}
\end{IEEEeqnarray}
is a commuting square for each $k$. Fix any $k\in\bbn$. Recall that the quadruple
\[
\begin{matrix}
B^u_{2k+2} &\subset & A_{2k+2}\\
\cup & & \cup\\
B^u_{2k} &\subset & A_{2k}\\
\end{matrix}
\]
is a co-commuting square by the construction of $R_u\subset R$. Hence, the quadruple
\[
\begin{matrix}
A_{2k}^\prime &\subset & (B^u_{2k})^\prime\\
\cup & & \cup\\
A_{2k+2}^\prime &\subset & (B^u_{2k+2})^\prime\\
\end{matrix}
\]
is a commuting square. Now, in \Cref{wq} we have the following,
\begin{IEEEeqnarray*}{lCl}
E^{(B^u_{2k})^\prime\cap(B^v_{2k})^\prime}_{A_{2k}^\prime}\big((B^u_{2k+2})^\prime\cap(B^v_{2k+2})^\prime\big) &=& E^{(B^u_{2k})^\prime\cap(B^v_{2k})^\prime}_{{A_{2k}}^\prime}E^{(B^u_{2k})^\prime}_{(B^u_{2k})^\prime\cap(B^v_{2k})^\prime}\big((B^u_{2k+2})^\prime\cap(B^v_{2k+2})^\prime\big)\\
&=& E^{(B^u_{2k})^\prime}_{A_{2k}^\prime}\big((B^u_{2k+2})^\prime\cap(B^v_{2k+2})^\prime\big)\\
&\subset& E^{(B^u_{2k})^\prime}_{A_{2k}^\prime}\big((B^u_{2k+2})^\prime\big)\\
&\subset& A_{2k+2}^\prime\,,
\end{IEEEeqnarray*}
which completes Step $1$.
\smallskip

\noindent\textbf{Step 2~:} $||\Lambda_{\mathcal{I}_{2k}}^{A_{2k}}||=||\Lambda^{\mathcal{I}_{2k}}_{B^u_{2k}}||=\sqrt{2}$ for all $k$.
\smallskip

\noindent Since $B^u_{2k}=\mbox{Ad}_{u_{2k}}(\Delta_4\otimes M_4^{(k)})$, we have $||\Lambda^{A_{2k}}_{B^u_{2k}}||=||\Lambda^{M_4\otimes M_4^{(k)}}_{\Delta_4\otimes M_4^{(k)}}||=2$. By the multiplicativity of the Watatani index, we have the equality $||\Lambda^{A_{2k}}_{B^u_{2k}}||^2=||\Lambda^{A_{2k}}_{\mathcal{I}_{2k}}||^2\,||\Lambda^{\mathcal{I}_{2k}}_{B^u_{2k}}||^2$. Therefore, we get that $2=||\Lambda^{A_{2k}}_{\mathcal{I}_{2k}}||\,||\Lambda^{\mathcal{I}_{2k}}_{B^u_{2k}}||$. By \Cref{proper} and \ref{proper1}, we have $||\Lambda^{A_{2k}}_{\mathcal{I}_{2k}}||>1$ and $||\Lambda^{\mathcal{I}_{2k}}_{B^u_{2k}}||>1$ respectively. Hence by Theorem $1.1.1$ in \cite{GHJ}, both $||\Lambda^{A_{2k}}_{\mathcal{I}_{2k}}||$ and $||\Lambda^{\mathcal{I}_{2k}}_{B^u_{2k}}||$ are greater than equal to $\sqrt{2}$. Thus, we get that $||\Lambda_{\mathcal{I}_{2k}}^{A_{2k}}||=||\Lambda^{\mathcal{I}_{2k}}_{B^u_{2k}}||=\sqrt{2}$, which finishes this step.
\smallskip

Finally, applying \Cref{proving com} together with Step $1$ and Step $2$, we obtain that the following quadruple
\[
\begin{matrix}
B^u_{2k+2}\vee B^v_{2k+2} &\subset & A_{2k+2}\\
\cup & & \cup\\
B^u_{2k}\vee B^v_{2k} &\subset & A_{2k}\\
\end{matrix}
\]
is a symmetric commuting square for all $k$. Since $R=\overline{\cup_k A_{2k}}^{\,\mbox{sot}}$ and $R_u\vee R_v=\big(\cup_k\mathcal{I}_{2k}\big)^{\prime\prime}\subset R$, we get that the following quadruple
\[
\begin{matrix}
R_u\vee R_v &\subset & R\\
\cup & & \cup\\
\mathcal{I}_0=B^u_0\vee B^v_0 &\subset & A_0\\
\end{matrix}
\]
is a symmetric commuting square. Therefore, $[R:R_u\vee R_v]=||\Lambda^{A_0}_{\mathcal{I}_0}||^2=2$ by Step $2$.\qed
\end{prf}


\subsubsection{The Sano-Watatani angle when $u\sim v$}

Let us begin at the beginning. We first compute the angle between $\mbox{Ad}_u(\Delta_4)\subset M_4$ and $\mbox{Ad}_v(\Delta_4)\subset M_4$, and show that the following quadruple
\begin{IEEEeqnarray}{lCl}\label{as}
\begin{matrix}
\mbox{Ad}_u(\Delta_4) &\subset & M_4\\
 \cup & & \cup\\
\mbox{Ad}_u(\Delta_4)\cap\mbox{Ad}_v(\Delta_4) &\subset & \mbox{Ad}_v(\Delta_4)\\
\end{matrix}
\end{IEEEeqnarray}
is not always a commuting square. This is another example where the floor (and consequently the roof due to \Cref{comm cube 2}) in \Cref{com1} need not be commuting square.

\begin{ppsn}\label{p4}
If $u\mbox{ and }v$ are parametrized by the circle parameters $a\mbox{ and }b$ respectively, then the cosine of the angle between $\mathrm{Ad}_u(\Delta_4)\mbox{ and }\mathrm{Ad}_v(\Delta_4)$ is the set $\{0,|\Re(\overline{a}b)|\}$, and the quadruple in \Cref{as} is a commuting square if and only if $b=\pm ia$.
\end{ppsn}
\begin{prf}
Consider the nonnegative operator $E_{u\Delta_4u^*}E_{v\Delta_4v^*}E_{u\Delta_4u^*}$, and observe that as an element of $\mathcal{B}(M_4,tr)=M_{16}$, this is same as the following operator
\[
\mbox{Ad}_uE^{M_4}_{\Delta_4}\mbox{Ad}_{u^*}\mbox{Ad}_vE^{M_4}_{\Delta_4}\mbox{Ad}_{v^*}\mbox{Ad}_uE^{M_4}_{\Delta_4}\mbox{Ad}_{u^*}=\mbox{Ad}_u\circ E^{M_4}_{\Delta_4}E^{M_4}_{\mbox{Ad}_{u^*v}(\Delta_4)}E^{M_4}_{\Delta_4}\circ\mbox{Ad}_{u^*}.
\]
First consider the operator $E^{M_4}_{\Delta_4}E^{M_4}_{\mbox{Ad}_{u^*v}(\Delta_4)}E^{M_4}_{\Delta_4}\in\mathcal{B}(M_4)$. Let $\,\theta=\overline{a}b$ and $\{E_{ij}:1\leq i,j\leq 4\}$ be the standard set of matrix units in $M_4$. Then, observe the following,
\begin{IEEEeqnarray}{lCl}\label{angle computation}
E^{M_4}_{\Delta_4}E^{M_4}_{u^*v\Delta_4v^*u}E^{M_4}_{\Delta_4}(E_{ij}) &=& E^{M_4}_{\Delta_4}E^{M_4}_{u^*v\Delta_4v^*u}(\delta_{ij}E_{ij})\nonumber\\
&=& \delta_{ij}E^{M_4}_{\Delta_4}\mbox{Ad}_{u^*v}\,E^{M_4}_{\Delta_4}\left(v^*uE_{ii}u^*v\right)\,.
\end{IEEEeqnarray}
It is easy to verify that $E^{M_4}_{\Delta_4}\left(v^*uE_{ii}u^*v\right)=E_{ii}$ for $i=1,3$, and
\begin{IEEEeqnarray*}{lCl}
E^{M_4}_{\Delta_4}\left(v^*uE_{22}u^*v\right) &=& \mbox{diag}\{0,(2+\theta+\overline{\theta})/4,0,(2-\theta-\overline{\theta})/4\}\,,\\
E^{M_4}_{\Delta_4}\left(v^*uE_{44}u^*v\right) &=& \mbox{diag}\{0,(2-\theta-\overline{\theta})/4,0,(2+\theta+\overline{\theta})/4\}\,.
\end{IEEEeqnarray*}
Therefore, we have the following,
\begin{IEEEeqnarray*}{lCl}
&  & E^{M_4}_{\Delta_4}\big(\mbox{Ad}_{v^*u}(E_{ii})\big)\\
&=& \delta_{i1}E_{11}+\delta_{i3}E_{33}+\delta_{i2}\left(\frac{2+\theta+\overline{\theta}}{4}E_{22}+\frac{2-\theta-\overline{\theta}}{4}E_{44}\right)+\delta_{i4}\left(\frac{2-\theta-\overline{\theta}}{4}E_{22}+\frac{2+\theta+\overline{\theta}}{4}E_{44}\right)\\
&=& \delta_{i1}E_{11}+\delta_{i3}E_{33}+\delta_{i2}\left(\frac{2+t}{4}E_{22}+\frac{2-t}{4}E_{44}\right)+\delta_{i4}\left(\frac{2-t}{4}E_{22}+\frac{2+t}{4}E_{44}\right)\,,
\end{IEEEeqnarray*}
where $t=\theta+\overline{\theta}=\overline{a}b+a\overline{b}\in\mathbb{R}$. Therefore, by \Cref{angle computation} we finally get the following
\begin{IEEEeqnarray}{lCl}\label{angle spectrum 4 by 4}
&  & E^{M_4}_{\Delta_4}E^{M_4}_{u^*v\Delta_4v^*u}E^{M_4}_{\Delta_4}(E_{ij})\nonumber\\
&=& \begin{cases}
\delta_{i1}E_{11}+\delta_{i3}E_{33}+\frac{1}{8}\big((4+t^2)\delta_{i2}+(4-t^2)\delta_{i4}\big)E_{22}\\
+\frac{1}{8}\big((4-t^2)\delta_{i2}+(4+t^2)\delta_{i4}\big)E_{44} & \mbox{ if }\,\,i=j, \cr
0 & \mbox{ if }\,\,i\neq j. \cr
\end{cases}
\end{IEEEeqnarray}
Now fix the ordered orthornormal basis $\mathfrak{B}:=\{2E_{11},2E_{22},2E_{33},2E_{44},2E_{k\ell}\,:\,1\leq k,\ell\leq 4\,,\,k\neq \ell\}$ of $M_4$. By \Cref{angle spectrum 4 by 4} we get that the matrix of the linear operator $E^{M_4}_{\Delta_4}E^{M_4}_{u^*v\Delta_4v^*u}E^{M_4}_{\Delta_4}$ acting on $M_4$, when viewed as a linear map from $\bbc^{16}$ to $\bbc^{16}$, with respect to the basis $\mathfrak{B}$ is a $16\times 16$ matrix having zero everywhere else except the following matrix
\[
B:=\left[{\begin{matrix}
1 & 0 & 0 & 0\\
0 & \frac{1}{8}(4+t^2) & 0 & \frac{1}{8}(4-t^2)\\
0 & 0 & 1 & 0\\
0 & \frac{1}{8}(4-t^2) & 0 & \frac{1}{8}(4+t^2)\\
\end{matrix}}\right]
\]
at the top left $4\times 4$ block.

Now, recall from \Cref{even intersection} that $\mathrm{Ad}_u(\Delta_4)\cap\mathrm{Ad}_v(\Delta_4)=\mathrm{Ad}_u(\bbc^3)$, where the emebedding $\bbc^3\xhookrightarrow{}\Delta_4\subseteq M_4$ is given by $(\lambda_1,\lambda_2,\lambda_3)\mapsto(\lambda_1,\lambda_2,\lambda_3,\lambda_2)$. The trace preserving conditional expectation $E^{M_4}_{\mathrm{Ad}_u(\Delta_4)\,\cap\,\mathrm{Ad}_v(\Delta_4)}$ is given by $\mbox{Ad}_u\circ E^{M_4}_{\Delta_4\,\cap\,\mathrm{Ad}_{u^*v}(\Delta_4)}\circ\mbox{Ad}_{u^*}=\mbox{Ad}_u\circ E^{M_4}_{\bbc^3}\circ \mbox{Ad}_{u^*}$, where $E^{M_4}_{\bbc^3}:(a_{ij})\longmapsto\mbox{diag}\big\{a_{11},\frac{1}{2}(a_{22}+a_{44}),a_{33},\frac{1}{2}(a_{22}+a_{44})\big\}$. Hence, the matrix of the linear operator $E^{M_4}_{\bbc^3}$ acting on $M_4$, when viewed as a linear map from $\bbc^{16}$ to $\bbc^{16}$, with respect to the basis $\mathfrak{B}$ is a $16\times 16$ matrix having zero everywhere else except the following matrix
\[
\widetilde{B}:=\left[{\begin{matrix}
1 & 0 & 0 & 0\\
0 & \frac{1}{2} & 0 & \frac{1}{2}\\
0 & 0 & 1 & 0\\
0 & \frac{1}{2} & 0 & \frac{1}{2}\\
\end{matrix}}\right]
\]
at the top left $4\times 4$ block.

Therefore, we see that the quadruple in \Cref{as} is a commuting square if and only if the matrix $B-\widetilde{B}=0$. This happens if and only if $t=0$, and consequently $b=\pm ia$.

Now, suppose that this condition is not satisfied. The spectrum of the following positive operator
\[
E^{M_4}_{\mbox{Ad}_u(\Delta_4)}E^{M_4}_{\mbox{Ad}_v(\Delta_4)}E^{M_4}_{\mbox{Ad}_u(\Delta_4)}-E^{M_4}_{\mathrm{Ad}_u(\Delta_4)\,\cap\,\mathrm{Ad}_v(\Delta_4)}
\]
acting on $M_4$, which is same as the spectrum of the operator $E^{M_4}_{\Delta_4}E^{M_4}_{\mbox{Ad}_{u^*v}(\Delta_4)}E^{M_4}_{\Delta_4}-E^{M_4}_{\bbc^3}$, simply becomes the spectrum of the $4\times 4$ matrix $B-\widetilde{B}$. This is because $\mbox{Ad}_u\in\mathcal{B}(M_4)$ is a unitary. Finally, observe that the spectrum of $B-\widetilde{B}$ is the union of $\{0\}$ and the spectrum of the $2\times 2$ matrix $\frac{t^2}{8}(I_2-\sigma_1)$, which is $\{0,t^2/4\}$. Therefore, the cosine of the angle between $\mbox{Ad}_u(\Delta_4)\mbox{ and }\mbox{Ad}_v(\Delta_4)$ is the set $\{0,|t|/2\}$. Since $t=\overline{a}b+a\overline{b}$, the proof is completed.\qed
\end{prf}

Therefore,  the quadruple $(\mbox{Ad}_u(\bbc^3)\subset \mbox{Ad}_u(\Delta_4),\mbox{Ad}_u(\Delta_4)\subset M_4)$ in \Cref{as} is not a commuting square when $b\neq\pm ia$. Using this, and the techniques of commuting cube introduced in \Cref{Sec 2}, we now obtain a complete characterization for the quadruple of von Neumann algebras $(R_u\cap R_v\subset R_u,R_v\subset R)$ being a commuting square. More precisely, we prove the following.

\begin{thm}\label{final com result}
(Characterization of commuting square) Given two $4\times 4$ Hadamard inequivalent complex Hadamard matrices $u(a)$ and $u(b)$ parametrized by the circle parameters $a\mbox{ and }b$, the quadruple of von Neumann algebras
\[
\begin{matrix}
R_{u(b)} &\subset & R\\
\cup & & \cup\\
R_{u(a)}\cap R_{u(b)} &\subset & R_{u(a)}\\
\end{matrix}
\]
is a commuting square if and only if $(\overline{a}b)^4=1$, that is $b=\pm ia$.
\end{thm}
\begin{prf}
We have the following commuting squares
\begin{IEEEeqnarray}{lCl}\label{las}
\begin{matrix}
R_u &\subset & R\\
\cup & & \cup \\
u\Delta_4u^* &\subset & M_4 \\
\end{matrix}\qquad\mbox{and}\qquad \begin{matrix}
R_v &\subset & R\\
\cup & & \cup \\
v\Delta_4v^* &\subset & M_4 \\
\end{matrix}\,.
 \end{IEEEeqnarray}
by the construction of $R_u,R_v\subset R$. Recall from \Cref{even intersection} that $u\Delta_4u^*\cap v\Delta_4v^*=\mbox{Ad}_u(\bbc^3)$. Obtain the cube as described in \Cref{com40},
\begin{figure}[!h]
\begin{center} \resizebox{7 cm}{!}{\includegraphics{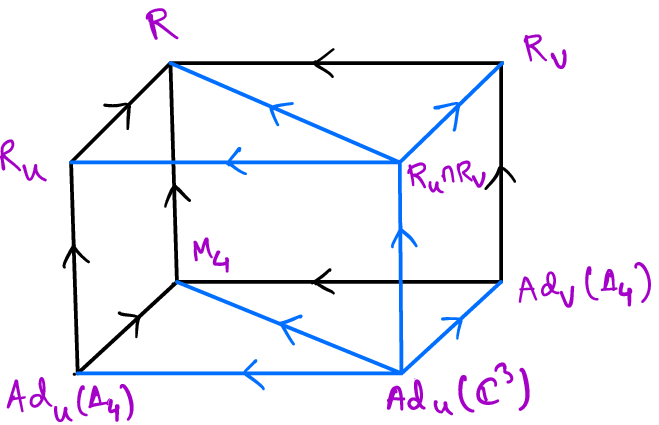}}\end{center}
\caption{Commuting cube in $4\times 4$}\label{com40}
\end{figure}
where the commuting squares in \Cref{las} are the adjacent faces. By \Cref{comm cube}, the slice
\[
\begin{matrix}
R_u\cap R_v &\subset & R\\
\cup & & \cup \\
\mbox{Ad}_u(\bbc^3) &\subset & M_4 \\
\end{matrix}
\]
is a commuting square. Hence, by \Cref{comm cube 2} we get that if the roof $(R_u\cap R_v\subset R_u,R_v\subset R)$ is a commuting square, then the floor $(\mbox{Ad}_u(\bbc^3)\subset u\Delta_4u^*,v\Delta_4v^*\subset M_4)$ must be also so. By \Cref{p4}, the quadruple $(\mbox{Ad}_u(\bbc^3)\subset u\Delta_4u^*,v\Delta_4v^*\subset M_4)$ is not a commuting square when $b\neq\pm ia$ (equivalently $(\overline{a}b)^4=1$), as there is a non-trivial angle between them. Therefore, when $b\neq\pm ia$ the quadruple $(R_u\cap R_v\subset R_u,R_v\subset R)$ is not a commuting square, or in other words, if the quadruple $(R_u\cap R_v\subset R_u,R_v\subset R)$ is a commuting square, then the condition $b=\pm ia$ is necessary. By \Cref{final cor 0}, we see that this is also a sufficient condition, which completes the proof.\qed
\end{prf}

\begin{thm}\label{dihedralangle}
Let $u=u(a)$ and $v=u(b)$ be distinct $4\times 4$ inequivalent complex Hadamard matrices parametrized by the circle parameters $a\mbox{ and }b$ such that $u\sim v$ with $(\overline{a}b)^m=1,\,m\in 2\bbn+2$. Then, we have the following.
\begin{enumerate}[$(i)$]
\item If $\,m=4$, then $\mathrm{Ang}_R(R_u,R_v)=\{\pi/2\};$
\item If $\,m\geq 6$, then $\mathrm{Ang}_{R}(R_u,R_v)=\left\{\frac{2k\pi}{m}\,:\,k=1,2,\ldots, \lfloor\frac{m-2}{4}\rfloor\right\}$.
\end{enumerate}
\end{thm}
\begin{prf}
Part $(i)$ is an immediate consequence of \Cref{final com result}. For part $(ii)$, first observe that $\text{Ang}_{R}(R_u,R_v)=\text{Ang}_{R_u\vee R_v}(R_u,R_v)$ (see Lemma $3.1$ in \cite{SW}). The rest follows from \Cref{intermediate} together with Theorem $6.1$ in \cite{SW} (putting $n=m/2$), since $[R_u\vee R_v:R_u\cap R_v]=m$ by Theorem (\ref{producing}, \ref{lop}).\qed
\end{prf}

\begin{rmrk}\rm
Although the subfactor $R_u\cap R_v\subset R$ becomes a composition of group-type subfactors $(R_u\vee R_v)^{D_n}\subset R_u\vee R_v\subset (R_u\vee R_v)\rtimes\bbz_2$ (see \cite{BH}), where $n=\frac{m}{2}$ and $D_n$ is the dihedral group of order $2n$, it does not help in characterization as it involves the unknown factor $R_u\vee R_v$ instead of $R$.
\end{rmrk}


\subsubsection{Connes-St{\o}rmer relative entropy}

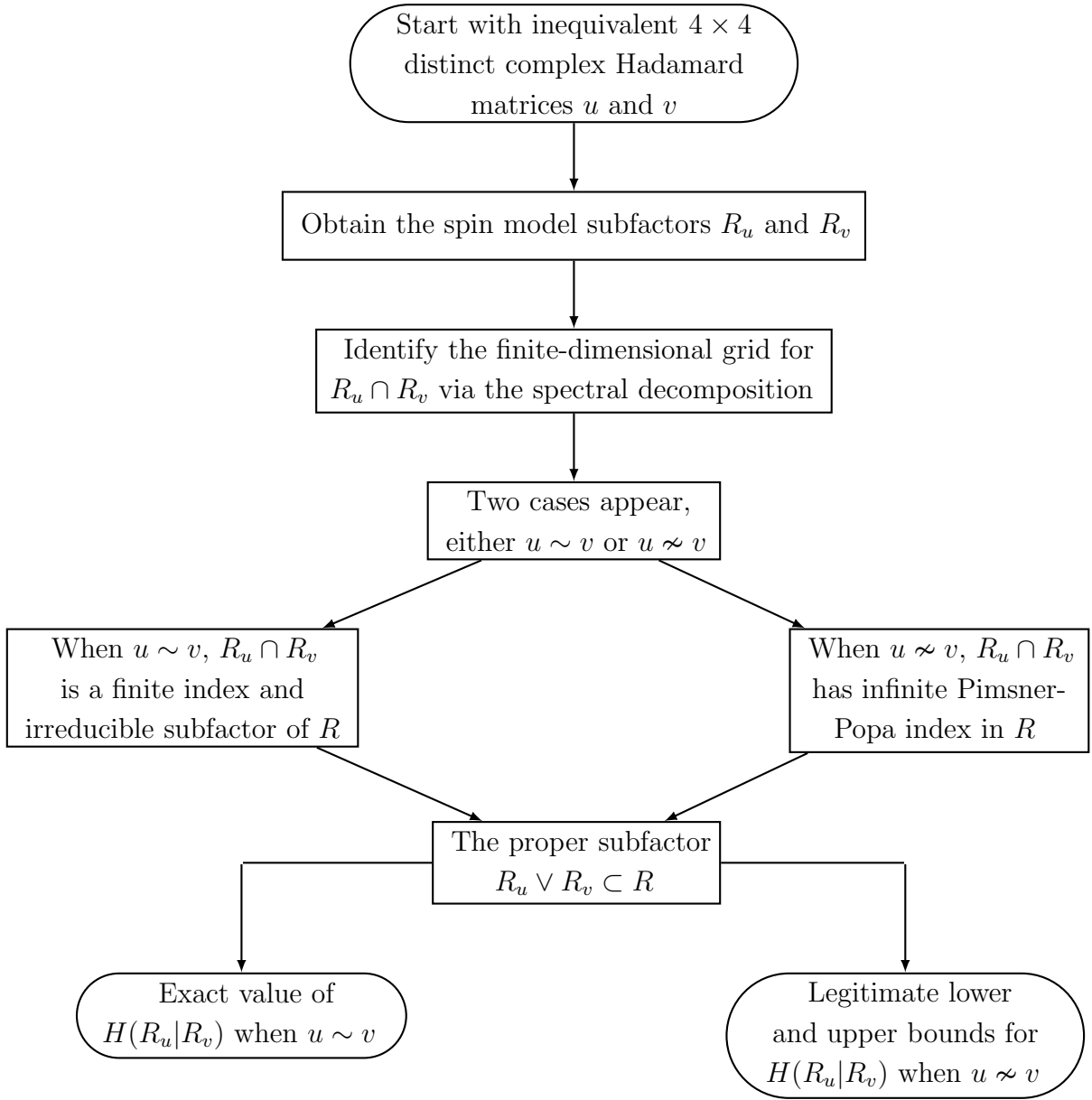
\begin{figure}
\centering
\begin{tikzpicture}[font=\large,thick]

\node[draw,
    align=center,
    rounded rectangle,
    minimum width=2.5cm,
    minimum height=1cm] (block0) {~Start with inequivalent $4\times 4$ ~\\
    distinct complex Hadamard\\
    matrices $u\mbox{ and }v$};
    
\node[draw,
    align=center,
    below=of block0,
    minimum width=2.5cm,
    minimum height=1cm] (block1) {~Obtain the spin model subfactors $R_u\mbox{ and }R_v\,$};
      
\node[draw,
    align=center,
    below=of block1,
    minimum width=3.5cm,
    minimum height=1cm] (block2) {~Identify the finite-dimensional grid for~\\
    $\,R_u\cap R_v$ via the spectral decomposition ~};
      
\node[draw,
    align=center,
    below=of block2,
    minimum width=3.5cm,
    minimum height=1cm] (block3) {~Two cases appear, \\
    ~either $u\sim v$ or $u\nsim v\,$};
    
\node[draw,
    align=center,
    below left=of block3,
    minimum width=3.5cm,
    minimum height=1cm] (block4) {~When $u\sim v,\,R_u\cap R_v\,$\\
    is a finite index and\\
    ~irreducible subfactor of $R$ ~};
    
\node[draw,
    align=center,
    below right=of block3,
    minimum width=3.5cm,
    minimum height=1cm] (block5) {~When $u\nsim v,\,R_u\cap R_v\,$\\
    ~has infinite Pimsner- ~\\
   Popa index in $R$ ~};

\node[draw,
    align=center,
    below right=of block4,
    below left=of block5,
    minimum width=3.5cm,
    minimum height=1cm] (block6) {~The proper subfactor\\ $R_u\vee R_v\subset R\,$};
    
\node[draw,
    rounded rectangle,
    align=center,
    below right=of block6,
    minimum width=3.5cm,
    minimum height=1cm] (block7) {~Legitimate lower \\
    and upper bounds for \\
    $H(R_u|R_v)$ when $u\nsim v\,$};
    
\node[draw,
    rounded rectangle,
    align=center,
    below left=of block6,
    minimum width=3.5cm,
    minimum height=1cm] (block8) {~Exact value of \\
    $H(R_u|R_v)$ when $u\sim v\,$};

\draw[-latex] (block0) edge (block1)
    (block1) edge (block2)
    (block2) edge (block3)
    (block3) edge (block4)
    (block3) edge (block5)
    (block4) edge (block6)
    (block5) edge (block6);
    
\draw[-latex] (block6) -| (block7)
    node[pos=0.5,fill=white,inner sep=0]{};   

\draw[-latex] (block6) -| (block8)
    node[pos=0.5,fill=white,inner sep=0]{};  
        
\end{tikzpicture}
\medskip

\caption{Roadmap for relative entropies}\label{entropy road map}
\end{figure}

In this section, we investigate various relative entropies for the following quadruples
\[
\begin{matrix}
R_u &\subset & R\\
\cup &  & \cup\\
R_u\cap R_v &\subset & R_v
\end{matrix}\quad\mbox{ if }u\nsim v\,;\qquad\qquad \begin{matrix}
 & & R_u &\subset & R\\
 & &\cup &  & \cup\\
 & & R^{(m)}_{u,v}=R_u\cap R_v &\subset & R_v
\end{matrix}\quad\mbox{ if }u\sim v.
\]
Most importantly, we obtain the exact value of $H(R_u|R_v)$ when $u\sim v$. Note that by \Cref{final com result}, it is clear that Theorem $6$ in \cite{Wa} is not applicable when $b\neq\pm ia$ to find the value of $H(R_u|R_v)$. However, it readily follows from \cite{PP} that $H(R|R_u)=H(R|R_v)=\log[R:R_u]=\log[R:R_v]=\log 4$, as spin model subfactors are irreducible. We give a pictorial roadmap in \Cref{entropy road map} to describe how we have reached at the stage of obtaining the exact value of relative entropy $H(R_u|R_v)$ when $u\sim v$.

\begin{ppsn}\label{555}
For the pair of spin model subfactors $R_u, R_v\subset R$, if $u\sim v$ with $(\overline{a}b)^m=1$, then we have the following,
\begin{enumerate}[$(i)$]
\item $H(R_u|R_u\cap R_v)=H(R_v|R_u\cap R_v)=\log m-\log 2\,,$
\item $H(R|R_u\cap R_v)=\log 2m$.
\end{enumerate}
\end{ppsn}
\begin{prf}
Let $u=u(a),\,v=u(b)$ and $\theta=\overline{a}b$. Since $u\sim v$, there exists $m\in2\bbn+2$ such that $\,\theta^{\,m}=1$, and we have the subfactor $R^{(m)}_{u,v}\subset R$ constructed in \Cref{producing} such that $[R:R^{(m)}_{u,v}]=2m$. Recall from Theorem $4.4$ and Corollary $4.5$ in \cite{PP} that for a finite index subfactor $N\subset M$, we have $H(M|N)=\log[M:N]$ if and only if the subfactor is extremal. Since $R^{(m)}_{u,v}\subset R$ is constructed through the basic construction, it is extremal. Therefore, $H(R|R^{(m)}_{u,v})=\log 2m$ and $H(R_u|R^{(m)}_{u,v})=\log[R_u:R^{(m)}_{u,v}]=\log\frac{2m}{4}$. Since $R^{(m)}_{u,v}=R_u\cap R_v$ by \Cref{lop}, we are done.\qed
\end{prf}

Now, recall the subfactor $R_u\vee R_v\subset R$ from \Cref{intermediate}. 

\begin{lmma}\label{zpb}
We have $H(R|R_u\vee R_v)=H(R_u\vee R_v|R_u)=H(R_u\vee R_v|R_v)=\log 2$. Moreover, when $u=u(a)$ and $v=u(b)$ are such that $(\overline{a}b)^m=1$ for some $m\in 2\bbn+2$ (that is, $u\sim v$), we have $H\big(R_u\vee R_v|R_u\cap R_v\big)=\log m$.
\end{lmma}
\begin{prf}
Recall that for any irreducible subfactor $N\subset M$, we have $H(M|N)=\log[M:N]$ by \cite{PP}. All the subfactors $R_u\vee R_v\subset R,\,R_u\subset R_u\vee R_v,\,R_v\subset R_u\vee R_v$ are of index $2$, and hence irreducible. This proves the first part. Now, given $u=u(a)\mbox{ and }v=u(b)$ recall that if $(\overline{a}b)^m=1$ for some $m\in 2\bbn+2$, then $R_u\cap R_v=R^{(m)}_{u,v}$ is an irreducible subfactor of $R$ with $[R:R_u\cap R_v]=2m$ (\Cref{producing}, \ref{interesting} and \ref{lop}). Hence, $[R_u\vee R_v:R_u\cap R_v]=m$. Since $R_u\cap R_v\subset R_u\vee R_v$ also becomes irreducible, we have $H\big(R_u\vee R_v|R^{(m)}_{u,v}=R_u\cap R_v\big)=\log m$.\qed
\end{prf}

\begin{thm}\label{entropy for 4 by 4}
For the pair of spin model subfactors $R_u, R_v\subset R$, we have the following.
\begin{enumerate}[$(i)$]
\item If $\,u\sim v$, then $H(R_u|R_v)=h(R_u|R_v)=\log 2;$
\item If $\,u\nsim v$, then
\[
0<\log 2+\frac{1}{8}\left(\eta|1+a\overline{b}|^2+\eta|1-a\overline{b}|^2\right)\leq H(R_u|R_v)\leq\log 2\,.
\]
\end{enumerate}
\end{thm}
\begin{prf}
For given $u\mbox{ and }v$, existence of the intermediate subfactor $R_u\vee R_v$ containing both $R_u\mbox{ and }R_v$ immediately gives us the following,
\begin{IEEEeqnarray}{lCl}\label{poq}
H(R_u|R_v)\leq H(R_u|R_u\vee R_v)+H(R_u\vee R_v|R_v)=H(R_u\vee R_v|R_v)=\log 2
\end{IEEEeqnarray}
by \Cref{zpb}.

If $\,u\sim v$ such that $\theta=\overline{a}b$ satisfies $\theta^{\,m}=1$ for some $m\in2\bbn+2$, then consider the quadruple $\big(R^{(m)}_{u,v}=R_u\cap R_v\subset R_u,R_v\subset R\big)$. By \Cref{interesting}, we have $R^{(m)}_{u,v}\subset R$ irreducible. Hence, by \Cref{B} we get that $h(R_u|R_v)=-\log\lambda(R_u,R_v)$ and $H(R_u|R_v)\geq-\log\lambda(R_u,R_v)$. By \Cref{Popa for 4 by 4}, we have $\lambda(R_u,R_v)=1/2$. Part $(i)$ now follows from \Cref{poq}.

For part $(ii)$, we only need to provide the lower bound. Note that $H(R_u|R_v)$ is limit of an increasing sequence where the first term is $H\left(\mbox{Ad}_u\Delta_4|\mbox{Ad}_v\Delta_4\right)$. We claim the following,
\begin{IEEEeqnarray}{lCl}\label{lower bound}
0<\log 2+\frac{1}{8}\left(\eta|1+a\overline{b}|^2+\eta|1-a\overline{b}|^2\right)\leq H\left(u\Delta_4u^*|v\Delta_4v^*\right)\,.
\end{IEEEeqnarray}
Consider the partition $\gamma=\{uE_{11}u^*,\ldots,uE_{44}u^*\}$ of $I_4\in M_4$ consisting of minimal projections. Then, by definition of relative entropy and using the fact that for $p\in(M_4)_+$ we have $\eta(p)=0$ if and only if $p$ is a projection (see \cite{PP}, for instance), we have the following,
\begin{IEEEeqnarray*}{lCl}
H_\gamma\left(u\Delta_4u^*|v\Delta_4v^*\right) &=& \sum_{j=1}^4\tau\eta E_{v\Delta_4v^*}^{M_4}\left(uE_{jj}u^*\right)-\tau\eta E_{u\Delta_4u^*}^{M_4}\left(uE_{jj}u^*\right)\\
&=& \sum_{j=1}^4\tau\eta\big(\mbox{Ad}_vE_{\Delta_4}^{M_4}\left(v^*uE_{jj}u^*v\right)\big)-\tau\eta\big(\mbox{Ad}_u E_{\Delta_4}^{M_4}\left(E_{jj}\right)\big)\\
&=& \sum_{j=1}^4\tau\eta\big(E_{\Delta_4}^{M_4}\left(v^*uE_{jj}u^*v\right)\big)\\
&=& \sum_{i,j=1}^4\frac{1}{4}\eta\left(|(v^*u)_{ij}|^2\right)\,,
\end{IEEEeqnarray*}
where $(v^*u)_{ij}$ denotes the $i,j$-entry of the matrix $v^*u$. Since $v^*u=p+q^*$, we get the following (recall few properties of $\eta$ from the beginning of section $3$ in \cite{PP}, for instance),
\begin{IEEEeqnarray*}{lCl}
H_\gamma\left(u\Delta_4u^*|v\Delta_4v^*\right) &=& \frac{1}{2}\left(\eta\left(\frac{1}{4}|1+a\overline{b}|^2\right)+\eta\left(\frac{1}{4}|1-a\overline{b}|^2\right)\right)\\
&=& \frac{1}{2}\eta\left(\frac{1}{4}\right)\left(|1+a\overline{b}|^2+|1-a\overline{b}|^2\right)+\frac{1}{8}\left(\eta|1+a\overline{b}|^2+\eta|1-a\overline{b}|^2\right)\\
&=& \log 2+\frac{1}{8}\left(\eta|1+a\overline{b}|^2+\eta|1-a\overline{b}|^2\right)\,.
\end{IEEEeqnarray*}
Finally, obeserve that
\begin{IEEEeqnarray*}{lCl}
\frac{1}{8}\left(\eta|1+a\overline{b}|^2+\eta|1-a\overline{b}|^2\right)&\geq& \frac{1}{8}\eta\left(|1+a\overline{b}|^2+|1-a\overline{b}|^2\right)=\frac{1}{8}\eta(4)=-\log 2\,.
\end{IEEEeqnarray*}
The inequality above is in fact strict. To see this, take $a\overline{b}=e^{i\gamma}$ for some $\gamma\in(-\pi,\pi)$, and observe that $\eta(1+\cos\gamma)+\eta(1-\cos\gamma)>\eta(2)$ (since in our situation $b\neq\pm a$). Since $H=\sup_\gamma\, H_\gamma$, this gives us the required estimate in \Cref{lower bound}.\qed
\end{prf}

Note that $\eta |1+a\overline{b}|^2+\eta|1-a\overline{b}|^2$ is always negative. To see this, take $a\overline{b}=e^{i\gamma}\in\mathbb{S}^1$ and observe that $\eta|1+a\overline{b}|^2+\eta|1-a\overline{b}|^2=-4\log 2+2\eta(1+\cos\gamma)+2\eta(1-\cos\gamma)$. On contrary, if we assume that there exists $\gamma$ such that $2\eta(1+\cos\gamma)+2\eta(1-\cos\gamma)\geq 4\log 2$, then one arrives at the stage $\eta\big(\cos^2\frac{\gamma}{2}\big)+\eta\big(\sin^2\frac{\gamma}{2}\big)\geq\log 4$. However, the global maxima for $\eta$ is $e^{-1}$, and therefore the left hand side can be at the best $2e^{-1}$, which is strictly smaller than $\log 4$. Hence, no such $\gamma$ can exist and the quantity $\eta |1+a\overline{b}|^2+\eta|1-a\overline{b}|^2$ is always negative.

\begin{crlre}\label{px}
If $u\sim v$, then for the pair of spin model subfactors $R_u, R_v\subset R$, we have $H(R_u|R_v)=-\log\lambda(R_u,R_v)$.
\end{crlre}
\begin{prf}
Combine \Cref{Popa for 4 by 4} and \Cref{entropy for 4 by 4}.\qed
\end{prf}

We summarize all the relative entropies $H(\,.\,|\,.\,)$ when $u\sim v$ (i,e. when $\theta\in\Gamma$) in a $5\times 5$ matrix form in \Cref{table2}. 

\begin{table}[ht]
\centering
\begin{tabular}{|c| c| c| c| c| c| }
\hline\hline
$\,H(\,.\,|\,.\,)\,\mbox{for }u\sim v\,$ & $\,\,R\,\,$ & $R_u$ & $R_v$ & $R^{(m)}_{u,v}=R_u\cap R_v$ & $R_u\vee R_v$ \\ [0.5ex]
\hline\hline
  $R$ & $0$ & $\log 4$ & $\log 4$ & $\log 2+\log m$ & $\log 2$ \\
  $R_u$  & 0 & $0$ & $\log 2$ & $\log m-\log 2$ & $0$ \\
  $R_v$ & 0 & $\log 2$ & $0$  & $\log m-\log 2$ & $0$ \\
  $R^{(m)}_{u,v}=R_u\cap R_v$  & $0$ & $0$ & $0$ & $0$ & $0$ \\
  $R_u\vee R_v$  & $0$ & $\log 2$ & $\log 2$ & $\log m$ & $0$ \\
\hline\hline
\end{tabular}
\vspace{2mm}
\caption{Relative entropies when $u\sim v$}\label{table2}
\end{table}


\section{Concluding remarks: open questions and perspective}\label{Sec 6}

In this concluding section, we briefly summarize our findings and discuss possible future directions. Consider a pair of subfactors of type $II_1$ factors :
\[
\begin{matrix}
P &\subset & M \cr
 & & \cup\cr
 &  & Q \end{matrix}
\]
with both $[M:P]\,,\,[M:Q]<\infty$. This work is guided by our attempt to answer  the following question~:
\smallskip

\noindent\textbf{Open Problem 1:} Given a pair of subfactors $P,Q\subset M$, provide formulae to compute the invariants $\lambda(P,Q),\,\mbox{Ang}_M(P,Q)$, and $H(P|Q)$.
\smallskip

This appears to be a  difficult problem. There are  two major obstacles in attacking this problem~:
\begin{enumerate}[$(i)$]
\item It may very well happen that $P\cap Q$ is not a factor.
\item The index of $[M:P\cap Q]$  may be infinite.
\end{enumerate}
To attack problem $1$ stated above, it seems reasonable to answer first the following particular case.
\smallskip

\noindent\textbf{Open Problem 2:} Consider a quadruple of $II_1$ factors $(N\subset P,Q\subset M)$ with $[M:N]<\infty$. How to compute the invariants $\,\lambda(P,Q),\,\mbox{Ang}_M(P,Q),\,\alpha^N_M(P,Q)$, $\beta^N_M(P,Q),\,H(P|Q)?$
\smallskip

One may expect that $H$ is related to $\lambda$ (for example, \Cref{px}). Only the formula for $\lambda(P,Q)$ in this irreducible situation is known due to \cite{B}. Even if $P\cap Q$ is a factor with $[M:P\cap Q]<\infty$, and $P\subset M,\,Q\subset M$ both are irreducible, then also $P\cap Q \subset M$ may not be irreducible (see \Cref{relativecommutant}, for instance). In view of the fact that $\lambda$ and $H$ are well-behaved under `controlled' limit of finite dimensional grid, it seems important to know the value of them in finite dimension. However, we stumbled upon the following problem.
\smallskip
 
\noindent\textbf{Open problem 3:} Given two finite-dimensional von Neumann-subalgebras $\mathcal{B}$ and $\mathcal{C}$ of a finite-dimensional von Neumann algebra $\mathcal{A}$, can we obtain formulae for $\lambda(\mathcal{B,C})$ and $H(\mathcal{B}|\mathcal{C})$ similar to Section 6 of \cite{PP}?
\smallskip

The following simpler version of this problem is open.
\smallskip

\noindent\textbf{Open Problem 4:} (M. Choda and D. Petz) What is the formula for $H(\Delta_n | u\Delta_n u^*)$, where $u\in M_n(\mathbb{C})$ is a unitary matrix?
\smallskip

To tackle Problem 1, it is evident that we need to first compute the invariants for important cases in order to gain insight. Motivated by this goal, in this article we have considered a pair $(u,v)$ of distinct complex Hadamard matrices, and attempted to compute these invariants for the corresponding pair of spin model subfactors $R_u,R_v\subset R$. A complete characterization of $R_u\neq R_v$ has been obtained. Given distinct $2\times 2$ complex Hadamard matrices $u$ and $v$, we have obtained a quadruple of $II_1$ factors which has been completely characterized, and the invariants are computed. The intersection of the corresponding spin model subfactors turns out to be a vertex model subfactor of index $4$. The relative position between the spin model subfactors in the $2\times 2$ situation is completely understood. However, the Hadamard inequivalent $4\times 4$ case is more involved and quite surprising.  Below we briefly explain the $4\times 4$ case in a more precise manner. 
\smallskip

\begin{figure}[h!]
\begin{center}
		\resizebox{4cm}{!}{\includegraphics{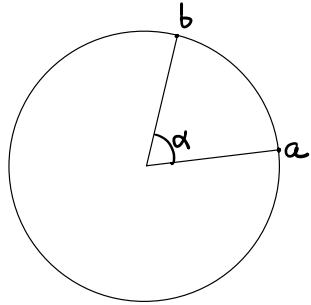}}
		\caption{Angle between circle parameters}\label{fig4}
		\end{center}
\end{figure}

Consider a pair $(u(a), u(b))$ of $4\times 4$ Hadamard inequivalent distinct complex Hadamard matrices parametrized by the circle parameters $a,b\in\mathbb{S}^1$ (semicircle to be precise). Note that if $a=e^{i\phi}$ and $b=e^{i\psi}$, then $\phi,\psi\in[0,\pi)$ (so that $\phi-\psi\in(-\pi,\pi)$). We obtain a pair of spin model subfactors $R_{u(a)}\subset R$ and $R_{u(b)}\subset R$. Let $\alpha$ be the angle between the two circle parameters $a$ and $b$ as described in \Cref{fig4}. Then, the following two cases arise.

\noindent  \textbf{Case I:} $\alpha=2k\pi/m$ for some $m\in 2\mathbb{N}+2$ and $1\leq k\leq m-1$ with $k\neq m/2$. In this case, we have $R_{u(a)}\cap R_{u(b)}$ is a factor and $[R: R_{u(a)}\cap R_{u(b)}]=2m$. Also, $H(R_{u(a)}|R_{u(b)})= H(R_{u(b)}| R_{u(a)})=\log 2$.
\smallskip

In this case, there are two major ingredients. The first one is proving that the intersection of $R_{u(a)}\cap R_{u(b)}$ is a factor, for which a bit of work is needed. We remark that we thereby obtain a series of potentially new finite index and irreducible subfactors of $R$. The second ingredient is proving that $R_{u(a)}\vee R_{u(b)}$ is a subfactor of $R$ with index $2$.
\smallskip

\noindent\textbf{Open problem 5:} What is the principal graph of the family of subfactors $R_{u(a)}\cap R_{u(b)}\subset R$? Can we describe the planar algebra?
\smallskip

\noindent \textbf{Case II:} $\alpha\neq 2k\pi/m$ for any $m\in 2\mathbb{N}+2$ with $1\leq k\leq m-1$. In this case $R_{u(a)}\cap R_{u(b)}$ has infinite Pimsner-Popa index in $R$.
\smallskip

\noindent\textbf{Open Problem 6:} For  Case II, the following problems seem to be of interest.
\begin{enumerate}[$(i)$]
\item Is it true that $\,R_{u(a)}\cap R_{u(b)}$ is a factor? 
\item What is the relative commutant ${(R_{u(a)}\cap R_{u(b)})}^{\prime}\cap R\,$? 
\item How to compute $\text{Ang}_R(R_{u(a)},R_{u(b)})$?
\item What are the values of relative entropies $H(R_{u(a)}|R_{u(b)})$ and $H(R|R_{u(a)}\cap R_{u(b)})?$ 
\end{enumerate}
 
We remark that we are only able to provide bounds for $H(R_{u(a)}|R_{u(b)})$ at present. To tackle Problem $6$, it seems to us that probably a vast generalization of the Ocneanu's compactness theorem beyond symmetric commuting squares is required, which is out of reach at present with the help of existing tools in the literature.
\smallskip
  
For a general  pair of $n\times n$ complex Hadamard matrices, we are not yet able to compute these invariants.  Our work in this article leads to many interesting open problems  that should lay the foundation for developing a general theory of two subfactors and non-commutaive entropy à la Connes-St{\o}rmer.  We hope to tackle some of these problems in years to come.


\section{Appendix}\label{Sec 7}

We provide a proof of the Jones' basic construction for the unital inclusion $\bbc\subset u\Delta_nu^*$ for $n=2\mbox{ and }4$. Although tower of spin model basic construction is well-known in the literature, our tower  is slightly different than the one in \cite{Jo3} as it is more convenient for our purpose. The Jones' projection for the inclusion $\bbc\subset\Delta_n$ will be denoted by $J_n=\frac{1}{n}\sum_{i,j}E_{ij}\in M_n$, and that for the inclusion $\Delta_n\subset M_n$ will be denoted by $e_1=\mbox{diag}\{E_{11},\ldots,E_{nn}\}\in\Delta_n\otimes M_n$.

\subsection{Basic construction for the \texorpdfstring{$2\times 2$}~ case}
Let $u=(u_{ij})_{1\leq i,j\leq 2}$ be an arbitrary complex Hadamard matrix in $M_2$ and consider the following matrices
\[
\eta=\left[{\begin{matrix}
\frac{u_{12}}{u_{11}} & 0\\
0 & \frac{u_{22}}{u_{21}}\\
\end{matrix}}\right]\quad \mbox{and}\quad \xi=\frac{1}{\sqrt{2}}\,\left[{\begin{matrix}
1 & \frac{u_{12}}{u_{11}}\\
1 & \frac{u_{22}}{u_{21}}\\
\end{matrix}}\right]\,,
\]
associated with $u$. Since $u$ is unitary, it is easy to check that $\,\eta\mbox{ and }\xi$ are unitary matrices.

\begin{lmma}\label{basis step 1}
Consider the Jones' projection $e_1=\mathrm{diag}\{E_{11},E_{22}\}\in\Delta_2\otimes M_2$ for the inclusion $\,\Delta_2\subset M_2$. The $\,4\times 4$ unitary matrix $\,u_1\,$ in $\Delta_2\otimes M_2$ satisfying $(u\Delta_2u^*)\,e_1\,(u\Delta_2u^*)=u_1M_2u_1^*$ is given by $u_1=\mbox{bl-diag}\{u,u\eta\}$.
\end{lmma}
\begin{prf}
Observe that $\Delta_2u^*\,e_1\,u\Delta_2=\mbox{bl-diag}\{\Delta_2u^*E_{11}u\Delta_2\,,\,\Delta_2u^*E_{22}u\Delta_2\}$. Using the fact that $u$ is a complex Hadamard matrix, i,e. all the entries of $u$ are of same modulus, and $\overline{u_{21}}u_{22}=\displaystyle{\frac{\overline{u_{11}u_{22}}\,u_{12}^2}{\overline{u_{21}}\,u_{11}}}$ because $u_{11}\overline{u_{21}}+u_{12}\overline{u_{22}}=0$, it is easy to check that $\eta u^*E_{11}u\eta^*=u^*E_{22}u\,$. Hence, $\Delta_2u^*E_{22}u\Delta_2=\eta\big(\Delta_2u^*E_{11}u\Delta_2\big)\eta^*$. This gives us the following,
\begin{center}
$(I_2\otimes u)(I_2\otimes\Delta_2)(I_2\otimes u^*)\,\mbox{bl-diag}\{E_{11},E_{22}\}\,(I_2\otimes u)(I_2\otimes\Delta_2)(I_2\otimes u^*)=u_1(I_2\otimes M_2)u_1^*$
\end{center}
with $u_1$ as the required matrix if we show that $\mbox{Alg}\,\big\{\Delta_2u^*E_{11}u\Delta_2\big\}=M_2\,$. However, this is obvious as $E_{ii}u^*E_{11}uE_{jj}$ gives the matrix unit $E_{ij}$ in $M_2$.\qed
\end{prf}

\begin{lmma}\label{basis step 2}
Consider the Jones' projection $e_2=J_2\otimes I_2\in M_2\otimes M_2$ for the inclusion $\bbc\otimes M_2\subset\Delta_2\otimes M_2$. The $4\times 4$ unitary matrix $\,u_2\,$ in $M_2\otimes M_2$ satisfying $(u_1 M_2u_1^*)\,e_2\,(u_1M_2u_1^*)=u_2(\Delta_2\otimes M_2)u_2^*$ is given by $\,u_2=u_1(\xi^*\otimes I_2)$.
\end{lmma}
\begin{prf}
First observe that
\[
u_1^*(J_2\otimes I_2)u_1=\left[{\begin{matrix}
I_2 & \eta\\
\eta^* & I_2\\
\end{matrix}}\right]\,.
\]
Then, using the fact that $\xi$ is unitary we have the following,
\begin{IEEEeqnarray*}{lCl}
&  & u_1(I_2\otimes M_2)u_1^*(J_2\otimes I_2)u_1(I_2\otimes M_2)u_1^*\\
&=& \mbox{bl-diag}\{u,u\eta\}(\xi^*\otimes I_2)(I_2\otimes A)(\xi\otimes I_2)\left[{\begin{matrix}
I_2 & \eta\\
\eta^* & I_2\\
\end{matrix}}\right](\xi^*\otimes I_2)(I_2\otimes B)(\xi\otimes I_2)\,\mbox{bl-diag}\{u^*,\eta^*u^*\}
\end{IEEEeqnarray*}
where $A,B\in M_2(\bbc)$. Now, it is easy to see that $I_2+\frac{u_{22}}{u_{21}}\eta^*+\displaystyle{\frac{\overline{u_{12}}}{\overline{u_{11}}}}\eta+\frac{u_{22}}{u_{21}}\displaystyle{\frac{\overline{u_{12}}}{\overline{u_{11}}}}I_2=0$ which gives us the following,
\begin{IEEEeqnarray*}{lCl}
(\xi\otimes I_2)\left[{\begin{matrix}
I_2 & \eta\\
\eta^* & I_2\\
\end{matrix}}\right](\xi^*\otimes I_2) &=& \frac{1}{4}\,\mbox{bl-diag}\left\{2I_2+\frac{u_{12}}{u_{11}}\eta^*+\frac{\overline{u_{12}}}{\overline{u_{11}}}\eta\,,\,2I_2+\frac{u_{22}}{u_{21}}\eta^*+\frac{\overline{u_{22}}}{\overline{u_{21}}}\eta\right\}\,.
\end{IEEEeqnarray*}
Since $\frac{1}{4}\big(2I_2+\frac{u_{12}}{u_{11}}\eta^*+\displaystyle{\frac{\overline{u_{12}}}{\overline{u_{11}}}}\eta\big)=E_{11}$ and $\frac{1}{4}\big(2I_2+\frac{u_{22}}{u_{21}}\eta^*+\displaystyle{\frac{\overline{u_{22}}}{\overline{u_{21}}}}\eta\big)=E_{22}$, and we have
\[
\mbox{Alg}\left\{\mbox{bl-diag}\{AE_{11}B,AE_{22}B\}\,:\,A,B\in M_2\right\}=\Delta_2\otimes M_2\,,
\]
the proof is concluded.\qed 
\end{prf}

\textbf{Proof of Theorem \ref{1st basic}~:} We prove by strong induction on $k$, where Lemma \ref{basis step 1} and \ref{basis step 2} are the basis step for part $(ii)$ and $(i)$ respectively. The Jones' projection for the inclusion $\bbc\otimes M_2^{(k)}\subset\Delta_2\otimes M_2^{(k)}$ is $e_{2k}=J_2\otimes I_2^{(k)}\in M_2\otimes M_2^{(k)}$, and that for $\Delta_2\otimes M_2^{(k)}\subset M_2\otimes M_2^{(k)}$ is $e_{2k+1}=\mbox{diag}\{E_{11},E_{22}\}\otimes I_2^{(k)}\in\Delta_2\otimes M_2\otimes M_2^{(k)}$. First we prove part $(i)$. Assume that the result is true up to $u_{2k-1}$ for some $k\geq 2$, and consider the Fourier matrix $F_2=\frac{1}{\sqrt{2}}\left[{\begin{smallmatrix}
1 & 1\\
1 & -1\\
\end{smallmatrix}}\right]$. Then, we have the following,
\begin{IEEEeqnarray*}{lCl}
&  & u_{2k-1}\left(I_2\otimes M_2^{(k)}\right)u_{2k-1}^*e_{2k}u_{2k-1}\left(I_2\otimes M_2^{(k)}\right)u_{2k-1}^*\\
&=& u_{2k-1}\big(I_2\otimes M_2^{(k)}\big)\,\mbox{bl-diag}\left\{I_2^{(k)},\eta^{(k-1)}\right\}\left(I_2\otimes u_{2k-2}^*\right)\left(J_2\otimes I_2^{(k)}\right)\left(I_2\otimes u_{2k-2}\right)\\
&  & \mbox{bl-diag}\left\{I_2^{(k)},\eta^{(k-1)}\right\}\left(I_2\otimes M_2^{(k)}\right)u_{2k-1}^*\\
&=& u_{2k-1}\left(I_2\otimes M_2^{(k)}\right)\mbox{bl-diag}\left\{I_2^{(k)},\eta^{(k-1)}\right\}\left(J_2\otimes I_2^{(k)}\right)\mbox{bl-diag}\left\{I_2^{(k)},\eta^{(k-1)}\right\}\left(I_2\otimes M_2^{(k)}\right)u_{2k-1}^*\\
&=& u_{2k-1}\left(I_2\otimes M_2^{(k)}\right)\left(\mbox{diag}\{1,1,1,-1\}\otimes I_2^{(k-1)}\right)\left(J_2\otimes I_2\otimes I_2^{(k-1)}\right)\left(\mbox{diag}\{1,1,1,-1\}\otimes I_2^{(k-1)}\right)\\
&  & \left(I_2\otimes M_2^{(k)}\right)u_{2k-1}^*\\
&=& u_{2k-1}\left(I_2\otimes M_2^{(k)}\right)\left(F_2\otimes I_2\otimes I_2^{(k-1)}\right)\left(\frac{1}{2}(I_2\otimes I_2+F_2\sigma_1F_2\otimes\sigma_3)\otimes I_2^{(k-1)}\right)\left(F_2\otimes I_2\otimes I_2^{(k-1)}\right)\\
&  & \left(I_2\otimes M_2^{(k)}\right)u_{2k-1}^*\\
&=& u_{2k-1}\left(F_2\otimes I_2\otimes I_2^{(k-1)}\right)\left(\left(F_2\otimes I_2\otimes I_2^{(k-1)}\right)\left(I_2\otimes M_2^{(k)}\right)\left(F_2\otimes I_2\otimes I_2^{(k-1)}\right)\right)\\
&  & \left(\frac{1}{2}(I_2\otimes I_2+F_2\sigma_1F_2\otimes\sigma_3)\otimes I_2^{(k-1)}\right)\left(\left(F_2\otimes I_2\otimes I_2^{(k-1)}\right)\left(I_2\otimes M_2^{(k)}\right)\left(F_2\otimes I_2\otimes I_2^{(k-1)}\right)\right)\\
&  & \left(F_2\otimes I_2\otimes I_2^{(k-1)}\right)u_{2k-1}^*\\
&=& u_{2k-1}\left(\xi_k\otimes I_2\otimes I_2^{(k-1)}\right)\left((F_2\otimes I_2)(I_2\otimes M_2)(F_2\otimes I_2)\otimes M_2^{(k-1)}\right)\left(\mbox{bl-diag}\{E_{11},E_{22}\}\otimes I_2^{(k-1)}\right)\\
&  & \left((F_2\otimes I_2)(I_2\otimes M_2)(F_2\otimes I_2)\otimes M_2^{(k-1)}\right)\left(\xi_k\otimes I_2\otimes I_2^{(k-1)}\right)^*u_{2k-1}^*\,.
\end{IEEEeqnarray*}
From here, part $(i)$ will follow if the following holds
\begin{IEEEeqnarray}{lCl}\label{xy}
\mbox{Alg}\,\big\{\mbox{Ad}_{F_2\otimes I_2}(I_2\otimes A)\,\mbox{bl-diag}\{E_{11},E_{22}\}\,\mbox{Ad}_{F_2\otimes I_2}(I_2\otimes B)\,:\,A,B\in M_2\big\}=\Delta_2\otimes M_2\,.\qquad
\end{IEEEeqnarray}
Now, the equality in \Cref{xy} holds because $\mbox{Ad}_{F_2\otimes I_2}(I_2\otimes A)=I_2\otimes A$, and 
\[
\mbox{Alg}\left\{\mbox{bl-diag}\{AE_{11}B,AE_{22}B\}:A,B\in M_2\right\}=\Delta_2\otimes M_2
\]
as in Lemma \ref{basis step 2}, which completes part $(i)$.
\smallskip

Now, we prove part $(ii)$. Assume that the result is true up to $u_{2k-1}$ for some $k\geq 2$. We also have $u_{2k}$ in our hand, since part $(i)$ is already proved. First observe that
\begin{IEEEeqnarray*}{lCl}
(I_2\otimes u_{2k}^*)e_{2k+1}(I_2\otimes u_{2k}) &=& \mbox{bl-diag}\left\{u_{2k}^*(E_{11}\otimes I_2^{(k)})u_{2k}\,,\,u_{2k}^*(E_{22}\otimes I_2^{(k)})u_{2k}\right\}\,.
\end{IEEEeqnarray*}
We claim that for $\eta_k=\mbox{diag}\{1,-1\}\otimes I_2^{(k)}\,,$ one has the following,
\begin{IEEEeqnarray}{lCl}\label{required eqn}
u_{2k}^*\left(E_{22}\otimes I_2^{(k)}\right)u_{2k}=\eta_k\left(u_{2k}^*(E_{11}\otimes I_2^{(k)})u_{2k}\right)\eta_k\,.
\end{IEEEeqnarray}
To prove this, we first write $u_{2k}=u_{2k-1}\big(\xi_k\otimes I_2^{(k)}\big)$ using part $(i)$, and then using the induction hypothesis we get that
\begin{IEEEeqnarray*}{lCl}
u_{2k} &=& \left(I_2\otimes u_{2k-2}\right)\left((E_{11}+E_{12})\otimes I_2\otimes I_2^{(k-1)}+(E_{21}-E_{22})\otimes\sigma_3\otimes I_2^{(k-1)}\right)\,.
\end{IEEEeqnarray*}
Now, it is a straightforward verification that Eqn. \ref{required eqn} holds. Therefore, we get the following,
\begin{IEEEeqnarray*}{lCl}
&  & (I_2\otimes u_{2k})\left(I_2\otimes\Delta_2\otimes M_2^{(k)}\right)(I_2\otimes u_{2k})^*e_{2k+1}(I_2\otimes u_{2k})\left(I_2\otimes\Delta_2\otimes M_2^{(k)}\right)(I_2\otimes u_{2k})^*\\
&=& (I_2\otimes u_{2k})\left(I_2\otimes\Delta_2\otimes M_2^{(k)}\right)\,\mbox{bl-diag}\left\{I_2^{(k+1)},\eta_k\right\}\left(I_2\otimes u_{2k}^*(E_{11}\otimes I_2^{(k)})u_{2k}\right)\\
&  & \mbox{bl-diag}\left\{I_2^{(k+1)},\eta_k\right\}\left(I_2\otimes\Delta_2\otimes M_2^{(k)}\right)(I_2\otimes u_{2k})^*\\
&=& (I_2\otimes u_{2k})\left(E_{11}\otimes I_2^{(k+1)}+E_{22}\otimes\eta_k\right)\left(I_2\otimes\left(\Delta_2\otimes M_2^{(k)}\right)u_{2k}^*\left(E_{11}\otimes I_2^{(k)}\right)u_{2k}\left(\Delta_2\otimes M_2^{(k)}\right)\right)\\
&  & \left(E_{11}\otimes I_2^{(k+1)}+E_{22}\otimes\eta_k\right)(I_2\otimes u_{2k})^*\,.
\end{IEEEeqnarray*}
This is equal to $\mbox{Ad}_{u_{2k+1}}(I_2\otimes M_2^{(k+1)})$, where $\,u_{2k+1}=\left(I_2\otimes u_{2k}\right)\left(E_{11}\otimes I_2^{(k+1)}+E_{22}\otimes\eta_k\right)$, if
\[
\mbox{Alg}\,\left\{\left(\Delta_2\otimes M_2^{(k)}\right)u_{2k}^*\left(E_{11}\otimes I_2^{(k)}\right)u_{2k}\left(\Delta_2\otimes M_2^{(k)}\right)\right\}=M_2\otimes M_2^{(k)}\,.
\]
We claim that $\,u_{2k}^*\left(E_{11}\otimes I_2^{(k)}\right)u_{2k}=J_2\otimes I_2^{(k)}$. To see this, observe that
\begin{IEEEeqnarray*}{lCl}
u_{2k} &=& u_{2k-1}\left(F_2\otimes I_2^{(k)}\right)\\
&=& \left(I_2\otimes u_{2k-2}\right)\big(\left((E_{11}+E_{12})\otimes I_2+(E_{21}-E_{22})\otimes\sigma_3\right)\otimes I_2^{(k-1)}\big)\,,
\end{IEEEeqnarray*}
and hence
\[
u_{2k}^*\left(E_{11}\otimes I_2^{(k)}\right)u_{2k}=J_2\otimes I_2\otimes I_2^{(k-1)}=J_2\otimes I_2^{(k)}\,.
\]
Since $\Delta_2 J_2\Delta_2=M_2$, proof of part $(ii)$ is now completed.\qed


\subsection{Basic construction for the \texorpdfstring{$4\times 4$}~ case}

Consider the $4\times 4$ complex Hadamard matrix
\[
u=\frac{1}{2}\,\left[{\begin{matrix}
1 & 1 & 1 & 1\\
1 & ia & -1 & -ia\\
1 & -1 & 1 & -1\\
1 & -ia & -1 & ia\\
\end{matrix}}\right]
\]
where $a\in\mathbb{S}^1$ (semicircle to be precise).
\begin{lmma}\label{basis step 12}
Let $e_1=\mbox{diag}\{E_{11},E_{22},E_{33},E_{44}\}\in\Delta_4\otimes M_4$, which is the Jones' projection for the basic construction of $\,\Delta_4\subset M_4$. The $16\times 16$ unitary matrix $\,u_1$ in $\Delta_4\otimes M_4$ such that $(u\Delta_4u^*)\,e_1\,(u\Delta_4u^*)=u_1M_4u_1^*$ is given by the block-diagonal matrix {\em bl-diag}$\big\{u,u\,W_1^{(a)},u\,W_2^{(a)},u\,W_3^{(a)}\big\}$.
\end{lmma}
\begin{prf}
It is easy to verify the following,
\begin{IEEEeqnarray*}{lCl}
u^*E_{22}u &=& \left(\sigma_3\otimes\mbox{diag}\{i,\overline{a}\}\right)u^*E_{11}u\left(\sigma_3\otimes\mbox{diag}\{-i,a\}\right)=W_1^{(a)}u^*E_{11}uW_3^{(\overline{a})}\,,\\
u^*E_{33}u &=& \left(I_2\otimes\sigma_3\right)u^*E_{11}u\left(I_2\otimes\sigma_3\right)=W_2^{(a)}u^*E_{11}uW_2^{(a)}\,,\\
u^*E_{44}u &=& \left(\sigma_3\otimes\mbox{diag}\{-i,\overline{a}\}\right)u^*E_{11}u\left(\sigma_3\otimes\mbox{diag}\{i,a\}\right)=W_3^{(a)}u^*E_{11}uW_1^{(\overline{a})}\,.
\end{IEEEeqnarray*}
Therefore,
\begin{IEEEeqnarray*}{lCl}
&  & u\Delta_4u^*e_1u\Delta_4u^*\\
&=& (I_4\otimes u)\,\mbox{bl-diag}\left\{I_4,W_1^{(a)},W_2^{(a)},W_3^{(a)}\right\}(I_4\otimes\Delta_4u^*E_{11}u\Delta_4)(I_4\otimes u^*)\,\mbox{bl-diag}\left\{I_4,W_3^{(\overline{a})},W_2^{(a)},W_1^{(\overline{a})}\right\}\\
&=& u_1(I_4\otimes M_4)u_1^*
\end{IEEEeqnarray*}
with $u_1$ in the required form if $\mbox{Alg}\,\big\{\Delta_4u^*E_{11}u\Delta_4\big\}=M_4$, i,e., $\mbox{Alg}\,\big\{\Delta_4(J_2\otimes J_2)\Delta_4\big\}=M_4$. But since $J_4$ is the Jones' projection for $\bbc\subseteq\Delta_4\subseteq M_4$, this is obvious.\qed
\end{prf}

\begin{lmma}\label{basis step 22}
Let $e_2=J_4\otimes I_4\in M_4\otimes M_4$, which is the Jones' projection for the basic construction of $\bbc\otimes M_4\subset\Delta_4\otimes M_4$. The $16\times 16$ unitary matrix $\,u_2\,$ in $M_4\otimes M_4$, where $(u_1 M_4u_1^*)\,e_2\,(u_1M_4u_1^*)=u_2(\Delta_4\otimes M_4)u_2^*\,,$ is given by $\,u_2=u_1(\xi\otimes I_4)$.
\end{lmma}
\begin{prf}
First observe that
\[ u_1^*e_2u_1=\mbox{bl-diag}\left\{u^*,\left(W_1^{(a)}\right)^*u^*,W_2^{(a)}u^*,\left(W_3^{(a)}\right)^*u^*\right\}(J_4\otimes I_4)\,\mbox{bl-diag}\left\{u,uW_1^{(a)},uW_2^{(a)},uW_3^{(a)}\right\}.
\]
and $(\xi\otimes I_4)^*u_1^*e_2u_1(\xi\otimes I_4)=16e_1\,,$ i,e. the unitary $\xi\otimes I_4$ diagonalize $u_1^*e_2u_1$. Therefore,
\begin{IEEEeqnarray*}{lCl}
u_1M_4u_1^*e_2u_1M_4u_1^*&=& u_1\left((I_4\otimes M_4)(\xi\otimes I_4)\right)\left((\xi\otimes I_4)^*u_1^*e_2 u_1(\xi\otimes I_4)\right)\left((\xi\otimes I_4)^*(I_4\otimes M_4)\right)u_1^*\\
&=& u_1(\xi\otimes I_4)\left((I_4\otimes M_4)e_1(I_4\otimes M_4)\right)(\xi\otimes I_4)^*u_1^*\\
&=& u_2\left(\Delta_4\otimes M_4\right)u_2^*
\end{IEEEeqnarray*}
with $u_2=u_1(\xi\otimes I_4)$, if we can show that $\mbox{Alg}\{(I_4\otimes M_4)e_1(I_4\otimes M_4)\}=\Delta_4\otimes M_4$. To conclude the following,
\begin{center}
$\mbox{Alg}\big\{\mbox{bl-diag}\left\{AE_{11}B,AE_{22}B,AE_{33}B,AE_{44}B\right\}\,:\,A,B\in M_4\big\}=\Delta_4\otimes M_4$
\end{center}
observe that taking $A=E_{ik}$ and $B=E_{kj}$ we obtain the matrix $E_{kk}\otimes E_{ij}$ in $\Delta_4\otimes M_4$ for $k,i,j\in\{1,2,3,4\}$.\qed 
\end{prf}

\textbf{Proof of Theorem \ref{tower in four by four}~:}
We prove by strong induction on $k$. Lemma \ref{basis step 12} and \ref{basis step 22} are the basis step for part $(i)$ and $(ii)$ respectively. The Jones' projection for the basic construction of $\bbc\otimes M_4^{(k)}\subset\Delta_4\otimes M_4^{(k)}$ is $e_{2k}=J_4\otimes I_4^{(k)}\in M_4\otimes M_4^{(k)}$, and that for $\Delta_4\otimes M_4^{(k)}\subset M_4\otimes M_4^{(k)}$ is $e_1\otimes I_4^{(k)}\in\Delta_4\otimes M_4\otimes M_4^{(k)}$.

First we prove part $(i)$. Assume that the result is true up to $u_{2k}$. Observe that $e_{2k+1}=e_1\otimes I_4^{(k)}$. Now, we have the following,
\begin{IEEEeqnarray*}{lCl}
&  & \left(I_4\otimes(\Delta_4\otimes M_4^{(k)})\right)(I_4\otimes u_{2k})^*(e_1\otimes I_4^{(k)})(I_4\otimes u_{2k})\left(I_4\otimes(\Delta_4\otimes M_4^{(k)})\right)\\
&=& \left(I_4\otimes(\Delta_4\otimes M_4^{(k)})\right)\left(I_4\otimes (\xi\otimes I_4^{(k)})^*u_{2k-1}^*\right)\Big(\sum_{j=1}^4E_{jj}\otimes(E_{jj}\otimes I_4^{(k)})\Big)\\
&  & \left(I_4\otimes u_{2k-1}(\xi\otimes I_4^{(k)})\right)\left(I_4\otimes(\Delta_4\otimes M_4^{(k)})\right)\\
&=& \sum_{j=1}^4E_{jj}\otimes\left((\Delta_4\otimes M_4^{(k)})(\xi^*\otimes I_4^{(k)})u_{2k-1}^*(E_{jj}\otimes I_4^{(k)})u_{2k-1}(\xi\otimes I_4^{(k)})(\Delta_4\otimes M_4^{(k)})\right)\\
&=& \sum_{j=1}^4E_{jj}\otimes\left((\Delta_4\otimes M_4^{(k)})(\xi^*\otimes I_4^{(k)})(E_{jj}\otimes I_4^{(k)})(\xi\otimes I_4^{(k)})(\Delta_4\otimes M_4^{(k)})\right)\\
&=& \sum_{j=1}^4\left(E_{jj}\otimes\Delta_4\xi^*E_{jj}\xi\Delta_4\right)\otimes M_4^{(k)}\,.
\end{IEEEeqnarray*}
Here, the second last equality follows from the following fact that
\[
u_{2k-1}^*\left(E_{jj}\otimes I_4^{(k)}\right)u_{2k-1}=E_{jj}\otimes I_4^{(k)}
\]
by using the induction hypothesis and the fact that $W_j^{(a)},\,j=1,2,3,$ are unitary matrices. Now, it is a straightforward verification that
\[
\xi^*E_{22}\xi=W_1^{(a)}\xi^*E_{11}\xi \big(W_1^{(a)}\big)^*\,,\,\xi^*E_{33}\xi=W_2^{(a)}\xi^*E_{11}\xi \big(W_2^{(a)}\big)^*\,,\,\xi^*E_{44}\xi=W_3^{(a)}\xi^*E_{11}\xi \big(W_3^{(a)}\big)^*\,.
\]
Therefore, if we set $W_0^{(a)}=I_4$, then we have the following,
\begin{IEEEeqnarray*}{lCl}
&  & (I_4\otimes u_{2k})\big(I_4\otimes(\Delta_4\otimes M_4^{(k)})\big)(I_4\otimes u_{2k})^*(e_1\otimes I_4^{(k)})(I_4\otimes u_{2k})\big(I_4\otimes(\Delta_4\otimes M_4^{(k)})\big)(I_4\otimes u_{2k})^*\\
&=& (I_4\otimes u_{2k})\Big(\sum_{j=1}^4\left((I_4\otimes W_j^{(a)})(E_{jj}\otimes\Delta_4\xi^*E_{11}\xi\Delta_4)(I
_4\otimes W_3^{(a)})^*\right)\otimes M_4^{(k)}\Big)(I_4\otimes u_{2k})^*\\
&=& \mbox{bl-diag}\left\{u_{2k},u_{2k}\left(W_1^{(a)}\otimes I_4^{(k)}\right),u_{2k}\left(W_2^{(a)}\otimes I_4^{(k)}\right),u_{2k}\left(W_3^{(a)}\otimes I_4^{(k)}\right)\right\}\big((I_4\otimes\Delta_4\xi^*E_{11}\xi\Delta_4)\otimes M_4^{(k)}\big)\\
&  & \mbox{bl-diag}\left\{u_{2k}^*, \left(W_1^{(a)}\otimes I_4^{(k)}\right)^*u_{2k}^*,\left(W_2^{(a)}\otimes I_4^{(k)}\right)u_{2k}^*,\left(W_3^{(a)}\otimes I_4^{(k)}\right)^*u_{2k}^*\right\}\\
&=& u_{2k+1}\big(I_4\otimes M_4\otimes M_4^{(k)}\big)u_{2k+1}^*
\end{IEEEeqnarray*}
with $u_{2k+1}$ in the required form if we show that $\Delta_4\big(\xi^*E_{11}\xi\big)\Delta_4=M_4\,.$ However, this is obvious since $\xi^*E_{11}\xi=J_4$, which is a direct verification, and $\Delta_4J_4\Delta_4=M_4$ because for any $i,j\in\{1,\ldots,4\}$, one has $E_{ii}J_4E_{jj}=E_{ij}$.
\smallskip

Now we prove part $(ii)$. Assume that the result is true up to $u_{2k}$. We also have $u_{2k+1}$ in our hand now since part (i) is already proved. Observe that $e_{2k+2}=J_4\otimes I_4^{(k+1)}$. Now, we have the following,
\begin{IEEEeqnarray*}{lCl}
&  & \left(I_4\otimes(M_4\otimes M_4^{(k)})\right)u_{2k+1}^*\left(J_4\otimes I_4^{(k+1)}\right)u_{2k+1}\left(I_4\otimes(M_4\otimes M_4^{(k)})\right)\\
&=& \left(\xi\otimes I_4^{(k+1)}\right)\left((\xi\otimes I_4^{(k+1)})^*(I_4\otimes M_4^{(k+1)})(\xi\otimes I_4^{(k+1)})\right)\Big((\xi\otimes I_4^{(k+1)})^*u_{2k+1}^*(J_4\otimes I_4^{(k+1)})\\
&  & u_{2k+1}(\xi\otimes I_4^{(k+1)})\Big)\left((\xi\otimes I_4^{(k+1)})^*(I_4\otimes M_4^{(k+1)})(\xi\otimes I_4^{(k+1)})\right)\left(\xi\otimes I_4^{(k+1)})\right)^*\\
&=& \left(\xi\otimes I_4^{(k+1)}\right)(I_4\otimes M_4^{(k+1)})\left((\xi\otimes I_4^{(k+1)})^*u_{2k+1}^*(J_4\otimes I_4^{(k+1)})u_{2k+1}(\xi\otimes I_4^{(k+1)})\right)\\
&  & (I_4\otimes M_4^{(k+1)})\left(\xi\otimes I_4^{(k+1)}\right)^*\,.
\end{IEEEeqnarray*}
Now, using part (i) for $u_{2k+1}$ and the induction hypothesis, we have the following,
\begin{IEEEeqnarray*}{lCl}
&  & \left(\xi\otimes I_4^{(k+1)}\right)^*u_{2k+1}^*\left(J_4\otimes I_4^{(k+1)}\right)u_{2k+1}\left(\xi\otimes I_4^{(k+1)}\right)\\
&=& \left(\xi\otimes I_4^{(k+1)}\right)^*\mbox{bl-diag}\left\{I_4\otimes I_4^{(k)},\left(W_1^{(a)}\right)^*\otimes I_4^{(k)},W_2^{(a)}\otimes I_4^{(k)},\left(W_3^{(a)}\right)^*\otimes I_4^{(k)}\right\}(I_4\otimes u_{2k})^*\\
&  & \left(J_4\otimes I_4^{(k+1)}\right)(I_4\otimes u_{2k})\,\mbox{bl-diag}\left\{I_4\otimes I_4^{(k)}, W_1^{(a)}\otimes I_4^{(k)},W_2^{(a)}\otimes I_4^{(k)},W_3^{(a)}\otimes I_4^{(k)}\right\}\left(\xi\otimes I_4^{(k+1)}\right)\\
&=& \left(\xi\otimes I_4^{(k+1)}\right)^*\mbox{bl-diag}\left\{\xi^*\otimes I_4^{(k)},\left(W_1^{(a)}\right)^*\xi^*\otimes I_4^{(k)},W_2^{(a)}\xi^*\otimes I_4^{(k)},\left(W_3^{(a)}\right)^*\xi^*\otimes I_4^{(k)}\right\}\\
&  & \left(J_4\otimes I_4^{(k+1)}\right)\,\mbox{bl-diag}\left\{\xi\otimes I_4^{(k)},\xi W_1^{(a)}\otimes I_4^{(k)},\xi W_2^{(a)}\otimes I_4^{(k)},\xi W_3^{(a)}\otimes I_4^{(k)}\right\}\left(\xi\otimes I_4^{(k+1)}\right)\\
&=& \left(\xi\otimes I_4^{(k+1)}\right)^*\Big(\Big(\mbox{bl-diag}\left\{\xi^*,\left(W_1^{(a)}\right)^*\xi^*,W_2^{(a)}\xi^*,\left(W_3^{(a)}\right)^*\xi^*\right\}(J_4\otimes I_4)\\
&  & \mbox{bl-diag}\left\{\xi,\xi W_1^{(a)},\xi W_2^{(a)},\xi W_3^{(a)}\right\}\Big)\otimes I_4^{(k)}\Big)\left(\xi\otimes I_4^{(k+1)}\right)\\
&=& \left(\xi\otimes I_4^{(k+1)}\right)^*\Big(\Big(\mbox{bl-diag}\left\{I_4,\left(W_1^{(a)}\right)^*,W_2^{(a)},\left(W_3^{(a)}\right)^*\right\}(J_4\otimes I_4)\\
&  & \mbox{bl-diag}\left\{I_4,W_1^{(a)},W_2^{(a)},W_3^{(a)}\right\}\Big)\otimes I_4^{(k)}\Big)\left(\xi\otimes I_4^{(k+1)}\right)\\
&=& \Big((\xi^*\otimes I_4)\,\mbox{bl-diag}\left\{I_4,\left(W_1^{(a)}\right)^*,W_2^{(a)},\left(W_3^{(a)}\right)^*\right\}(J_4\otimes I_4)\\
&  & \mbox{bl-diag}\left\{I_4,W_1^{(a)},W_2^{(a)},W_3^{(a)}\right\}(\xi\otimes I_4)\Big)\otimes I_4^{(k)}\\
&=& \mbox{bl-diag}\left\{E_{11},E_{22},E_{33},E_{44}\right\}\otimes I_4^{(k)}
\end{IEEEeqnarray*}
in $\Delta_4\otimes\Delta_4\otimes M_4^{(k)}$. Here, the last line follows from the proof of Lemma \ref{basis step 22}. We finally observe that
\begin{center}
$\mbox{Alg}\left\{(I_4\otimes A)\mbox{bl-diag}\left\{E_{11},E_{22},E_{33},E_{44}\right\}(I_4\otimes B):A,B\in M_4(\bbc)\right\}=\Delta_2\otimes M_4(\bbc)$
\end{center}
which concludes the proof of part (ii).\qed

\bigskip

\bigskip

\noindent{\sc Keshab Chandra Bakshi} (\texttt{keshab@iitk.ac.in, bakshi209@gmail.com})\\
         {\footnotesize Department of Mathematics and Statistics,\\
         Indian Institute of Technology, Kanpur,\\
         Uttar Pradesh 208016, India}
\bigskip

\noindent{\sc Satyajit Guin} (\texttt{sguin@iitk.ac.in})\\
         {\footnotesize Department of Mathematics and Statistics,\\
         Indian Institute of Technology, Kanpur,\\
         Uttar Pradesh 208016, India}

\end{document}